\documentclass{elsarticle}
  
\usepackage[margin=1in]{geometry}

\usepackage{bbm}
\usepackage{graphicx}
\usepackage{pstool}
\usepackage{wrapfig}
\usepackage{caption}
\usepackage{subcaption}
\usepackage[section]{placeins} 

\usepackage{fancyhdr} 
\usepackage{lastpage} 
\usepackage{extramarks} 

\usepackage{xcolor} 
\usepackage{enumerate} 
\usepackage{paralist} 
\usepackage{amsmath, amsthm, amssymb, mathtools}
\usepackage{mathabx, pifont, stmaryrd} 
\usepackage[explicit]{titlesec} 
\usepackage{etoolbox} 
\usepackage{bibentry} 
\makeatletter\let\saved@bibitem\@bibitem\makeatother 
\usepackage[colorlinks, bookmarksopen, bookmarksnumbered,
            citecolor=red,urlcolor=red]{hyperref} 
\makeatletter\let\@bibitem\saved@bibitem\makeatother 
\usepackage{rotating}

\usepackage{algorithm}
\usepackage{algorithmic}

\newtheorem{remark}{Remark}

\theoremstyle{definition}

\usepackage{bm} 
\usepackage{amsfonts} 

\DeclareMathOperator*{\argmin}{arg\,min}

\newcommand{\ds}[1]{\ensuremath{\displaystyle{#1}}}
\newcommand{\func}[3]{\ensuremath{#1 : #2 \rightarrow #3}}

\newcommand{\norm}[1]{\ensuremath{\left\| #1 \right\|}}
\newcommand{\jump}[1]{\ensuremath{\left\llbracket #1 \right\rrbracket}}
\newcommand{\suchthat}{\mathrel{}\middle|\mathrel{}}

\newcommand{\optconOne}[3]{
\begin{aligned}
& \underset{#1}{\text{minimize}}
& & #2 \\
& \text{subject to} & & #3
\end{aligned}}

\newcommand{\pder}[2]{\ensuremath{\frac{\partial #1}{\partial #2}}}

\newcommand{\Bcal}{\ensuremath{\mathcal{B}}}

\newcommand{\Ecal}{\ensuremath{\mathcal{E}}}
\newcommand{\Fcal}{\ensuremath{\mathcal{F}}}
\newcommand{\Gcal}{\ensuremath{\mathcal{G}}}
\newcommand{\Hcal}{\ensuremath{\mathcal{H}}}
\newcommand{\Ical}{\ensuremath{\mathcal{I}}}
\newcommand{\Jcal}{\ensuremath{\mathcal{J}}}

\newcommand{\Lcal}{\ensuremath{\mathcal{L}}}

\newcommand{\Ncal}{\ensuremath{\mathcal{N}}}
\newcommand{\Ocal}{\ensuremath{\mathcal{O}}}
\newcommand{\Pcal}{\ensuremath{\mathcal{P}}}

\newcommand{\Scal}{\ensuremath{\mathcal{S}}}
\newcommand{\Tcal}{\ensuremath{\mathcal{T}}}

\newcommand{\Vcal}{\ensuremath{\mathcal{V}}}
\newcommand{\Wcal}{\ensuremath{\mathcal{W}}}

\newcommand{\Gbb}{\ensuremath{\mathbb{G}}}

\newcommand{\Nbb}{\ensuremath{\mathbb{N} }}

\newcommand{\Rbb}{\ensuremath{\mathbb{R} }}

\newcommand\Abm{{\ensuremath{\bm{A}}}}
\newcommand\Bbm{{\ensuremath{\bm{B}}}}

\newcommand\Dbm{{\ensuremath{\bm{D}}}}

\newcommand\Fbm{{\ensuremath{\bm{F}}}}

\newcommand\Jbm{{\ensuremath{\bm{J}}}}

\newcommand\Rbm{{\ensuremath{\bm{R}}}}

\newcommand\Xbm{{\ensuremath{\bm{X}}}}

\newcommand\bbm{{\ensuremath{\bm{b}}}}
\newcommand\cbm{{\ensuremath{\bm{c}}}}

\newcommand\gbm{{\ensuremath{\bm{g}}}}

\newcommand\rbm{{\ensuremath{\bm{r}}}}

\newcommand\ubm{{\ensuremath{\bm{u}}}}
\newcommand\vbm{{\ensuremath{\bm{v}}}}

\newcommand\xbm{{\ensuremath{\bm{x}}}}
\newcommand\ybm{{\ensuremath{\bm{y}}}}
\newcommand\zbm{{\ensuremath{\bm{z}}}}

\newcommand\lambdabold{{\ensuremath{\boldsymbol{\lambda}}}}

\newcommand\etabold{{\ensuremath{\boldsymbol{\eta}}}}

\newcommand\phibold{{\ensuremath{\boldsymbol{\phi}}}}

\newcommand\zerobold{\ensuremath{\mathbf{0}}}

\usepackage{tikz}
\usepackage{pgfplots}
\usepackage{pgfplotstable, filecontents, booktabs}
\pgfplotsset{compat=1.9}

\usetikzlibrary{pgfplots.groupplots}
\usepgfplotslibrary{fillbetween}
\usetikzlibrary{calc,fit,matrix,arrows,automata,positioning,shapes}
\usetikzlibrary{arrows.meta}

\pgfplotsset{select coords between index/.style 2 args={
    x filter/.code={
        \ifnum\coordindex<#1\fi
        \ifnum\coordindex>#2\fi
    }
}}

\tikzset{
 invisible/.style={opacity=0},
 visible on/.style={alt={#1{}{invisible}}},
 alt/.code args={<#1>#2#3}{%
   \alt<#1>{\pgfkeysalso{#2}}{\pgfkeysalso{#3}}
 },
}

\newcommand{\colorbarMatlabParula}[5]{
\begin{tikzpicture}
\begin{axis}[
   hide axis, scale only axis,
   height=0pt, width=0pt,
   colormap={parula}{rgb255=(62,38,168) rgb255=(62,39,172) rgb255=(63,40,175) rgb255=(63,41,178) rgb255=(64,42,180) rgb255=(64,43,183) rgb255=(65,44,186) rgb255=(65,45,189) rgb255=(66,46,191) rgb255=(66,47,194) rgb255=(67,48,197) rgb255=(67,49,200) rgb255=(67,50,202) rgb255=(68,51,205) rgb255=(68,52,208) rgb255=(69,53,210) rgb255=(69,55,213) rgb255=(69,56,215) rgb255=(70,57,217) rgb255=(70,58,220) rgb255=(70,59,222) rgb255=(70,61,224) rgb255=(71,62,225) rgb255=(71,63,227) rgb255=(71,65,229) rgb255=(71,66,230) rgb255=(71,68,232) rgb255=(71,69,233) rgb255=(71,70,235) rgb255=(72,72,236) rgb255=(72,73,237) rgb255=(72,75,238) rgb255=(72,76,240) rgb255=(72,78,241) rgb255=(72,79,242) rgb255=(72,80,243) rgb255=(72,82,244) rgb255=(72,83,245) rgb255=(72,84,246) rgb255=(71,86,247) rgb255=(71,87,247) rgb255=(71,89,248) rgb255=(71,90,249) rgb255=(71,91,250) rgb255=(71,93,250) rgb255=(70,94,251) rgb255=(70,96,251) rgb255=(70,97,252) rgb255=(69,98,252) rgb255=(69,100,253) rgb255=(68,101,253) rgb255=(67,103,253) rgb255=(67,104,254) rgb255=(66,106,254) rgb255=(65,107,254) rgb255=(64,109,254) rgb255=(63,110,255) rgb255=(62,112,255) rgb255=(60,113,255) rgb255=(59,115,255) rgb255=(57,116,255) rgb255=(56,118,254) rgb255=(54,119,254) rgb255=(53,121,253) rgb255=(51,122,253) rgb255=(50,124,252) rgb255=(49,125,252) rgb255=(48,127,251) rgb255=(47,128,250) rgb255=(47,130,250) rgb255=(46,131,249) rgb255=(46,132,248) rgb255=(46,134,248) rgb255=(46,135,247) rgb255=(45,136,246) rgb255=(45,138,245) rgb255=(45,139,244) rgb255=(45,140,243) rgb255=(45,142,242) rgb255=(44,143,241) rgb255=(44,144,240) rgb255=(43,145,239) rgb255=(42,147,238) rgb255=(41,148,237) rgb255=(40,149,236) rgb255=(39,151,235) rgb255=(39,152,234) rgb255=(38,153,233) rgb255=(38,154,232) rgb255=(37,155,232) rgb255=(37,156,231) rgb255=(36,158,230) rgb255=(36,159,229) rgb255=(35,160,229) rgb255=(35,161,228) rgb255=(34,162,228) rgb255=(33,163,227) rgb255=(32,165,227) rgb255=(31,166,226) rgb255=(30,167,225) rgb255=(29,168,225) rgb255=(29,169,224) rgb255=(28,170,223) rgb255=(27,171,222) rgb255=(26,172,221) rgb255=(25,173,220) rgb255=(23,174,218) rgb255=(22,175,217) rgb255=(20,176,216) rgb255=(18,177,214) rgb255=(16,178,213) rgb255=(14,179,212) rgb255=(11,179,210) rgb255=(8,180,209) rgb255=(6,181,207) rgb255=(4,182,206) rgb255=(2,183,204) rgb255=(1,183,202) rgb255=(0,184,201) rgb255=(0,185,199) rgb255=(0,186,198) rgb255=(1,186,196) rgb255=(2,187,194) rgb255=(4,187,193) rgb255=(6,188,191) rgb255=(9,189,189) rgb255=(13,189,188) rgb255=(16,190,186) rgb255=(20,190,184) rgb255=(23,191,182) rgb255=(26,192,181) rgb255=(29,192,179) rgb255=(32,193,177) rgb255=(35,193,175) rgb255=(37,194,174) rgb255=(39,194,172) rgb255=(41,195,170) rgb255=(43,195,168) rgb255=(44,196,166) rgb255=(46,196,165) rgb255=(47,197,163) rgb255=(49,197,161) rgb255=(50,198,159) rgb255=(51,199,157) rgb255=(53,199,155) rgb255=(54,200,153) rgb255=(56,200,150) rgb255=(57,201,148) rgb255=(59,201,146) rgb255=(61,202,144) rgb255=(64,202,141) rgb255=(66,202,139) rgb255=(69,203,137) rgb255=(72,203,134) rgb255=(75,203,132) rgb255=(78,204,129) rgb255=(81,204,127) rgb255=(84,204,124) rgb255=(87,204,122) rgb255=(90,204,119) rgb255=(94,205,116) rgb255=(97,205,114) rgb255=(100,205,111) rgb255=(103,205,108) rgb255=(107,205,105) rgb255=(110,205,102) rgb255=(114,205,100) rgb255=(118,204,97) rgb255=(121,204,94) rgb255=(125,204,91) rgb255=(129,204,89) rgb255=(132,204,86) rgb255=(136,203,83) rgb255=(139,203,81) rgb255=(143,203,78) rgb255=(147,202,75) rgb255=(150,202,72) rgb255=(154,201,70) rgb255=(157,201,67) rgb255=(161,200,64) rgb255=(164,200,62) rgb255=(167,199,59) rgb255=(171,199,57) rgb255=(174,198,55) rgb255=(178,198,53) rgb255=(181,197,51) rgb255=(184,196,49) rgb255=(187,196,47) rgb255=(190,195,45) rgb255=(194,195,44) rgb255=(197,194,42) rgb255=(200,193,41) rgb255=(203,193,40) rgb255=(206,192,39) rgb255=(208,191,39) rgb255=(211,191,39) rgb255=(214,190,39) rgb255=(217,190,40) rgb255=(219,189,40) rgb255=(222,188,41) rgb255=(225,188,42) rgb255=(227,188,43) rgb255=(230,187,45) rgb255=(232,187,46) rgb255=(234,186,48) rgb255=(236,186,50) rgb255=(239,186,53) rgb255=(241,186,55) rgb255=(243,186,57) rgb255=(245,186,59) rgb255=(247,186,61) rgb255=(249,186,62) rgb255=(251,187,62) rgb255=(252,188,62) rgb255=(254,189,61) rgb255=(254,190,60) rgb255=(254,192,59) rgb255=(254,193,58) rgb255=(254,194,57) rgb255=(254,196,56) rgb255=(254,197,55) rgb255=(254,199,53) rgb255=(254,200,52) rgb255=(254,202,51) rgb255=(253,203,50) rgb255=(253,205,49) rgb255=(253,206,49) rgb255=(252,208,48) rgb255=(251,210,47) rgb255=(251,211,46) rgb255=(250,213,46) rgb255=(249,214,45) rgb255=(249,216,44) rgb255=(248,217,43) rgb255=(247,219,42) rgb255=(247,221,42) rgb255=(246,222,41) rgb255=(246,224,40) rgb255=(245,225,40) rgb255=(245,227,39) rgb255=(245,229,38) rgb255=(245,230,38) rgb255=(245,232,37) rgb255=(245,233,36) rgb255=(245,235,35) rgb255=(245,236,34) rgb255=(245,238,33) rgb255=(246,239,32) rgb255=(246,241,31) rgb255=(246,242,30) rgb255=(247,244,28) rgb255=(247,245,27) rgb255=(248,247,26) rgb255=(248,248,24) rgb255=(249,249,22) rgb255=(249,251,21) },
   colorbar horizontal,
   point meta min=#1, point meta max=#5,
   colorbar style={width=10cm, xtick={#1,#2,#3,#4,#5}}
]
\addplot [draw=none] coordinates {(0,0)};
\end{axis}
\end{tikzpicture}
}

\usepackage{setspace}

\newbool{fastcompile}
\setbool{fastcompile}{false}

\begin{document}
\title{A robust, high-order implicit shock tracking method for simulation of complex, high-speed flows}

\author[rvt1]{Tianci Huang\fnref{fn1}}
\ead{thuang5@nd.edu}

\author[rvt1]{Matthew J. Zahr\fnref{fn2}\corref{cor1}}
\ead{mzahr@nd.edu}

\address[rvt1]{Department of Aerospace and Mechanical Engineering, University
               of Notre Dame, Notre Dame, IN 46556, United States}
\cortext[cor1]{Corresponding author}

\fntext[fn1]{Graduate Student, Department of Aerospace and Mechanical
             Engineering, University of Notre Dame}
\fntext[fn2]{Assistant Professor, Department of Aerospace and Mechanical
             Engineering, University of Notre Dame}

\begin{keyword} 
Shock tracking, shock fitting, high-order methods, discontinuous Galerkin, numerical optimization,
high-speed flows
\end{keyword}

\begin{abstract}
High-order implicit shock tracking (fitting) is a new class of numerical methods to
approximate solutions of conservation laws with non-smooth features, e.g., contact lines,
shock waves, and rarefactions. These methods align elements of the computational mesh
with non-smooth features to represent them perfectly, allowing high-order basis functions
to approximate smooth regions of the solution without the need for nonlinear stabilization, which
leads to accurate approximations on traditionally coarse meshes. The hallmark of these
methods is the underlying optimization formulation whose solution is a feature-aligned
mesh and the corresponding high-order approximation to the flow; the key challenge is
robustly solving the central optimization problem.
In this work, we develop a robust optimization solver
for high-order implicit shock tracking methods so they can be reliably used to simulate complex,
high-speed, compressible flows in multiple dimensions. The proposed method integrates
practical robustness measures into a sequential quadratic programming method that employs
a Levenberg-Marquardt approximation of the Hessian and line-search globalization.
The robustness measures include dimension- and order-independent simplex element
collapses, mesh smoothing, and element-wise solution re-initialization, which prove to be 
necessary to reliably track complex discontinuity surfaces, such as curved and reflecting
shocks, shock formation, and shock-shock interaction. A series of nine numerical experiments---%
including two- and three-dimensional compressible flows with complex discontinuity
surfaces---are used to demonstrate:
\begin{inparaenum}[1)]
 \item the robustness of the solver,
 \item the meshes produced are high-quality and track continuous, non-smooth features
  in addition to discontinuities,
 \item the method achieves the optimal convergence rate of the underlying discretization even
  for flows containing discontinuities, and
 \item the method produces highly accurate solutions on extremely coarse meshes relative to
  approaches based on shock capturing.
\end{inparaenum}
\end{abstract}

\maketitle

\section{Introduction}
\label{sec:intro}
High-order methods such as discontinuous Galerkin (DG) methods
\cite{cockburn_rungekutta_2001, hesthaven_nodal_2008} offer a
number of advantages for high-fidelity simulation of fluid flow
including high accuracy per degree of freedom, low dissipation,
geometric flexibility, and a high degree of parallel scalability.
However, for high-speed flows that contain shocks and other
discontinuities, spurious oscillations arise that degrade
the approximation accuracy and usually lead to failure of the
simulation. These nonlinear instabilities become increasingly
problematic as the strength of the shock increases. The
development of robust and accurate approaches to stabilize
high-order approximations near shocks is critical to make
them competitive for real-world problems \cite{wang_high-order_2013}.

Several approaches have been proposed to stabilize shocks, most of which
are based on shock capturing, where the numerical discretization accounts
for discontinuities on a fixed computational grid. Limiters, used to
limit the gradient of the solution in the vicinity of shocks, are
commonly used with second-order finite volume methods
\cite{van1979towards} and high-order DG methods
\cite{cockburn_rungekutta_2001} to yield total variation
diminishing schemes; these methods are widely accepted and
commonly used to simulate real-world flows \cite{candler_advances_2015}.
Weighted essentially non-oscillatory (WENO) methods
\cite{harten_uniformly_1987, liu_weighted_1994, jiang_efficient_1996}
use a high-order reconstruction with a stencil tailored to the
flow solution to mitigate spurious oscillations in the solutions
near shocks. These methods lead to crisp shocks, although they
require a large stencil and have not been demonstrated on
unstructured meshes for practical problems.
For high-order methods, artificial viscosity
has also proven to be competitive, and has emerged as the
preferred method for finite-element-based methods
\cite{persson_sub-cell_2006,barter_shock_2010,
      fernandez_physics-based_2018,ching_shock_2019}
because it can smoothly resolve steep gradients with
sub-cell accuracy. A recent comparative study of artificial
viscosity models \cite{yu_artificial_2020} discusses their relative
merits, but notes they all suffer from a relatively strong dependency
on a large number of empirical parameters that must be tuned.
The main problem with all these approaches is that they are first-order
near the shock, which translates into a globally first-order accurate
scheme. This can be remedied by using local mesh refinement around the
shock ($h$-adaptivity) \cite{dervieux_about_2003}, although the anisotropic
elements that are required for efficiency are difficult to generate and
extremely fine elements are needed near the shock.

An alternative approach is so-called \emph{shock tracking} or \emph{shock fitting}, 
where the computational mesh is moved to align faces of mesh elements with the 
solution discontinuities, representing them perfectly with the inter-element
jump in the solution basis without requiring additional stabilization.
However, it is a difficult meshing problem since it essentially requires
generating a fitted mesh to the (unknown) shock surface. Most of these
methods employ specialized formulations and solvers which are
dimension-dependent and do not easily generalize
\cite{harten1983self, glimm2003conservative, bell1982fully}
and/or are limited to relatively simple problems
\cite{shubin1981steady, shubin1982steady, rosendale1994floating}.
In addition, early approaches to shock tracking have been applied to
low-order schemes where the relative advantage over shock capturing is
smaller than for high-order methods
\cite{trepanier1996conservative, baines2002multidimensional}.
One popular class of methods that we will call explicit shock tracking
is surveyed in \cite{moretti2002thirty, salas2009shock}.
These methods largely consist of explicitly identifying the
shock and using the Rankine-Hugoniot conditions
to compute its motion and states upstream and downstream of the shock.
More recent developments in explicit shock tracking \cite{rawat2010high}
use more sophisticated methods to compute shock velocities and
discretize the flow equations; however, they ultimately
still require a specialized strategy to explicitly track
the shock separately from the remainder of the
flow. These methods are not easily applicable to
discontinuities whose topologies are not known \textit{a priori}.
While interest in shock tracking/fitting has seen
somewhat of a resurgence in recent years
\cite{ciallella_extrapolated_2020, bonfiglioli_unsteady_2016, geisenhofer_extended_2020, daquila_novel_2021},
shock tracking is not widely used for practical applications.

A new class of high-order numerical methods has recently emerged,
\textit{implicit shock tracking}, that includes the High-Order
Implicit Shock Tracking (HOIST) method
\cite{zahr_optimization-based_2018,zahr_implicit_2020, zahr_high-order_2020}
and the Moving Discontinous Galerkin Method with Interface Condition Enforcement
(MDG-ICE) \cite{corrigan_moving_2019,kercher_least-squares_2020,kercher_moving_2021}.
Like traditional shock tracking, these methods
align element faces with discontinuity surfaces to represent them
perfectly and the high-order basis functions approximate the smooth
solution away from shocks. The key difference is the implicit
shock tracking methods do not attempt to explicitly generate
a mesh fitted to the unknown discontinuity surface, e.g., by treating
the discontinuity as an interior boundary. Rather, they discretize
the conservation law on a mesh without knowledge of the flow field
and pose an optimization problem over the discrete flow variables
and nodal coordinates of the mesh whose solution is a discontinuity-aligned
mesh and the corresponding flow solution. That is, discontinuity tracking
is achieved implicitly through the solution of an optimization problem.
This is the critical innovation that has overcome the limitations of
explicit shock tracking and led to a general approach
that is not tailored to the governing equations or specific
flow problem (equation- and problem-independent), only
requires one simple topological mesh operation (element removal),
and is able to handle intricate shock structures (curved and
reflecting shocks, shock formation, and shock-shock interaction).
As such, implicit shock tracking has been used to solve steady and
unsteady, viscous and inviscid, inert and reacting flows of varying
degrees of complexity. While these methods have shown considerable
promise, they are still lacking a robust solver for the underlying
optimization problem that converges efficiently and reliably to a
discontinuity-aligned mesh with high-quality elements for high Mach
flows with complex discontinuity structures. This is the gap this
work aims to fill.

We propose a series of algorithmic developments to the
HOIST method to improve its automation and robustness,
although many of the developments will apply to MDG-ICE
with minor modifications. The HOIST method is based on a
high-order DG discretization of the governing equations and
formulates implicit tracking as an optimization problem constrained
by the DG residual to endow the method with the desirable
properties of DG: consistency, conservation, and stability.
The objective function penalizes violation of the DG residual
in an enriched test space; it 
is a surrogate for violation of the infinite-dimensional weak
formulation of the conservation law, which endows the method with
$r$-adaptive behavior. The optimization problem is solved using a
sequential quadratic programming (SQP) method with a Levenberg-Marquardt
Hessian approximation that simultaneously converges the mesh and flow solution
to their optimal values, which never requires the fully converged DG solution
on a non-aligned mesh and therefore does not require nonlinear stabilization.
The contributions of this paper are summarized as follows.
\begin{itemize}
 \item[-] We develop a general, automated procedure to guarantee
  planar boundaries and their intersections will be preserved as the
  computational mesh moves to align with shocks and detail its implementation.
 \item[-] We develop a new solver for the implicit tracking optimization problem
  based on the SQP solver in \cite{zahr_implicit_2020}.
The new solver features
  a new merit function penalty parameter, an adaptive penalty parameter
  for the mesh distortion term in the objective function, and, most
  importantly, a number of practical robustness measures. The most
  critical robustness measures are dimension- and order-independent
  simplex element removal via edge collapse that preserves the
  boundaries of the domain and the tracked shocks (to high-order) and
  element-wise solution re-initialization that resets the solution
  in oscillatory elements, identified with standard shock sensors
  \cite{persson_sub-cell_2006}, to a constant value to promote high-quality
  SQP steps.
 \item[-] We demonstrate implicit shock tracking, when equipped with robust solvers,
  is practical for two-dimensional, compressible flows with intricate shock
  structures (curved and reflecting shocks, shock formation, and shock-shock
  interaction); three-dimensional, compressible flows with a curved bow shock;
  and high Mach flows.
\end{itemize}
The remainder of the paper is organized as follows.
Section~\ref{sec:govern} introduces the governing system of inviscid
conservation laws, its reformulation on a fixed reference domain,
and its discretization using a discontinuous Galerkin method.
Section~\ref{sec:dommap} constructs the space of admissible domain mappings
(i.e., mesh motion) using nodal basis functions associated with a mesh
of the reference domain and a parametrization that guarantees the
mapping will approximate the boundary of the physical domain to
high-order accuracy; a general, automated procedure to construct
the parametrization is introduced for the special case where the
boundary is piecewise planar.
Section~\ref{sec:ist-form} recalls the implicit shock tracking formulation
originally proposed in \cite{zahr_implicit_2020} that incorporates the
boundary-preserving parametrization of the nodal coordinates.
Section~\ref{sec:ist_solve} introduces the SQP solver with the new
robustness measures. 
Finally, Section~\ref{sec:numexp} presents a series of nine increasingly
difficult numerical experiments that demonstrate the robustness of the
solver, the high accuracy per degree of freedom of the method, and ability
of the method to handle two- and three-dimensional compressible flows with
complex features (curved and reflecting shocks, shock formation, and 
shock-shock interaction).

\section{Governing equations and high-order discretization}
\label{sec:govern}
In this section, we introduce the governing partial differential
equations (steady or space-time, inviscid conservation law)
(Section~\ref{sec:govern:claw}), its transformation to a reference
domain so that domain deformations appear explicitly in the governing
equations (Section~\ref{sec:govern:transf}), and its discretization
via a high-order DG method (high-order with respect to both the
solution and geometry) (Section~\ref{sec:govern:dg}).

\subsection{System of conservation laws}
\label{sec:govern:claw}
Consider a general system of $m$ inviscid conservation laws, defined on the
fixed domain $\Omega \subset \Rbb^d$ and subject to appropriate boundary
conditions,
\begin{equation} \label{eqn:claw-phys}
 \nabla\cdot F(U) = S(U) \quad \text{in}~~\Omega,
\end{equation}
where $\func{U}{\Omega}{\Rbb^m}$ is the solution of the system of
conservation laws, $\func{F}{\Rbb^m}{\Rbb^{m\times d}}$ is the flux
function, $\func{S}{\Rbb^m}{\Rbb^m}$ is the source term,
$\ds{\nabla \coloneqq (\partial_{x_1},\dots,\partial_{x_d})}$
is the gradient operator in the physical domain, and the boundary of
the domain $\partial\Omega$ has outward unit normal
$\func{n}{\partial\Omega}{\Rbb^d}$. The formulation of the conservation
law in (\ref{eqn:claw-phys}) is sufficiently general to encapsulate
steady conservation laws in a $d$-dimensional spatial domain or
time-dependent conservation laws in a $(d-1)$-dimensional domain, i.e.,
a $d$-dimensional space-time domain. In general, the solution $U(x)$ may
contain discontinuities, in which case, the conservation law
(\ref{eqn:claw-phys}) holds away from the discontinuities
and the Rankine-Hugoniot conditions \cite{majda2012compressible}
hold at discontinuities.

\subsection{Transformed system of conservation laws on a fixed reference domain}
\label{sec:govern:transf}
Before discretizing (\ref{eqn:claw-phys}), it is convenient
to explicitly treat deformations to the domain of the conservation law
$\Omega$---such deformations will eventually be induced by deformation
to the mesh as nodal coordinates are moved to track discontinuities---by
transforming to a fixed reference domain $\Omega_0\subset\Rbb^d$.
Let $\Gbb$ be the collection of diffeomorphisms from the reference domain
$\Omega_0$ to the physical domain $\Omega$, i.e., for any $\Gcal\in\Gbb$,
we have
\begin{equation} \label{eqn:dom-map}
 \func{\Gcal}{\Omega_0}{\Omega}, \quad
 \Gcal : X \mapsto \Gcal(X).
\end{equation}
Following the approach in \cite{zahr_optimization-based_2018},
for any $\Gcal\in\Gbb$, the conservation law on the physical
domain $\Omega$ is transformed to a conservation law on
the reference domain $\Omega_0$ as
\begin{equation} \label{eqn:claw-ref}
 \bar\nabla\cdot\bar{F}(\bar{U};G) =
 \bar{S}(\bar{U};g) \quad \text{in}~~\Omega_0,
\end{equation}
where $\func{\bar{U}}{\Omega_0}{\Rbb^m}$ is the solution of the transformed
conservation law, $\func{\bar{F}}{\Rbb^m\times\Rbb^{d\times d}}{\Rbb^{m\times d}}$ is
the transformed flux function, $\bar\nabla \coloneqq (\partial_{X_1},\dots,\partial_{X_d})$ is the gradient operator on the reference domain, and the
deformation gradient $\func{G}{\Omega_0}{\Rbb^{d\times d}}$ and mapping
Jacobian $\func{g}{\Omega_0}{\Rbb}$ are defined as
\begin{equation}
 G = \bar\nabla \Gcal, \qquad g = \det G.
\end{equation}
The unit outward normal to the reference domain is denoted
$\func{N}{\partial\Omega_0}{\Rbb^d}$  and related to the unit
normal in the physical domain by
\begin{equation} \label{eqn:transf-normal}
 n \circ \Gcal = \frac{g G^{-T}N}{\norm{g G^{-T}N}}.
\end{equation}
For any $X\in\Omega_0$, the transformed and physical solution are related as
\begin{equation}
 \bar{U}(X) = U(\Gcal(X))
\end{equation}
and the transformed flux and source term are defined as
\begin{equation}
 \bar{F} : (\bar{W}; \Theta) \mapsto (\det\Theta) F(\bar{W}) \Theta^{-T},
 \qquad
 \bar{S} : (\bar{W}; q) \mapsto q S(\bar{W}).
\end{equation}
\begin{remark} \label{rem:refphys}
In general, the reference domain can be defined such that it maps to the
physical domain under the action of a smooth, invertible mapping
$\func{\hat{\Gcal}}{\Rbb^d}{\Rbb^d}$, i.e., $\Omega_0 = \hat{\Gcal}^{-1}(\Omega)$.
In this work, we take the reference and physical domains to be the same
set, i.e., $\hat{\Gcal}=\mathrm{Id}$.
\end{remark}

\subsection{Discontinuous Galerkin discretization of
            the transformed conservation law}
\label{sec:govern:dg}
We use a nodal discontinuous Galerkin method
\cite{cockburn_rungekutta_2001, hesthaven_nodal_2008}
to discretize the transformed conservation law (\ref{eqn:claw-ref}).
Let $\Ecal_h$ represent a discretization of the reference domain
$\Omega_0$ into non-overlapping, potentially curved, computational elements.
To establish the finite-dimensional DG formulation, we introduce
the DG approximation (trial) space of discontinuous piecewise
polynomials associated with the mesh $\Ecal_h$
\begin{equation}
 \Vcal_h^p = \left\{v \in [L^2(\Omega_0)]^m \suchthat
         \left.v\right|_K \in [\Pcal_p(K)]^m,
         ~\forall K \in \Ecal_h\right\},
\end{equation}
where $\Pcal_p(K)$ is the space of polynomial functions of degree at most
$p \geq 1$ on the element $K$. Furthermore, we define the space of globally
continuous piecewise polynomials of degree $q$ associated with the mesh
$\Ecal_h$ as
\begin{equation} \label{eqn:gfcnsp}
  \Wcal_h = \left\{v \in C^0(\Omega_0) \suchthat
         \left.v\right|_K \in \Pcal_q(K),~\forall K \in \Ecal_h\right\}
\end{equation}
and discretize the domain mapping with the corresponding vector-valued
space $[\Wcal_h]^d$.

Taking the DG test space to be $\Vcal_h^{p'}$, where $p'\geq p$, the
DG formulation is: given $\Gcal_h\in[\Wcal_h]^d$, find
$\bar{U}_h\in\Vcal_h^p$ such that for all $\bar\psi_h\in\Vcal_h^{p'}$, we have
\begin{equation} \label{eqn:claw-weak-elem}
 \int_{\partial K} \bar\psi_h^+ \cdot
   \bar\Hcal(\bar{U}_h^+,\bar{U}_h^-,N_h;\bar\nabla\Gcal_h) \, dS -
 \int_K \bar{F}(\bar{U}_h; \bar\nabla\Gcal_h):\bar\nabla \bar\psi_h \, dV =
 \int_K \bar\psi_h \cdot \bar{S}(\bar{U}_h; \det(\bar\nabla\Gcal_h)) \, dV,
\end{equation}
where $\func{N_h}{\partial K}{\Rbb^d}$ is the unit outward normal to
element $K\in\Ecal_h$, $\bar{W}_h^+$ ($\bar{W}_h^-$) denotes the interior
(exterior) trace of $\bar{W}_h$ to the element $K$ for $\bar{W}_h\in\Vcal_h^s$
for any $s\in\Nbb$ (for $X\in\partial K\cap\partial\Omega_0$,
$\bar{U}_h^-$ is a boundary state constructed to enforce the appropriate
boundary condition). Furthermore, 
$\func{\bar\Hcal}{\Rbb^m\times\Rbb^m\times\Rbb^d\times\Rbb^{d\times d}}{\Rbb^m}$
is the numerical flux function associated with the reference inviscid flux
$\bar{F}$, which ensures the boundary integral is single-valued and can be
constructed to ensure the DG discretization is consistent, conservative,
and stable \cite{hesthaven_nodal_2008}. An expression for the reference numerical
flux function can be obtained from the standard physical numerical flux
function \cite{zahr_implicit_2020}. The residual form of the
DG equation in (\ref{eqn:claw-weak-elem}) is given by
$\func{r_h^{p',p}}{\Vcal_h^{p'}\times\Vcal_h^p\times[\Wcal_h]^d}{\Rbb}$
\begin{equation}
 r_h^{p',p} : (\bar\psi_h,\bar{W}_h,\Gcal_h) \mapsto
 \sum_{K\in\Ecal_h} r_K^{p',p}(\bar\psi_h,\bar{W}_h,\Gcal_h),
\end{equation}
where the elemental DG form is given by
$\func{r_K^{p',p}}{\Vcal_h^{p'}\times\Vcal_h^p\times[\Wcal_h]^d}{\Rbb}$
\begin{equation}
\begin{aligned}
 r_K^{p',p} : (\bar\psi_h,\bar{W}_h,\Gcal_h) \mapsto
 &\int_{\partial K} \bar\psi_h^+ \cdot
   \bar\Hcal(\bar{W}_h^+,\bar{W}_h^-,N_h;\bar\nabla\Gcal_h) \, dS \\
 &-\int_K \bar{F}(\bar{W}_h; \bar\nabla\Gcal_h):\bar\nabla \bar\psi_h \, dV \\
 &-\int_K \bar\psi_h \cdot \bar{S}(\bar{W}_h; \det(\bar\nabla\Gcal_h)) \, dV.
\end{aligned}
\end{equation}

Next, we introduce a (nodal) basis for the test space ($\Vcal_h^{p'}$),
trial space ($\Vcal_h^p$), and domain mapping space ($[\Wcal_h]^d$)
to reduce the weak formulation in residual form to a system of
nonlinear algebraic equations in residual form. In the case where
$p'=p$, we denote the algebraic residual
\begin{equation}
 \func{\rbm}{\Rbb^{N_\ubm}\times\Rbb^{N_\xbm}}{\Rbb^{N_\ubm}}, \qquad
 \rbm : (\ubm,\xbm) \mapsto \rbm(\ubm,\xbm)
\end{equation}
where $N_\ubm=\dim \Vcal_h^p$ and $N_\xbm=\dim([\Wcal_h]^d)$.
In this notation, a standard DG
discretization (algebraic form) reads: given $\xbm\in\Rbb^{N_\xbm}$,
find $\ubm\in\Rbb^{N_\ubm}$ such that $\rbm(\ubm,\xbm)=\zerobold$,
where $\ubm$ are the DG solution coefficients and $\xbm$ are the
coefficients of the domain mapping (nodal coordinates). Typically,
$\xbm$ is known
(mesh generation) and fixed; however, in this work, it will be
determined via optimization such that the mesh tracks (aligns
element faces with) all discontinuities in the flow.
Furthermore, we define the algebraic \textit{enriched residual}
\begin{equation}
 \func{\Rbm}{\Rbb^{N_\ubm}\times\Rbb^{N_\xbm}}{\Rbb^{N_\ubm'}}, \qquad
 \Rbm : (\ubm,\xbm) \mapsto \Rbm(\ubm,\xbm)
\end{equation}
associated with a trial space of degree $p'$, where $N_\ubm'=\dim\Vcal_h^{p'}$,
which will be used to construct the implicit shock tracking objective function.
In this work, we take $p'=p+1$.
Finally, to maintain a connection between the algebraic and functional
representation of the DG solution, we define the operator
$\func{\Xi}{\Rbb^{N_\ubm}}{\Vcal_h^p}$ that maps
$\vbm\in\Rbb^{N_\ubm}$ to its representation as a function
over the reference domain $\Omega_0$, $V_h = \Xi(\vbm)$,
i.e., $\vbm$ is the encoding of $V_h\in\Vcal_h^p$ as a
vector in $\Rbb^{N_\ubm}$.

\section{Construction of admissible domain mappings}
\label{sec:dommap}
In this section, we discretize the domain mapping on the computational
mesh $\Ecal_h$ such that the mapped mesh approximates the boundaries of
the physical domain to high-order, even as nodes are moved to track
flow features (Section~\ref{sec:dommap:disc}). Then, we construct a general,
automated approach for the special case where all boundaries are planar and
intersect to form non-smooth edges and corners (Section~\ref{sec:dommap:planar}).

\subsection{Domain mapping discretization and parametrization}
\label{sec:dommap:disc}
Let $\{\hat{X}_I\}_{I=1}^{N_\mathrm{v}}$ denote (ordered) nodes associated
with the mesh $\Ecal_h$ (including high-order nodes). Then, any element
of $[\Wcal_h]^d$ is uniquely determined by its action on the nodes
$\hat{X}_I$. That is, let $\{\Psi_I\}_{I=1}^{N_\mathrm{v}}$ be a nodal
basis of $\Wcal_h$ associated with the nodes $\{\hat{X}_I\}_{I=1}^{N_\mathrm{v}}$,
i.e., $\Psi_I(\hat{X}_J)=\delta_{IJ}$, and define
\begin{equation}
 \func{\Gcal_h(\,\cdot\,;\xbm)}{\Omega_0}{\Rbb^d}, \qquad
 \Gcal_h(\,\cdot\,;\xbm): X \mapsto
 \sum_{I=1}^{N_\mathrm{v}} \hat{x}_I \Psi_I(X),
\end{equation}
where the coefficients $\hat{x}_I\in\Rbb^d$ can be interpreted
as the physical coordinates of the reference nodes $\hat{X}_I$
because they coincide with the action of the mapping at $\hat{X}_I$
due to the choice of nodal basis and $\xbm\in\Rbb^{N_\xbm}$
($N_\xbm=dN_\mathrm{v}$) is the concatenation of $\{\hat{x}_I\}_{I=1}^{N_\mathrm{v}}$.
To ensure that $\Gcal_h(\,\cdot\,;\xbm)$ is a bijection from
$\Omega_0$ to $\Omega$, $\xbm\in\Rbb^{N_\xbm}$ must be
defined such that $\Gcal_h(\Ecal_h;\xbm)$ is a valid
mesh of $\Omega$, i.e., a partition of $\Omega$ into
non-overlapping and non-inverted elements. Because
$\Gcal_h(\,\cdot\,;\xbm)$ is continuous for any
$\xbm\in\Rbb^{N_\xbm}$, the mesh $\Gcal_h(\Ecal_h;\xbm)$
will inherit air-tightness (no gaps between elements) from
the reference mesh $\Ecal_h$. Therefore, $\xbm$ only needs to
be restricted to ensure the elements of $\Gcal_h(\Ecal_h;\xbm)$:
\begin{inparaenum}[1)]
 \item are not inverted (for injectivity) and
 \item conform to the boundary $\partial\Omega$ (for surjectivity with
  respect to $\Omega$).
\end{inparaenum}
The first condition is difficult to impose explicitly so it will
be accounted for in the implicit tracking optimization problem,
whereas the second condition will be explicitly baked into the
definition of the admissible mappings. To this end, we will
introduce a parametrization of the physical nodes
\begin{equation} \label{eqn:mapparam}
 \func{\phibold}{\Rbb^{N_\ybm}}{\Rbb^{N_\xbm}}, \qquad
 \phibold: \ybm \mapsto \phibold(\ybm)
\end{equation}
such that $\Gcal_h(\Ecal_h;\phibold(\ybm))$ conforms to $\partial\Omega$
for any $\ybm\in\Rbb^{N_\ybm}$ that does not cause element inversion.
References \cite{corrigan_convergence_2019,zahr_radaptive_2020,zahr_implicit_2020}
introduced an approach to directly construct a boundary-preserving parametrization
of the mesh motion from the analytical representation of each boundary;
however, the approach is cumbersome when multiple boundaries intersect
(e.g., at edges or corners). In the next section, we introduce a systematic
procedure to define $\phibold$ in the special case where all boundaries
are planar; future work will extend the approach to curved boundaries.
\begin{remark}
For domains with curved boundaries, in particular boundaries that
are non-polynomial functions or polynomials of degree greater than
$q$, $\Gcal_h(\Ecal_h;\xbm)$ will not partition $\Omega$ exactly
because the boundaries will not be perfectly represented by functions
in $[\Wcal_h]^d$ (they will be approximated by polynomials of degree $q$).
\end{remark}

\subsection{Enforcement of planar physical boundaries}
\label{sec:dommap:planar}
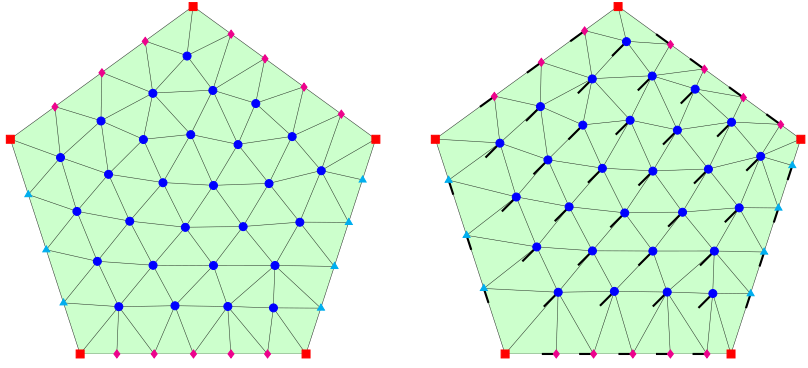
\begin{figure}
\centering
\begin{tikzpicture}
\begin{axis}[
axis equal image,
width=0.45\textwidth,
ymin=-0.849016994375,
ymax=1.04,
xmax=0.991056516295,
axis line style={draw=none},
xmin=-0.991056516295,
ticks=none]
\addplot []
graphics [xmin=-0.951056516295,xmax=0.951056516295,ymin=-0.809016994375,ymax=1.0] { 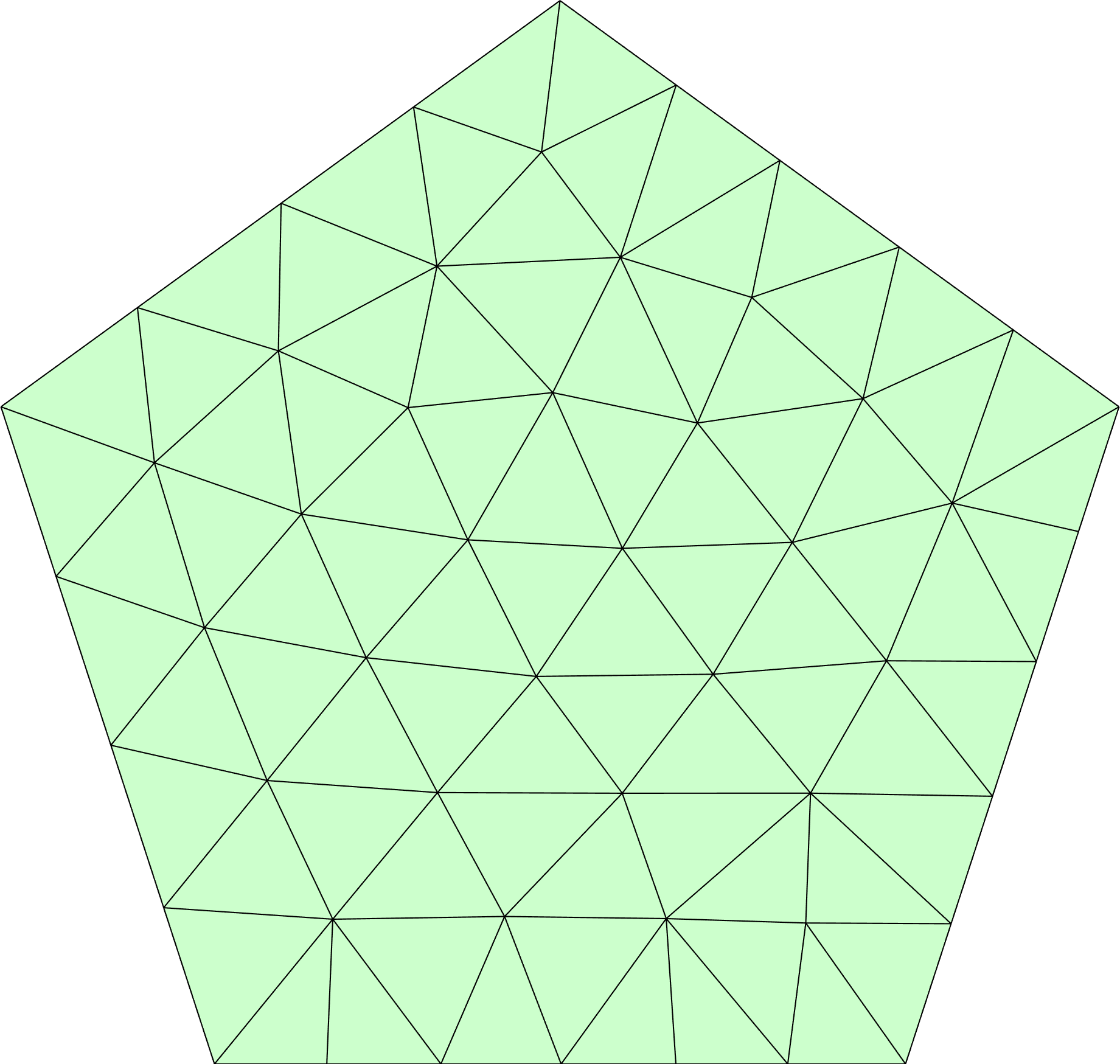};

\addplot [mark options={solid, thin}, mark=triangle*, mark size=1.5, cyan, only marks]
coordinates {
(-8.57312537e-01,  2.05026931e-02)
(-7.63974754e-01, -2.66761465e-01)
(-6.74174773e-01, -5.43137390e-01)
( 6.65324023e-01, -5.70377195e-01)
( 7.35771025e-01, -3.53563616e-01)
( 8.10243357e-01, -1.24361349e-01)
( 8.82108941e-01,  9.68181755e-02)};\label{line:nd_fix_x}

\addplot [mark options={solid, thin}, mark=diamond*, mark size=1.5, magenta, only marks]
coordinates {
(-7.18803077e-01,  4.77758995e-01)
(-4.74426180e-01,  6.55309204e-01)
(-3.96714897e-01, -8.09016994e-01)
(-2.49134478e-01,  8.18993207e-01)
(-2.02567007e-01, -8.09016994e-01)
( 1.91324110e-03, -8.09016994e-01)
( 1.97141037e-01, -8.09016994e-01)
( 1.97742290e-01,  8.56331817e-01)
( 3.74112067e-01,  7.28191673e-01)
( 3.87404805e-01, -8.09016994e-01)
( 5.77127422e-01,  5.80692383e-01)
( 7.71179643e-01,  4.39705193e-01)};\label{line:nd_fix_y}

\addplot [mark options={solid, thin}, mark=square*, mark size=1.5, red, only marks]
coordinates {
(-9.51056516e-01,  3.09016994e-01)
(-5.87785252e-01, -8.09016994e-01)
( 6.12323400e-17,  1.00000000e+00)
( 5.87785252e-01, -8.09016994e-01)
( 9.51056516e-01,  3.09016994e-01)};\label{line:nd_fix_all}

\addplot [mark options={solid, thin}, mark=*, mark size=1.5, blue, only marks]
coordinates {
(-6.89908017e-01,  2.13747076e-01)
(-6.04725624e-01, -6.67534301e-02)
(-4.98208203e-01, -3.26773006e-01)
(-4.79063635e-01,  4.04193294e-01)
(-4.40010111e-01,  1.26457090e-01)
(-3.86581444e-01, -5.62708248e-01)
(-3.29729311e-01, -1.17808715e-01)
(-2.58504538e-01,  3.07614211e-01)
(-2.09270854e-01,  5.48237631e-01)
(-2.08360699e-01, -3.47315837e-01)
(-1.56719591e-01,  8.26603085e-02)
(-9.42616161e-02, -5.58066966e-01)
(-4.04021774e-02, -1.49714397e-01)
(-3.13933914e-02,  7.42466753e-01)
(-1.23262891e-02,  3.33141353e-01)
( 1.02653672e-01,  5.63312279e-01)
( 1.05933018e-01, -3.48409636e-01)
( 1.06077204e-01,  6.82422614e-02)
( 1.80712893e-01, -5.61662480e-01)
( 2.33765057e-01,  2.81416115e-01)
( 2.60485362e-01, -1.45971916e-01)
( 3.26171860e-01,  4.95034299e-01)
( 3.95092892e-01,  7.83467552e-02)
( 4.18122088e-01, -5.69166177e-01)
( 4.26175800e-01, -3.48263155e-01)
( 5.15485366e-01,  3.22883872e-01)
( 5.55407562e-01, -1.23072447e-01)
( 6.67102293e-01,  1.45208561e-01)};\label{line:nd_fix_none}

\end{axis}
\end{tikzpicture} \quad
\begin{tikzpicture}
\begin{axis}[
axis equal image,
width=0.45\textwidth,
ymin=-0.849016994375,
ymax=1.04,
xmax=0.991056516295,
axis line style={draw=none},
xmin=-0.991056516295,
ticks=none]
\addplot []
graphics [xmin=-0.951056516295,xmax=0.951056516295,ymin=-0.809016994375,ymax=1.0] { 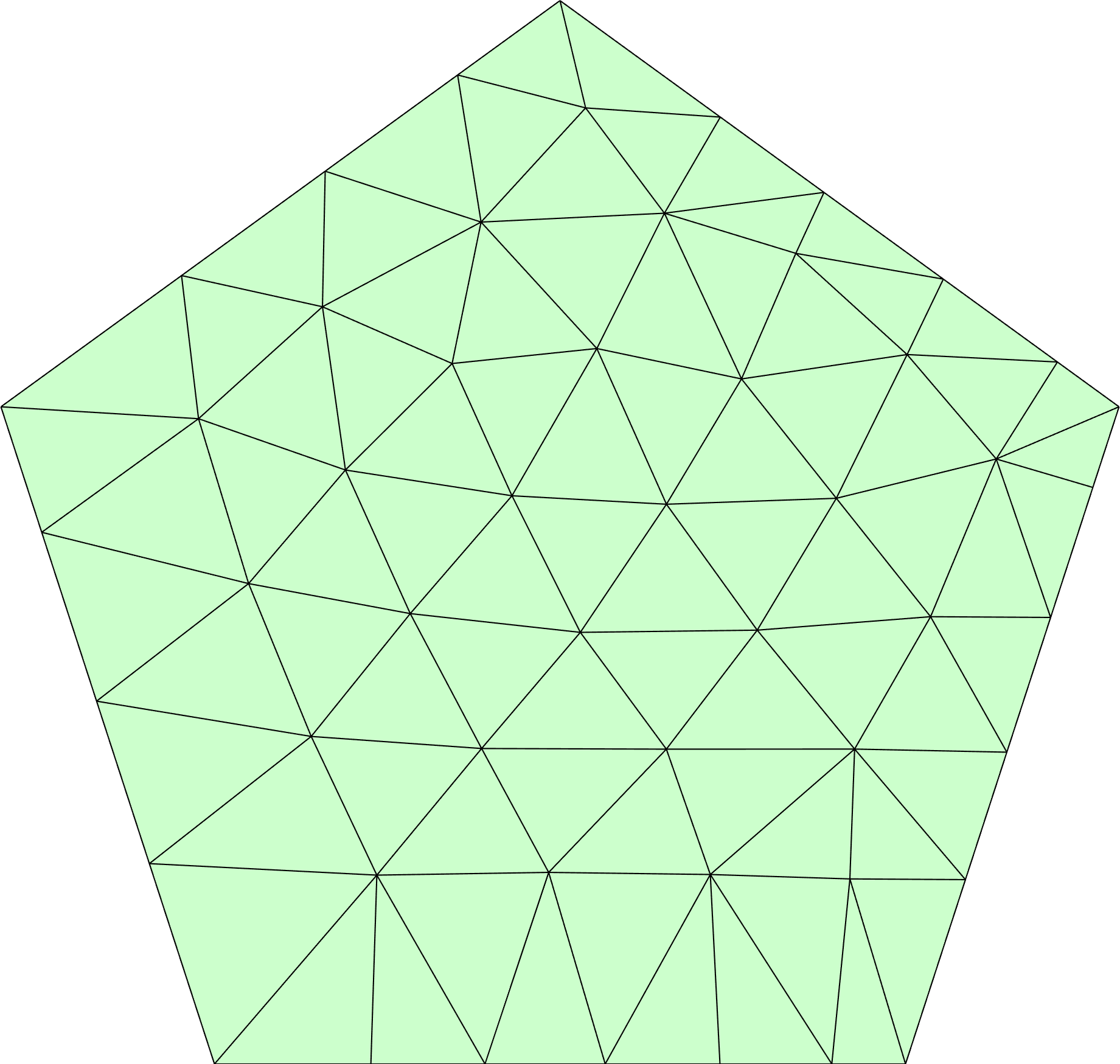};

\addplot [mark options={solid, thin}, mark=triangle*, mark size=1.5, cyan, only marks]
coordinates {
(-8.81681514e-01,  9.55026931e-02)
(-7.88343731e-01, -1.91761465e-01)
(-6.98543750e-01, -4.68137390e-01)
( 6.89693000e-01, -4.95377195e-01)
( 7.60140003e-01, -2.78563616e-01)
( 8.34612334e-01, -4.93613486e-02)
( 9.06477918e-01,  1.71818175e-01)};\label{line:nd_fix_x}

\addplot [mark options={solid, thin}, mark=diamond*, mark size=1.5, magenta, only marks]
coordinates {
(-6.43803077e-01,  5.32249685e-01)
(-3.99426180e-01,  7.09799894e-01)
(-3.21714897e-01, -8.09016994e-01)
(-1.74134478e-01,  8.73483896e-01)
(-1.27567007e-01, -8.09016994e-01)
( 7.69132411e-02, -8.09016994e-01)
( 2.72141037e-01, -8.09016994e-01)
( 2.72742290e-01,  8.01841127e-01)
( 4.49112067e-01,  6.73700983e-01)
( 4.62404805e-01, -8.09016994e-01)
( 6.52127422e-01,  5.26201694e-01)
( 8.46179643e-01,  3.85214503e-01)};\label{line:nd_fix_y}

\addplot [mark options={solid, thin}, mark=square*, mark size=1.5, red, only marks]
coordinates {
(-9.51056516e-01,  3.09016994e-01)
(-5.87785252e-01, -8.09016994e-01)
( 0.00000000e+00,  1.00000000e+00)
( 5.87785252e-01, -8.09016994e-01)
( 9.51056516e-01,  3.09016994e-01)};\label{line:nd_fix_all}

\addplot [mark options={solid, thin}, mark=*, mark size=1.5, blue, only marks]
coordinates {
(-6.14908017e-01,  2.88747076e-01)
(-5.29725624e-01,  8.24656991e-03)
(-4.23208203e-01, -2.51773006e-01)
(-4.04063635e-01,  4.79193294e-01)
(-3.65010111e-01,  2.01457090e-01)
(-3.11581444e-01, -4.87708248e-01)
(-2.54729311e-01, -4.28087149e-02)
(-1.83504538e-01,  3.82614211e-01)
(-1.34270854e-01,  6.23237631e-01)
(-1.33360699e-01, -2.72315837e-01)
(-8.17195912e-02,  1.57660308e-01)
(-1.92616161e-02, -4.83066966e-01)
( 3.45978226e-02, -7.47143965e-02)
( 4.36066086e-02,  8.17466753e-01)
( 6.26737109e-02,  4.08141353e-01)
( 1.77653672e-01,  6.38312279e-01)
( 1.80933018e-01, -2.73409636e-01)
( 1.81077204e-01,  1.43242261e-01)
( 2.55712893e-01, -4.86662480e-01)
( 3.08765057e-01,  3.56416115e-01)
( 3.35485362e-01, -7.09719162e-02)
( 4.01171860e-01,  5.70034299e-01)
( 4.70092892e-01,  1.53346755e-01)
( 4.93122088e-01, -4.94166177e-01)
( 5.01175800e-01, -2.73263155e-01)
( 5.90485366e-01,  3.97883872e-01)
( 6.30407562e-01, -4.80724474e-02)
( 7.42102293e-01,  2.20208561e-01)};\label{line:nd_fix_none}

\addplot [black, solid, thick]
coordinates {
(-9.51056516e-01,  3.09016994e-01)
(-9.51056516e-01,  3.09016994e-01)};\label{line:transl}

\addplot [black, solid, thick]
coordinates {
(-8.57312537e-01,  2.05026931e-02)
(-8.81681514e-01,  9.55026931e-02)};\label{line:transl}

\addplot [black, solid, thick]
coordinates {
(-7.63974754e-01, -2.66761465e-01)
(-7.88343731e-01, -1.91761465e-01)};\label{line:transl}

\addplot [black, solid, thick]
coordinates {
(-7.18803077e-01,  4.77758995e-01)
(-6.43803077e-01,  5.32249685e-01)};\label{line:transl}

\addplot [black, solid, thick]
coordinates {
(-6.89908017e-01,  2.13747076e-01)
(-6.14908017e-01,  2.88747076e-01)};\label{line:transl}

\addplot [black, solid, thick]
coordinates {
(-6.74174773e-01, -5.43137390e-01)
(-6.98543750e-01, -4.68137390e-01)};\label{line:transl}

\addplot [black, solid, thick]
coordinates {
(-6.04725624e-01, -6.67534301e-02)
(-5.29725624e-01,  8.24656991e-03)};\label{line:transl}

\addplot [black, solid, thick]
coordinates {
(-5.87785252e-01, -8.09016994e-01)
(-5.87785252e-01, -8.09016994e-01)};\label{line:transl}

\addplot [black, solid, thick]
coordinates {
(-4.98208203e-01, -3.26773006e-01)
(-4.23208203e-01, -2.51773006e-01)};\label{line:transl}

\addplot [black, solid, thick]
coordinates {
(-4.79063635e-01,  4.04193294e-01)
(-4.04063635e-01,  4.79193294e-01)};\label{line:transl}

\addplot [black, solid, thick]
coordinates {
(-4.74426180e-01,  6.55309204e-01)
(-3.99426180e-01,  7.09799894e-01)};\label{line:transl}

\addplot [black, solid, thick]
coordinates {
(-4.40010111e-01,  1.26457090e-01)
(-3.65010111e-01,  2.01457090e-01)};\label{line:transl}

\addplot [black, solid, thick]
coordinates {
(-3.96714897e-01, -8.09016994e-01)
(-3.21714897e-01, -8.09016994e-01)};\label{line:transl}

\addplot [black, solid, thick]
coordinates {
(-3.86581444e-01, -5.62708248e-01)
(-3.11581444e-01, -4.87708248e-01)};\label{line:transl}

\addplot [black, solid, thick]
coordinates {
(-3.29729311e-01, -1.17808715e-01)
(-2.54729311e-01, -4.28087149e-02)};\label{line:transl}

\addplot [black, solid, thick]
coordinates {
(-2.58504538e-01,  3.07614211e-01)
(-1.83504538e-01,  3.82614211e-01)};\label{line:transl}

\addplot [black, solid, thick]
coordinates {
(-2.49134478e-01,  8.18993207e-01)
(-1.74134478e-01,  8.73483896e-01)};\label{line:transl}

\addplot [black, solid, thick]
coordinates {
(-2.09270854e-01,  5.48237631e-01)
(-1.34270854e-01,  6.23237631e-01)};\label{line:transl}

\addplot [black, solid, thick]
coordinates {
(-2.08360699e-01, -3.47315837e-01)
(-1.33360699e-01, -2.72315837e-01)};\label{line:transl}

\addplot [black, solid, thick]
coordinates {
(-2.02567007e-01, -8.09016994e-01)
(-1.27567007e-01, -8.09016994e-01)};\label{line:transl}

\addplot [black, solid, thick]
coordinates {
(-1.56719591e-01,  8.26603085e-02)
(-8.17195912e-02,  1.57660308e-01)};\label{line:transl}

\addplot [black, solid, thick]
coordinates {
(-9.42616161e-02, -5.58066966e-01)
(-1.92616161e-02, -4.83066966e-01)};\label{line:transl}

\addplot [black, solid, thick]
coordinates {
(-4.04021774e-02, -1.49714397e-01)
( 3.45978226e-02, -7.47143965e-02)};\label{line:transl}

\addplot [black, solid, thick]
coordinates {
(-3.13933914e-02,  7.42466753e-01)
( 4.36066086e-02,  8.17466753e-01)};\label{line:transl}

\addplot [black, solid, thick]
coordinates {
(-1.23262891e-02,  3.33141353e-01)
( 6.26737109e-02,  4.08141353e-01)};\label{line:transl}

\addplot [black, solid, thick]
coordinates {
( 6.12323400e-17,  1.00000000e+00)
( 0.00000000e+00,  1.00000000e+00)};\label{line:transl}

\addplot [black, solid, thick]
coordinates {
( 1.91324110e-03, -8.09016994e-01)
( 7.69132411e-02, -8.09016994e-01)};\label{line:transl}

\addplot [black, solid, thick]
coordinates {
( 1.02653672e-01,  5.63312279e-01)
( 1.77653672e-01,  6.38312279e-01)};\label{line:transl}

\addplot [black, solid, thick]
coordinates {
( 1.05933018e-01, -3.48409636e-01)
( 1.80933018e-01, -2.73409636e-01)};\label{line:transl}

\addplot [black, solid, thick]
coordinates {
( 1.06077204e-01,  6.82422614e-02)
( 1.81077204e-01,  1.43242261e-01)};\label{line:transl}

\addplot [black, solid, thick]
coordinates {
( 1.80712893e-01, -5.61662480e-01)
( 2.55712893e-01, -4.86662480e-01)};\label{line:transl}

\addplot [black, solid, thick]
coordinates {
( 1.97141037e-01, -8.09016994e-01)
( 2.72141037e-01, -8.09016994e-01)};\label{line:transl}

\addplot [black, solid, thick]
coordinates {
( 1.97742290e-01,  8.56331817e-01)
( 2.72742290e-01,  8.01841127e-01)};\label{line:transl}

\addplot [black, solid, thick]
coordinates {
( 2.33765057e-01,  2.81416115e-01)
( 3.08765057e-01,  3.56416115e-01)};\label{line:transl}

\addplot [black, solid, thick]
coordinates {
( 2.60485362e-01, -1.45971916e-01)
( 3.35485362e-01, -7.09719162e-02)};\label{line:transl}

\addplot [black, solid, thick]
coordinates {
( 3.26171860e-01,  4.95034299e-01)
( 4.01171860e-01,  5.70034299e-01)};\label{line:transl}

\addplot [black, solid, thick]
coordinates {
( 3.74112067e-01,  7.28191673e-01)
( 4.49112067e-01,  6.73700983e-01)};\label{line:transl}

\addplot [black, solid, thick]
coordinates {
( 3.87404805e-01, -8.09016994e-01)
( 4.62404805e-01, -8.09016994e-01)};\label{line:transl}

\addplot [black, solid, thick]
coordinates {
( 3.95092892e-01,  7.83467552e-02)
( 4.70092892e-01,  1.53346755e-01)};\label{line:transl}

\addplot [black, solid, thick]
coordinates {
( 4.18122088e-01, -5.69166177e-01)
( 4.93122088e-01, -4.94166177e-01)};\label{line:transl}

\addplot [black, solid, thick]
coordinates {
( 4.26175800e-01, -3.48263155e-01)
( 5.01175800e-01, -2.73263155e-01)};\label{line:transl}

\addplot [black, solid, thick]
coordinates {
( 5.15485366e-01,  3.22883872e-01)
( 5.90485366e-01,  3.97883872e-01)};\label{line:transl}

\addplot [black, solid, thick]
coordinates {
( 5.55407562e-01, -1.23072447e-01)
( 6.30407562e-01, -4.80724474e-02)};\label{line:transl}

\addplot [black, solid, thick]
coordinates {
( 5.77127422e-01,  5.80692383e-01)
( 6.52127422e-01,  5.26201694e-01)};\label{line:transl}

\addplot [black, solid, thick]
coordinates {
( 5.87785252e-01, -8.09016994e-01)
( 5.87785252e-01, -8.09016994e-01)};\label{line:transl}

\addplot [black, solid, thick]
coordinates {
( 6.65324023e-01, -5.70377195e-01)
( 6.89693000e-01, -4.95377195e-01)};\label{line:transl}

\addplot [black, solid, thick]
coordinates {
( 6.67102293e-01,  1.45208561e-01)
( 7.42102293e-01,  2.20208561e-01)};\label{line:transl}

\addplot [black, solid, thick]
coordinates {
( 7.35771025e-01, -3.53563616e-01)
( 7.60140003e-01, -2.78563616e-01)};\label{line:transl}

\addplot [black, solid, thick]
coordinates {
( 7.71179643e-01,  4.39705193e-01)
( 8.46179643e-01,  3.85214503e-01)};\label{line:transl}

\addplot [black, solid, thick]
coordinates {
( 8.10243357e-01, -1.24361349e-01)
( 8.34612334e-01, -4.93613486e-02)};\label{line:transl}

\addplot [black, solid, thick]
coordinates {
( 8.82108941e-01,  9.68181755e-02)
( 9.06477918e-01,  1.71818175e-01)};\label{line:transl}

\addplot [black, solid, thick]
coordinates {
( 9.51056516e-01,  3.09016994e-01)
( 9.51056516e-01,  3.09016994e-01)};\label{line:transl}

\end{axis}
\end{tikzpicture}
\caption{Reference domain $\Omega_0$ (regular pentagon) and mesh
$\Ecal_h$ with the reference nodes $\hat{X}_I$ indicated with markers
(\textit{left}) and physical domain $\Omega$ and mapped mesh
obtained by boundary-preserving mapping $\Gcal_h(\,\cdot\,;\phibold(\ybm))$ 
with the physical nodes $\hat{x}_I$ indicated with markers (\textit{right});
the node movement from its original position is indicated by
(\ref{line:transl}). The nodes are partitioned into four
categories based on the degree of freedom partitioning in
(\ref{eqn:partdof0})-(\ref{eqn:partdof1}): 
$\Ical_I^\mathrm{u}=\{1,2\}$ ($\hat{x}_I$ moves freely) (\ref{line:nd_fix_none}),
$\Ical_I^\mathrm{u}=\{1\}$ ($(\hat{x}_I)_1$ moves freely, $(\hat{x}_I)_2$ constrained) (\ref{line:nd_fix_y}),
$\Ical_I^\mathrm{u}=\{2\}$ ($(\hat{x}_I)_2$ moves freely, $(\hat{x}_I)_1$ constrained) (\ref{line:nd_fix_x}),
$\Ical_I^\mathrm{u}=\emptyset$ ($\hat{x}_I$ fixed) (\ref{line:nd_fix_all}).}
\label{fig:demo_mshmot}
\end{figure}

In the special case where all boundaries are planar, we
define an affine parametrization of the physical nodes
\begin{equation} \label{eqn:dommap_affine}
\phibold : \ybm \mapsto \Abm \ybm + \bbm,
\end{equation}
where $\Abm\in\Rbb^{N_\xbm\times N_\ybm}$ and $\bbm\in\Rbb^{N_\xbm}$
are given analytically in terms of the surface normal vectors
and the reference nodes. Because we take the reference and physical
domains to be the same set (Remark~\ref{rem:refphys}), the
reference-to-physical mapping can be interpreted as moving
nodes throughout the domain. Suppose the
boundary of $\partial\Omega$ is the union of $N_b$ planar surfaces
$\{\partial\Omega_i\}_{i=1}^{N_b}$, i.e.,
$\partial\Omega = \bigcup_{i=1}^{N_b} \partial\Omega_i$,
that intersect to form non-smooth features, e.g., edges and corners.
For each planar surface $\partial\Omega_i$,
we denote its unit normal
(oriented outward with respect to $\Omega$) as $\eta_i\in\Rbb^d$, which
implies any two points $x,\bar{x}\in\partial\Omega_i$ must satisfy
$\eta_i^T(x-\bar{x})=0$. We assume there are no redundancies
in the specification of these $N_b$ boundaries, i.e., the normal vectors
of all boundaries passing through a given point are linearly independent.
We will require $\hat{x}_I$ to lie
on all boundaries on which $\hat{X}_I$ lies, i.e., $\hat{x}_I$ can be
obtained by sliding $\hat{X}_I$ along the boundaries on which it lies.
This condition is sufficient to ensure the mapped mesh
$\Gcal_h(\Ecal_h;\xbm)$ conforms to all boundaries
$\{\partial\Omega_i\}_{i=1}^{N_b}$, provided none of the
elements in $\Gcal_h(\Ecal_h;\xbm)$ are inverted, because
the same element faces/edges used to represent each boundary/intersection
in the reference mesh $\Ecal_h$ will be used to represent the
same boundary/intersection in the mapped mesh $\Gcal_h(\Ecal_h;\xbm)$
(Figure~\ref{fig:demo_mshmot}).

Let $\Bcal_I\subset\{1,\dots,N_b\}$ ($N_I^\mathrm{c}\coloneqq|\Bcal_I|$)
be the set of boundaries that pass through $\hat{X}_I$
and $B_I\in\Rbb^{N_I^\mathrm{c}\times d}$ be the corresponding collection
of boundary normals
\begin{equation}
 \Bcal_I = \{ i\in\{1,\dots,N_b\} \mid \hat{X}_I \in \partial\Omega_i\},
 \quad
B_I = \begin{bmatrix}
        \eta_{(\Bcal_I)_1} &
        \cdots &
        \eta_{(\Bcal_I)_{N_I^\mathrm{c}}}
      \end{bmatrix}^T.
\end{equation}
We implicitly assume the set $\Bcal_I$ is ordered with $(\Bcal_I)_i$
denoting its $i$th entry, although the specific ordering of the
entries is unimportant, and define $\Bcal_I$ such that its entries
are unique, i.e., $(\Bcal_I)_i=(\Bcal_I)_j$ implies $i=j$.
Owing to the assumption that the normal vectors of all boundaries
passing through a given point are linearly independent, we have
$N_I^\mathrm{c}\leq d$ and $B_I$ is full rank. In order for the
corresponding physical node $\hat{x}_I$ to remain on all the boundaries
in $\Bcal_I$, it must satisfy
\begin{equation} \label{eqn:domdef_lincon_nodal}
 B_I (\hat{x}_I-\hat{X}_I) = 0
\end{equation}
which defines $N_I^\mathrm{c}$ unique constraints.
We expand the coordinate of the physical node as
\begin{equation} \label{eqn:physnode_expand}
 \hat{x}_I = Y_I \hat{x}_I^\mathrm{c} + Z_I \hat{x}_I^\mathrm{u},
\end{equation}
where $Y_I\in\Rbb^{d \times N_I^\mathrm{c}}$ is an orthogonal matrix such
that $B_I Y_I$ is non-singular, $Z_I\in\Rbb^{d\times(d-N_I^\mathrm{c})}$
is an orthogonal matrix such that
$\mathrm{Ran}(Y_I)\oplus\mathrm{Ran}(Z_I) = \Rbb^d$
($\mathrm{Ran}(A)$ denotes the range of the matrix $A$),
and $\hat{x}_I^\mathrm{c}\in\Rbb^{N_I^\mathrm{c}}$ and
$\hat{x}_I^\mathrm{u}\in\Rbb^{d-N_I^\mathrm{c}}$ are the
coordinates of $\hat{x}_I$ in the respective basis
that will constitute the constrained and unconstrained
degrees of freedom of $\hat{x}_I$, respectively.
Substituting (\ref{eqn:physnode_expand}) into the constraints
(\ref{eqn:domdef_lincon_nodal}) leads to
the following relationship between $\hat{x}_I^\mathrm{c}$
and $\hat{x}_I^\mathrm{u}$ that must hold for $\hat{x}_I$
to remain on the same boundaries as $\hat{X}_I$
\begin{equation} \label{eqn:physnode_constr_dof}
 \hat{x}_I^\mathrm{c} =
 \left[B_I Y_I\right]^{-1}B_I\left(\hat{X}_I-Z_I\hat{x}_I^\mathrm{u}\right).
\end{equation}
This leads to an affine relationship between the physical node
coordinate and its unconstrained degrees of freedom
\begin{equation}
 \hat{x}_I = A_I \hat{x}_I^\mathrm{u} + b_I, \quad
 A_I \coloneqq Z_I - Y_I\left[B_IY_I\right]^{-1}B_IZ_I, \quad
 b_I \coloneqq Y_I\left[B_IY_I\right]^{-1}B_I\hat{X}_I,
\end{equation}
which we concatenate over all nodes $I=1,\dots,N_\mathrm{v}$ to define
the affine parametrization in (\ref{eqn:dommap_affine}) with
\begin{equation}
 \Abm = \begin{bmatrix} A_1 & & \\ & \ddots & \\ & & A_{N_\mathrm{v}} \end{bmatrix},
 \quad
 \ybm = \begin{bmatrix} \hat{x}_1^\mathrm{u} \\ \vdots \\ \hat{x}_{N_\mathrm{v}}^\mathrm{u} \end{bmatrix},
 \quad
 \bbm = \begin{bmatrix} b_1 \\ \vdots \\ b_{N_\mathrm{v}} \end{bmatrix}
\end{equation}
and $N_\ybm = N_\xbm-\sum_{I=1}^{N_\mathrm{v}} N_I^\mathrm{c}$.

We choose $Y_I$ and $Z_I$ to be non-overlapping subsets of the columns
of the $d\times d$ identity matrix
\begin{equation} \label{eqn:basis0}
 Y_I =
 \begin{bmatrix}
  e_{(\Ical_I^\mathrm{c})_1} & \cdots & e_{(\Ical_I^\mathrm{c})_{N_I^\mathrm{c}}}
 \end{bmatrix}, \quad
 Z_I =
 \begin{bmatrix}
  e_{(\Ical_I^\mathrm{u})_1} & \cdots & e_{(\Ical_I^\mathrm{u})_{d-N_I^\mathrm{c}}}
 \end{bmatrix},
\end{equation}
where $e_i\in\Rbb^d$ is the $i$th canonical unit vector,
$\Ical_I^\mathrm{u}\subset\{1,\dots,d\}$
($|\Ical_I^\mathrm{u}|=d-N_I^\mathrm{c}$)
denotes the unconstrained degrees of freedom, and
$\Ical_I^\mathrm{c}=\{1,\dots,d\}\setminus\Ical_I^\mathrm{u}$
denotes the constrained degrees of freedom; similar to $\Bcal_I$,
$\Ical_I^\mathrm{u}$ and $\Ical_I^\mathrm{c}$ are
ordered sets with unique elements. This means
$\hat{x}_I^\mathrm{u}$ and $\hat{x}_I^\mathrm{c}$ are subsets
of the entries of $\hat{x}_I$
\begin{equation} \label{eqn:partdof0}
 (\hat{x}_I^\mathrm{u})_i = (\hat{x}_I)_{(\Ical_I^\mathrm{u})_i}, \quad
 (\hat{x}_I^\mathrm{c})_j = (\hat{x}_I)_{(\Ical_I^\mathrm{c})_j}
\end{equation}
for $i=1,\dots,d-N_I^\mathrm{c}$ and $j=1,\dots,N_I^\mathrm{c}$,
i.e., $\hat{x}_I^\mathrm{u}$ are the components of $\hat{x}_I$
allowed to move freely and $\hat{x}_I^\mathrm{c}$ are the
components uniquely determined to ensure $\hat{x}_I$ lies
on the appropriate boundaries. 
The decomposition of $\{1,\dots,d\}$ into unconstrained and
constrained degrees of freedom is not unique; the only condition
is $B_IY_I$ be non-singular because
the $\mathrm{Ran}(Y_I)\oplus\mathrm{Ran}(Z_I)=\Rbb^d$ condition is
guaranteed by construction in (\ref{eqn:basis0}). We choose the unconstrained
degrees of freedom $\hat{x}_I^\mathrm{u}$ to correspond to the coordinate
directions closest to the null space of $B_I$, denoted $\mathrm{Null}(B_I)$, because
this will require the least action on the constrained degrees of freedom.
In the extreme case where $d-N_I^\mathrm{c}$ coordinate directions lie in
$\mathrm{Null}(B_I)$, this choice implies the columns of $Z_I$ will form
a basis for $\mathrm{Null}(B_I)$, which implies $B_I Z_I = 0$ and the
expression for $\hat{x}_I^\mathrm{c}$ in (\ref{eqn:physnode_constr_dof})
is independent of $\hat{x}_I^\mathrm{u}$. 
To this end, let $V_I \in \Rbb^{d\times (d-N_I^\mathrm{c})}$ be a
matrix whose columns form an orthonormal basis for $\mathrm{Null}(B_I)$,
then the set $\Ical_I^\mathrm{u}$ is defined as
\begin{equation} \label{eqn:partdof1}
 (\Ical_I^\mathrm{u})_i \coloneqq
 \argmin_{j\in\Jcal_I^{(i)}} \norm{e_j-V_I V_I^T e_j}, \quad
 \Jcal_I^{(i)} = \{1,\dots,d\}\setminus\{(\Ical_I^\mathrm{u})_1,\dots,(\Ical_I^\mathrm{u})_{i-1}\}
\end{equation}
for $i=1,\dots,d-N_I^\mathrm{c}$ and
$\Ical_I^\mathrm{c}=\{1,\dots,d\}\setminus\Ical_I^\mathrm{u}$.
The partitioning of the degrees of freedom of each node in a 
mesh of a regular pentagon using this procedure is illustrated
in Figure~\ref{fig:demo_mshmot}.

\begin{remark}
\label{rem:dommap:null}
There are other many valid choices for $Y_I$ and $Z_I$ that effectively
decompose $\hat{x}_I$ into constrained and unconstrained parts.
One choice that mimics approaches for enforcing linear equality
constraints in optimization solvers \cite{gill_practical_1981}
chooses the columns of $Z_I$ to be an orthonormal basis for
$\mathrm{Null}(B_I)$ and $Y_I$ to be an orthonormal basis for the range
space of $B_I^T$. This guarantees that $B_I Y_I$ is non-singular
and $B_I Z_I = 0$, which simplifies the expressions in
(\ref{eqn:physnode_constr_dof}) to $A_I = Z_I$ and
$\hat{x}_I^\mathrm{c} = [B_IY_I]^{-1}B_I\hat{X}_I$
(independent of $\hat{x}_I^\mathrm{u}$). While this approach is
clean and elegant, the constrained and unconstrained parts of
$\hat{x}_I$ do not have an intuitive physical interpretation
as components of $\hat{x}_I$.
\end{remark}

\section{High-order implicit shock tracking formulation}
\label{sec:ist-form}
The high-order implicit shock tracking method introduced in
\cite{zahr_optimization-based_2018,zahr_implicit_2020}
simultaneously computes the discrete solution of the conservation law
and the nodal coordinates of the mesh that causes element
faces to align with discontinuities. This is achieved through a fully
discrete, full space PDE-constrained optimization formulation with the
optimization variables taken to be the discrete flow solution and nodal
coordinates of the mesh. 
With the boundary-preserving parametrization of the mesh motion in
(\ref{eqn:mapparam}), the HOIST method is formulated as
\begin{equation} \label{eqn:pde-opt}
 (\ubm^\star,\ybm^\star) \coloneqq
 \argmin_{\ubm\in\Rbb^{N_\ubm},\ybm\in\Rbb^{N_\ybm}} f(\ubm,\ybm) \quad \text{subject to:} \quad \rbm(\ubm,\phibold(\ybm)) = \zerobold,
\end{equation}
where $\func{f}{\Rbb^{N_\ubm}\times\Rbb^{N_\ybm}}{\Rbb}$ is the objective
function and the nodal coordinates of the aligned mesh are
$\xbm^\star = \phibold(\ybm^\star)$. The objective function is composed
of two terms as
\begin{equation}\label{eqn:obj0}
 f : (\ubm,\ybm) \mapsto f_\text{err}(\ubm,\ybm) + \kappa^2 f_\text{msh}(\ybm),
\end{equation}
which balances alignment of the mesh with non-smooth features and the quality
of the elements and $\kappa\in\Rbb_{\geq 0}$ is the mesh penalty parameter.
The mesh alignment term,
$\func{f_\text{err}}{\Rbb^{N_\ubm}\times\Rbb^{N_\ybm}}{\Rbb}$, is taken
to be the norm of the enriched DG residual
\begin{equation}\label{eqn:obj1}
 f_\text{err} : (\ubm,\ybm) \mapsto \frac{1}{2}\norm{\Rbm(\ubm,\phibold(\ybm))}_2^2.
\end{equation}
To ensure elements of the
discontinuity-aligned mesh are high-quality, we define the mesh distortion
term, $\func{f_\text{msh}}{\Rbb^{N_\ybm}}{\Rbb}$, as
\begin{equation}\label{eqn:obj2}
 f_\text{msh} : \ybm \mapsto \frac{1}{2}\norm{\Rbm_\text{msh}(\phibold(\ybm))}_2^2,
\end{equation}
where $\func{\Rbm_\text{msh}}{\Rbb^{N_\ybm}}{\Rbb^{|\Ecal_h|}}$ is the
element-wise mesh distortion with respect to an ideal element
\cite{zahr_implicit_2020,knupp_algebraic_2001,roca_defining_2012}.
Concretely, assume the elements of $\Ecal_h$
are ordered and let $\Omega_{0,e}\in\Ecal_h$ ($\Omega_{0,e}\subset\Omega_0$)
for $e=1,\dots,|\Ecal_h|$ denote the $e$th element. Furthermore, define
the corresponding physical element as $\Omega_e(\xbm)=\Gcal_h(\Omega_{0,e};\xbm)$
for a given instance of the physical nodes $\xbm\in\Rbb^{N_\xbm}$ and
let $K_{\star,e}\subset\Rbb^d$ denote the ideal shape of element $e$.
Then, the element-wise mesh distortion is defined as
\begin{equation}
 \left[\Rbm_\text{msh}(\xbm)\right]_e =
 \frac{1}{|K_{\star,e}|}
 \int_{K_{\star,e}} \left(\frac{\norm{\nabla\Fcal_{\star,e}(\,\cdot\,;\xbm)}_F^2}{d \det(\nabla\Fcal_{\star,e}(\,\cdot\,;\xbm))_+^{2/d}}\right)^2\, dV,
\end{equation}
where $\func{\Fcal_{\star,e}(\,\cdot\,;\xbm)}{K_{\star,e}}{\Omega_e(\xbm)}$ is a
diffeomorphism between the ideal and physical element. We consider two choices for
the ideal element ($K_{\star,e}$), which lead to a mesh distortion term with
different properties: 1) $\Omega_{0,e}$ (the reference element) and
2) $K_\star\subset\Rbb^d$ (the regular $d$-simplex). The first choice
implies $\Fcal_{\star,e}(\,\cdot\,;\xbm) = \left.\Gcal_h(\,\cdot\,;\xbm)\right|_{\Omega_{0,e}}$
and the resulting mesh distortion---identical to the one used in our previous
work \cite{zahr_implicit_2020}---penalizes deviation of the physical mesh from
the reference mesh, whereas the second choice drives elements toward the
regular $d$-simplex.

To obtain the first-order optimality system of the implicit shock
tracking formulation in (\ref{eqn:pde-opt}), we introduce the corresponding
Lagrangian,
$\func{\Lcal}{\Rbb^{N_\ubm}\times\Rbb^{N_\ybm}\times\Rbb^{N_\ubm}}{\Rbb}$,
defined as
\begin{equation}
 \Lcal : (\ubm,\ybm,\lambdabold) \mapsto f(\ubm,\ybm)-\lambdabold^T\rbm(\ubm,\phibold(\ybm)).
\end{equation}
Then, the first-order optimality, or Karush-Kuhn-Tucker (KKT), conditions state
that $(\ubm^\star,\ybm^\star)\in\Rbb^{N_\ubm}\times\Rbb^{N_\ybm}$ is a
first-order solution of the optimization problem in (\ref{eqn:pde-opt}) if there
exists $\lambdabold^\star\in\Rbb^{N_\ubm}$ such that the Lagrangian is
stationary, which leads to the following conditions
\begin{equation}\label{eqn:kkt0}
\begin{aligned}
 \pder{f}{\ubm}(\ubm^\star,\ybm^\star)^T-\pder{\rbm}{\ubm}(\ubm^\star,\phibold(\ybm^\star))^T\lambdabold^\star &= \zerobold \\
 \pder{f}{\ybm}(\ubm^\star,\ybm^\star)^T-\pder{\phibold}{\ybm}(\ybm^\star)^T\pder{\rbm}{\xbm}(\ubm^\star,\phibold(\ybm^\star))^T\lambdabold^\star &= \zerobold \\
 \rbm(\ubm^\star,\phibold(\ybm^\star)) &= \zerobold.
\end{aligned}
\end{equation}
Following \cite{zahr_implicit_2020}, we define the Lagrange multiplier estimate,
$\func{\hat\lambdabold}{\Rbb^{N_\ubm}\times\Rbb^{N_\ybm}}{\Rbb^{N_\ubm}}$, as
\begin{equation} \label{eqn:lagmul-est}
 \hat\lambdabold : (\ubm, \ybm) \mapsto
 \left[\pder{\rbm}{\ubm}(\ubm,\phibold(\ybm))\right]^{-T}
 \pder{f}{\ubm}(\ubm,\ybm)^T,
\end{equation}
which ensures the first equation in (\ref{eqn:kkt0}) (adjoint equation) is
satisfied for any $\ubm\in\Rbb^{N_\ubm}$ and $\ybm\in\Rbb^{N_\ybm}$.
With this estimate of the Lagrange multipliers, the optimality conditions
reduce to
\begin{equation}
 \cbm(\ubm^\star,\ybm^\star) = \zerobold, \qquad
 \rbm(\ubm^\star,\phibold(\ybm^\star)) = \zerobold,
\end{equation}
where $\func{\cbm}{\Rbb^{N_\ubm}\times\Rbb^{N_\ybm}}{\Rbb^{N_\ubm}}$
is defined as
\begin{equation}
 \cbm : (\ubm, \ybm) \mapsto
 \pder{f}{\ybm}(\ubm,\ybm)^T-\pder{\phibold}{\ybm}(\ybm)^T\pder{\rbm}{\xbm}(\ubm,\phibold(\ybm))^T\hat\lambdabold(\ubm,\ybm).
\end{equation}

\begin{remark}
Intuitively, the enriched DG residual is an improved approximation of the
infinite-dimensional weak formulation of the conservation law than the
standard DG residual and, as a result, is sensitive to non-physical
oscillations that result from approximating discontinuities on non-aligned
meshes. Thus, minimizing this term promotes alignment of the mesh with
discontinuities, as demonstrated in \cite{zahr_implicit_2020}.
\end{remark}

\begin{remark}
Because the DG equations ($\rbm$) with standard upwind numerical
flux are enforced as constraints, the properties of the DG method are directly 
inherited by the HOIST method independent of the objective function. Therefore,
there is flexibility in the choice of the numerical flux function for the enriched
DG residual ($\Rbm$). In this work, we use a central flux to define the enriched
DG residual because it is smooth and conservative \cite{zahr_implicit_2020} and
sufficient to reliably track shocks (Section~\ref{sec:numexp}).
\end{remark}

\begin{remark}
The reference element is a reasonable choice for the ideal element
if a high-quality reference mesh $\Ecal_h$ is used; however, the quality
of the reference elements degrade as elements are collapsed during the
tracking iterations (\cite{zahr_implicit_2020},
Section~\ref{sec:ist_solve:mod:clps})---particularly prevalent
in problems with complex discontinuity surfaces---which has the adverse
effect of the mesh distortion term favoring some poor-quality elements.
On the other hand, taking the regular $d$-simplex to be the ideal element
for all $e=1,\dots,|\Ecal_h|$ ensures high-quality elements are always
favored by the mesh distortion term; however, it can also lead to substantial
movement of the physical mesh away from the reference mesh even in regions
away from the shock. In practice, we usually use the mesh distortion term
based on the $d$-regular simplex ideal element, unless the large mesh motion
causes many unwanted element collapses.
\end{remark}

\begin{remark}
In the case where the ideal element is the regular $d$-simplex, i.e.,
$K_{\star,e}=K_\star$ for $e=1,\dots,|\Ecal_h|$, the mapping $\Fcal_{\star,e}$
is described in detail in \cite{roca_defining_2012}; we briefly outline
its construction. Let $\{\psi_i\}_{i=1}^n$ be a nodal basis of $\Pcal_q(K_\star)$
associated with well-spaced nodes $\{\hat\xi_i\}_{i=1}^n\subset K_\star$, where
$n=\dim\Pcal_q(K_\star)={{q+d}\choose{d}}$. Then, 
$\Fcal_{\star,e}:\xi\mapsto\sum_{i=1}^n \hat{x}_i^e \psi_i(\xi)$,
where $\{\hat{x}_i^e\}_{i=1}^n\subset\{\hat{x}_I\}_{I=1}^{N_\mathrm{v}}$
are the physical nodes associated with element $e$.
\end{remark}

\begin{remark}
The objective function in (\ref{eqn:obj0})-(\ref{eqn:obj2}) can be written
as the norm of a residual function,
$\func{\Fbm}{\Rbb^{N_\ubm}\times\Rbb^{N_\ybm}}{\Rbb^{N_\ubm'+|\Ecal_h|}}$, as
\begin{equation}
 f(\ubm,\ybm) = \frac{1}{2} \Fbm(\ubm,\ybm)^T\Fbm(\ubm,\ybm), \qquad
 \Fbm : (\ubm,\ybm) \mapsto \begin{bmatrix}\Rbm(\ubm,\phibold(\ybm)) \\ \kappa \Rbm_\text{msh}(\phibold(\ybm)) \end{bmatrix}.
\end{equation}
The implication is the HOIST optimization problem in (\ref{eqn:pde-opt})
is a constrained residual minimization problem, which can be exploited
in the development of the solver; see \cite{zahr_implicit_2020} and
Section~\ref{sec:ist_solve:sqp:lmhess}.
\end{remark}

\section{Robust implicit shock tracking solver}
\label{sec:ist_solve}
In this section, we introduce an iterative solver for the optimization
problem in (\ref{eqn:pde-opt}) based on the SQP method in
\cite{zahr_implicit_2020} with improved robustness and capability.
Because the DG system cannot be solved robustly on non-aligned meshes without
additional measures to suppress numerical oscillations, nested approaches
to solve the discrete PDE-constrained optimization problem in
(\ref{eqn:pde-opt}) are not appropriate; instead, we use a full-space
approach that aims to simultaneously converge the DG solution and the
mesh to their optimal values simultaneously. To this end, we define a
new variable $\zbm\in\Rbb^{N_\zbm}$ ($N_\zbm=N_\ubm+N_\ybm$) that combines
the DG solution $\ubm$ and unconstrained mesh coordinates $\ybm$ as
\begin{equation}
 \zbm = (\ubm, \ybm),
\end{equation}
and use $\zbm$ interchangeably with the tuple $(\ubm,\ybm)$.
For brevity, we introduce the following notation for the
derivatives of the objective function,
$\func{\gbm}{\Rbb^{N_\zbm}}{\Rbb^{N_\zbm}}$, and the DG residual,
$\func{\Jbm}{\Rbb^{N_\zbm}}{\Rbb^{N_\ubm}\times\Rbb^{N_\zbm}}$, as
\begin{equation}
 \gbm : \zbm \mapsto
 \begin{bmatrix}
   \ds{\pder{f}{\ubm}(\ubm,\ybm)^T} \\[1em]
   \ds{\pder{f}{\ybm}(\ubm,\ybm)^T}
 \end{bmatrix}, \quad
 \Jbm : \zbm \mapsto
 \begin{bmatrix}
  \ds{\pder{\rbm}{\ubm}(\ubm,\phibold(\ybm))} &
  \ds{\pder{\rbm}{\xbm}(\ubm,\phibold(\ybm))\pder{\phibold}{\ybm}(\ybm)}
 \end{bmatrix}.
\end{equation}

\subsection{Sequential quadratic programming method}
\label{sec:ist_solve:sqp}
The SQP method in \cite{zahr_implicit_2020} produces a sequence of iterates
$\{\zbm_k\}_{k=0}^\infty$ such that
$\zbm_k=(\ubm_k,\ybm_k)\rightarrow\zbm^\star=(\ubm^\star,\ybm^\star)$,
where $(\ubm^\star,\ybm^\star)$ satisfies the first-order
optimality conditions in (\ref{eqn:kkt0}). The sequence of iterates
is generated as
\begin{equation} \label{eqn:sqp_step0}
 \zbm_{k+1} = \zbm_k + \alpha_k\Delta \zbm_k,
\end{equation}
where the search direction $\Delta\zbm_k\in\Rbb^{N_\zbm}$
is computed as the solution of the following quadratic program
\begin{equation}\label{eqn:qp0}
\optconOne{\Delta\zbm \in \Rbb^{N_\zbm}}
{\gbm_k^T\Delta\zbm+\frac{1}{2}\Delta\zbm^T\Bbm_k\Delta\zbm}
{\rbm_k+\Jbm_k\Delta\zbm = \zerobold},
\end{equation}
$\gbm_k\in\Rbb^{N_\zbm}$, $\rbm_k\in\Rbb^{N_\ubm}$, and
$\Jbm_k\in\Rbb^{N_\ubm\times N_\zbm}$ are the objective gradient,
residual, and residual Jacobian, respectively, evaluated at $\zbm_k$
\begin{equation}
 \rbm_k \coloneqq \rbm(\ubm_k,\phibold(\ybm_k)), \qquad
 \gbm_k \coloneqq \gbm(\zbm_k), \qquad
 \Jbm_k \coloneqq \Jbm(\zbm_k),
\end{equation}
$\Bbm_k\in\Rbb^{N_\zbm \times N_\zbm}$ is a symmetric positive
definite (SPD) approximation to the Hessian of the Lagrangian at $\zbm_k$,
and $\alpha_k\in\Rbb_{>0}$ is the step length. The first-order
optimality conditions of the quadratic program leads to the linear system
of equations that must be solved at each iteration $k$
\begin{equation} \label{eqn:sqp_sys0}
 \begin{bmatrix}
  \Bbm_k & \Jbm_k^T \\
  \Jbm_k & \zerobold
 \end{bmatrix}
 \begin{bmatrix}
  \Delta\zbm_k \\ \etabold_k
 \end{bmatrix}
 =
 -\begin{bmatrix}
   \gbm_k \\ \rbm_k
  \end{bmatrix},
\end{equation}
where $\etabold_k\in\Rbb^{N_\ubm}$ are the Lagrange multipliers associated
with the linearized constraint in (\ref{eqn:qp0}).

\subsubsection{Line search globalization}
\label{sec:ist_solve:sqp:lsrch}
To ensure the sequence converges to a first-order critical point of
(\ref{eqn:pde-opt}) from any initial guess, the step length, $\alpha_k$,
is computed to guarantee sufficient decrease of a merit function,
$\func{\varphi_k}{\Rbb}{\Rbb}$, that combines the objective function
and constraint violation. In this work, $\alpha_k$ is computed using
backtracking \cite{nocedal2006numerical} such that
\begin{equation} \label{eqn:steep_descent}
 \varphi_k(\alpha_k) \leq \varphi_k(0) + c\alpha_k\varphi_k'(0),
\end{equation}
where $c\in(0,1)$. We use the $\ell_1$ merit function, defined as
\begin{equation}
\varphi_k : \alpha \mapsto f(\zbm+\alpha\Delta\zbm_k) + \mu_k \norm{\rbm(\zbm_k+\alpha\Delta\zbm_k)}_1;
\end{equation}
with the penalty parameter, $\mu_k\in\Rbb_{>0}$, taken as \cite{nocedal2006numerical}
\begin{equation} \label{eqn:pnlty0}
 \mu_0 = 0, \qquad
 \mu_k = \min\{\max\{\varpi \bar{\mu}_k,\mu_{k-1}\}, \mu_{\max}\}, \qquad
 \bar\mu_k = \frac{\gbm_k^T\Delta\zbm_k+(1/2)\Delta\zbm_k^T\Bbm_k\Delta\zbm_k}{(1-\rho)\norm{\rbm_k}_1},
\end{equation}
where $\varpi>1$ and $\rho\in(0,1)$ are parameters taken to be $\varpi=1.2$
and $\rho=0.95$ in this work.

\begin{remark} \label{rem:lsrch_pnlty}
The choice of penalty parameter in (\ref{eqn:pnlty0}) for the
$\ell_1$ merit function differs from the choice in
\cite{zahr_implicit_2020} that defines
\begin{equation}
 \mu_k = 2\norm{\hat\lambdabold(\zbm_k)}_\infty
\end{equation}
because this guarantees, under some assumptions, the SQP search direction
will be a direction of descent for the $\ell_1$ merit function
\cite{nocedal2006numerical}.
However, we observed this choice can lead to a scaling mismatch
between the objective and constraint causing many unnecessary
line search iterations and interfering with progress of the
iteration, which is consistent with other observations
\cite{nocedal2006numerical}.
\end{remark}

\subsubsection{Levenberg-Marquardt Hessian approximation}
\label{sec:ist_solve:sqp:lmhess}
We use the modified Levenberg-Marquardt Hessian approximation
introduced in \cite{corrigan_moving_2019, zahr_implicit_2020}
to define $\Bbm_k$. To this end, we expand $\Bbm_k$ as
\begin{equation} \label{eqn:hess0}
 \Bbm_k =
 \begin{bmatrix}
  \Bbm_k^{\ubm\ubm} & \Bbm_k^{\ubm\ybm} \\[0.5em]
  (\Bbm_k^{\ubm\ybm})^T & \Bbm_k^{\ybm\ybm}
 \end{bmatrix},
\end{equation}
where the individual components
$\Bbm_k^{\ubm\ubm}\in\Rbb^{N_\ubm\times N_\ubm}$,
$\Bbm_k^{\ubm\ybm}\in\Rbb^{N_\ubm\times N_\ybm}$, and
$\Bbm_k^{\ybm\ybm}\in\Rbb^{N_\ybm\times N_\ybm}$
are defined as
\begin{equation} \label{eqn:hess1}
\begin{aligned}
\Bbm_k^{\ubm\ubm} &\coloneqq \pder{\Fbm}{\ubm}(\zbm_k)^T\pder{\Fbm}{\ubm}(\zbm_k) \\
\Bbm_k^{\ubm\ybm} &\coloneqq \pder{\Fbm}{\ubm}(\zbm_k)^T\pder{\Fbm}{\ybm}(\zbm_k) \\
\Bbm_k^{\ybm\ybm} &\coloneqq \pder{\Fbm}{\ybm}(\zbm_k)^T\pder{\Fbm}{\ybm}(\zbm_k)+
 \gamma_k\pder{\phibold}{\ybm}(\ybm_k)^T\Dbm_k\pder{\phibold}{\ybm}(\ybm_k)
\end{aligned}
\end{equation}
and $\Dbm_k\in\Rbb^{N_\xbm\times N_\xbm}$ is a SPD matrix constructed to
regularize the mesh motion and $\gamma_k\in\Rbb_{\geq 0}$ is a regularization
parameter. The regularization matrix, $\Dbm_k$, is taken to be the linear
elasticity (isotropic) stiffness matrix with elasticity modulus inversely
proportional to the volume of the elements in the reference mesh
\cite{zahr_implicit_2020}.
The regularization parameter, $\gamma_k$, is difficult to choose
\emph{a priori} because the performance of the solver is sensitive
to its value: small values can lead to an ill-conditioned Hessian
approximation that results in poor-quality search directions and
large mesh motion, whereas large values over-regularize the Hessian
approximation, which can lead to small steps with little information
from the actual Hessian.
Following \cite{corrigan_moving_2019, zahr_implicit_2020},
we reduce the sensitivity of the SQP method
to the choice of the regularization parameter by choosing it
adaptively such that the magnitude of the mesh motion lies
in a reasonable range. That is, if the magnitude of
$\Delta\xbm_k\coloneqq\phibold(\ybm_k)-\phibold(\ybm_{k-1})$
is too large (small), we increase (decrease) the regularization
parameter for the next iteration as
\begin{equation}
\gamma_{k+1} = \max \{ \bar{\gamma}_{k+1}, \gamma_{\min} \}, \qquad 
\bar{\gamma}_{k+1} = \begin{cases}
\tau^{-1}\gamma_k &\text{if }\norm{\Delta\xbm_k}<\sigma_1 L \\
\tau \gamma_k &\text{if }\norm{\Delta\xbm_k}>\sigma_2 L \\
\gamma_k &\text{otherwise},
\end{cases}
\end{equation}
where $L\in\Rbb_{>0}$ is a reference length for the domain,
$\tau\in\Rbb_{>0}$ is a parameter that controls the aggressiveness
of the adaptation, $\sigma_1,\sigma_2\in\Rbb_{>0}$ are parameters
that define acceptable magnitudes of the mesh motion relative to the
reference length scale, $\gamma_{\min}\in\Rbb_{>0}$ is a minimum
value for the regularization parameter, and $\gamma_0\in\Rbb_{\geq 0}$
is the starting value for the regularization parameter.
In this work, we typically choose $\tau\in(1,2]$, $\sigma_1 = 10^{-2}$, 
$\sigma_2 = 10^{-1}$, $\gamma_{\min}$ and $\gamma_0$ are problem-specific,
and $L$ is determined from the extents of the domain.
\begin{remark} \label{rem:hess_approx}
The choice of Hessian approximation in (\ref{eqn:hess0})-(\ref{eqn:hess1})
is identical to that introduced in \cite{zahr_implicit_2020} with the sole
exception being the definition of the regularization matrix, $\Dbm_k$,
that was previously defined as the stiffness matrix of independent
Poisson equations for each coordinate direction. Both
choices lead to similar performance; however, isotropic
linear elasticity is a more natural and intuitive choice
because it leads to the interpretation of the mesh as a
pseudo-structure, which is common in settings that involve
mesh motion, e.g., fluid-structure interaction
\cite{lesoinne_linearized_2001}.
\end{remark}

\subsubsection{Mesh quality parameter adaptation}
\label{sec:ist_solve:sqp:mshqual}
In addition to the choice of Hessian approximation
(Section~\ref{sec:ist_solve:sqp:lmhess}), the mesh quality parameter
($\kappa$) in (\ref{eqn:obj0}) plays a critical role in
obtaining high-quality search directions; however, for complex problems,
selecting a single value of $\kappa$ \emph{a priori} as in
\cite{zahr_implicit_2020}
does not work well in practice. If $\kappa$ is taken too small, it
can lead to search directions that severely degrade the quality of
the mesh, whereas excessively large values of $\kappa$ will prioritize
mesh quality over alignment with discontinuities, potentially
resulting in convergence to a non-aligned mesh. This is
complicated by the fact that the meaning of ``small'' and
``large'' values of $\kappa$ changes throughout the
optimization process as $f_\mathrm{err}(\zbm_k)$ and
$f_\mathrm{msh}(\ybm_k)$ change relative to each other.
We circumvent this issue by adapting the value of $\kappa$
during the optimization procedure, i.e., at iteration $k$ of
the optimization procedure, $\kappa$ in the definition of
the objective function (\ref{eqn:obj0}) is replaced with $\kappa_k$,
where $\{\kappa_k\}_{k=0}^\infty$ is a sequence of mesh
quality parameters. Because modifying the penalty parameter
changes the definition of the objective function, we aim to
adapt it infrequently, i.e., $\kappa_k=\kappa_{k-1}$ for many
iterations, and for only a finite number of iterations, i.e.,
$\kappa_k = \bar{\kappa}$ for $k>M$, where $M>0$ is a constant.
The latter condition ensures the definition of the objective
function is fixed in the asymptotic regime to not hinder or
stall convergence.

Our adaptation procedure aims to keep the contributions from
$f_\mathrm{err}$ and $\kappa^2 f_\mathrm{msh}$ to the objective
function relatively balanced with $f_\mathrm{err}$ slightly
dominate to ensure tracking discontinuities is prioritized.
To this end, we take
\begin{equation} \label{eqn:kappa_update}
\kappa_k = \max \{ \upsilon_k \kappa_{k-1}, \kappa_{\min} \},
\qquad 
\upsilon_k =
 \begin{cases}
   \upsilon &\text{if }
   f_\text{err}(\zbm_k) < \xi \kappa_{k-1}^2 f_\mathrm{msh}(\ybm_k) \text{ and } k \leq M \\
   1 &\text{otherwise},
\end{cases}
\end{equation}
where $\xi\in\Rbb_{>0}$ is a constant that defines the target
ratio between $f_\mathrm{err}$ and $\kappa^2f_\mathrm{msh}$,
$\upsilon\in(0,1)$ is the rate at which $\kappa$ decreases,
$\kappa_{\min}\in\Rbb_{\geq 0}$ is a lower bound on the mesh
quality parameter, and $\kappa_0$ is the starting point for
the mesh quality parameter, usually chosen as
$\kappa_0 = \sqrt{f_\mathrm{err}(\zbm_0)/f_\mathrm{msh}(\ybm_0)}$.

\begin{remark}
Based on the updated formula in (\ref{eqn:kappa_update}), the value
of $\kappa$ is only modified periodically. To see this, let $k'\ll M$
be the last iteration in which the value of $\kappa$ was modified,
which implies
\begin{equation}\label{eqn:kappa_remark0}
 f_\mathrm{err}(\zbm_{k'}) < \xi \kappa_{k'-1}^2 f_\mathrm{msh}(\ybm_{k'}),
\end{equation}
and suppose
$\upsilon\kappa_{k'-1}>\kappa_\mathrm{min}$, then
$\kappa_{k'}=\upsilon\kappa_{k'-1}$. Furthermore, we assume
\begin{equation} \label{eqn:kappa_remark1}
 f_\mathrm{err}(\zbm_{k'})\geq
 \xi(\upsilon\kappa_{k'-1})^2 f_\mathrm{msh}(\ybm_{k'});
\end{equation}
although this is not explicitly enforced algorithmically, it will
eventually hold because $\kappa_k$ will decrease geometrically if
$f_\mathrm{err}(\zbm_k)$ is dominated by
$\kappa_k^2 f_\mathrm{msh}(\ybm_k)$ due to the adaptation
criteria in (\ref{eqn:kappa_update}). Let $k''>k'$ be the
next iteration in which the value of $\kappa$ is updated,
i.e., the smallest value such that
\begin{equation} \label{eqn:kappa_remark2}
 f_\text{err}(\zbm_{k''}) <
 \xi\kappa_{k''-1}^2f_\text{msh}(\ybm_{k''}) =
 \xi(\upsilon\kappa_{k'-1})^2f_\text{msh}(\ybm_{k''}),
\end{equation}
where we used $\kappa_k=\kappa_{k'}=\upsilon \kappa_{k'-1}$
for $k=k',\dots,k''-1$ (definition of $k''$).
In the early iterations of the optimization procedure ($k'$ small),
the values of $f_\mathrm{err}(\zbm_k)$ and $f_\mathrm{msh}(\ybm_k)$
will change significantly between iterations, which may cause frequent
updates depending on the choice of $\xi$ and $\upsilon$. However,
in later iterations ($k'$ large), the values of $f_\mathrm{err}(\zbm_k)$
and $f_\mathrm{msh}(\ubm_k)$ stabilize close to their optimal values.
If $f_\mathrm{err}(\zbm_k) \approx f_\mathrm{err}(\zbm_{k'})$ and
$f_\mathrm{msh}(\ybm_k) \approx f_\mathrm{msh}(\ybm_{k'})$ for
$k>k'$, there may not exist a value of $k''$ such that
(\ref{eqn:kappa_remark2}) is satisfied in light of (\ref{eqn:kappa_remark1}),
which means the final update to $\kappa$ occurs at iteration $k'$.
Because of this, in practice, the mesh penalty
updates cease naturally and it is not necessary to terminate them after $M$
iterations. In practice, we initially observe numerous successive 
iterations where updates occur to calibrate the value
of $\kappa_k$ until both (\ref{eqn:kappa_remark0}) and
(\ref{eqn:kappa_remark1}) hold, after which updates occur
only periodically.
\end{remark}

\begin{remark}
Adapting penalty parameters is a common practice in numerical (constrained)
optimization. For example, penalty, augmented Lagrangian, and interior point
methods construct unconstrained subproblems whose objective function is the
original objective penalized by the constraint violation weighted by a penalty
parameter \cite{gill_practical_1981,nocedal2006numerical}.
The penalty parameter is fixed while an approximate
local minima to the unconstrained problem is sought and then adapted. The
approach in (\ref{eqn:kappa_update}) uses the structure of the objective
function to update the penalty parameter rather than stationarity of
the subproblems.
\end{remark}


\subsubsection{Incorporation of robustness measures}
\label{sec:ist_solve:sqp:robust}
The SQP solver outlined to this point is identical to the one
introduced in \cite{zahr_implicit_2020}, barring a few exceptions
(Remarks~\ref{rem:lsrch_pnlty}-\ref{rem:hess_approx},
 Section~\ref{sec:ist_solve:sqp:mshqual}).
Numerical experimentation with this baseline SQP solver revealed
it is not sufficiently robust to handle more complex problems than
those explored in \cite{zahr_implicit_2020}, particularly high Mach flows
and flows with complex discontinuity (shock and contacts) surfaces. In this
work, we propose a number of modifications to the SQP step with the
intention of improving its overall robustness for complex flows. To
this end, the iteration update in (\ref{eqn:sqp_step0}) will be changed to
\begin{equation} \label{eqn:sqp_step1}
 \tilde\zbm_{k+1} = \zbm_k + \alpha_k\Delta \zbm_k, \qquad
 \zbm_{k+1} = \Upsilon_{k+1}(\tilde\zbm_{k+1}),
\end{equation}
where $\tilde\zbm_{k+1}$ is the step that comes directly from SQP,
which is then modified using the mapping
$\func{\Upsilon_k}{\Rbb^{N_\zbm}}{\Rbb^{N_\zbm}}$ to
incorporate robustness measures (Section~\ref{sec:ist_solve:mod}).
To ensure this perturbation of the SQP step does not hinder or
stall convergence, we require $\Upsilon_k = \mathrm{Id}$ for
$k>M$, where $M>0$ is a constant. This ensures the modified
update in (\ref{eqn:sqp_step1}) reduces to the standard SQP
update in (\ref{eqn:sqp_step0}) after $M$ iterations. The
remainder of this section provides an algorithmic construction
of $\Upsilon_k$ that incorporates critical measures to ensure
the modified SQP method can robustly converge, even for complex flows.


\subsection{Modifications to SQP step for improved robustness}
\label{sec:ist_solve:mod}
We introduce three robustness measures that collectively define
the mapping $\Upsilon_k$ that are critical for the SQP solver
in the previous section to reliably converge for complex flows:
\begin{inparaenum}[1)]
 \item boundary-preserving, shock-aware element removal,
 \item removal of geometric curvature from inverted or ill-conditioned elements, and
 \item element-wise solution re-initialization.
\end{inparaenum}

\subsubsection{Element removal via edge collapse}
\label{sec:ist_solve:mod:clps}
For flows with complex shock surfaces, e.g., reflecting and interacting
shocks
(e.g., Figures~\ref{fig:euler:diamond:sltn}-\ref{fig:euler:scramjet:coarse}),
it is not reasonable to expect smooth deformations to a shock-agnostic mesh
(fixed mesh topology) to cause element faces to track all discontinuity
surfaces. This intuition was confirmed by numerical experiments that showed,
despite measures to keep the mesh well-conditioned, minimization of the
enriched residual causes elements to be crushed near discontinuity surfaces.
Rather than attempting to prevent this behavior, e.g., by increasing the
mesh quality and regularization parameters, the best option is to simply
remove the crushed elements. To accomplish this, we remove appropriate
elements by collapsing their shortest edge, an operation that is well-defined
for simplical meshes of any polynomial degree in any dimension. As such,
in this work, we only consider simplex meshes. Edge collapses are
the only topological mesh operation required for the proposed
implicit shock tracking method; more complex operations such as
edge flips, face swaps, refinement, and complete re-meshing are
not required, all of which become more complicated as the polynomial
degree and dimension of the mesh increase. This is a substantial
distinction of the proposed approach from explicit shock tracking
approaches \cite{moretti2002thirty, salas2009shock} that require
sophisticated mesh operations, which limits their utility
in three dimensions, particularly for complex shock structures.
While the concept of element removal via edge collapse is not new
in the context of implicit shock tracking
\cite{corrigan_moving_2019,zahr_implicit_2020},
we provide implementation details and identify two important
considerations, namely, preservation of boundaries and shocks.

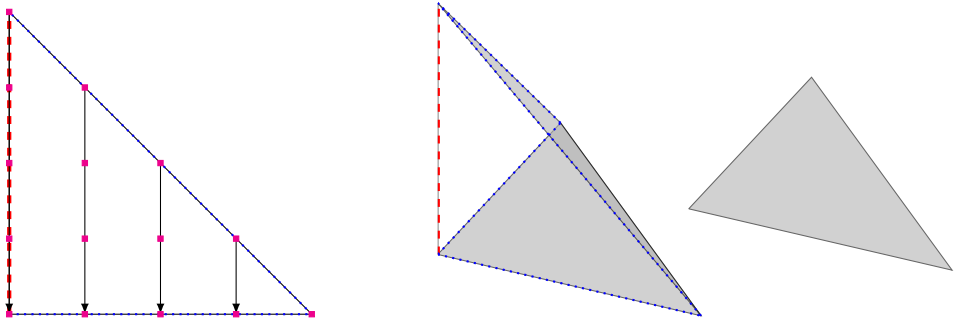
\begin{figure}
\centering
\begin{tikzpicture}
\begin{axis}[
axis equal image,
axis lines=none,
width=0.45\textwidth]
\addplot [, forget plot]
coordinates {
( 0.00000000e+00,  0.00000000e+00)
( 1.00000000e+00,  0.00000000e+00)
( 0.00000000e+00,  1.00000000e+00)
( 0.00000000e+00,  0.00000000e+00)};

\addplot [red, dashed, ultra thick]
coordinates {
( 0.00000000e+00,  0.00000000e+00)
( 0.00000000e+00,  1.00000000e+00)};\label{line:mshop_demo0:collapse}

\addplot [blue, dotted, thick]
coordinates {
( 0.00000000e+00,  0.00000000e+00)
( 1.00000000e+00,  0.00000000e+00)};\label{line:mshop_demo0:align}

\addplot [blue, dotted, thick, forget plot]
coordinates {
( 1.00000000e+00,  0.00000000e+00)
( 0.00000000e+00,  1.00000000e+00)};

\addplot [magenta, mark options={solid, thin}, mark=square*, mark size=1, only marks]
coordinates {
( 0.00000000e+00,  0.00000000e+00)
( 2.50000000e-01,  0.00000000e+00)
( 5.00000000e-01,  0.00000000e+00)
( 7.50000000e-01,  0.00000000e+00)
( 1.00000000e+00,  0.00000000e+00)};\label{line:mshop_demo0:node}

\addplot [magenta, mark options={solid, thin}, mark=square*, mark size=1, only marks, forget plot]
coordinates {
( 0.00000000e+00,  2.50000000e-01)
( 2.50000000e-01,  2.50000000e-01)
( 5.00000000e-01,  2.50000000e-01)
( 7.50000000e-01,  2.50000000e-01)};

\addplot [magenta, mark options={solid, thin}, mark=square*, mark size=1, only marks, forget plot]
coordinates {
( 0.00000000e+00,  5.00000000e-01)
( 2.50000000e-01,  5.00000000e-01)
( 5.00000000e-01,  5.00000000e-01)};

\addplot [magenta, mark options={solid, thin}, mark=square*, mark size=1, only marks, forget plot]
coordinates {
( 0.00000000e+00,  7.50000000e-01)
( 2.50000000e-01,  7.50000000e-01)};

\addplot [magenta, mark options={solid, thin}, mark=square*, mark size=1, only marks, forget plot]
coordinates {
( 0.00000000e+00,  1.00000000e+00)};

\draw [->, >=latex] plot [] coordinates {(axis cs:0.0, 1.0) (axis cs:0.0, 0.0)
};

\draw [->, >=latex] plot [] coordinates {(axis cs:0.25, 0.75) (axis cs:0.25, 0.0)
};

\draw [->, >=latex] plot [] coordinates {(axis cs:0.5, 0.5) (axis cs:0.5, 0.0)
};

\draw [->, >=latex] plot [] coordinates {(axis cs:0.75, 0.25) (axis cs:0.75, 0.0)
};

\end{axis}
\end{tikzpicture} \qquad \qquad
\begin{tikzpicture}
\begin{axis}[
axis equal image,
axis lines=none,
width=0.8\textwidth]
\addplot3 [opacity=0.6, fill=black!30!white, opacity=0.6, forget plot]
coordinates {
( 0.00000000e+00,  0.00000000e+00,  0.00000000e+00)
( 1.00000000e+00,  0.00000000e+00,  0.00000000e+00)
( 0.00000000e+00,  1.00000000e+00,  0.00000000e+00)
( 0.00000000e+00,  0.00000000e+00,  0.00000000e+00)};

\addplot3 [opacity=0.6, fill=black!30!white, opacity=0.6, forget plot]
coordinates {
( 1.00000000e+00,  0.00000000e+00,  0.00000000e+00)
( 0.00000000e+00,  1.00000000e+00,  0.00000000e+00)
( 0.00000000e+00,  0.00000000e+00,  1.00000000e+00)
( 1.00000000e+00,  0.00000000e+00,  0.00000000e+00)};

\addplot3 [opacity=0.6, fill=black!30!white, opacity=0.6, forget plot]
coordinates {
( 6.50000000e-01,  6.50000000e-01,  0.00000000e+00)
( 1.65000000e+00,  6.50000000e-01,  0.00000000e+00)
( 6.50000000e-01,  1.65000000e+00,  0.00000000e+00)
( 6.50000000e-01,  6.50000000e-01,  0.00000000e+00)};

\addplot3 [black, solid, forget plot]
coordinates {
( 0.00000000e+00,  0.00000000e+00,  0.00000000e+00)
( 0.00000000e+00,  0.00000000e+00,  1.00000000e+00)};

\addplot3 [red, dashed, ultra thick, forget plot]
coordinates {
( 0.00000000e+00,  0.00000000e+00,  0.00000000e+00)
( 0.00000000e+00,  0.00000000e+00,  1.00000000e+00)};

\addplot3 [blue, dotted, thick, forget plot]
coordinates {
( 0.00000000e+00,  0.00000000e+00,  0.00000000e+00)
( 1.00000000e+00,  0.00000000e+00,  0.00000000e+00)};

\addplot3 [blue, dotted, thick, forget plot]
coordinates {
( 0.00000000e+00,  0.00000000e+00,  0.00000000e+00)
( 0.00000000e+00,  1.00000000e+00,  0.00000000e+00)};

\addplot3 [blue, dotted, thick, forget plot]
coordinates {
( 1.00000000e+00,  0.00000000e+00,  0.00000000e+00)
( 0.00000000e+00,  0.00000000e+00,  1.00000000e+00)};

\addplot3 [blue, dotted, thick, forget plot]
coordinates {
( 0.00000000e+00,  1.00000000e+00,  0.00000000e+00)
( 0.00000000e+00,  0.00000000e+00,  1.00000000e+00)};

\end{axis}
\end{tikzpicture}
\caption{Mechanics of edge collapse for two-dimensional (\textit{left})
  and three-dimensional (\textit{middle}) simplex, including the collapsing
  edge (\ref{line:mshop_demo0:collapse}), edges/faces that will align after
  the collapse (\ref{line:mshop_demo0:align} and gray shaded), and the
  high-order nodes of the element (\ref{line:mshop_demo0:node}). In the
  three-dimensional case, the nodes are excluded for clarity and the
  remaining zero-volume element is included (\textit{right}).}
\label{fig:mshop_demo0}
\end{figure}

\begin{figure}
\centering
\begin{tikzpicture}
\begin{axis}[
axis equal image,
axis lines=none,
width=0.45\textwidth,
ymax=2,
xmax=3,
xmin=0,
ymin=0]
\addplot []
graphics [xmin=0,xmax=3,ymin=0,ymax=2] { 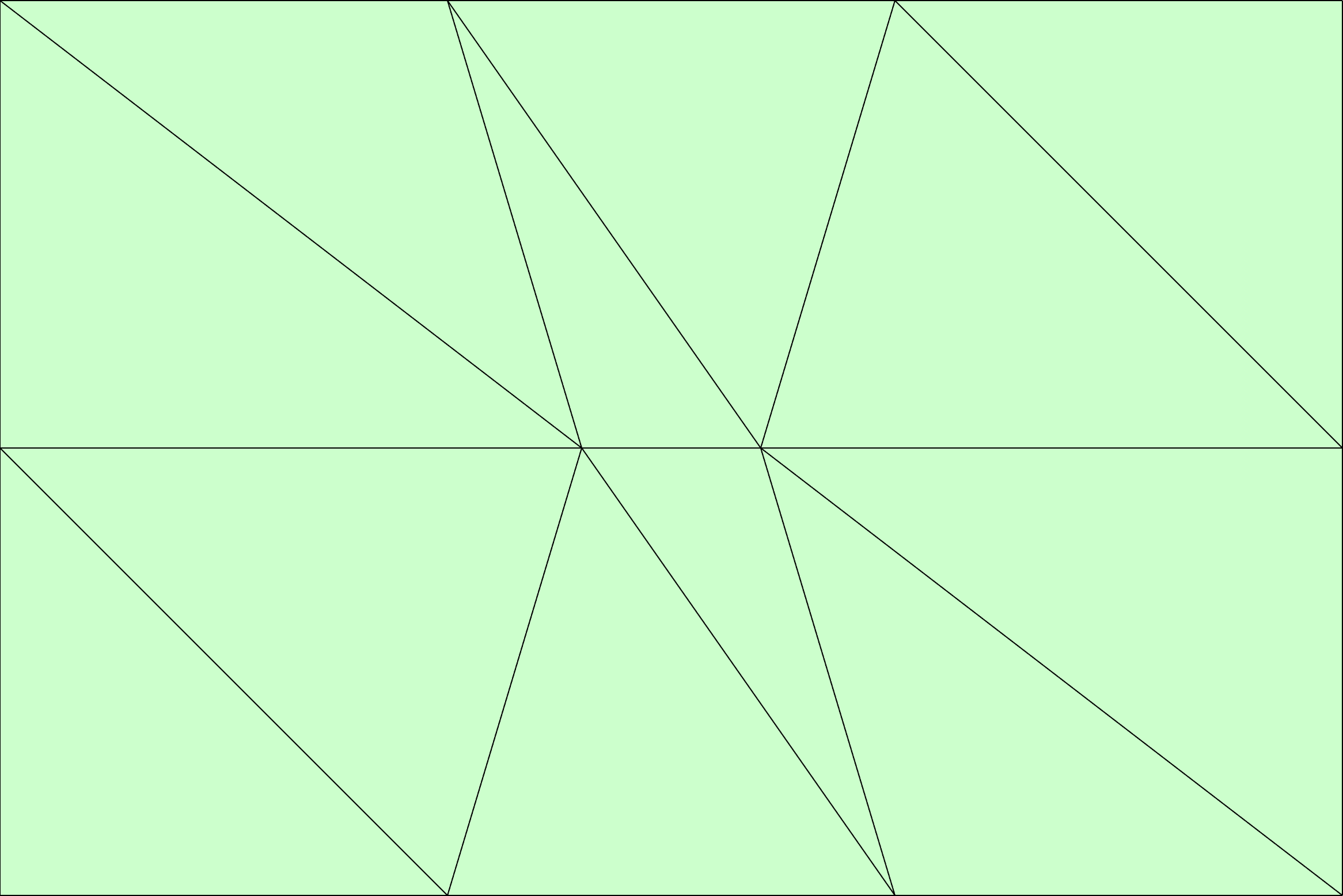};

\addplot [red, dashed, thick]
coordinates {
( 1.30000000e+00,  1.00000000e+00)
( 1.70000000e+00,  1.00000000e+00)};\label{line:mshop_demo2a:edge}

\end{axis}
\end{tikzpicture} \qquad
\begin{tikzpicture}
\begin{axis}[
axis equal image,
axis lines=none,
width=0.45\textwidth,
ymax=2,
xmax=3,
xmin=0,
ymin=0]
\addplot []
graphics [xmin=0,xmax=3,ymin=0,ymax=2] { 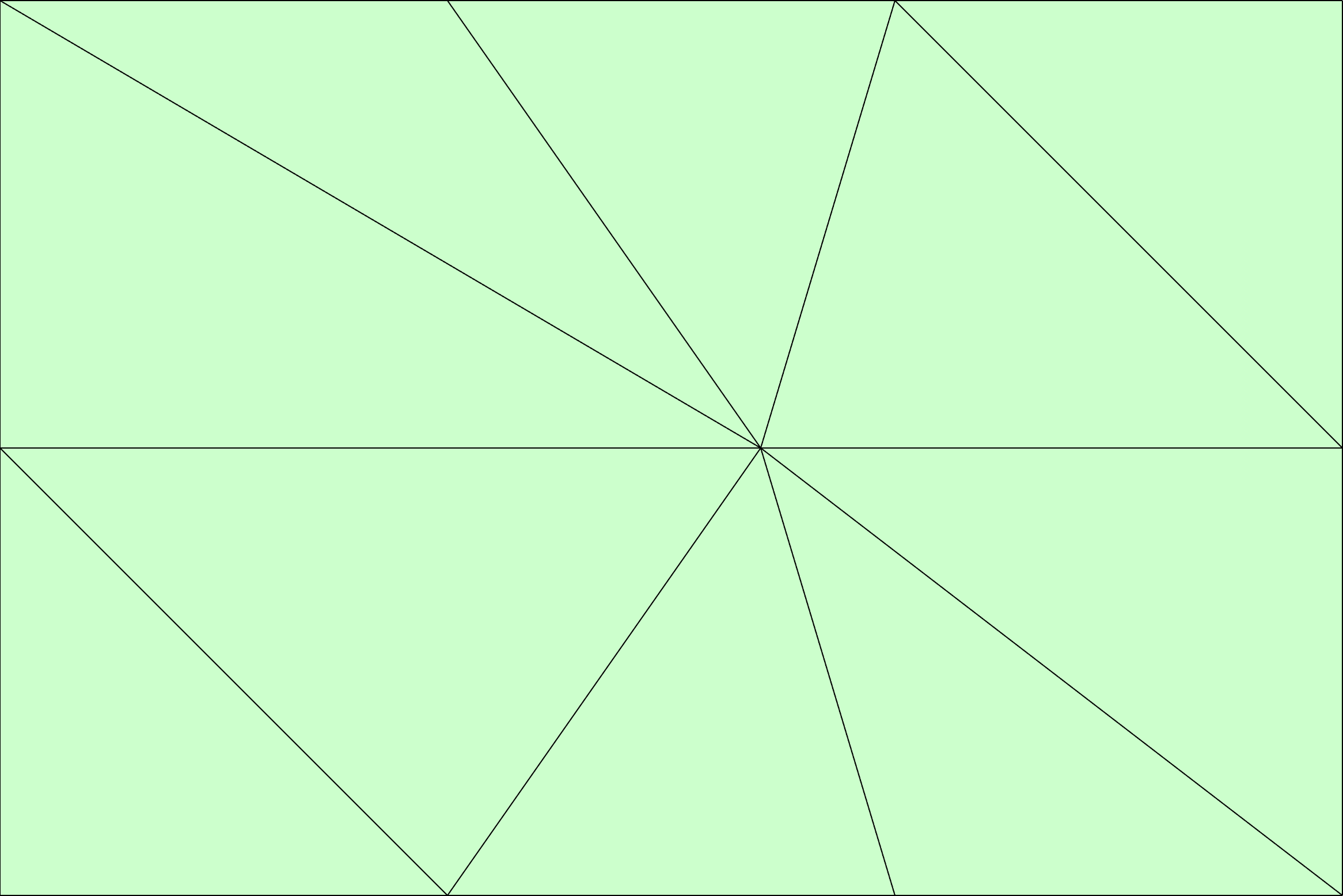};

\end{axis}
\end{tikzpicture}
\caption{Mesh prior to (\textit{left}) and after (\textit{right})
  collapse of the highlighted edge (\ref{line:mshop_demo2a:edge}) and
  removal of zero-volume elements.}
\label{fig:mshop_demo2}
\end{figure}

\emph{Mechanics of simplex collapse.}
From the definition of a simplex, collapsing any of its edges causes
the volume of the element to go to zero; this is not true for other
element geometries and the fundamental reason why element removal
via edge collapse is simple and straightforward for simplex element,
but not necessarily other geometries. Furthermore, the edge collapse
causes the surface area of $d-1$ faces of the simplex (those containing
the collapsed edge) to go to zero and the other $2$ faces (those the
collapsed edge connects) to overlap (Figure~\ref{fig:mshop_demo0}).
In the context of a mesh comprised of (potentially high-order) simplex
elements, collapsing an edge causes the volume of all elements containing
that edge to go to zero and the faces connected by that edge to overlap;
the latter condition ensures there are no gaps in the mesh. Subsequently,
these zero-volume elements are removed and the connectivity is adjusted
using the information in the newly overlapped faces
(Figure~\ref{fig:mshop_demo2}).

\begin{figure}
\centering
\begin{subfigure}[b]{0.25\textwidth}
 \includegraphics[width=\textwidth]{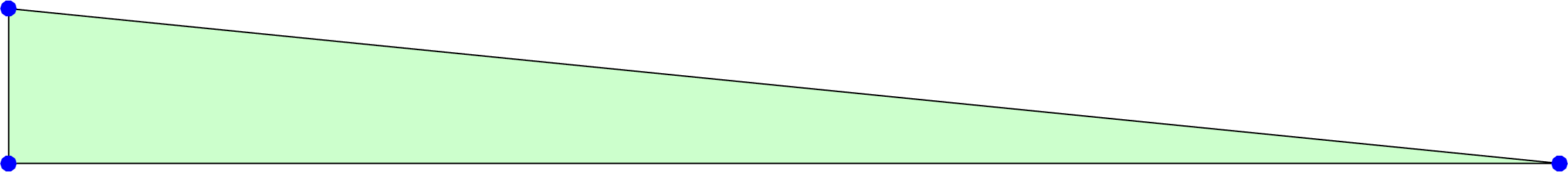}
 \caption{} \label{fig:mshop_demo3a}
\end{subfigure}
\quad
\begin{subfigure}[b]{0.15\textwidth}
 \includegraphics[width=\textwidth]{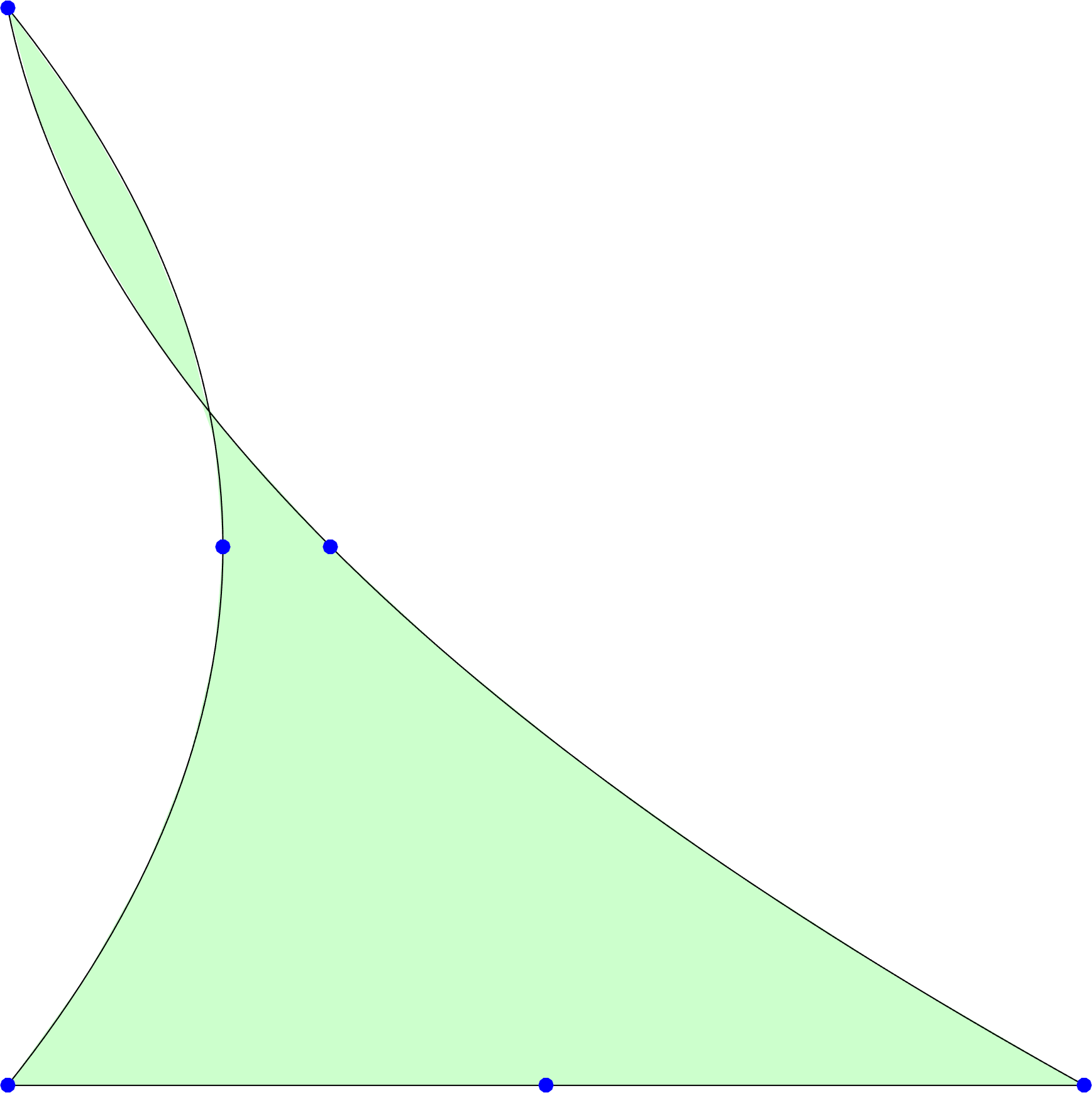}
 \caption{} \label{fig:mshop_demo3c}
\end{subfigure}
\quad
\begin{subfigure}[b]{0.076\textwidth}
 \includegraphics[width=\textwidth]{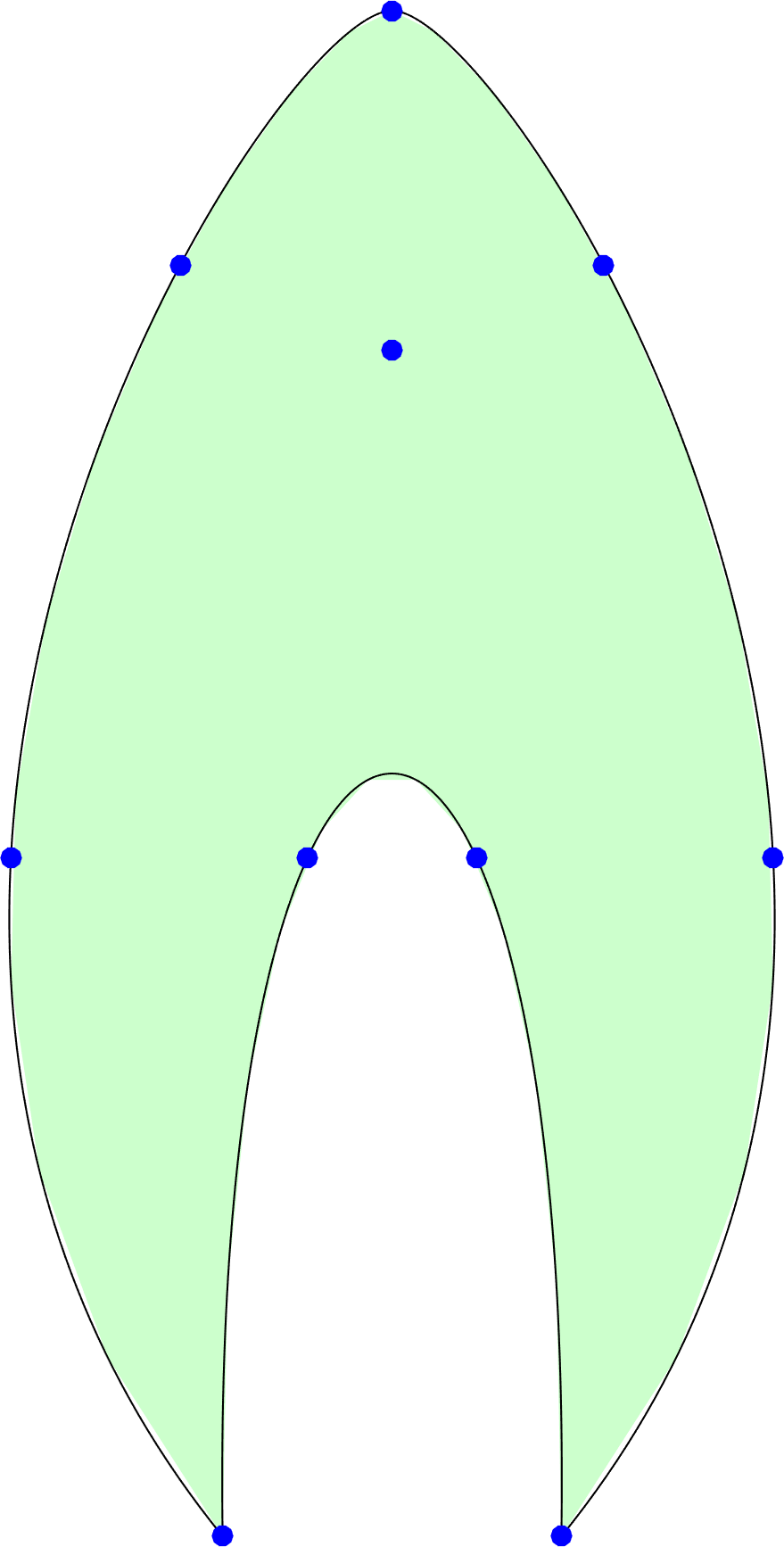}
 \caption{} \label{fig:mshop_demo3b}
\end{subfigure}
\quad
\begin{subfigure}[b]{0.15\textwidth}
 \includegraphics[width=\textwidth]{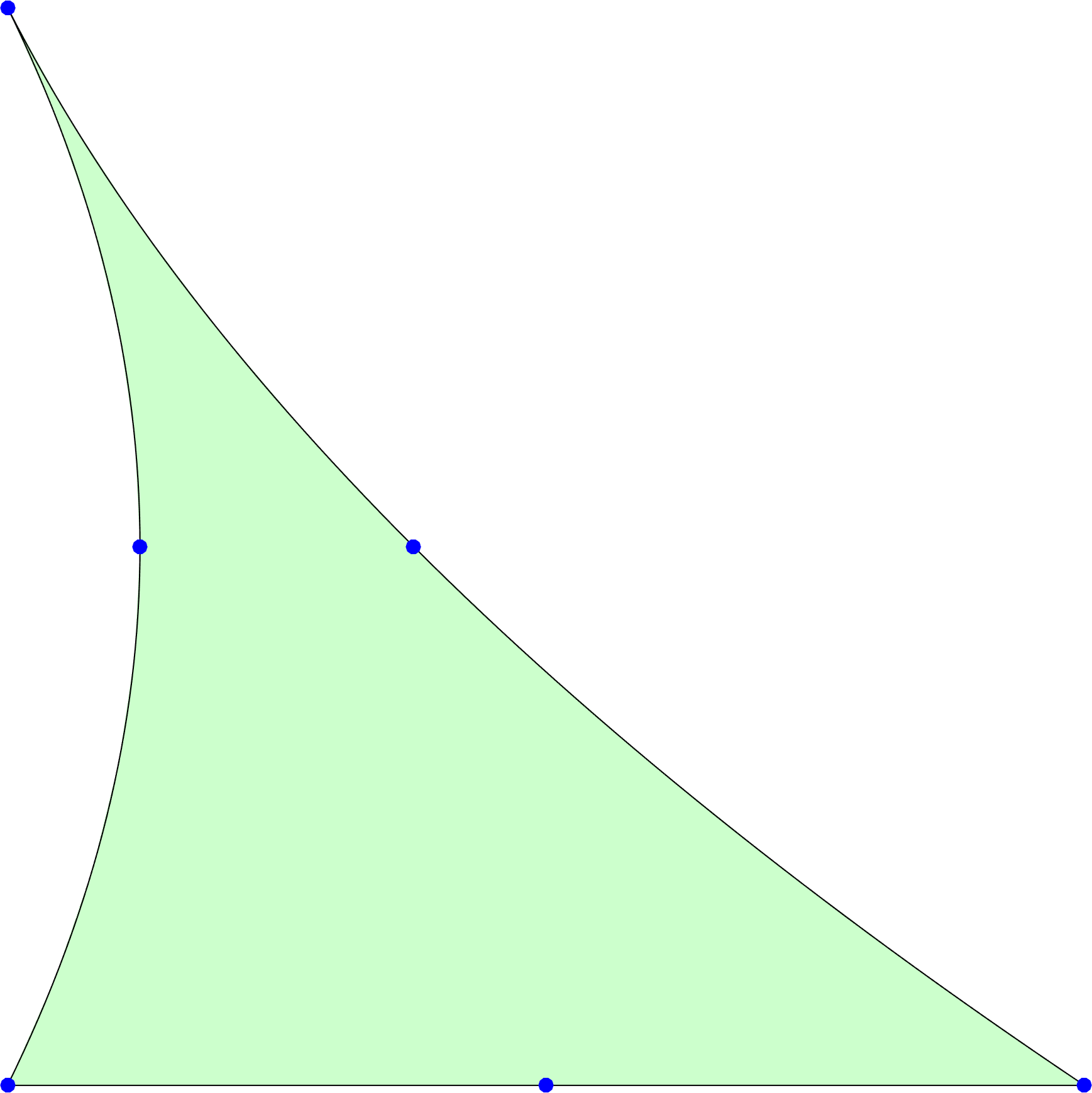}
 \caption{} \label{fig:mshop_demo3d}
\end{subfigure}
\quad
\begin{subfigure}[b]{0.15\textwidth}
 \includegraphics[width=\textwidth]{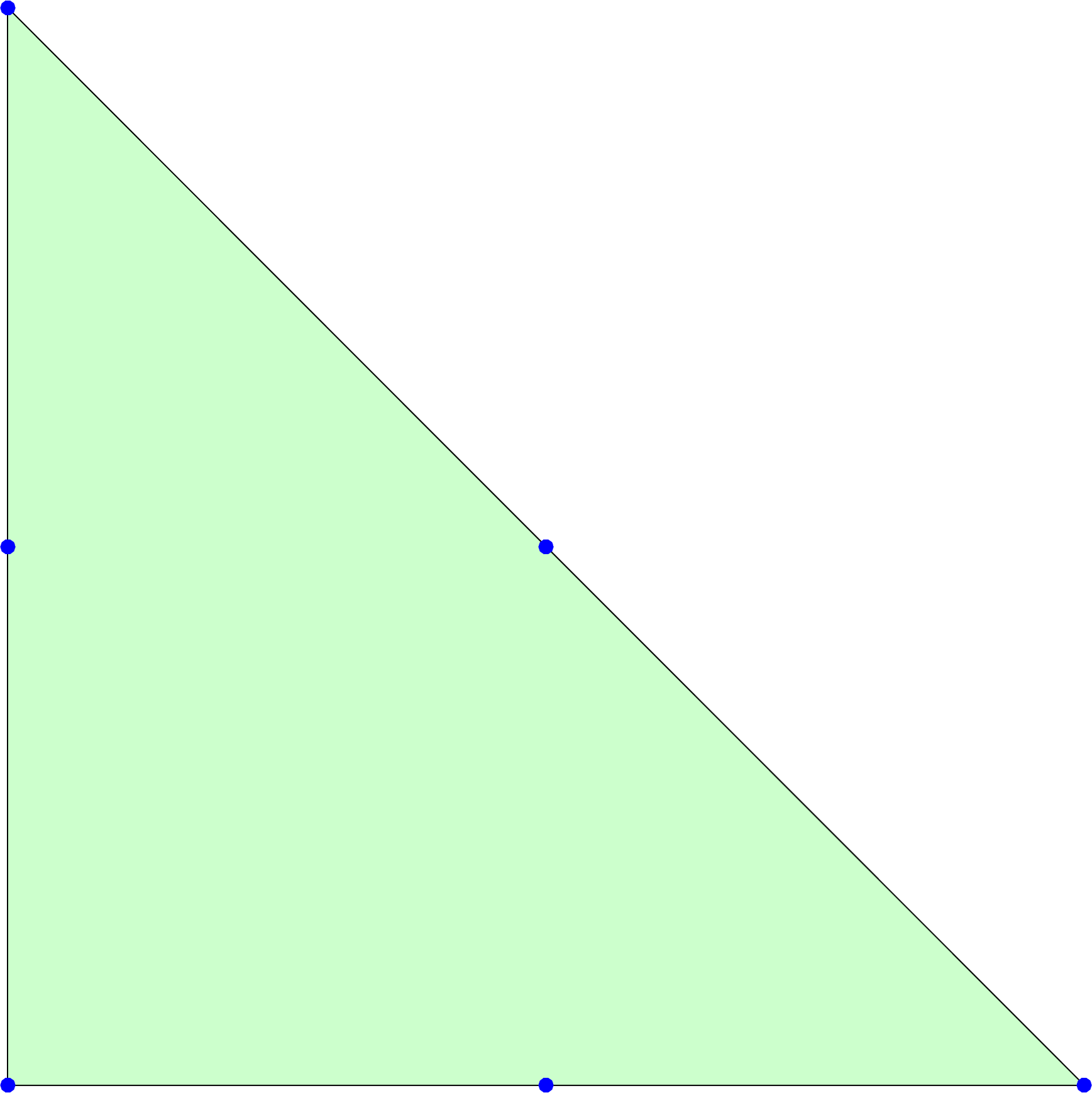}
 \caption{} \label{fig:mshop_demo3e}
\end{subfigure}
\caption{Illustration of problematic two-dimensional simplex elements.
   (a) An element with a small volume and short edge.
   (b) An inverted element ($g_{e,\mathrm{inf}}<0$).
   (c) A skewed element with a large volume, edges of similar lengths,
       and the bottom edge with a short distance between its endpoints,
       which justifies use of distance between edge endpoints rather
       than edge length in the definition of $\Ecal_{\mathrm{rmv},3}$.
   (d) An ill-conditioned element with $g_{e,\mathrm{inf}}<(1/60)g_{e,\mathrm{sup}}$.
   (e) The well-conditioned element that results from straightening the
       ill-conditioned element in (d).}
\label{fig:mshop_demo3}
\end{figure}

\emph{Identification of elements for removal.}
With the mechanics of simplical edge collapses established above, it
remains to decide which edges to collapse in the mesh. Because the
need to collapse elements arises from the desire to remove poor-quality
elements that emerge during the implicit shock tracking iterations, we
begin by identifying problematic elements and remove each one by collapsing
its shortest edge. In this work, problematic elements are elements with a
small volume (either in the reference or physical mesh), with a short minimum
edge length, or that are inverted or nearly inverted. To this end, define
$v_{0,e}\in\Rbb_{\geq 0}$ and $\func{v_e}{\Rbb^{N_\xbm}}{\Rbb_{\geq 0}}$
as the volume of the $e$th element of the reference and physical domain,
respectively, i.e.,
\begin{equation}
 v_{0,e} \coloneqq \int_{\Omega_{0,e}} \, dV, \qquad
 v_e : \xbm \mapsto \int_{\Gcal_h(\Omega_{0,e};\xbm)} \, dv.
\end{equation}
Furthermore, define $\func{\ell_{e,j}}{\Rbb^{N_\xbm}}{\Rbb_{\geq 0}}$,
where $\ell_{e,j}:\xbm\mapsto\ell_{e,j}(\xbm)$ denotes the distance
between the endpoints of the $j$th edge of the $e$th physical element,
$\Gcal_h(\Omega_{0,e};\xbm)$; this will be the notion of ``edge length''
used in this work rather than the true arc length of the edge
(Remark~\ref{rem:clps_edgelen}). We define the minimum and maximum
edge lengths of the $e$th element,
$\func{\ell_{e,\mathrm{min}}}{\Rbb^{N_\xbm}}{\Rbb_{\geq0}}$ and
$\func{\ell_{e,\mathrm{max}}}{\Rbb^{N_\xbm}}{\Rbb_{\geq0}}$, respectively,
as
\begin{equation}
 \ell_{e,\mathrm{min}} : \xbm \mapsto
 \min_{j\in\{1,\dots,N_\mathrm{ed}\}} \ell_{e,j}(\xbm), \qquad
 \ell_{e,\mathrm{max}} : \xbm \mapsto
 \max_{j\in\{1,\dots,N_\mathrm{ed}\}} \ell_{e,j}(\xbm),
\end{equation}
where $N_\mathrm{ed}=d(d+1)/2$ is the number of edges in a simplex.
Finally, we define the smallest and largest values of the domain
mapping Jacobian of the $e$th element,
$\func{g_{e,\mathrm{inf}}}{\Rbb^{N_\xbm}}{\Rbb}$ and
$\func{g_{e,\mathrm{sup}}}{\Rbb^{N_\xbm}}{\Rbb}$, respectively,
as
\begin{equation}
 g_{e,\mathrm{inf}} : \xbm \mapsto
 \inf_{X\in\Omega_{0,e}} \det(\bar\nabla\Gcal_h(X;\xbm)), \qquad
 g_{e,\mathrm{sup}} : \xbm \mapsto
 \sup_{X\in\Omega_{0,e}} \det(\bar\nabla\Gcal_h(X;\xbm)).
\end{equation}
With these definitions, we define the collection of problematic
elements as the set-valued function,
$\func{\Ecal_\mathrm{rmv}}{\Rbb^{N_\xbm}}{\Scal_h}$, defined as
\begin{equation}
 \Ecal_\mathrm{rmv} : \xbm \mapsto \Ecal_{\mathrm{rmv},1}(\xbm) \cup \Ecal_{\mathrm{rmv},2}(\xbm) \cup \Ecal_{\mathrm{rmv},3}(\xbm) \cup \Ecal_{\mathrm{rmv},4}(\xbm),
\end{equation}
where $\Scal_h$ is the collection of all subsets of $\Ecal_h$ and
the set-valued functions
$\func{\Ecal_{\mathrm{rmv},1}}{\Rbb^{N_\xbm}}{\Scal_h}$,
$\func{\Ecal_{\mathrm{rmv},2}}{\Rbb^{N_\xbm}}{\Scal_h}$,
$\func{\Ecal_{\mathrm{rmv},3}}{\Rbb^{N_\xbm}}{\Scal_h}$,
$\func{\Ecal_{\mathrm{rmv},4}}{\Rbb^{N_\xbm}}{\Scal_h}$ are
\begin{equation}
 \begin{aligned}
  \Ecal_{\mathrm{rmv},1}: \xbm &\mapsto \left\{\Omega_{0,e}\in\Ecal_h \suchthat v_e(\xbm) \leq c_1 v_{0,e} \right\} \\
  \Ecal_{\mathrm{rmv},2}: \xbm &\mapsto \left\{\Omega_{0,e}\in\Ecal_h \suchthat \min\{v_{0,e},v_e(\xbm)\} \leq c_2 \right\} \\
  \Ecal_{\mathrm{rmv},3}: \xbm &\mapsto \left\{\Omega_{0,e}\in\Ecal_h \suchthat \ell_{e,\mathrm{min}}(\xbm) \leq c_3 \ell_{e,\mathrm{max}}(\xbm) \right\} \\
  \Ecal_{\mathrm{rmv},4}: \xbm &\mapsto \left\{\Omega_{0,e}\in\Ecal_h \suchthat g_{e,\mathrm{inf}}(\xbm) \leq c_4 g_{e,\mathrm{sup}}(\xbm) \right\},
 \end{aligned}
\end{equation}
and $c_1,c_2,c_3,c_4\in\Rbb_{\geq 0}$ are tolerances.
That is, we identify elements for removal if its:
\begin{inparaenum}[1)]
 \item physical volume is small relative to that of the corresponding
   reference element ($\Ecal_{\mathrm{rmv},1}$),
 \item absolute volume in either the reference or physical domain is
   small ($\Ecal_{\mathrm{rmv},2}$),
 \item shortest edge is small relative to its longest edge in the
   physical domain ($\Ecal_{\mathrm{rmv},3}$), or
 \item mapping is nearly inverted in a relative sense
   ($\Ecal_{\mathrm{rmv},4}$).
\end{inparaenum}
In this work, we usually take $c_1 = 0.2$, $c_2 = 10^{-10}$, $c_3 = 0.2$, $c_4 = 0$,
with a few exceptions (Section~\ref{sec:numexp}). Commonly occurring
problematic elements are illustrated in Figure~\ref{fig:mshop_demo3}.

\begin{figure}
\centering
\begin{tikzpicture}
\begin{axis}[
axis equal image,
axis lines=none,
width=0.8\textwidth,
ymax=2.1,
xmax=5.1,
xmin=-0.1,
ymin=-0.1]
\addplot [white, forget plot]
coordinates {
(-1.00000000e-01, -1.00000000e-01)
( 5.10000000e+00, -1.00000000e-01)
( 5.10000000e+00,  2.10000000e+00)
(-1.00000000e-01,  2.10000000e+00)
(-1.00000000e-01, -1.00000000e-01)};

\addplot []
graphics [xmin=0,xmax=5,ymin=0,ymax=2] { 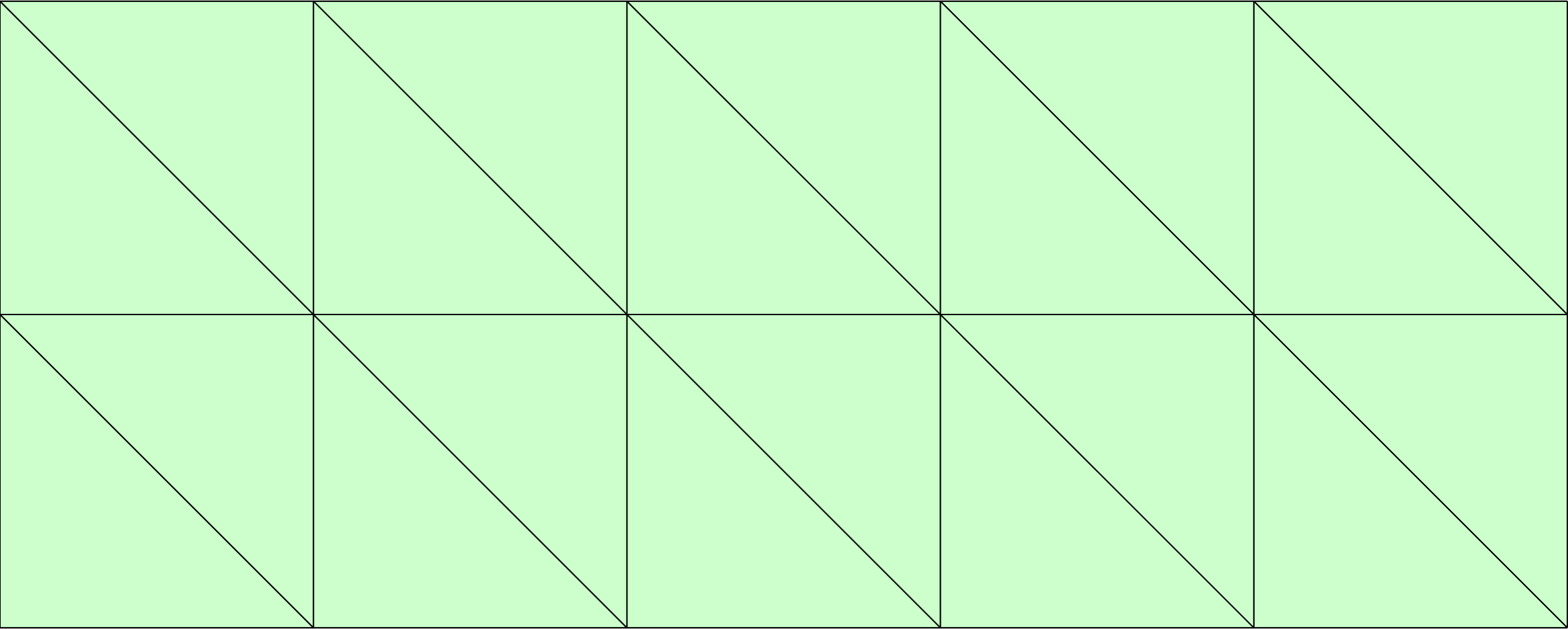};

\addplot [mark options={solid, thin}, mark=*, mark size=1.5, blue, only marks]
coordinates {
( 0.00000000e+00,  1.00000000e+00)
( 1.00000000e+00,  0.00000000e+00)
( 1.00000000e+00,  1.00000000e+00)
( 2.00000000e+00,  0.00000000e+00)
( 2.00000000e+00,  1.00000000e+00)
( 3.00000000e+00,  1.00000000e+00)
( 4.00000000e+00,  1.00000000e+00)};\label{line:edge_endpt_movable}

\addplot [mark options={solid, thin}, mark=square*, mark size=1.5, red, only marks]
coordinates {
( 0.00000000e+00,  0.00000000e+00)
( 0.00000000e+00,  2.00000000e+00)
( 5.00000000e+00,  1.00000000e+00)};\label{line:edge_endpt_fixed}

\draw [<-, >=latex] plot [] coordinates {(axis cs:0.0, 0.0) (axis cs:0.0, 1.0)
};

\draw [<-, >=latex] plot [] coordinates {(axis cs:0.0, 2.0) (axis cs:1.0, 1.0)
};

\draw [<->, >=latex] plot [] coordinates {(axis cs:1.0, 0.0) (axis cs:2.0, 0.0)
};

\draw [<->, >=latex] plot [] coordinates {(axis cs:2.0, 1.0) (axis cs:3.0, 1.0)
};

\draw [->, >=latex] plot [] coordinates {(axis cs:4.0, 1.0) (axis cs:5.0, 1.0)
};

\end{axis}
\end{tikzpicture}
\caption{Demonstration of movable (\ref{line:edge_endpt_movable}) and
  fixed (\ref{line:edge_endpt_fixed}) endpoints for potential edge
  collapses throughout a mesh. A point must be fixed if it lies at
  the intersection of more boundaries than the opposite endpoint of
  the edge.}
\label{fig:mshop_demo4}
\end{figure}
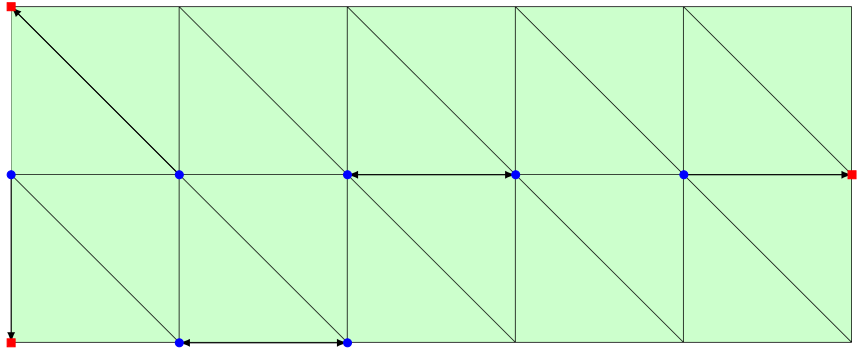

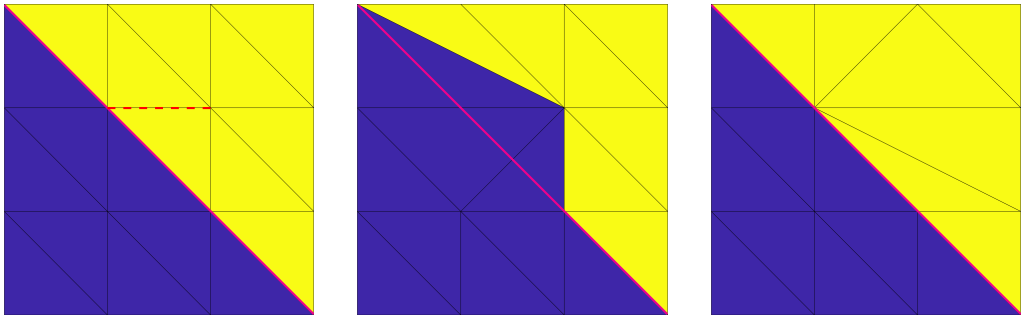
\begin{figure}[h]
\centering
\begin{tikzpicture}
\begin{axis}[
axis equal image,
axis lines=none,
width=0.4\textwidth,
ymax=1,
xmax=1,
xmin=0,
ymin=0]
\addplot []
graphics [xmin=0,xmax=1,ymin=0,ymax=1] { 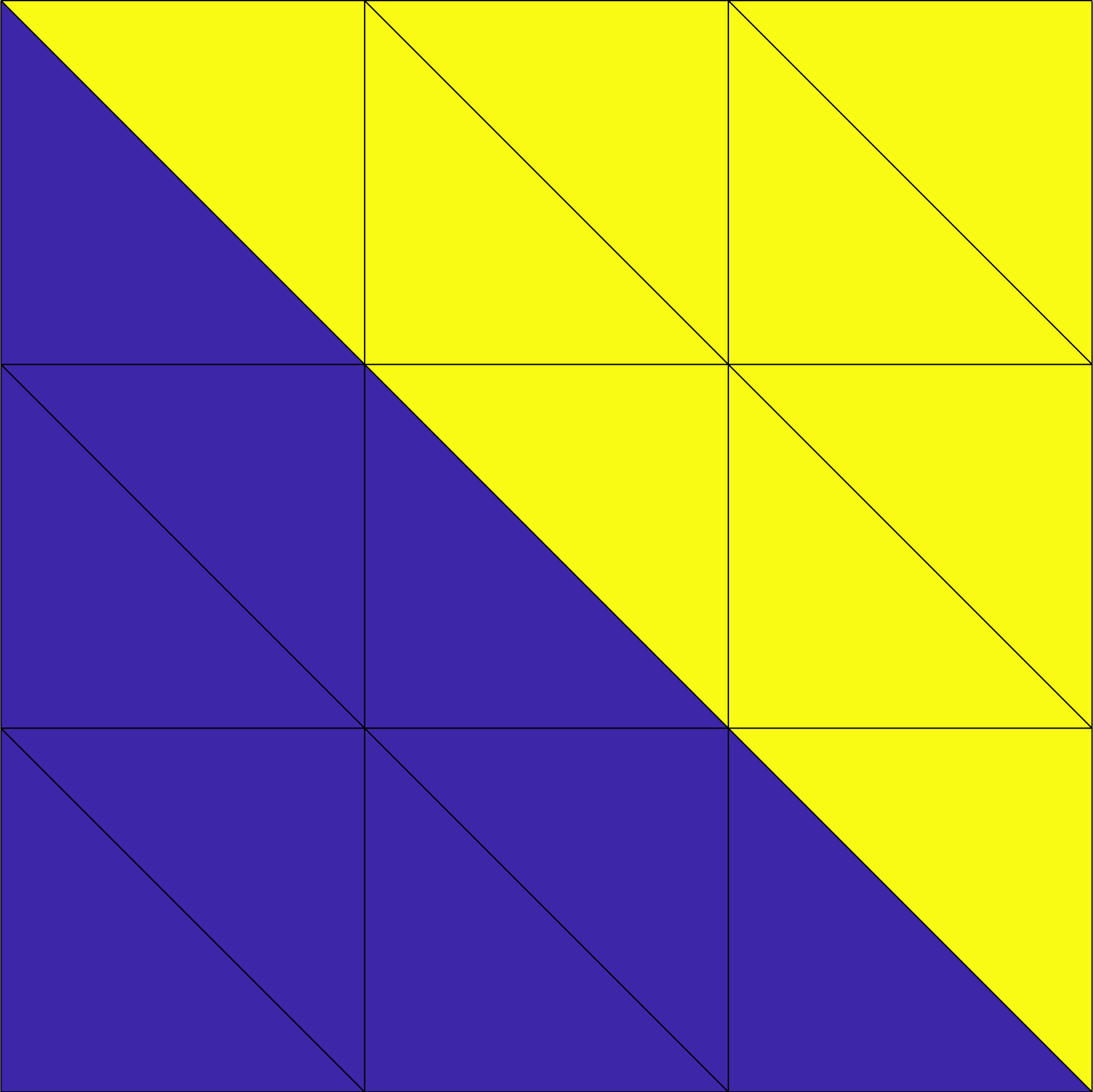};

\addplot [magenta, solid, thick]
coordinates {
( 0.00000000e+00,  1.00000000e+00)
( 1.00000000e+00,  0.00000000e+00)};\label{line:shk_preserve:shock}

\addplot [red, dashed, thick]
coordinates {
( 3.33333333e-01,  6.66666667e-01)
( 6.66666667e-01,  6.66666667e-01)};\label{line:shk_preserve:edge}

\end{axis}
\end{tikzpicture} \quad
\begin{tikzpicture}
\begin{axis}[
axis equal image,
axis lines=none,
width=0.4\textwidth,
ymax=1,
xmax=1,
xmin=0,
ymin=0]
\addplot []
graphics [xmin=0,xmax=1,ymin=0,ymax=1] { 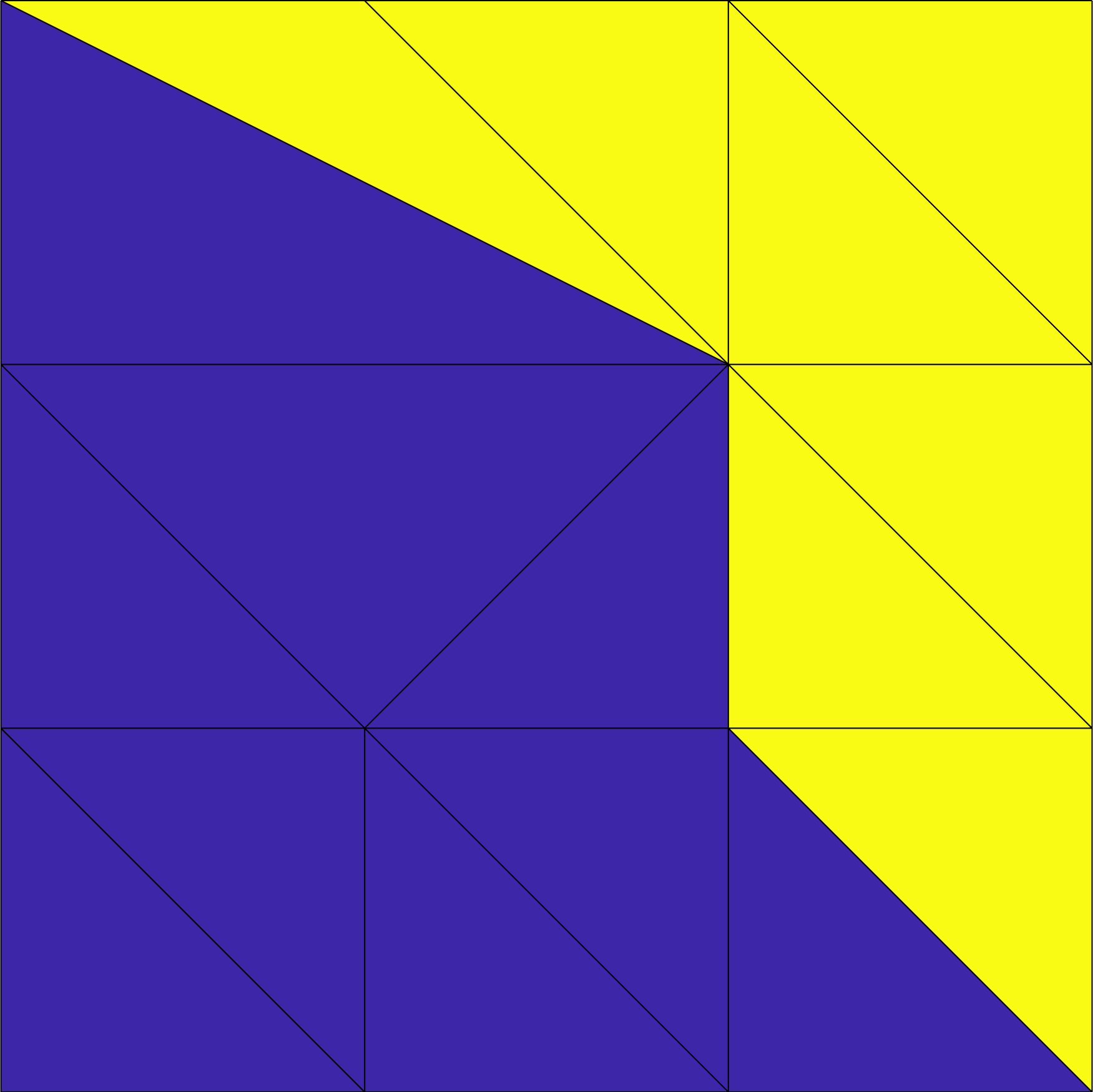};

\addplot [magenta, solid, thick]
coordinates {
( 0.00000000e+00,  1.00000000e+00)
( 1.00000000e+00,  0.00000000e+00)};\label{line:shk_preserve:shock}

\end{axis}
\end{tikzpicture} \quad
\begin{tikzpicture}
\begin{axis}[
axis equal image,
axis lines=none,
width=0.4\textwidth,
ymax=1,
xmax=1,
xmin=0,
ymin=0]
\addplot []
graphics [xmin=0,xmax=1,ymin=0,ymax=1] { 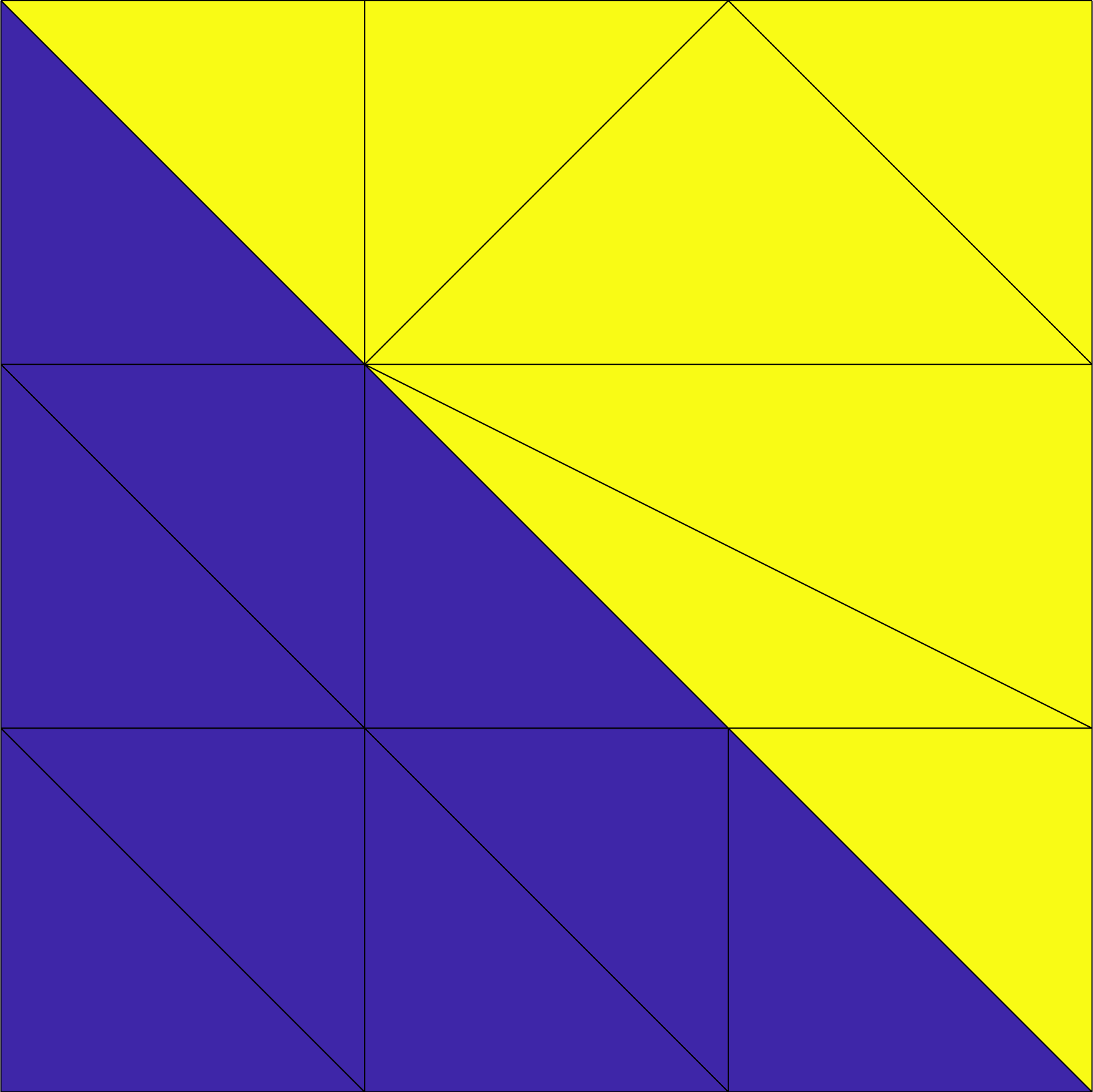};

\addplot [magenta, solid, thick]
coordinates {
( 0.00000000e+00,  1.00000000e+00)
( 1.00000000e+00,  0.00000000e+00)};\label{line:shk_preserve:shock}

\end{axis}
\end{tikzpicture}
\caption{A shock-aligned mesh (true shock given by
  \ref{line:shk_preserve:shock}) and corresponding solution (\textit{left})
  and two meshes that result from degenerating the highlighted edge
  (\ref{line:shk_preserve:edge}) at its right endpoint (\textit{middle})
  and its left endpoint (\textit{right}).
  Because the jump in the solution at the left endpoint is larger, we choose
  to degenerate the edge at that point; as a result, the shock is preserved.}
\label{fig:shk-aware-clps}
\end{figure}

\emph{Boundary-preserving and shock-aware edge collapse.}
With the elements identified for removal, we remove them one-by-one
by collapsing their shortest edge (edge with smallest distance between
endpoints) and updating the connectivity of the resulting mesh as
discussed previously. For the edge collapse to be well-defined,
the position of the degenerate edge (point) must be defined; the
two endpoints of the original edge are natural choices, although
any point in the vicinity of the original edge are candidates.
In this work, we position the collapsed edge at one of the
endpoints of the original edge with constraints on which
endpoint is chosen to ensure boundaries are preserved
(to high-order) and tracked shocks are maintained to the
extent possible. The first condition (boundary preservation)
is a strict requirement to ensure the PDE domain does not
change, whereas the second condition (shock preservation)
enhances the robustness of the SQP solver. To preserve
the domain boundary (to high-order), we collapse the
edge at whichever of the original endpoints lies at
the intersection of more boundaries (constraints).
This ensures remaining mesh nodes will not move off
boundaries due to an edge collapse, rather they
will slide along boundaries (Figure~\ref{fig:mshop_demo4}). In the
common case where both endpoints lie at the intersection
of the same boundaries, the position of the
degenerate edge will be determined from the value
of the DG solution $\Xi(\tilde\ubm_{k+1})$ to maintain the
tracked shock. Because the DG solution is multi-valued at
edge endpoints, we degenerate the edge at the endpoint at
which the range of values in $\Xi(\tilde\ubm_{k+1})$ is
largest. Because endpoints that lie on shocks will possess
a large range of values (equal to the magnitude of the
jump in the discontinuity), this effectively ensures nodes
will not be moved off a shock surface. The opposite choice
could degrade the representation of the shock, which could
either destroy or significantly slow convergence
(Figure~\ref{fig:shk-aware-clps}).

\textit{Element removal and solution transfer.}
Because the reference and physical meshes share the same
connectivity, edge collapses and subsequent element removal
are applied to the reference mesh 
and inherited by the physical mesh. 
Once all problematic elements are removed, the DG solution is
trivially transferred to the new mesh by removing degrees
of freedom corresponding to removed elements. In the
context of the SQP iterations,
$\tilde\zbm_{k+1}=(\tilde\ubm_{k+1},\tilde\ybm_{k+1})$
denotes the SQP update at iteration $k$ in (\ref{eqn:sqp_step1}).
Then, the problematic elements $\Ecal_\mathrm{rmv}(\phibold(\tilde\ybm_{k+1}))$
are removed via edge collapses to produce a new mesh, the parametrization
in (\ref{eqn:mapparam}) is reconstructed for the new mesh and the
unconstrained degrees of freedom are extracted, and the DG solution
is transferred to the new mesh to produce the state
$\hat\zbm_{k+1}=(\hat\ubm_{k+1},\hat\ybm_{k+1})$.

\begin{remark} \label{rem:clps_edgelen}
We use the distance between endpoints rather than edge length to
avoid favoring edges that have an artificially large length despite
its end nodes being close together (Figure~\ref{fig:mshop_demo3b});
in our experience, these edges are advantageous to collapse because
they are indicative of low-quality elements and degenerating it at either
of its endpoints induces a relatively small motion to the remaining nodes.
\end{remark}

\begin{remark}
Even though the reference element volume is independent of the
physical nodes of the mesh ($\xbm$), they will change every time
an element is collapsed. We have observed some pathological cases
where an element can be small in the reference mesh but not in the
physical mesh, leading to an ill-conditioned SQP system
(\ref{eqn:sqp_sys0}).
A robust algorithm must remove them or use another strategy to
prevent them from arising; we opt for the former case by including
an absolute tolerance on the reference element volume as a
condition for element removal ($\Ecal_{\mathrm{rmv},2}$).
\end{remark}

\begin{figure}
\centering
\begin{tikzpicture}
\begin{axis}[
axis equal image,
axis lines=none,
width=0.7\textwidth]
\addplot [opacity=0.6, fill=black!30!white, opacity=0.6, forget plot]
coordinates {
( 0.00000000e+00,  0.00000000e+00)
( 1.00000000e+00,  0.00000000e+00)
( 1.00000000e+01,  0.00000000e+00)
( 1.00000000e+01,  1.00000000e+00)
( 1.00000000e+00,  0.00000000e+00)};

\addplot [black, solid, thick]
coordinates {
( 0.00000000e+00,  0.00000000e+00)
( 1.00000000e+01,  0.00000000e+00)};\label{line:mshop_demo5:bnd}

\addplot [blue, solid, thick]
coordinates {
( 1.00000000e+00,  0.00000000e+00)
( 1.00000000e+01,  1.00000000e+00)};\label{line:mshop_demo5:shock}

\addplot [opacity=1.0, fill={rgb,1:red,0.8000000000;green,1.0000000000;blue,0.8000000000}, opacity=1.0, forget plot]
coordinates {
( 1.00000000e+00,  0.00000000e+00)
( 5.50000000e+00,  0.00000000e+00)
( 5.50000000e+00,  5.00000000e-01)
( 1.00000000e+00,  0.00000000e+00)};

\addplot [red, solid, thick]
coordinates {
( 5.50000000e+00,  0.00000000e+00)
( 5.50000000e+00,  5.00000000e-01)};\label{line:mshop_demo5:edge}

\addplot [mark options={solid, thin}, mark=square*, mark size=1.5, red, only marks]
coordinates {
( 5.50000000e+00,  0.00000000e+00)};\label{line:mshop_demo5:fixed}

\end{axis}
\end{tikzpicture} 
\caption{A situation that arises when a shock (\ref{line:mshop_demo5:shock})
  makes an acute angle with a boundary (\ref{line:mshop_demo5:bnd}). If the
  indicated element is targeted for removal, its shortest edge
  (\ref{line:mshop_demo5:edge}) will be collapsed, which makes
  it impossible to preserve both the shock and boundary. Because
  preservation of the boundary is a hard constraint, the edge is
  degenerated at its bottom endpoint (\ref{line:mshop_demo5:fixed}).}
\label{fig:mshop_demo5}
\end{figure}
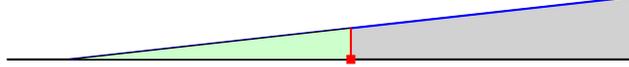

\begin{remark}
Because domain preservation is prioritized over shock
preservation, some pathological cases where the position
of the degenerate edge cannot be chosen to maintain the
shock, e.g., shocks that form sharp edges with boundaries
(Figure~\ref{fig:mshop_demo5}). In these cases, additional
SQP iterations are needed to fix the shock position.
Such situations are less likely to occur if the mesh
is finer in these regions or the element removal
parameters, particularly $c_3$, are weakened.
\end{remark}

\begin{remark}
In the situation where the shortest edge of an element identified
for removal cannot be collapsed without changing the domain boundary,
e.g., its endpoints lie on different boundaries, we collapse the
shortest admissible edge (i.e., edges whose collapse will preserve
the domain boundary) and use the same logic to identify the endpoint
at which to degenerate the edge to preserve the boundary and shock.
\end{remark}

\begin{remark}
In the physical
domain, the edge collapse is a relatively benign operation because
of the choice to collapse the shortest edge of highly skewed elements;
however, in the reference domain, the mesh motion is more severe because
the edge being collapsed is not necessarily short. If only the nodes
along the degenerated edge are moved as described above, the resulting
elements will be distorted (Figure~\ref{fig:elem-smooth}). This is
particularly problematic is the ideal element is taken to be the
reference element, i.e., $K_{\star,e}=\Omega_{0,e}$
(Section~\ref{sec:ist-form}),
because the objective function will drive the corresponding
physical element toward this distorted element. In this
work, we avoid this situation by removing the curvature
from the elements in the reference domain
(Figure~\ref{fig:elem-smooth}); this approach
is simpler than global mesh smoothing and
works well in practice.
\end{remark}
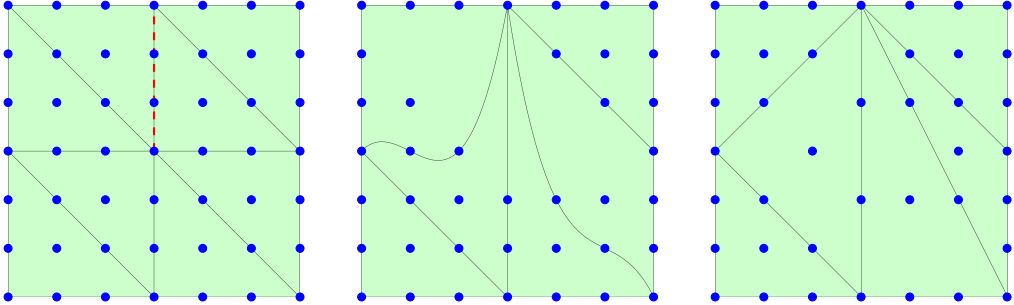
\begin{figure}
\centering
\begin{tikzpicture}
\begin{axis}[
axis equal image,
axis lines=none,
width=0.4\textwidth,
ymax=0.686666666667,
xmax=0.353333333333,
xmin=-0.353333333333,
ymin=-0.02]
\addplot [white, forget plot]
coordinates {
(-3.53333333e-01, -2.00000000e-02)
( 3.53333333e-01, -2.00000000e-02)
( 3.53333333e-01,  6.86666667e-01)
(-3.53333333e-01,  6.86666667e-01)
(-3.53333333e-01, -2.00000000e-02)};

\addplot []
graphics [xmin=-0.333333333333,xmax=0.333333333333,ymin=0.0,ymax=0.666666666667] { 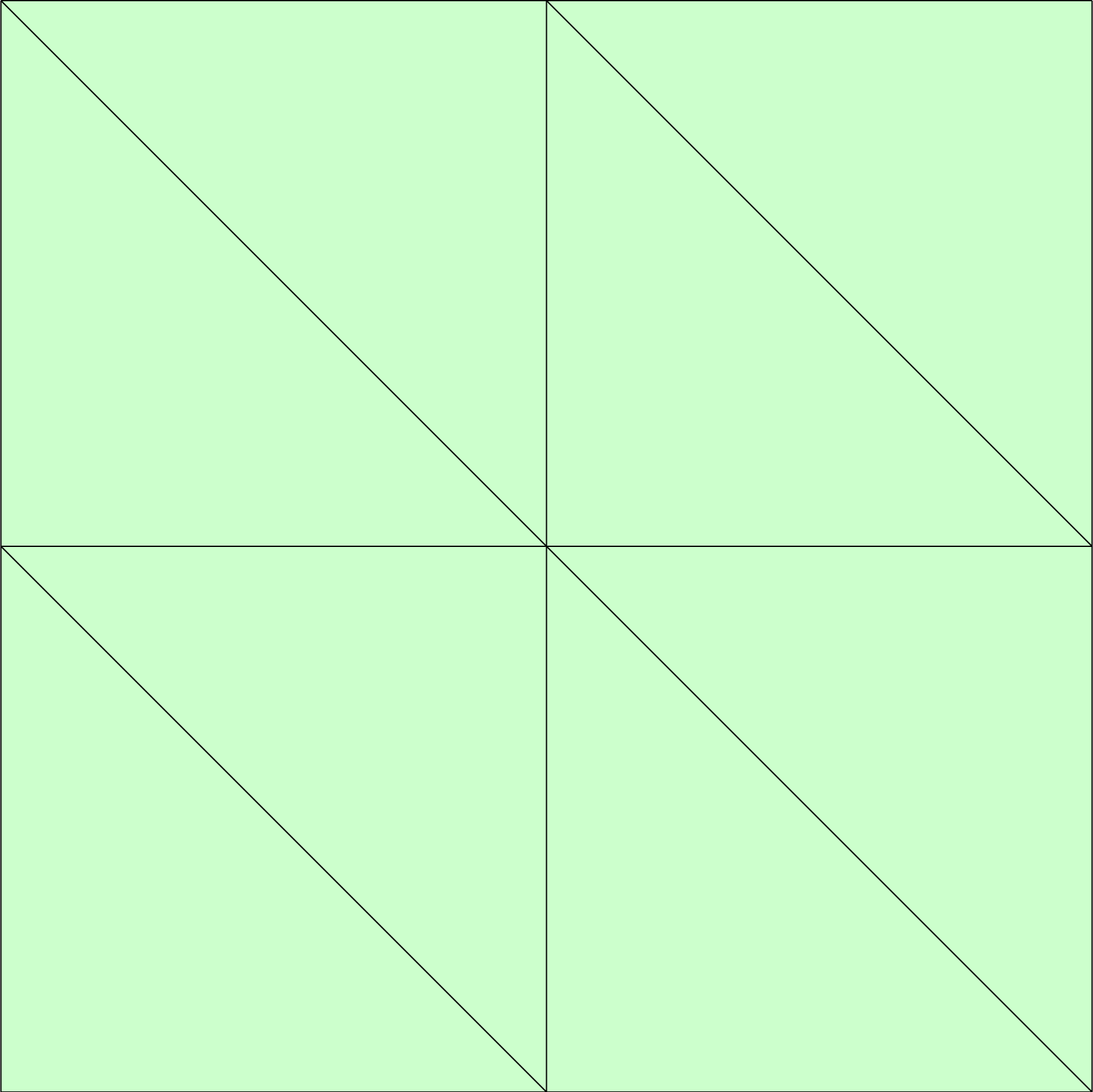};

\addplot [blue, mark options={solid, thin}, mark=*, mark size=1.5, only marks]
coordinates {
(-3.33333333e-01,  0.00000000e+00)
(-3.33333333e-01,  2.22222222e-01)
(-3.33333333e-01,  3.33333333e-01)
(-3.33333333e-01,  5.55555556e-01)
(-3.33333333e-01,  6.66666667e-01)
(-3.33333333e-01,  1.11111111e-01)
(-3.33333333e-01,  4.44444444e-01)
(-2.22222222e-01,  0.00000000e+00)
(-2.22222222e-01,  2.22222222e-01)
(-2.22222222e-01,  3.33333333e-01)
(-2.22222222e-01,  5.55555556e-01)
(-2.22222222e-01,  1.11111111e-01)
(-2.22222222e-01,  4.44444444e-01)
(-2.22222222e-01,  6.66666667e-01)
(-1.11111111e-01,  1.11111111e-01)
(-1.11111111e-01,  4.44444444e-01)
(-1.11111111e-01,  0.00000000e+00)
(-1.11111111e-01,  2.22222222e-01)
(-1.11111111e-01,  3.33333333e-01)
(-1.11111111e-01,  5.55555556e-01)
(-1.11111111e-01,  6.66666667e-01)
( 0.00000000e+00,  0.00000000e+00)
( 0.00000000e+00,  1.11111111e-01)
( 0.00000000e+00,  2.22222222e-01)
( 0.00000000e+00,  3.33333333e-01)
( 0.00000000e+00,  4.44444444e-01)
( 0.00000000e+00,  5.55555556e-01)
( 0.00000000e+00,  6.66666667e-01)
( 1.11111111e-01,  1.11111111e-01)
( 1.11111111e-01,  4.44444444e-01)
( 1.11111111e-01,  0.00000000e+00)
( 1.11111111e-01,  2.22222222e-01)
( 1.11111111e-01,  3.33333333e-01)
( 1.11111111e-01,  5.55555556e-01)
( 1.11111111e-01,  6.66666667e-01)
( 2.22222222e-01,  0.00000000e+00)
( 2.22222222e-01,  1.11111111e-01)
( 2.22222222e-01,  3.33333333e-01)
( 2.22222222e-01,  4.44444444e-01)
( 2.22222222e-01,  2.22222222e-01)
( 2.22222222e-01,  5.55555556e-01)
( 2.22222222e-01,  6.66666667e-01)
( 3.33333333e-01,  0.00000000e+00)
( 3.33333333e-01,  1.11111111e-01)
( 3.33333333e-01,  3.33333333e-01)
( 3.33333333e-01,  4.44444444e-01)
( 3.33333333e-01,  6.66666667e-01)
( 3.33333333e-01,  2.22222222e-01)
( 3.33333333e-01,  5.55555556e-01)};\label{line:smooth_ref:node}

\addplot [red, dashed, thick]
coordinates {
( 0.00000000e+00,  3.33333333e-01)
( 0.00000000e+00,  6.66666667e-01)};\label{line:smooth_ref:edge}

\end{axis}
\end{tikzpicture} \quad
\begin{tikzpicture}
\begin{axis}[
axis equal image,
axis lines=none,
width=0.4\textwidth,
ymax=0.686666666667,
xmax=0.353333333333,
xmin=-0.353333333333,
ymin=-0.02]
\addplot [white, forget plot]
coordinates {
(-3.53333333e-01, -2.00000000e-02)
( 3.53333333e-01, -2.00000000e-02)
( 3.53333333e-01,  6.86666667e-01)
(-3.53333333e-01,  6.86666667e-01)
(-3.53333333e-01, -2.00000000e-02)};

\addplot []
graphics [xmin=-0.333333333333,xmax=0.333333333333,ymin=0.0,ymax=0.666666666667] { 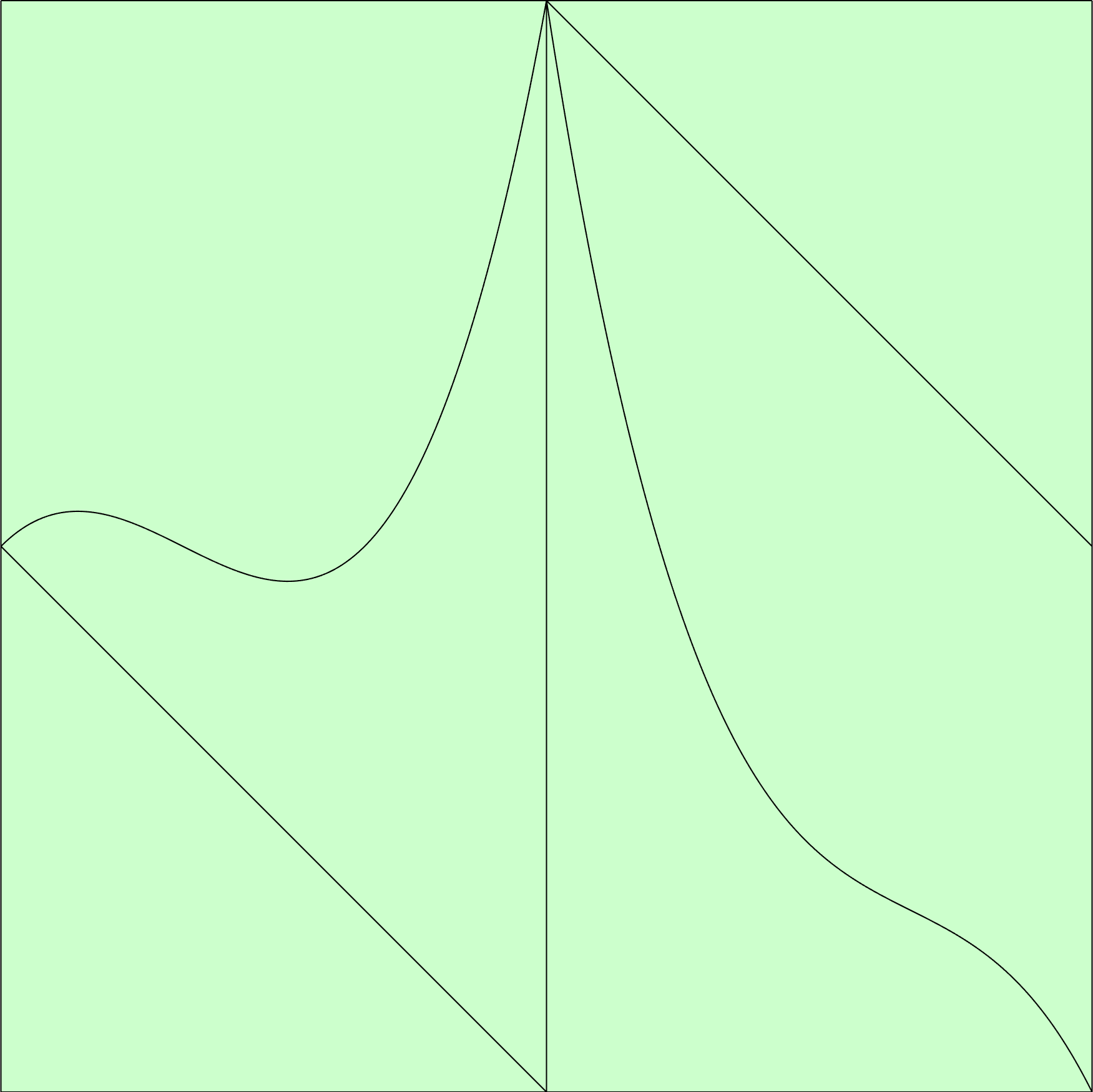};

\addplot [blue, mark options={solid, thin}, mark=*, mark size=1.5, only marks]
coordinates {
(-3.33333333e-01,  0.00000000e+00)
(-3.33333333e-01,  2.22222222e-01)
(-3.33333333e-01,  3.33333333e-01)
(-3.33333333e-01,  5.55555556e-01)
(-3.33333333e-01,  6.66666667e-01)
(-3.33333333e-01,  1.11111111e-01)
(-3.33333333e-01,  4.44444444e-01)
(-2.22222222e-01,  0.00000000e+00)
(-2.22222222e-01,  2.22222222e-01)
(-2.22222222e-01,  3.33333333e-01)
(-2.22222222e-01,  1.11111111e-01)
(-2.22222222e-01,  4.44444444e-01)
(-2.22222222e-01,  6.66666667e-01)
(-1.11111111e-01,  1.11111111e-01)
(-1.11111111e-01,  0.00000000e+00)
(-1.11111111e-01,  2.22222222e-01)
(-1.11111111e-01,  3.33333333e-01)
(-1.11111111e-01,  6.66666667e-01)
( 0.00000000e+00,  0.00000000e+00)
( 0.00000000e+00,  1.11111111e-01)
( 0.00000000e+00,  2.22222222e-01)
( 0.00000000e+00,  6.66666667e-01)
( 1.11111111e-01,  1.11111111e-01)
( 1.11111111e-01,  0.00000000e+00)
( 1.11111111e-01,  2.22222222e-01)
( 1.11111111e-01,  5.55555556e-01)
( 1.11111111e-01,  6.66666667e-01)
( 2.22222222e-01,  0.00000000e+00)
( 2.22222222e-01,  1.11111111e-01)
( 2.22222222e-01,  4.44444444e-01)
( 2.22222222e-01,  2.22222222e-01)
( 2.22222222e-01,  5.55555556e-01)
( 2.22222222e-01,  6.66666667e-01)
( 3.33333333e-01,  0.00000000e+00)
( 3.33333333e-01,  1.11111111e-01)
( 3.33333333e-01,  3.33333333e-01)
( 3.33333333e-01,  4.44444444e-01)
( 3.33333333e-01,  6.66666667e-01)
( 3.33333333e-01,  2.22222222e-01)
( 3.33333333e-01,  5.55555556e-01)};\label{line:smooth_ref:node}

\end{axis}
\end{tikzpicture} \quad
\begin{tikzpicture}
\begin{axis}[
axis equal image,
axis lines=none,
width=0.4\textwidth,
ymax=0.686666666667,
xmax=0.353333333333,
xmin=-0.353333333333,
ymin=-0.02]
\addplot [white, forget plot]
coordinates {
(-3.53333333e-01, -2.00000000e-02)
( 3.53333333e-01, -2.00000000e-02)
( 3.53333333e-01,  6.86666667e-01)
(-3.53333333e-01,  6.86666667e-01)
(-3.53333333e-01, -2.00000000e-02)};

\addplot []
graphics [xmin=-0.333333333333,xmax=0.333333333333,ymin=0.0,ymax=0.666666666667] { 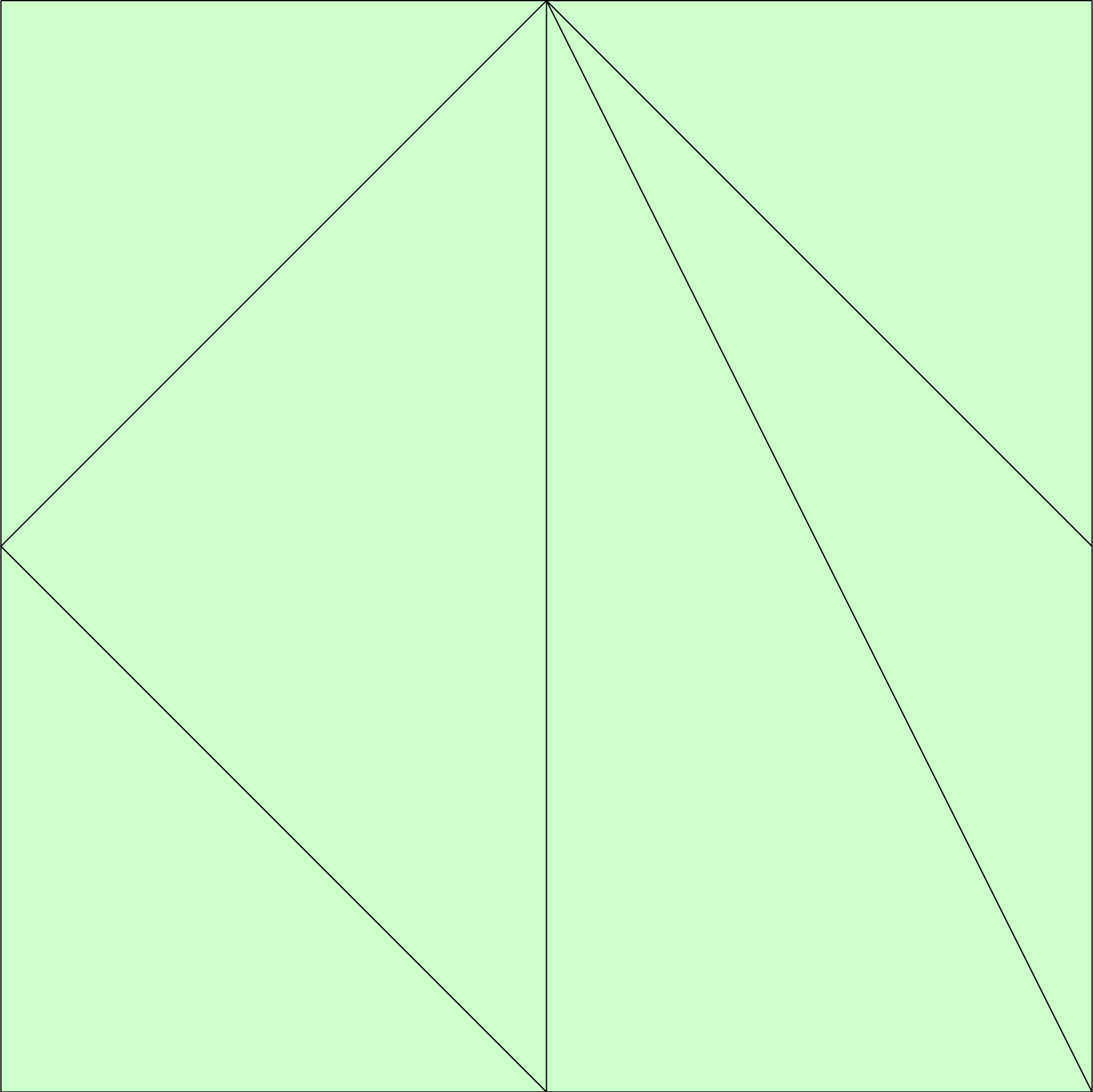};

\addplot [blue, mark options={solid, thin}, mark=*, mark size=1.5, only marks]
coordinates {
(-3.33333333e-01,  0.00000000e+00)
(-3.33333333e-01,  2.22222222e-01)
(-3.33333333e-01,  3.33333333e-01)
(-3.33333333e-01,  5.55555556e-01)
(-3.33333333e-01,  6.66666667e-01)
(-3.33333333e-01,  1.11111111e-01)
(-3.33333333e-01,  4.44444444e-01)
(-2.22222222e-01,  0.00000000e+00)
(-2.22222222e-01,  2.22222222e-01)
(-2.22222222e-01,  4.44444444e-01)
(-2.22222222e-01,  1.11111111e-01)
(-2.22222222e-01,  5.55555556e-01)
(-2.22222222e-01,  6.66666667e-01)
(-1.11111111e-01,  1.11111111e-01)
(-1.11111111e-01,  0.00000000e+00)
(-1.11111111e-01,  3.33333333e-01)
(-1.11111111e-01,  5.55555556e-01)
(-1.11111111e-01,  6.66666667e-01)
( 0.00000000e+00,  0.00000000e+00)
( 0.00000000e+00,  2.22222222e-01)
( 0.00000000e+00,  4.44444444e-01)
( 0.00000000e+00,  6.66666667e-01)
( 1.11111111e-01,  2.22222222e-01)
( 1.11111111e-01,  0.00000000e+00)
( 1.11111111e-01,  4.44444444e-01)
( 1.11111111e-01,  5.55555556e-01)
( 1.11111111e-01,  6.66666667e-01)
( 2.22222222e-01,  0.00000000e+00)
( 2.22222222e-01,  2.22222222e-01)
( 2.22222222e-01,  4.44444444e-01)
( 2.22222222e-01,  3.33333333e-01)
( 2.22222222e-01,  5.55555556e-01)
( 2.22222222e-01,  6.66666667e-01)
( 3.33333333e-01,  0.00000000e+00)
( 3.33333333e-01,  1.11111111e-01)
( 3.33333333e-01,  3.33333333e-01)
( 3.33333333e-01,  4.44444444e-01)
( 3.33333333e-01,  6.66666667e-01)
( 3.33333333e-01,  2.22222222e-01)
( 3.33333333e-01,  5.55555556e-01)};\label{line:smooth_ref:node}

\end{axis}
\end{tikzpicture}
\caption{An edge collapse (\ref{line:smooth_ref:edge}) in the reference
  mesh (\textit{left}) can lead to highly distorted, artificially curved
  elements if only the nodes (\ref{line:smooth_ref:node}) along the
  collapsed edge are moved (\textit{middle}). To avoid this situation,
  we reset all elements whose nodes were affected by the edge collapse
  to be straight-sided (\textit{right}).}
\label{fig:elem-smooth}
\end{figure}

\begin{remark}
Element removal constitutes a fundamental change to the
optimization problem as it alters the definition of the
objective function and constraints and changes the number
of optimization variables. Such operations are inconsistent
with the spirit of optimization solvers that aim to find
the solution of a given optimization problem and poses
the risk of never converging to a meaningful solution
if the underlying optimization problem is continually
changing. In practice, all element collapses occur in
the early optimization iterations meaning the optimization
problem is fixed in the asymptotic regime and this issue
does not arise.
\end{remark}

\subsubsection{Mesh straightening}
\label{sec:ist_solve:mod:smth}
For high-order meshes, it is not necessarily advantageous to remove
ill-conditioned elements because such elements can be large
(Figure~\ref{fig:mshop_demo3d}) so the associated edge collapse
could be a severe mesh operation. This is not an issue for $q=1$
(straight-sided) meshes because ill-conditioning is
akin to a short edge whose collapse will be a benign
mesh operation. Therefore, our approach is to only
remove elements if they are inverted (Figure~\ref{fig:mshop_demo3c}),
e.g., take $c_4 = 0$, and reset ill-conditioned elements to
be straight-sided (from Figure~\ref{fig:mshop_demo3d} to
Figure~\ref{fig:mshop_demo3e}).
We define ill-conditioned elements to be elements
in the set $\Ecal_\mathrm{ill}(\xbm)$,
where $\func{\Ecal_\mathrm{ill}}{\Rbb^{N_\xbm}}{\Scal_h}$
is defined as
\begin{equation}
 \Ecal_\mathrm{ill}: \xbm \mapsto \left\{\Omega_{0,e}\in\Ecal_h \suchthat g_{e,\mathrm{inf}}(\xbm) \leq c_4' g_{e,\mathrm{sup}}(\xbm) \right\}
\end{equation}
and $c_4'\in\Rbb_{\geq 0}$; $c_4' = 0.05$ in this work. In the SQP
setting, this is applied after the element collapses in
Section~\ref{sec:ist_solve:mod:clps}, i.e., we remove
the high-order information from the elements in
$\Ecal_\mathrm{ill}(\phibold(\hat\ybm_{k+1}))$ to produce
the new unconstrained degrees of freedom $\ybm_{k+1}$.

\subsubsection{Solution re-initialization}
\label{sec:ist_solve:mod:reinit}
\begin{figure}
 \centering
 \begin{subfigure}[b]{0.22\textwidth}
  \includegraphics[width=\textwidth]{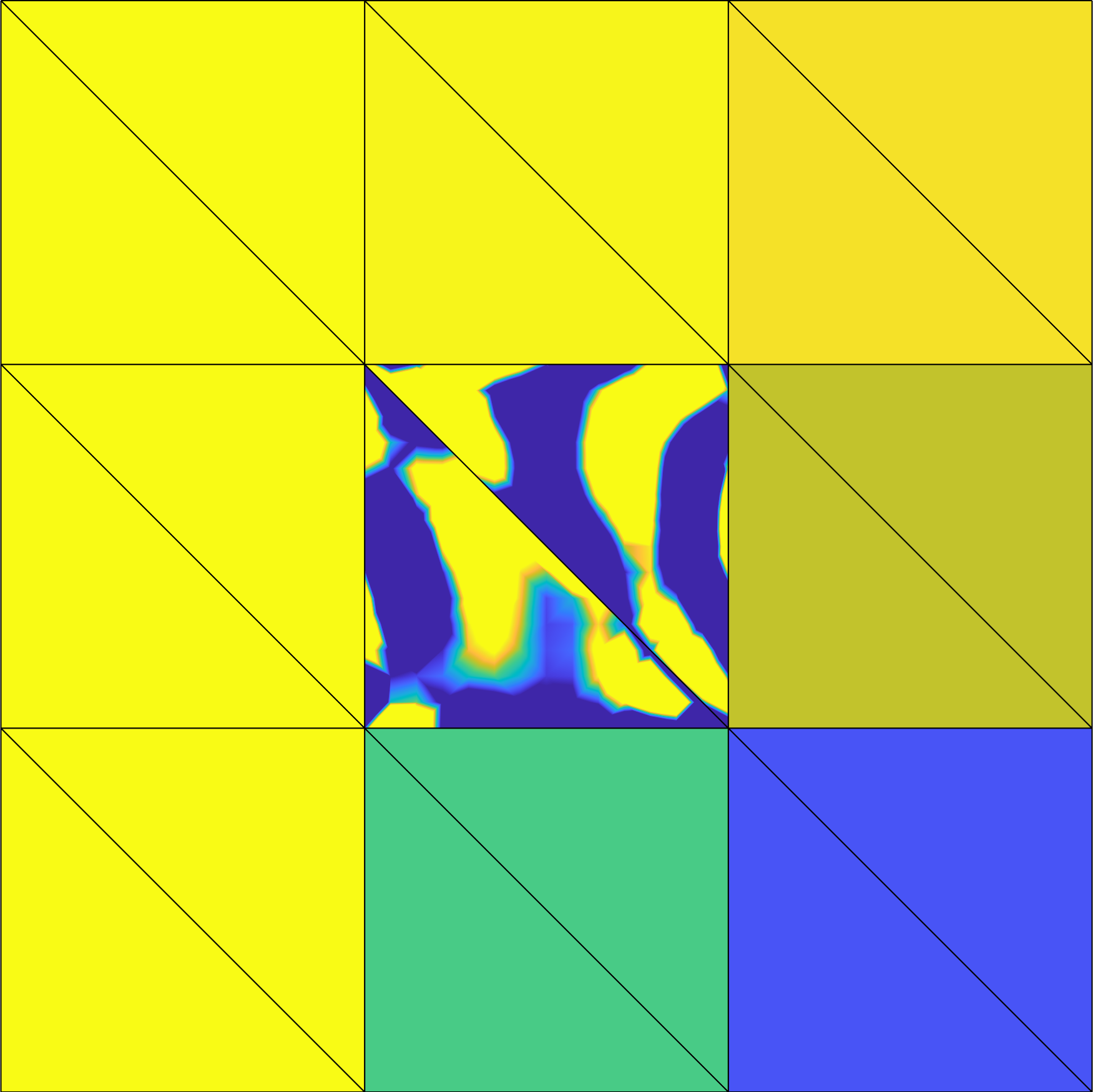}
  \caption{} \label{fig:reinita}
 \end{subfigure} \quad
 \begin{subfigure}[b]{0.22\textwidth}
  \includegraphics[width=\textwidth]{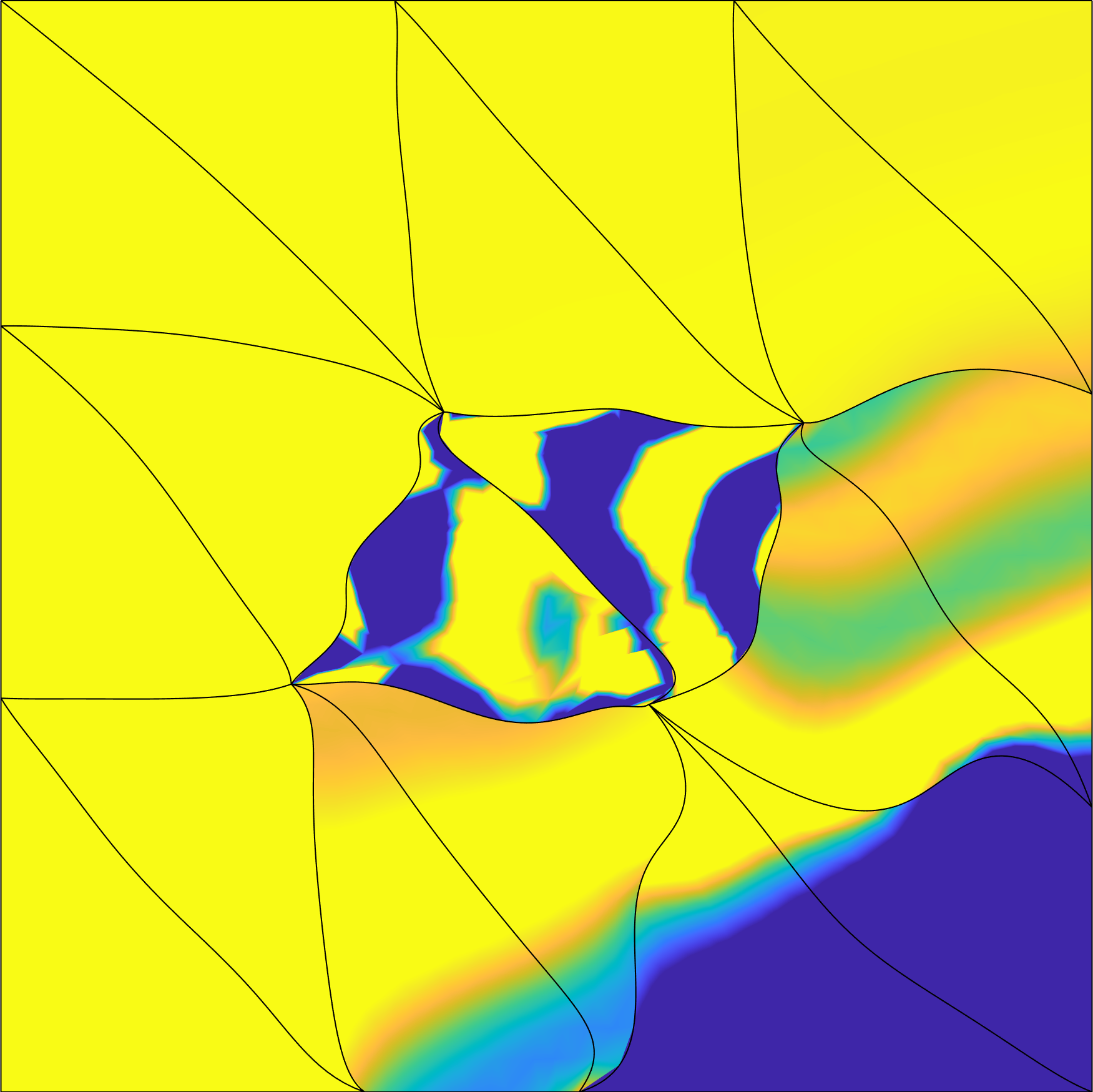}
  \caption{} \label{fig:reinitb}
 \end{subfigure} \quad
 \begin{subfigure}[b]{0.22\textwidth}
  \includegraphics[width=\textwidth]{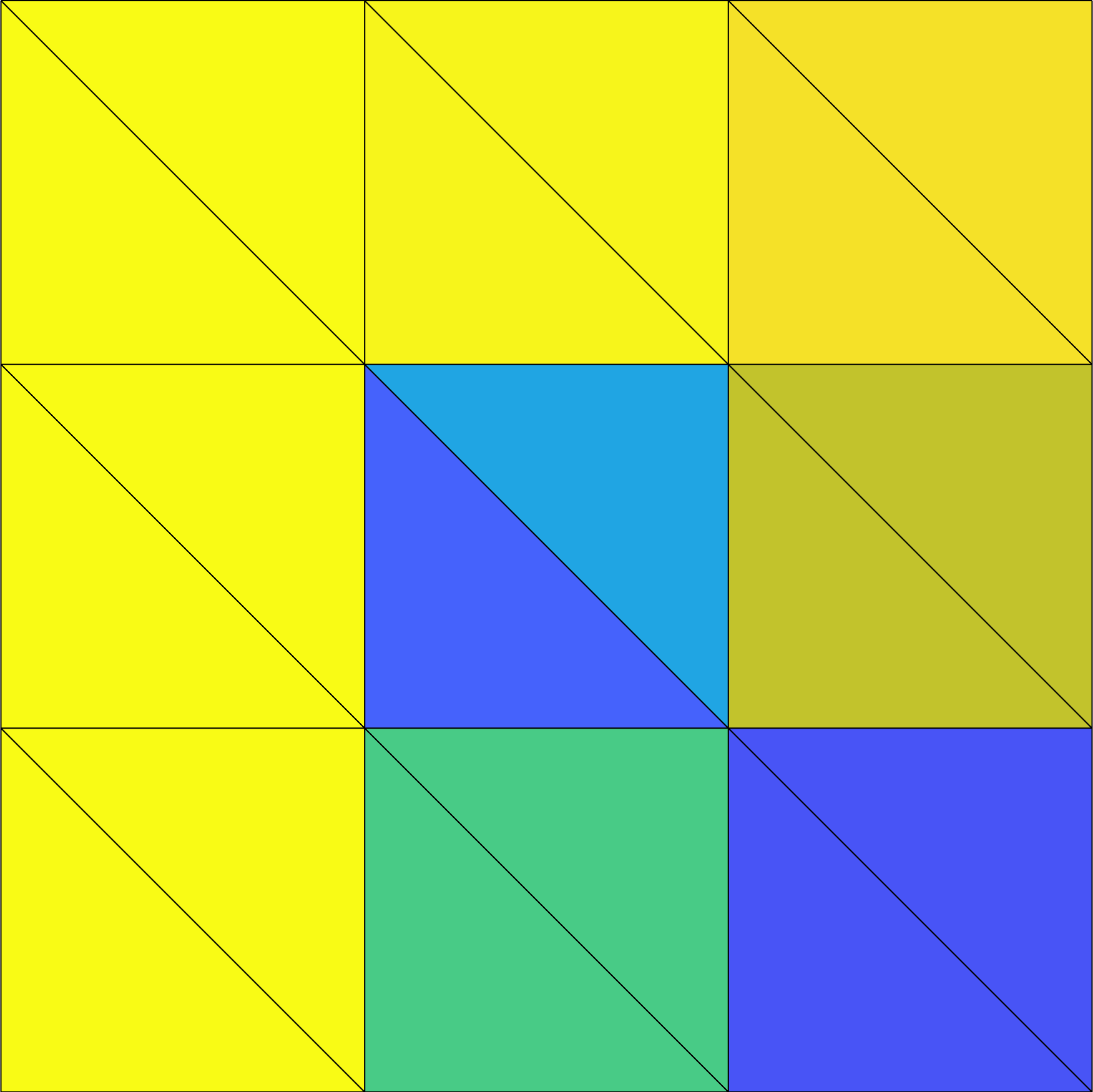}
  \caption{} \label{fig:reinitc}
 \end{subfigure} \quad
 \begin{subfigure}[b]{0.22\textwidth}
  \includegraphics[width=\textwidth]{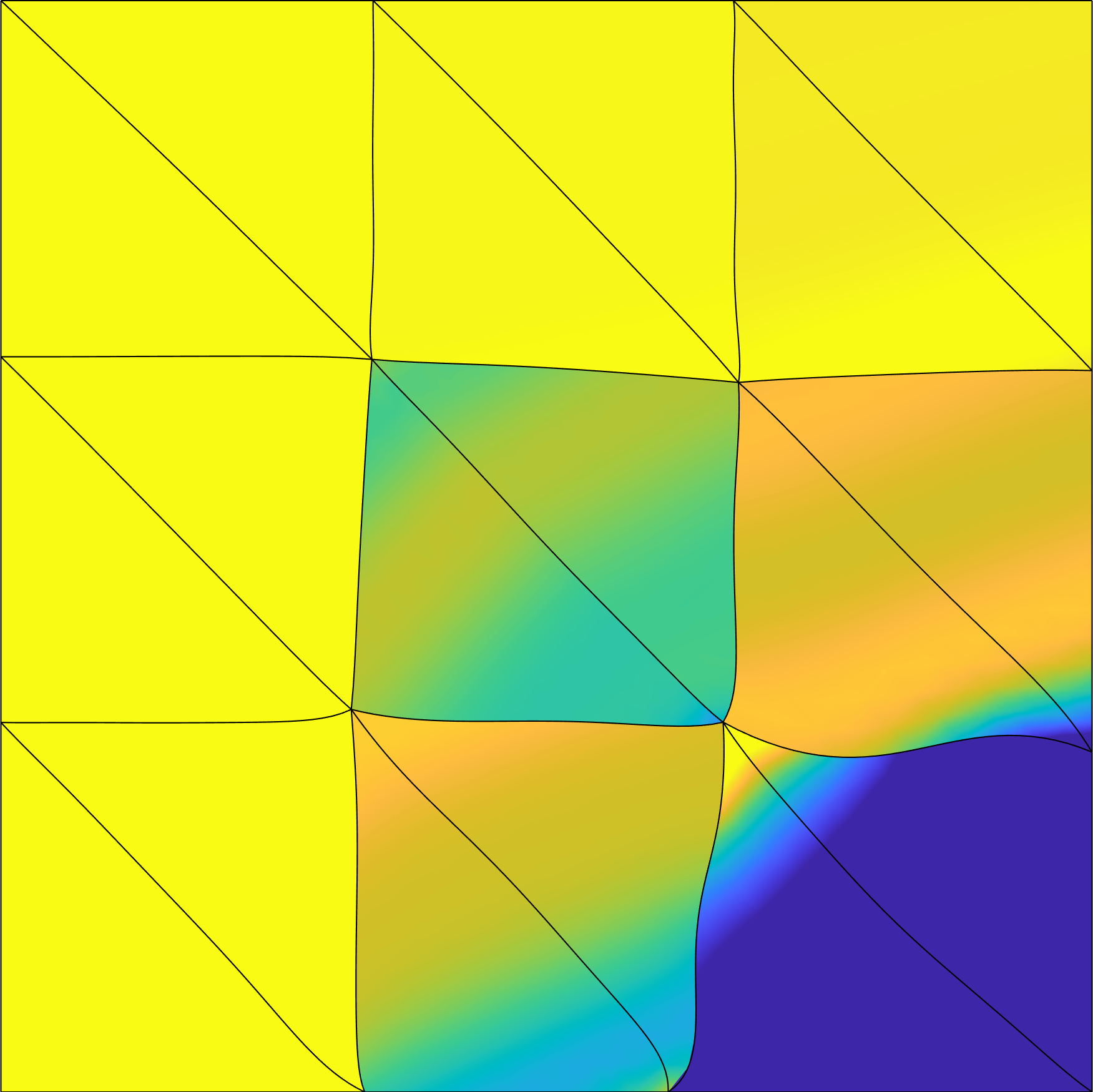}
  \caption{} \label{fig:reinitd}
 \end{subfigure}
 \caption{(a) An oscillatory intermediate shock tracking solution and
  mesh that can arise and (b) the corresponding SQP update (including
  line search). (c) The intermediate solution in (a) after re-initialization
  and (d) the corresponding SQP update (including line search). Oscillations
  in the DG solution lead to an SQP step that degrades the quality of the
  mesh and spreads the oscillations to neighboring elements, whereas
  the piecewise constant DG solution has a less severe effect.}
 \label{fig:reinit}
\end{figure}

We have observed that non-physical oscillations often arise in the DG
solution at intermediate SQP steps, which lead to poor SQP search
directions in subsequent iterations that require many line search
iterations and substantially degrade both the mesh quality and DG
solution. Numerical experiments have shown that resetting the DG
solution in elements where it is oscillatory to a constant value
promotes high-quality SQP search directions for both the mesh and DG solution
(illustration in Figure~\ref{fig:reinit}, full numerical experiment in
 Figure~\ref{fig:iburg:nrl0:reinit-demo}).
Our approach is to re-initialize the DG solution in an oscillatory
element, identified using the Persson-Peraire shock sensor, with the
average of the DG solution over a patch of elements in its vicinity.

\textit{Identification of oscillatory elements.}
To identify oscillatory elements, we turn to the shock capturing
literature and use the popular shock sensor of Persson and
Peraire \cite{persson_sub-cell_2006},
$\func{s_e}{\Vcal_h^p}{\Rbb}$ for $e=1,\dots,|\Ecal_h|$, defined as
\begin{equation}
 s_e : \bar{W}_h \mapsto \log_{10}\left(\sqrt{\frac{\int_{\Omega_{0,e}} (\chi(\bar{W}_h)-\chi(\Pi_{p-1}\bar{W}_h))^2 \, dV}{\int_{\Omega_{0,e}} \chi(\bar{W}_h)^2 \, dV}}\right),
\end{equation}
where $\func{\Pi_{p-1}}{\Vcal_h^p}{\Vcal_h^{p-1}}$ is a projection
onto the space of piecewise polynomials of total degree $p-1$ and
$\func{\chi}{\Rbb^m}{\Rbb}$ is a suitable scalar of the state. In this
work, we choose $\chi : (W_1,\dots,W_m) \mapsto W_1$. Then, we define the
set of oscillatory elements $\func{\Ecal_\mathrm{osc}}{\Rbb^{N_\ubm}}{\Scal_h}$
as all elements where the shock sensor exceeds a threshold, i.e.,
\begin{equation} \label{eqn:Eosc}
 \Ecal_\mathrm{osc} : \ubm \mapsto
 \left\{\Omega_{0,e}\in\Ecal_h \suchthat s_e(\Xi(\ubm)) \geq c_5 \right\},
\end{equation}
where $c_5\in\Rbb_{\geq 0}$. Numerical experiments have suggested it
is advantageous to re-initialize not only elements in
$\Ecal_\mathrm{osc}(\ubm)$ but also neighboring elements. 
To this end, define $\func{\Ncal}{\Scal_h}{\Scal_h}$ as
$\Ncal: A \mapsto \Ncal(A)$, where for any $K\in\Ncal(A)$,
there exists $K'\in A$ ($K\neq K'$) such that $K$ and
$K'$ share a face (i.e., neighboring elements). Then
the set of elements to be re-initialized,
$\func{\Ecal_\mathrm{reinit}}{\Rbb^{N_\ubm}}{\Scal_h}$,
is defined as
\begin{equation}
 \Ecal_\mathrm{reinit} : \ubm \mapsto
 \Ecal_\mathrm{osc}(\ubm) \cup
 \Ncal(\Ecal_\mathrm{osc}(\ubm)).
\end{equation}
In the case where 
the number of line search iterations is excessive
(defined as more than $5$ in this work), we force
re-initialization for a fraction of the most oscillatory
elements despite not meeting the criteria in (\ref{eqn:Eosc}).
That is, we re-define $\Ecal_\mathrm{reinit}$ as
\begin{equation}
 \Ecal_\mathrm{reinit} : \ubm \mapsto
 \left\{\Omega_{0,e}\in\Ecal_h \suchthat s_e(\Xi(\ubm)) \geq c_6 \max_{e'=1,\dots,|\Ecal_h|} s_{e'}(\Xi(\ubm)) \right\},
\end{equation}
where $c_6\in\Rbb_{\geq 0}$; in this work, we typically take $c_6 = 0.01$.

\textit{Re-initialization from patch of elements.} With the oscillatory
elements identified, we re-initialize them by resetting the DG solution
in each of these elements to a constant value. The constant value assigned
to element $K\in\Ecal_\mathrm{init}(\ubm)$ is given by the average of the
current DG solution ($\ubm$) over a patch of elements neighboring $K$
\textit{without crossing the shock}. To this end, define
$\func{a}{\Ecal_h\times\Ecal_h\times\Vcal_h^p}{\Rbb}$
as the average jump between two elements, i.e.,
\begin{equation}
 a : (K,K',\bar{W}_h) \mapsto
 \frac{1}{|\partial K\cap \partial K'|}\int_{\partial K\cap \partial K'} \jump{\chi(\bar{W}_h)} \, dS,
 \qquad |\partial K\cap \partial K'| \coloneqq \int_{\partial K\cap \partial K'}\,dS,
\end{equation}
where $\jump{\cdot}$ denotes the jump operator. Next, we define
$\func{\tilde\Ncal}{\Ecal_h\times\Rbb^{N_\ubm}}{\Scal_h}$ as
\begin{equation}
 \tilde\Ncal : (K, \ubm) \mapsto
 \left\{K'\in\Ncal(K) \suchthat |a(K,K',\Xi(\ubm))|\leq c_7\right\},
\end{equation}
where $c_7\in\Rbb_{\geq 0}$, i.e., for $K\in\Ecal_\mathrm{reinit}(\ubm)$,
$\tilde\Ncal(K,\ubm)\subset\Ncal(K)$ is the set of elements that neighbor
$K$ that do not cross the shock (Figure~\ref{fig:mshop_demo6}).
Finally, the DG solution in each element $K\in\Ecal_\mathrm{reinit}(\ubm)$
is reset to the constant value $\func{\zeta}{\Ecal_h\times\Rbb^{N_\ubm}}{\Rbb}$
given by
\begin{equation} \label{eqn:reinit0}
 \zeta : (K, \ubm) \mapsto
 \frac{1}{|K\cup\tilde{N}(K,\ubm)|}
 \int_{K\cup\tilde{N}(K,\ubm)} \Xi(\ubm) \, dV.
\end{equation}
In the SQP setting, elements are re-initialized after the mesh operations
in Section~\ref{sec:ist_solve:mod:clps}-\ref{sec:ist_solve:mod:smth}, i.e.,
the elements in $\Ecal_\mathrm{reinit}(\hat\ubm_{k+1})$ are re-initialized,
and $\ubm_{k+1}$ is defined as
\begin{equation} \label{eqn:reinit1}
 \left.\Xi(\ubm_{k+1})\right|_K =
 \begin{cases}
  \zeta(K,\hat\ubm_{k+1}) &\text{if } K \in \Ecal_\mathrm{reinit}(\hat\ubm_{k+1}) \\[0.8em]
  \left.\Xi(\hat\ubm_{k+1})\right|_K &\text{otherwise.}
 \end{cases}
\end{equation}
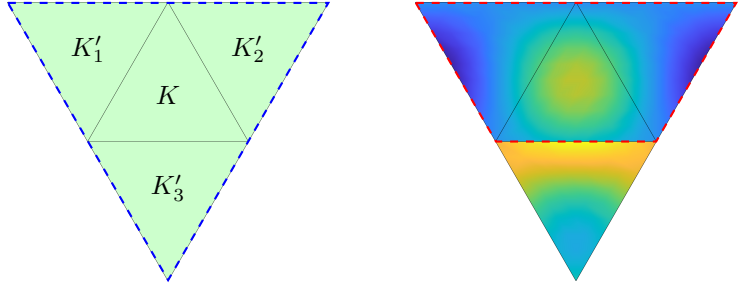
\begin{figure}
\centering
\begin{tikzpicture}
\begin{axis}[
axis equal image,
axis lines=none,
width=0.4\textwidth,
ymax=0.1,
xmax=2.1,
xmin=-0.1,
ymin=-1.8320508075688773]
\addplot [white, forget plot]
coordinates {
(-1.00000000e-01, -1.83205081e+00)
( 1.10000000e+00, -1.83205081e+00)
( 1.10000000e+00,  1.00000000e-01)
(-1.00000000e-01,  1.00000000e-01)
(-1.00000000e-01, -1.83205081e+00)};

\addplot []
graphics [xmin=0,xmax=2,ymin=-1.7320508075688772,ymax=0] { 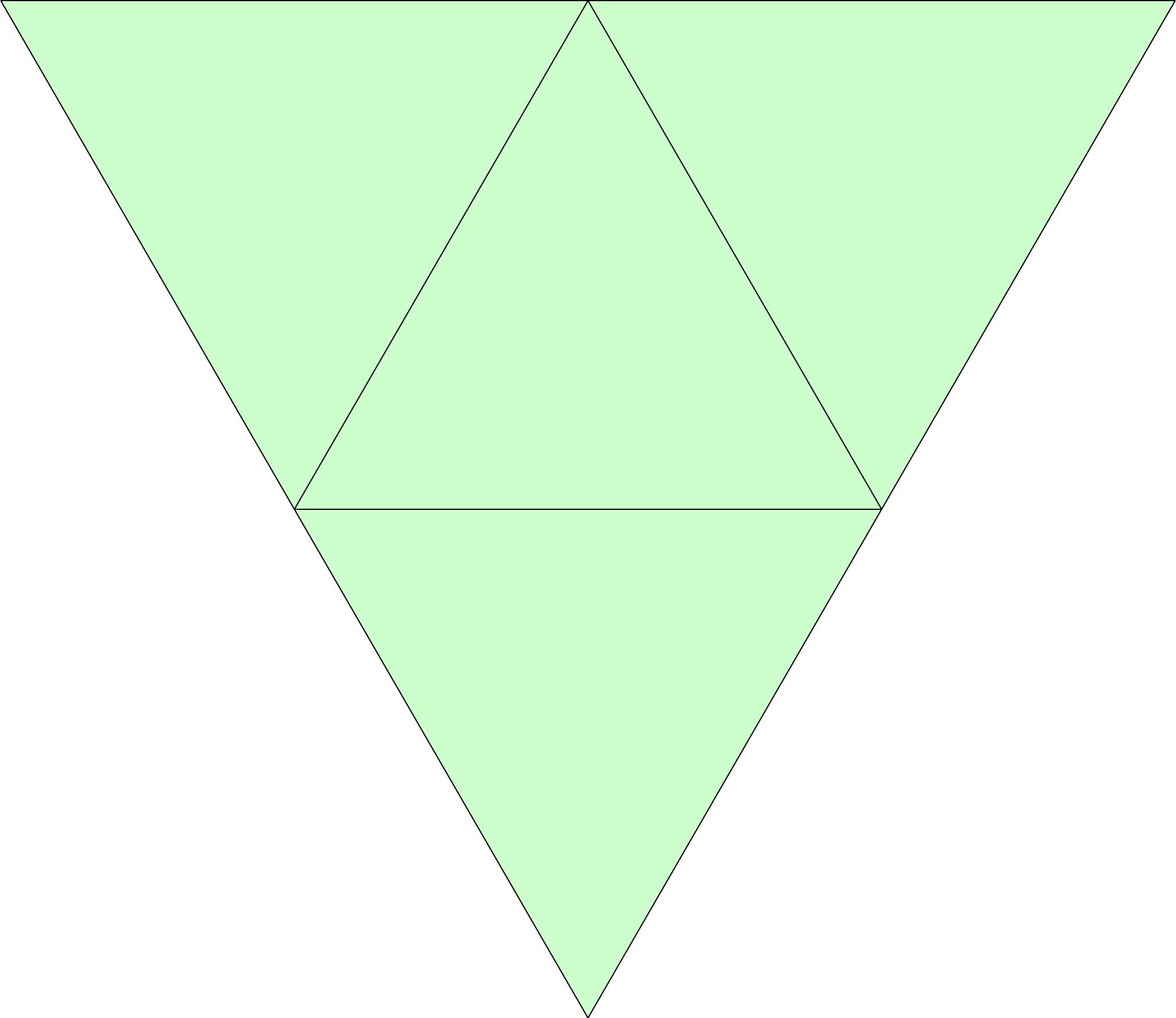};

\node[]    at    (axis cs:0.5, -0.28867513459481287) {$K'_1$};
\node[]    at    (axis cs:1.0, -0.5773502691896257) {$K$};
\node[]    at    (axis cs:1.5, -0.28867513459481287) {$K'_2$};
\node[]    at    (axis cs:1.0, -1.1547005383792515) {$K'_3$};
\addplot [blue, dashed, thick]
coordinates {
( 0.00000000e+00,  0.00000000e+00)
( 2.00000000e+00,  0.00000000e+00)
( 1.00000000e+00, -1.73205081e+00)
( 0.00000000e+00,  0.00000000e+00)};\label{line:mshop_demo6:neigh}

\end{axis}
\end{tikzpicture} \qquad
\begin{tikzpicture}
\begin{axis}[
axis equal image,
axis lines=none,
width=0.4\textwidth,
ymax=0.1,
xmax=2.1,
xmin=-0.1,
ymin=-1.8320508075688773]
\addplot [white, forget plot]
coordinates {
(-1.00000000e-01, -1.83205081e+00)
( 1.10000000e+00, -1.83205081e+00)
( 1.10000000e+00,  1.00000000e-01)
(-1.00000000e-01,  1.00000000e-01)
(-1.00000000e-01, -1.83205081e+00)};

\addplot []
graphics [xmin=0,xmax=2,ymin=-1.7320508075688772,ymax=0] { 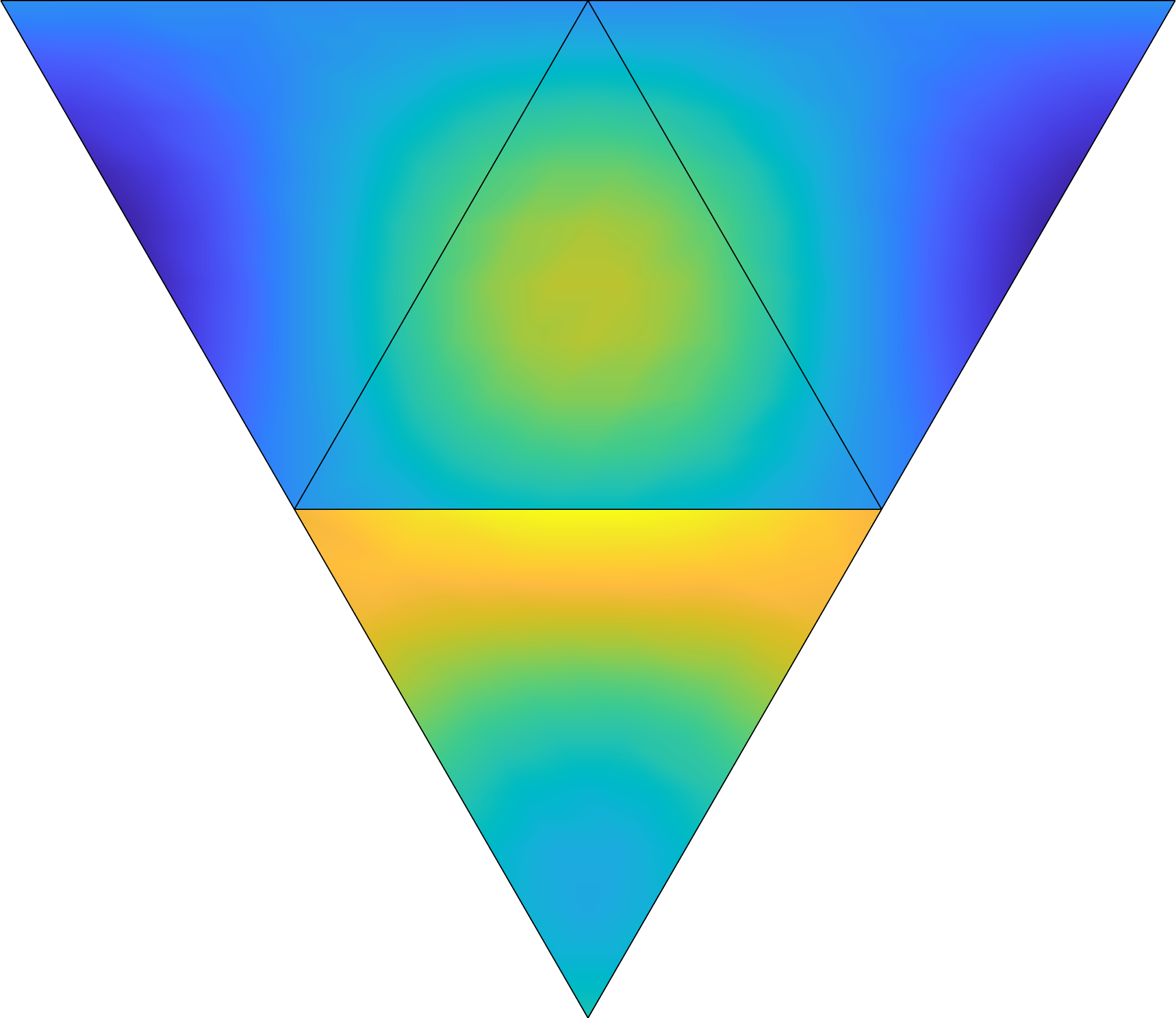};

\addplot [red, dashed, thick]
coordinates {
( 0.00000000e+00,  0.00000000e+00)
( 2.00000000e+00,  0.00000000e+00)
( 1.50000000e+00, -8.66025404e-01)
( 5.00000000e-01, -8.66025404e-01)
( 0.00000000e+00,  0.00000000e+00)};\label{line:mshop_demo6:neighshk}

\end{axis}
\end{tikzpicture}
\caption{A mesh consisting of element $K$ and its neighbors
 $\Ncal(K)=\{K_1',K_2',K_3'\}$ (\textit{left}) and a function
 $W_h\in\Vcal_h^p$ ($p=4$) with a strong jump between $K$ and
 $K_3'$ (\textit{right}), which leads to the shock-aware
 neighbors $\tilde\Ncal(K,W_h)=\{K_1',K_2'\}$. The set
 $K\cup\Ncal(K)$ is enclosed in (\ref{line:mshop_demo6:neigh})
 and the set $K\cup\tilde\Ncal(K)$ is enclosed in
 (\ref{line:mshop_demo6:neighshk}).}
\label{fig:mshop_demo6}
\end{figure}

\textit{Termination of solution re-initialization.}
To ensure element re-initialization does not continue indefinitely,
we only allow re-initialization for iterations $k$ such that
\begin{equation}
\norm{\rbm_k} > c_8, \qquad k \leq M,
\end{equation}
where $c_8\in\Rbb_{\geq 0}$ and $M$ is defined in
Section~\ref{sec:ist_solve:sqp:robust}.
The first condition prevents re-initialization when the constraints
are nearly satisfied because this condition usually only happens
near convergence where re-initialization of the DG solution would
interfere with clean and rapid convergence. This condition is
usually sufficient to ensure re-initialization terminates
after a finite number of iterations. The second condition
is a safeguard to ensure $\Upsilon_k = \mathrm{Id}$ for
$k > M$.

\begin{remark} \label{rem:reinit_smth}
After solution re-initialization, it can be advantageous to also
straighten the elements in $\Ecal_\mathrm{reinit}$, i.e., remove
high-order geometry information, as discussed in
Section~\ref{sec:ist_solve:mod:smth} to remove highly curved
edges from intermediate configurations. This combination is especially
helpful for problems with complex shock structures, but not necessary
for simpler problems.
\end{remark}

\begin{remark}
A major concern in the context of compressible flows is whether
oscillations will lead to catastrophic breakdown of the solution,
which usually manifests as negative pressures and densities
throughout the domain. While the oscillations that arise
during the optimization iterations adversely impact the
performance of the SQP solver, negative pressures and
densities do not arise at \emph{accepted} iterations
because the line search ensures sufficient decrease
of the merit function. That is, suppose iteration
$\zbm_k$ is free of negative pressures and densities
and recall $\zbm_{k+1}= \zbm_k + \alpha_k\Delta\zbm_k$.
The backtracking line search will choose $\alpha_k$
sufficiently small that $\zbm_{k+1}$ leads to sufficient
decrease of the merit function, which cannot be a state
with negative pressures and densities. Therefore, the
line search is sufficient to ensure the sequence
$\{\zbm_k\}_{k=0}^\infty$ is free of negative pressures
and densities; however, the re-initialization procedure
is necessary to avoid steps becoming prohibitively small
as the DG solution at intermediate iterations becomes
oscillatory.
\end{remark}

\begin{remark}
Both element removal or solution re-initialization are operations that
cause abrupt changes to the objective function and constraints, which
interferes with the value of $\kappa$ that has been calibrated as the
optimization iterations progress. Suppose either an element is removed
or the DG solution is re-initialized at iteration $k$, then we take
\begin{equation}
 \kappa_{k+1} =
 \max\left\{\upsilon\sqrt{\frac{f_\text{err}(\zbm_{k+1})}{f_\text{msh}(\ybm_{k+1})}},
       \kappa_\mathrm{min}\right\}
\end{equation}
to balance the tracking and mesh quality objectives.
\end{remark}


\subsection{SQP initialization}
\label{sec:ist_solve:init}
The HOIST optimization problem in (\ref{eqn:pde-opt}) is a nonlinear,
non-convex program, which makes the starting point for the
optimization solver important. Let $\Xbm\in\Rbb^{N_\xbm}$
be the nodal coordinates of the reference mesh $\Ecal_h$,
i.e., concatenation of $\{\hat{X}_I\}_{I=1}^{N_\mathrm{v}}$ into a
vector, which is produced through a mesh generation procedure
agnostic to the shock location. Then, we take the initial
physical mesh to be the reference mesh, $\xbm_0 = \Xbm$,
and the corresponding unconstrained degrees of freedom
are determined as
\begin{equation} \label{eqn:y0_planar}
 \ybm_0 = \Abm^\dagger(\xbm_0-\bbm),
\end{equation}
i.e., the least-squares fit to the physical mesh, in the
case where all domain boundaries are planar; otherwise,
we solve the nonlinear least-square problem:
$\ybm_0 = \argmin_{\ybm\in\Rbb^{N_\ybm}} \norm{\phibold(\ybm)-\xbm_0}$.
Because $\xbm_0 = \Xbm$ lies in the range of $\Abm$ by construction
(Section~\ref{sec:dommap:planar}), perfect inversion (to
the tolerances of the linear solver) of the parametrization
is guaranteed ($\xbm_0 = \Abm\ybm_0+\bbm$). The corresponding
DG solution is initialized from the first-order finite volume or
$p=0$ DG solution on the shock-agnostic initial mesh
$\Gcal_h(\Ecal_h;\xbm_0)$ because it is simple, fast,
and does not require stabilization. Numerical investigations
into more sophisticated initialization of the DG
solution---$L^2$ projection of the exact solution
and shock capturing solutions based on artificial
viscosity---did not offer noticeable advantages
over this simple choice.

\begin{remark}
The approach described in this section is always used
to initialize the SQP solver, regardless of polynomial
degrees ($p$ and $q$). This is in stark contrast to the
predecessor to this work \cite{zahr_implicit_2020} that
required a homotopy procedure: the strategy described in this
section was used to initialize a $p=q=1$ simulation, which was
in turn used to initialize $p>1$, $q>1$ simulations. The robustness
measures in Section~\ref{sec:ist_solve:mod} are the key
innovations of this work that enable reliable convergence
for complex flows from a simple, generic starting point (shock
agnostic mesh and corresponding first-order finite volume solution)
for any polynomial degrees.
\end{remark}


\subsection{Complete algorithm}
\label{sec:ist_solve:alg}
The complete HOIST method described in Section~\ref{sec:ist_solve} is
summarized in Algorithm~\ref{alg:hoist} for domains with planar boundaries.
The only generalization required for non-planar boundaries is the
definition of $\phibold$ for curved surfaces
\cite{corrigan_convergence_2019,zahr_radaptive_2020,zahr_implicit_2020}
and definition of $\ybm_0$ via nonlinear least-squares
(Section~\ref{sec:ist_solve:init}).
\begin{algorithm}
 \caption{HOIST method}
 \label{alg:hoist}
 \begin{algorithmic}[1]
  \REQUIRE Reference mesh $\Ecal_h$ with nodal coordinates $\Xbm$, the
  parameters in Table~\ref{tab:params}, and tolerances $\epsilon_1$, $\epsilon_2$
  \ENSURE Shock-aligned mesh $\phibold(\ybm^\star)$ and corresponding
  DG solution $\ubm^\star$
  \STATE \textbf{Construct parametrization}: Construct $\Abm$ and $\bbm$
   according to Section~\ref{sec:dommap:planar} and define $\phibold$
   from (\ref{eqn:dommap_affine})
  \STATE \textbf{Solution initialization}: Define $\ybm_0$ from
   (\ref{eqn:y0_planar}) and $\ubm_0$ as the DG($p=0$) solution
   on the reference mesh
  \FOR{$k=0,1,2,\dots$}
   \STATE \textbf{SQP search direction}: Compute $\Delta\zbm_k$ from
    (\ref{eqn:sqp_sys0})
   \STATE \textbf{Line search}: Determine $\alpha_k$ that satisfies
    (\ref{eqn:steep_descent}) by backtracking
   \STATE \textbf{SQP update}: Compute $\tilde\zbm_{k+1}=\zbm_k+\alpha_k\Delta\zbm_k$
   \STATE \textbf{Check convergence}: Terminate iterations with
    $(\ubm^\star,\ybm^\star)=\tilde\zbm_{k+1}$ if
    $\norm{\rbm(\tilde\zbm_{k+1})}\leq\epsilon_1$,
    $\norm{\cbm(\tilde\zbm_{k+1})}\leq\epsilon_2$
   \STATE \textbf{Modify SQP step}: Compute $\zbm_{k+1}=\Upsilon_{k+1}(\tilde\zbm_{k+1})$
   \begin{itemize}
    \item \textbf{Remove elements} in set $\Ecal_\mathrm{rmv}(\phibold(\tilde\ybm_{k+1}))$
    \begin{itemize}
     \item Fix connectivity to define new mesh ($\Ecal_h$) with physical
         nodes $\hat\xbm_{k+1}$
     \item Reconstruct parametrization $\phibold$ and unconstrained
         degrees of freedom $\hat\ybm_{k+1} = \Abm^\dagger(\hat\xbm_{k+1}-\bbm)$
     \item Transfer solution to new mesh $\hat\ubm_{k+1}$
    \end{itemize}
    \item \textbf{Straighten elements} in set
     $\Ecal_\mathrm{ill}(\phibold(\hat\ybm_{k+1}))$ to yield
     $\xbm_{k+1}$ and $\ybm_{k+1} = \Abm^\dagger(\xbm_{k+1}-\bbm)$
    \item \textbf{Re-initialize solution} in elements in set
     $\Ecal_\mathrm{reinit}(\hat\ubm_{k+1})$ to yield $\ubm_{k+1}$
     from (\ref{eqn:reinit1}), provided $\norm{\rbm(\tilde\zbm_{k+1})}>c_8$
     and $k\leq M$
   \end{itemize}
  \ENDFOR
\end{algorithmic}
\end{algorithm}

\section{Numerical experiments}
\label{sec:numexp}
In this section, we demonstrate the HOIST framework on a sequence of
increasingly difficult inviscid conservation laws with discontinuous
solutions in $d=1,2,3$. In particular, we show the HOIST method achieves
optimal $\Ocal(h^{p+1})$ convergence rates and the proposed solver is
able to robustly track complex discontinuity surfaces (e.g., contact
discontinuities, curved and reflecting shocks, shock formation, and
shock-shock interaction) as well as surfaces along which the solution
is continuous but its derivative is discontinuous (e.g., head and tail
of rarefactions).
\subsection{Linear advection}
\label{sec:numexp:advec}
Consider steady, linear advection of a scalar quantity
$\func{U}{\Omega}{\Rbb}$ through a domain $\Omega\subset\Rbb^d$
\begin{equation} \label{eqn:advec}
 \nabla \cdot (U(x)\beta(x)^T) = 0 \quad \text{for}~x\in\Omega, \qquad
 U(x) = U_\infty(x) \quad \text{for}~x\in\Gamma_\mathrm{in},
\end{equation}
where $\func{\beta}{\Omega}{\Rbb^d}$ is the local flow direction,
$\Gamma_\mathrm{in}\coloneqq\{x\in\partial\Omega\mid\beta(x)\cdot n(x)<0\}$
is the inflow boundary, and $\func{U_\infty}{\Omega}{\Rbb}$ is
the inflow boundary condition. For the DG discretization, we use
the smoothed upwind flux described in
\cite{zahr_implicit_2020,zahr_high-order_2020}
as the inviscid numerical flux function. 
\subsubsection{Straight shock, piecewise constant solution}
\label{sec:numexp:advec:planar}
First, we consider a three-dimensional problem ($d=3$) with
a planar shock surface and piecewise constant solution
(test case: \texttt{advec-planar}) to demonstrate
the HOIST method in a simple three-dimensional setting.
Consider the domain $\Omega\coloneqq(-1, 1)\times(0, 1)\times (0, 1)$
with constant advection field and piecewise constant boundary condition
\begin{equation}
 \beta: x \mapsto (-1.25, 1, 0), \qquad
 U_\infty: x \mapsto H(x_1),
\end{equation}
where $\func{H}{\Rbb}{\{0,1\}}$ is the Heaviside function; this is the
three-dimensional extrusion of the problem considered in \cite{zahr_implicit_2020}.
We label the six planar boundaries $\{\partial\Omega_i\}_{i=1}^6$, where
$\partial\Omega_i$ has outward unit normal $\eta_i=-e_i$ for $i=1,2,3$
and $\partial\Omega_{i+d}$ has outward unit normal $\eta_i=e_i$ for
$i=1,2,3$. The choice of advection field implies
$\Gamma_\mathrm{in}=\partial\Omega_2\cup\partial\Omega_4$.
To ensure the Heaviside function is accurately represented and
integrated in the weak form along boundary $\partial\Omega_2$,
we require the reference mesh to have element edges aligned
with the line $\{(0, 0, s) | s\in[0,1]\}$
and do not allow the corresponding nodes to move throughout
the HOIST iterations (Figure~\ref{fig:advec3d_geom}).
\ifbool{fastcompile}{}{
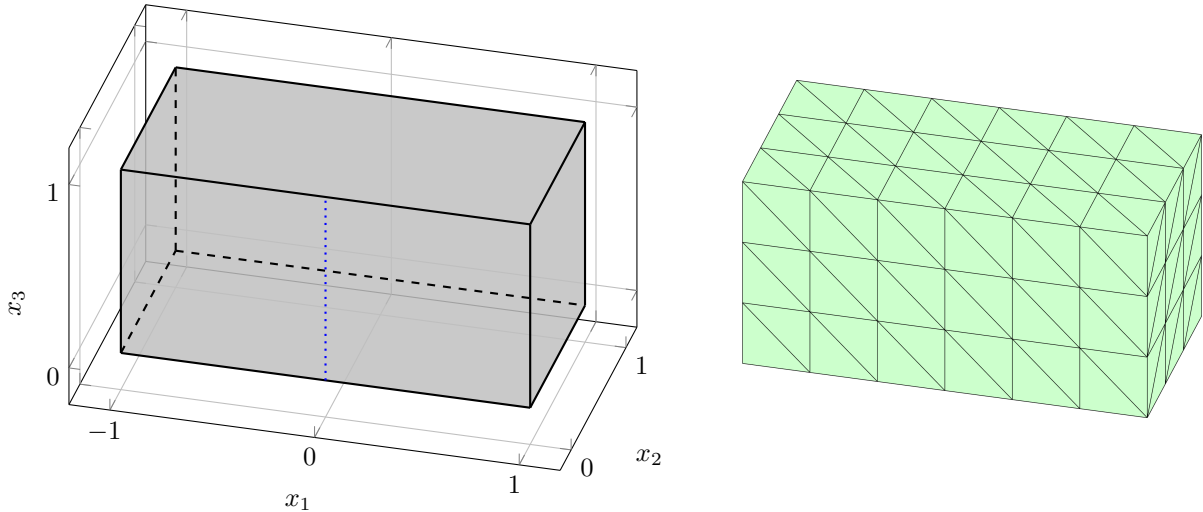
\begin{figure}[!htbp]
\centering
\begin{tikzpicture}
\begin{axis}[
axis equal image,
grid,
zmax=1.2,
ymax=1.2,
zlabel={$x_3$},
ztick={0, 1},
width=0.7\textwidth,
xtick={-1, 0, 1},
ytick={0, 1},
xlabel={$x_1$},
xmax=1.2,
ylabel={$x_2$},
xmin=-1.2,
ymin=-0.2,
zmin=-0.2,
view={15}{30}]
\addplot3 [opacity=0.6, draw=none, fill=lightgray, forget plot]
coordinates {
(-1.00000000e+00,  0.00000000e+00,  0.00000000e+00)
( 1.00000000e+00,  0.00000000e+00,  0.00000000e+00)
( 1.00000000e+00,  1.00000000e+00,  0.00000000e+00)
(-1.00000000e+00,  1.00000000e+00,  0.00000000e+00)
(-1.00000000e+00,  0.00000000e+00,  0.00000000e+00)};

\addplot3 [opacity=0.6, draw=none, fill=lightgray, forget plot]
coordinates {
(-1.00000000e+00,  0.00000000e+00,  1.00000000e+00)
( 1.00000000e+00,  0.00000000e+00,  1.00000000e+00)
( 1.00000000e+00,  1.00000000e+00,  1.00000000e+00)
(-1.00000000e+00,  1.00000000e+00,  1.00000000e+00)
(-1.00000000e+00,  0.00000000e+00,  1.00000000e+00)};

\addplot3 [opacity=0.6, draw=none, fill=lightgray, forget plot]
coordinates {
(-1.00000000e+00,  0.00000000e+00,  0.00000000e+00)
( 1.00000000e+00,  0.00000000e+00,  0.00000000e+00)
( 1.00000000e+00,  0.00000000e+00,  1.00000000e+00)
(-1.00000000e+00,  0.00000000e+00,  1.00000000e+00)
(-1.00000000e+00,  0.00000000e+00,  0.00000000e+00)};

\addplot3 [opacity=0.6, draw=none, fill=lightgray, forget plot]
coordinates {
(-1.00000000e+00,  1.00000000e+00,  0.00000000e+00)
( 1.00000000e+00,  1.00000000e+00,  0.00000000e+00)
( 1.00000000e+00,  1.00000000e+00,  1.00000000e+00)
(-1.00000000e+00,  1.00000000e+00,  1.00000000e+00)
(-1.00000000e+00,  1.00000000e+00,  0.00000000e+00)};

\addplot3 [opacity=0.6, draw=none, fill=lightgray, forget plot]
coordinates {
(-1.00000000e+00,  0.00000000e+00,  0.00000000e+00)
(-1.00000000e+00,  1.00000000e+00,  0.00000000e+00)
(-1.00000000e+00,  1.00000000e+00,  1.00000000e+00)
(-1.00000000e+00,  0.00000000e+00,  1.00000000e+00)
(-1.00000000e+00,  0.00000000e+00,  0.00000000e+00)};

\addplot3 [opacity=0.6, draw=none, fill=lightgray, forget plot]
coordinates {
( 1.00000000e+00,  0.00000000e+00,  0.00000000e+00)
( 1.00000000e+00,  1.00000000e+00,  0.00000000e+00)
( 1.00000000e+00,  1.00000000e+00,  1.00000000e+00)
( 1.00000000e+00,  0.00000000e+00,  1.00000000e+00)
( 1.00000000e+00,  0.00000000e+00,  0.00000000e+00)};

\addplot3 [black, solid, thick, forget plot]
coordinates {
(-1.00000000e+00,  0.00000000e+00,  0.00000000e+00)
( 1.00000000e+00,  0.00000000e+00,  0.00000000e+00)};

\addplot3 [black, dashed, thick, forget plot]
coordinates {
(-1.00000000e+00,  1.00000000e+00,  0.00000000e+00)
( 1.00000000e+00,  1.00000000e+00,  0.00000000e+00)};

\addplot3 [black, solid, thick, forget plot]
coordinates {
(-1.00000000e+00,  0.00000000e+00,  1.00000000e+00)
( 1.00000000e+00,  0.00000000e+00,  1.00000000e+00)};

\addplot3 [black, solid, thick, forget plot]
coordinates {
(-1.00000000e+00,  1.00000000e+00,  1.00000000e+00)
( 1.00000000e+00,  1.00000000e+00,  1.00000000e+00)};

\addplot3 [black, dashed, thick, forget plot]
coordinates {
(-1.00000000e+00,  0.00000000e+00,  0.00000000e+00)
(-1.00000000e+00,  1.00000000e+00,  0.00000000e+00)};

\addplot3 [black, solid, thick, forget plot]
coordinates {
( 1.00000000e+00,  0.00000000e+00,  0.00000000e+00)
( 1.00000000e+00,  1.00000000e+00,  0.00000000e+00)};

\addplot3 [black, solid, thick, forget plot]
coordinates {
(-1.00000000e+00,  0.00000000e+00,  1.00000000e+00)
(-1.00000000e+00,  1.00000000e+00,  1.00000000e+00)};

\addplot3 [black, solid, thick, forget plot]
coordinates {
( 1.00000000e+00,  0.00000000e+00,  1.00000000e+00)
( 1.00000000e+00,  1.00000000e+00,  1.00000000e+00)};

\addplot3 [black, solid, thick, forget plot]
coordinates {
(-1.00000000e+00,  0.00000000e+00,  0.00000000e+00)
(-1.00000000e+00,  0.00000000e+00,  1.00000000e+00)};

\addplot3 [black, solid, thick, forget plot]
coordinates {
( 1.00000000e+00,  0.00000000e+00,  0.00000000e+00)
( 1.00000000e+00,  0.00000000e+00,  1.00000000e+00)};

\addplot3 [black, dashed, thick, forget plot]
coordinates {
(-1.00000000e+00,  1.00000000e+00,  0.00000000e+00)
(-1.00000000e+00,  1.00000000e+00,  1.00000000e+00)};

\addplot3 [black, solid, thick, forget plot]
coordinates {
( 1.00000000e+00,  1.00000000e+00,  0.00000000e+00)
( 1.00000000e+00,  1.00000000e+00,  1.00000000e+00)};

\addplot3 [blue, dotted, thick]
coordinates {
( 0.00000000e+00,  0.00000000e+00,  0.00000000e+00)
( 0.00000000e+00,  0.00000000e+00,  1.00000000e+00)};\label{line:advec_geom:fix}

\end{axis}
\end{tikzpicture} \qquad
\raisebox{0.3\height}{\includegraphics[width=0.37\textwidth]{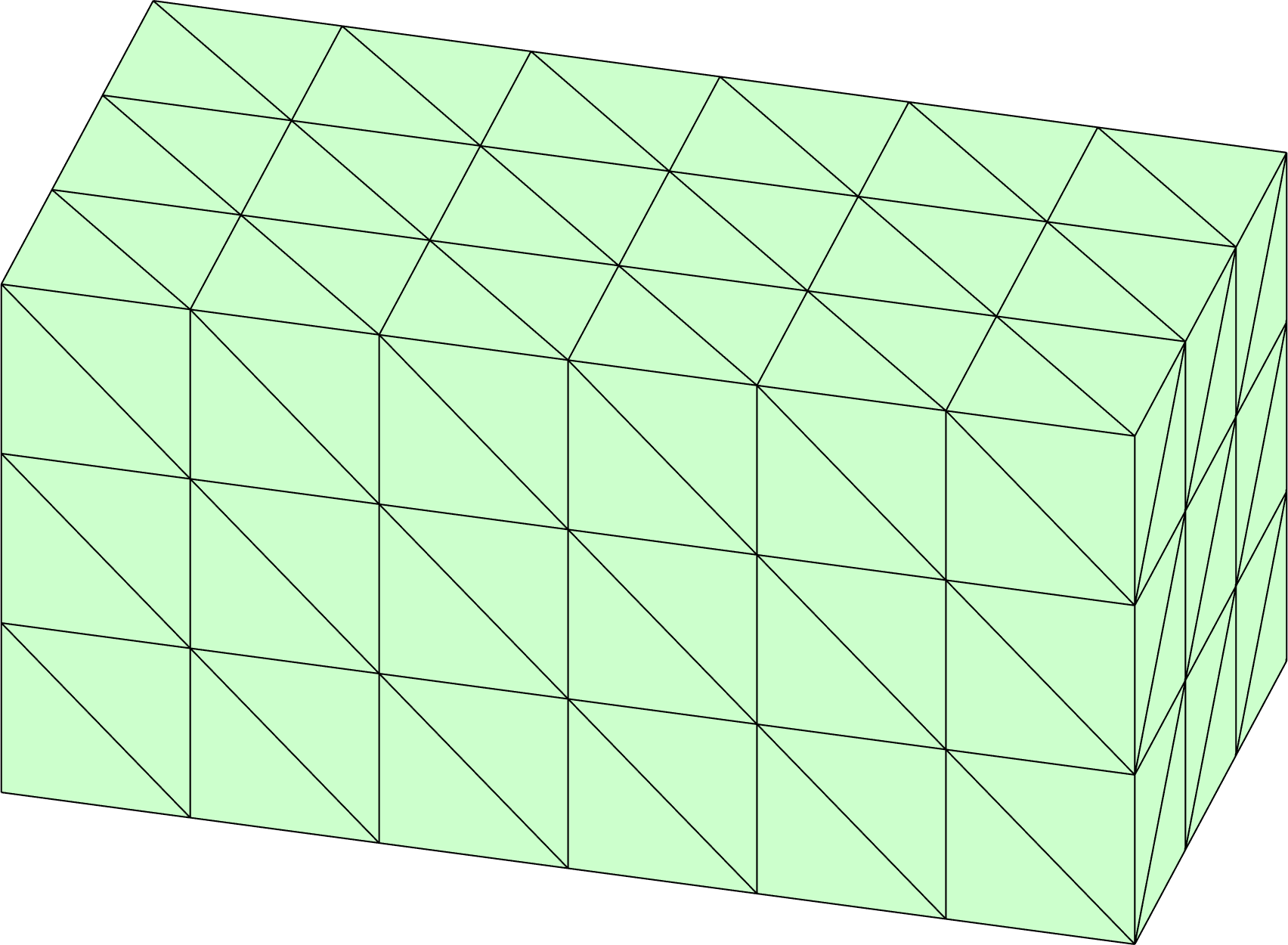}}
\caption{The domain for the advection experiments in
Sections~\ref{sec:numexp:advec:planar} and 
\ref{sec:numexp:advec:curved}; the reference mesh is
constructed such that element edges lie along (\ref{line:advec_geom:fix})
and the corresponding nodes are constrained to slide along
(\ref{line:advec_geom:fix}) to ensure the Heaviside function is
integrated accurately.}
\label{fig:advec3d_geom}
\end{figure}
}

We discretize the domain using a structured mesh consisting
of 324 uniform right tetrahedra.
Because the analytical solution is piecewise constant with a
straight-sided discontinuity surface, we choose a $p=0$ solution
and $q=1$ mesh approximation.
The HOIST parameters are defined in Table~\ref{tab:params}.
Because of the piecewise constant
($p=0$) DG solution, we do not perform re-initialization
and the mesh quality term in the objective function is
not necessary (the enriched residual induces minute
mesh motion away from the shock).

\begin{figure}[!htbp]
\centering
\ifbool{fastcompile}{}{
\includegraphics[width=0.35\textwidth]{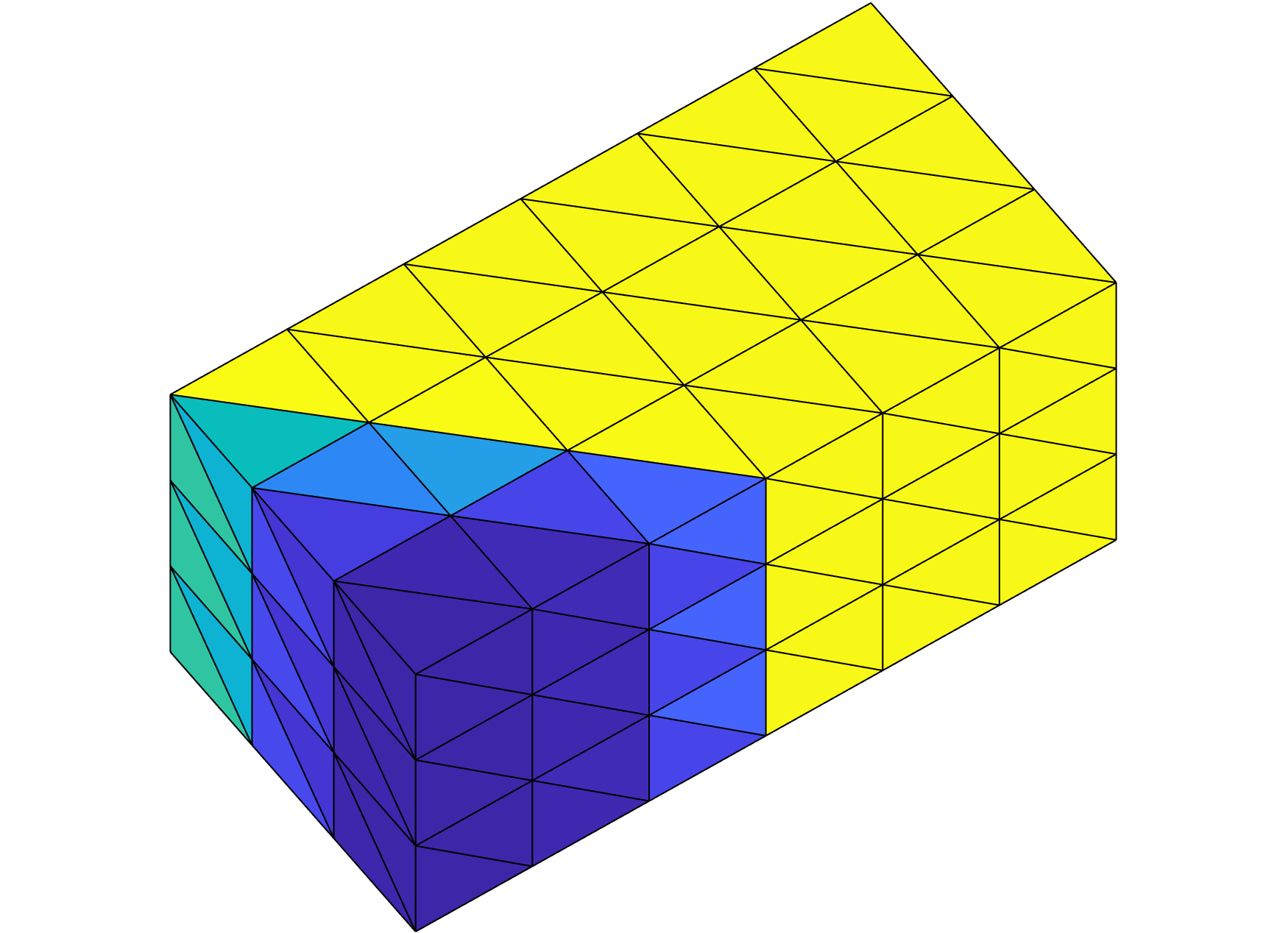}
\includegraphics[width=0.35\textwidth]{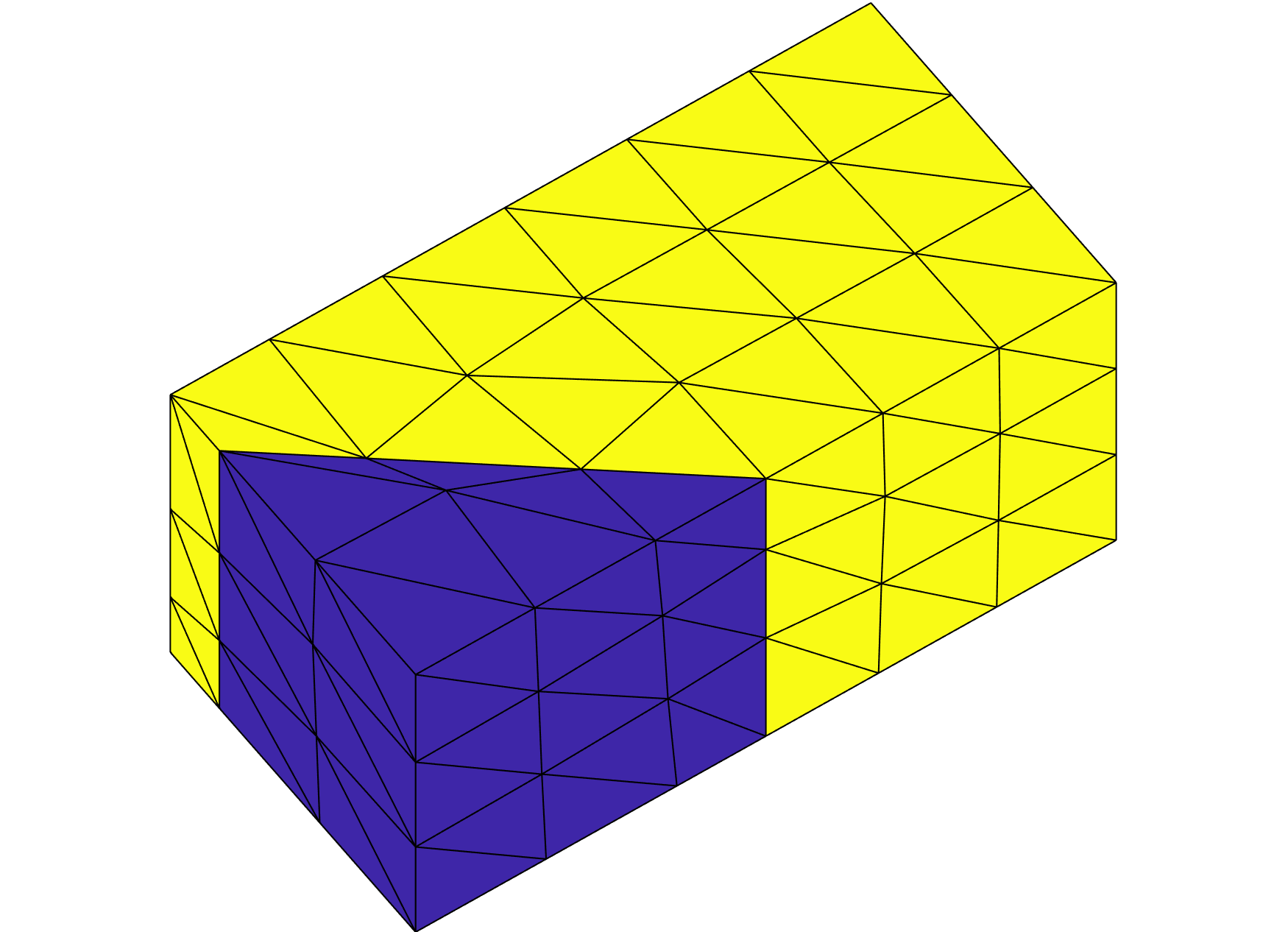}
\includegraphics[width=0.25\textwidth]{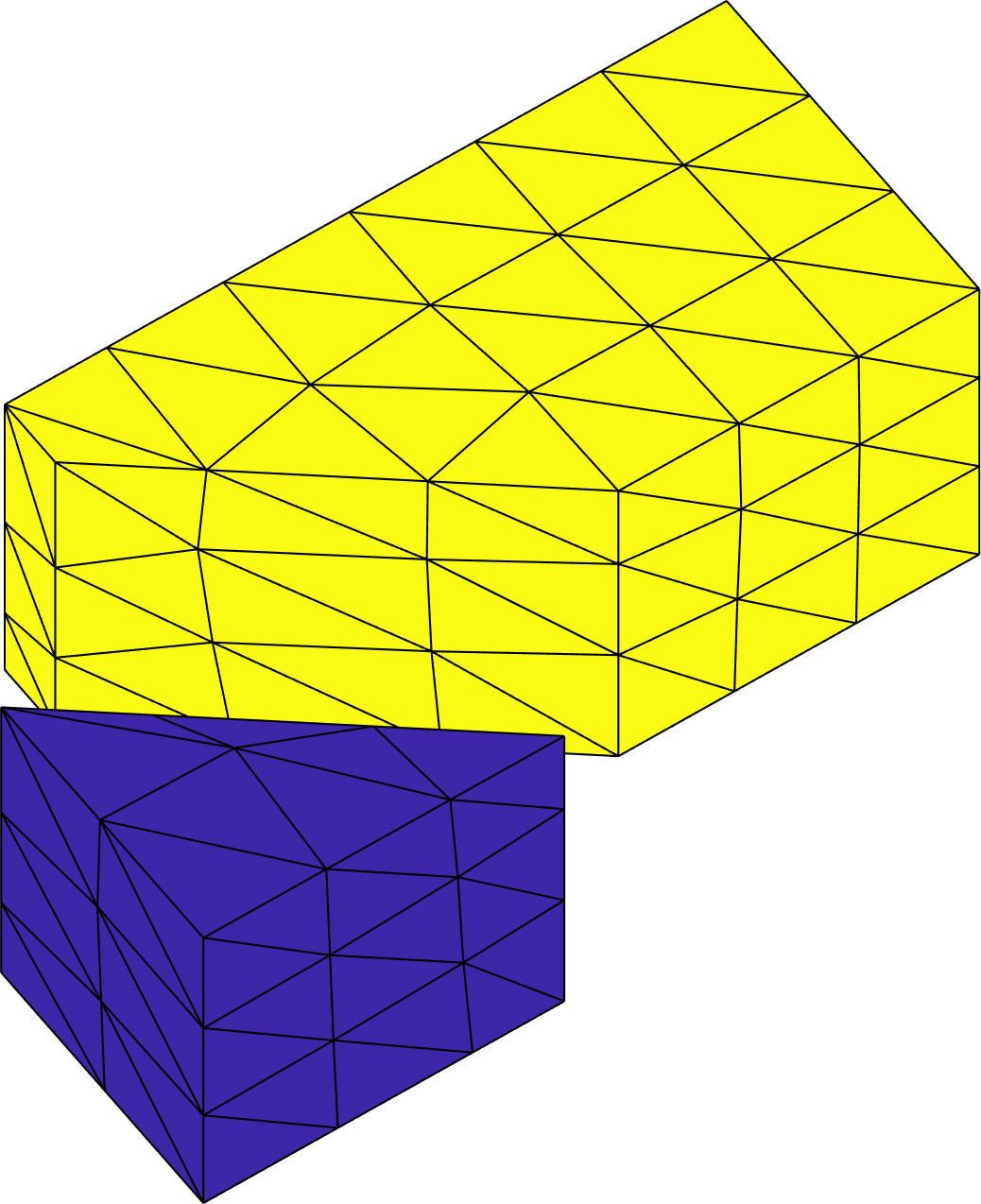}
}

\colorbarMatlabParula{0}{0.25}{0.5}{0.75}{1}
\caption{The starting point (\textit{left}) and the tracked configuration
 (\textit{middle}) for the \texttt{advec-planar} test case.
 The mesh is separate along the discontinuity surface
 for a clearer view of the tracked mesh (\textit{right}).}
\label{fig:advec3d-case0}
\end{figure}

\begin{figure}[!htbp]
\centering
\ifbool{fastcompile}{}{
\begin{tikzpicture}
\begin{groupplot} [
group style={group size = 3 by 1, horizontal sep = 1.4cm, vertical sep = 1cm}]
\nextgroupplot[width=0.33\textwidth, xtick={0,15,30}, xlabel={Iteration ($k$)}, ymax=30, xmax=30, xmin=-1, ymode=log, ymin=5e-11]
\addplot [black, thick, mark options={solid, thin}, mark=*, mark size=1.5, mark repeat=4]
coordinates {
( 0.00000000e+00,  5.62288107e-08)
( 1.00000000e+00,  2.18111680e-07)
( 2.00000000e+00,  8.21686071e-07)
( 3.00000000e+00,  2.93433331e-06)
( 4.00000000e+00,  9.55297623e-06)
( 5.00000000e+00,  2.68221448e-05)
( 6.00000000e+00,  6.03826022e-05)
( 7.00000000e+00,  1.01145922e-04)
( 8.00000000e+00,  1.34020182e-04)
( 9.00000000e+00,  6.31797541e-05)
( 1.00000000e+01,  5.22174139e-05)
( 1.10000000e+01,  6.08093379e-05)
( 1.20000000e+01,  7.88871801e-05)
( 1.30000000e+01,  1.12250023e-04)
( 1.40000000e+01,  1.75841917e-04)
( 1.50000000e+01,  3.08558972e-04)
( 1.60000000e+01,  6.13659030e-04)
( 1.70000000e+01,  1.17348224e-03)
( 1.80000000e+01,  3.71649424e-03)
( 1.90000000e+01,  2.69748244e-03)
( 2.00000000e+01,  2.77444964e-03)
( 2.10000000e+01,  1.12462079e-03)
( 2.20000000e+01,  2.20351564e-03)
( 2.30000000e+01,  4.01227298e-04)
( 2.40000000e+01,  7.87545324e-05)
( 2.50000000e+01,  8.92363240e-06)
( 2.60000000e+01,  9.56466174e-07)
( 2.70000000e+01,  3.72044232e-08)
( 2.80000000e+01,  4.62696015e-10)};\label{line:R0}

\addplot [red, thick, mark options={solid, thin}, mark=triangle*, mark size=1.5, mark repeat=4]
coordinates {
( 0.00000000e+00,  1.00000000e+00)
( 1.00000000e+00,  9.96978726e-01)
( 2.00000000e+00,  9.91015104e-01)
( 3.00000000e+00,  9.79392648e-01)
( 4.00000000e+00,  9.57269799e-01)
( 5.00000000e+00,  9.16905021e-01)
( 6.00000000e+00,  8.48437199e-01)
( 7.00000000e+00,  7.45873536e-01)
( 8.00000000e+00,  6.22373419e-01)
( 9.00000000e+00,  5.22876684e-01)
( 1.00000000e+01,  4.87039331e-01)
( 1.10000000e+01,  4.85449532e-01)
( 1.20000000e+01,  5.03602963e-01)
( 1.30000000e+01,  5.37760532e-01)
( 1.40000000e+01,  5.88403729e-01)
( 1.50000000e+01,  6.58425548e-01)
( 1.60000000e+01,  7.53044858e-01)
( 1.70000000e+01,  8.74515779e-01)
( 1.80000000e+01,  9.53298621e-01)
( 1.90000000e+01,  1.06214352e+00)
( 2.00000000e+01,  7.61769460e-01)
( 2.10000000e+01,  5.35502722e-01)
( 2.20000000e+01,  5.01258161e-01)
( 2.30000000e+01,  3.86457310e-01)
( 2.40000000e+01,  8.53287898e-02)
( 2.50000000e+01,  1.82193203e-02)
( 2.60000000e+01,  3.89955185e-03)
( 2.70000000e+01,  6.33869350e-04)
( 2.80000000e+01,  6.96945801e-05)};\label{line:dLdY}

\nextgroupplot[width=0.33\textwidth, xtick={0,15,30}, xlabel={Iteration ($k$)}, ymax=0.2, xmax=30, xmin=-1, ymode=log, ymin=8e-09]
\addplot [blue, thick, mark options={solid, thin}, mark=square*, mark size=1.5, mark repeat=4]
coordinates {
( 0.00000000e+00,  1.97759580e-02)
( 1.00000000e+00,  1.97571038e-02)
( 2.00000000e+00,  1.97204736e-02)
( 3.00000000e+00,  1.96512169e-02)
( 4.00000000e+00,  1.95266600e-02)
( 5.00000000e+00,  1.93212142e-02)
( 6.00000000e+00,  1.90245237e-02)
( 7.00000000e+00,  1.86543665e-02)
( 8.00000000e+00,  1.81928743e-02)
( 9.00000000e+00,  1.78013842e-02)
( 1.00000000e+01,  1.74091429e-02)
( 1.10000000e+01,  1.69887086e-02)
( 1.20000000e+01,  1.65197324e-02)
( 1.30000000e+01,  1.59745048e-02)
( 1.40000000e+01,  1.53090737e-02)
( 1.50000000e+01,  1.44486231e-02)
( 1.60000000e+01,  1.32699333e-02)
( 1.70000000e+01,  1.16890364e-02)
( 1.80000000e+01,  9.95908926e-03)
( 1.90000000e+01,  7.78509868e-03)
( 2.00000000e+01,  6.24444603e-03)
( 2.10000000e+01,  4.62963390e-03)
( 2.20000000e+01,  2.28158357e-03)
( 2.30000000e+01,  6.33457405e-04)
( 2.40000000e+01,  2.02485368e-04)
( 2.50000000e+01,  5.46902240e-05)
( 2.60000000e+01,  9.82711351e-06)
( 2.70000000e+01,  1.06377564e-06)
( 2.80000000e+01,  8.30996480e-08)};\label{line:R1err}

\addplot [magenta, thick, mark options={solid, thin}, mark=pentagon*, mark size=1.5, mark repeat=4]
coordinates {
( 0.00000000e+00,  2.15445249e-04)
( 1.00000000e+00,  2.15448798e-04)
( 2.00000000e+00,  2.15455851e-04)
( 3.00000000e+00,  2.15469768e-04)
( 4.00000000e+00,  2.15496820e-04)
( 5.00000000e+00,  2.15547878e-04)
( 6.00000000e+00,  2.15639747e-04)
( 7.00000000e+00,  2.15798553e-04)
( 8.00000000e+00,  2.16092681e-04)
( 9.00000000e+00,  2.16426541e-04)
( 1.00000000e+01,  2.16813134e-04)
( 1.10000000e+01,  2.17264890e-04)
( 1.20000000e+01,  2.17791446e-04)
( 1.30000000e+01,  2.18401313e-04)
( 1.40000000e+01,  2.19106001e-04)
( 1.50000000e+01,  2.19923325e-04)
( 1.60000000e+01,  2.20872616e-04)
( 1.70000000e+01,  2.21925540e-04)
( 1.80000000e+01,  2.22937258e-04)
( 1.90000000e+01,  2.23781170e-04)
( 2.00000000e+01,  2.24274235e-04)
( 2.10000000e+01,  2.24714577e-04)
( 2.20000000e+01,  2.25160544e-04)
( 2.30000000e+01,  2.25297366e-04)
( 2.40000000e+01,  2.25335042e-04)
( 2.50000000e+01,  2.25343550e-04)
( 2.60000000e+01,  2.25344061e-04)
( 2.70000000e+01,  2.25341022e-04)
( 2.80000000e+01,  2.25334277e-04)};\label{line:R1msh}

\nextgroupplot[width=0.33\textwidth, xtick={0,15,30}, xlabel={Iteration ($k$)}, xmax=30, xmin=-1, ymode=log]
\addplot [black, solid]
coordinates {
( 0.00000000e+00,  0.00000000e+00)
( 1.00000000e+00,  1.00000000e+00)
( 2.00000000e+00,  1.00000000e+00)
( 3.00000000e+00,  1.00000000e+00)
( 4.00000000e+00,  1.00000000e+00)
( 5.00000000e+00,  1.00000000e+00)
( 6.00000000e+00,  1.00000000e+00)
( 7.00000000e+00,  1.00000000e+00)
( 8.00000000e+00,  1.00000000e+00)
( 9.00000000e+00,  1.00000000e+00)
( 1.00000000e+01,  1.00000000e+00)
( 1.10000000e+01,  1.00000000e+00)
( 1.20000000e+01,  1.00000000e+00)
( 1.30000000e+01,  1.00000000e+00)
( 1.40000000e+01,  1.00000000e+00)
( 1.50000000e+01,  1.00000000e+00)
( 1.60000000e+01,  1.00000000e+00)
( 1.70000000e+01,  1.00000000e+00)
( 1.80000000e+01,  1.00000000e+00)
( 1.90000000e+01,  1.00000000e+00)
( 2.00000000e+01,  1.00000000e+00)
( 2.10000000e+01,  1.00000000e+00)
( 2.20000000e+01,  1.00000000e+00)
( 2.30000000e+01,  1.00000000e+00)
( 2.40000000e+01,  1.00000000e+00)
( 2.50000000e+01,  1.00000000e+00)
( 2.60000000e+01,  1.00000000e+00)
( 2.70000000e+01,  1.00000000e+00)
( 2.80000000e+01,  1.00000000e+00)};\label{line:alpha}

\addplot [blue, dashed]
coordinates {
( 0.00000000e+00,  1.00000000e-02)
( 1.00000000e+00,  1.00000000e-02)
( 2.00000000e+00,  5.00000000e-03)
( 3.00000000e+00,  2.50000000e-03)
( 4.00000000e+00,  1.25000000e-03)
( 5.00000000e+00,  6.25000000e-04)
( 6.00000000e+00,  3.12500000e-04)
( 7.00000000e+00,  1.56250000e-04)
( 8.00000000e+00,  7.81250000e-05)
( 9.00000000e+00,  3.90625000e-05)
( 1.00000000e+01,  3.90625000e-05)
( 1.10000000e+01,  3.90625000e-05)
( 1.20000000e+01,  3.90625000e-05)
( 1.30000000e+01,  3.90625000e-05)
( 1.40000000e+01,  3.90625000e-05)
( 1.50000000e+01,  3.90625000e-05)
( 1.60000000e+01,  3.90625000e-05)
( 1.70000000e+01,  3.90625000e-05)
( 1.80000000e+01,  3.90625000e-05)
( 1.90000000e+01,  3.90625000e-05)
( 2.00000000e+01,  3.90625000e-05)
( 2.10000000e+01,  3.90625000e-05)
( 2.20000000e+01,  3.90625000e-05)
( 2.30000000e+01,  3.90625000e-05)
( 2.40000000e+01,  3.90625000e-05)
( 2.50000000e+01,  1.95312500e-05)
( 2.60000000e+01,  9.76562500e-06)
( 2.70000000e+01,  4.88281250e-06)
( 2.80000000e+01,  2.44140625e-06)};\label{line:lam}

\addplot [red, dash dot]
coordinates {
( 0.00000000e+00,  1.00000000e-10)
( 1.00000000e+00,  1.00000000e-10)
( 2.00000000e+00,  1.00000000e-10)
( 3.00000000e+00,  1.00000000e-10)
( 4.00000000e+00,  1.00000000e-10)
( 5.00000000e+00,  1.00000000e-10)
( 6.00000000e+00,  1.00000000e-10)
( 7.00000000e+00,  1.00000000e-10)
( 8.00000000e+00,  1.00000000e-10)
( 9.00000000e+00,  1.00000000e-10)
( 1.00000000e+01,  1.00000000e-10)
( 1.10000000e+01,  1.00000000e-10)
( 1.20000000e+01,  1.00000000e-10)
( 1.30000000e+01,  1.00000000e-10)
( 1.40000000e+01,  1.00000000e-10)
( 1.50000000e+01,  1.00000000e-10)
( 1.60000000e+01,  1.00000000e-10)
( 1.70000000e+01,  1.00000000e-10)
( 1.80000000e+01,  1.00000000e-10)
( 1.90000000e+01,  1.00000000e-10)
( 2.00000000e+01,  1.00000000e-10)
( 2.10000000e+01,  1.00000000e-10)
( 2.20000000e+01,  1.00000000e-10)
( 2.30000000e+01,  1.00000000e-10)
( 2.40000000e+01,  1.00000000e-10)
( 2.50000000e+01,  1.00000000e-10)
( 2.60000000e+01,  1.00000000e-10)
( 2.70000000e+01,  1.00000000e-10)
( 2.80000000e+01,  1.00000000e-10)};\label{line:kappa}

\end{groupplot}\end{tikzpicture}
}
\caption{SQP convergence history and behavior of adaptive parameters for
 the \texttt{advec-planar} test case (legend in Table~\ref{tab:legend}).}
\label{fig:advec3d_case0_hist}
\end{figure}
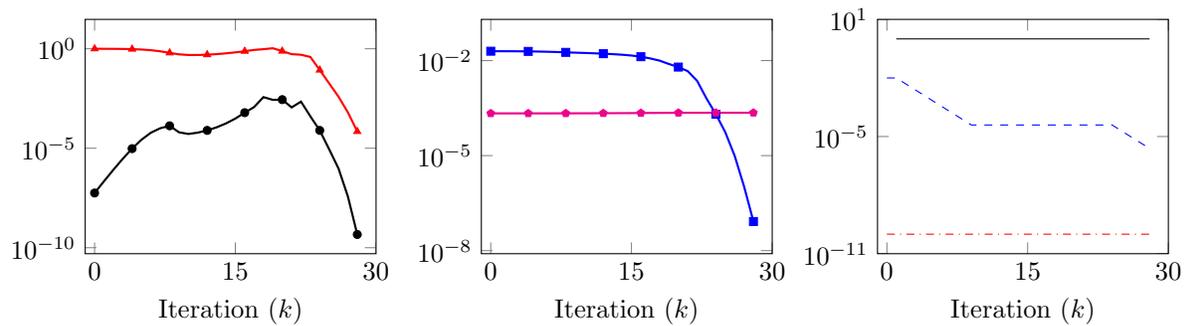

Starting from a shock-agnostic (structured) mesh and
first-order finite volume solution, only 29 SQP iterations
are required to track the discontinuity surface with
the mesh (Figure~\ref{fig:advec3d-case0}) and reduce the
standard DG residual to a magnitude of $\Ocal(10^{-10})$
and the enriched residual to $\Ocal(10^{-8})$ with no elements
being removed during the solution procedure. The SQP convergence
history and behavior of the adaptive parameters are shown
in Figure~\ref{fig:advec3d_case0_hist}. The mesh tracks the
shock nearly perfectly at iteration 23 and exhibits
Newton-like convergence after that point
(Figure~\ref{fig:advec3d_case0_hist}).

\begin{table}
\centering
\caption{Legend for SQP convergence history and parameter adaption plots,
 including the DG residual ($\norm{\rbm(\ubm_k,\phibold(\ybm_k))}$),
 enriched DG residual ($\norm{\Rbm(\ubm_k,\phibold(\ybm_k))}$),
 optimality condition ($\norm{\cbm(\ubm_k,\ybm_k)}$) (\ref{line:dLdY}),
 mesh distortion ($\norm{\kappa_k\Rbm_\text{msh}(\phibold(\ybm_k))}$),
 step size ($\alpha_k$,
 mesh regularization parameter ($\gamma_k$),
 and mesh distortion parameter ($\kappa_k$).}
\label{tab:legend}
\begin{tabular}{c|c|c|c|c|c|c}
$\norm{\rbm(\ubm_k,\phibold(\ybm_k))}$ &
$\norm{\Rbm(\ubm_k,\phibold(\ybm_k))}$ &
$\norm{\cbm(\ubm_k,\ybm_k)}$ &
$\norm{\kappa_k\Rbm_\text{msh}(\phibold(\ybm_k))}$ &
$\alpha_k$ & $\gamma_k$ & $\kappa_k$ \\\hline
(\ref{line:R0}) & (\ref{line:R1err}) & (\ref{line:dLdY}) &
(\ref{line:R1msh}) & (\ref{line:alpha}) & (\ref{line:lam}) &
(\ref{line:kappa})
\end{tabular}
\end{table}
\subsubsection{Curved shock, piecewise constant solution}
\label{sec:numexp:advec:curved}
Next, we consider a three-dimensional problem ($d=3$) with
a curved (trigonometric) shock surface
(test case: \texttt{advec-trig}) to demonstrate the
HOIST method is able to represent curved discontinuity surfaces
in three dimensions. Consider the domain
$\Omega\coloneqq(-1, 1)\times(0, 1)\times (0, 1)$ with trigonometric
advection field and piecewise constant boundary condition
\begin{equation}
 \beta: x \mapsto (-\sin(\pi x_2), 1, 0), \qquad
 U_\infty: x \mapsto H(x_1);
\end{equation}
this is the three-dimensional extrusion of the problem considered in
\cite{zahr_implicit_2020}.
Similar to the previous section, we require the reference mesh to have element
edges aligned with the line $\{(0, 0, s) | s\in[0,1]\}$ to ensure the Heaviside
function is accurately represented and integrated along boundary
$\partial\Omega_2$ (same boundary labeling as
Section~\ref{sec:numexp:advec:planar}) and discretize
the domain using a structured mesh consisting of 324 uniform right tetrahedra
(Figure~\ref{fig:advec3d_geom}).
Because the analytical solution is piecewise constant, we again use a
$p=0$ DG solution; however, we use a $q=2$ mesh to approximate the shock
surface to third-order accuracy. The HOIST parameters are defined in
Table~\ref{tab:params}, which are identical to the planar case 
except we take $c_4'=0.05$ to prevent high-order elements from becoming singular.

Starting from a shock-agnostic (structured) mesh and
first-order finite volume solution, the HOIST method
leads to an accurate solution and a curved mesh that
is well-aligned with the discontinuity surface
(Figure~\ref{fig:advec3d-case1}) and the magnitude of the standard and
enriched DG residual are reduced $\Ocal(10^{-8})$ and
$\Ocal(10^{-4})$, respectively; 12 element collapses are
required. Unlike the planar case, we do not see Newton-like
convergence; however, we observe a rapid drop in the residual
from $\Ocal(10^{-2})$ to $\Ocal(10^{-6})$ in just 11 iterations
($k=33$ to $k=45$) once the discontinuity surface is nearly tracked
(Figure~\ref{fig:advec3d_case1_hist}).
The main difference between the planar and trigonometric shock case
is that the analytical solution lies in the nonlinear approximation space
(defined as the composition of the DG solution with the domain mapping)
in the planar case, whereas it does not in the trigonometric case. As a
result, the enriched DG residual also converged to a near-zero optimal
value $\Ocal(10^{-8})$) in the planar case, whereas it was much larger
in the trigonometric case $\Ocal(10^{-4})$, which implies the
Levenberg-Marquardt Hessian in (\ref{eqn:hess0})-(\ref{eqn:hess1}) is
a better approximation to the true Hessian of the objective function in
the planar case and likely the reason for the different convergence rates.
\ifbool{fastcompile}{}{
\begin{figure}[!htbp]
\centering
\includegraphics[width=0.35\textwidth]{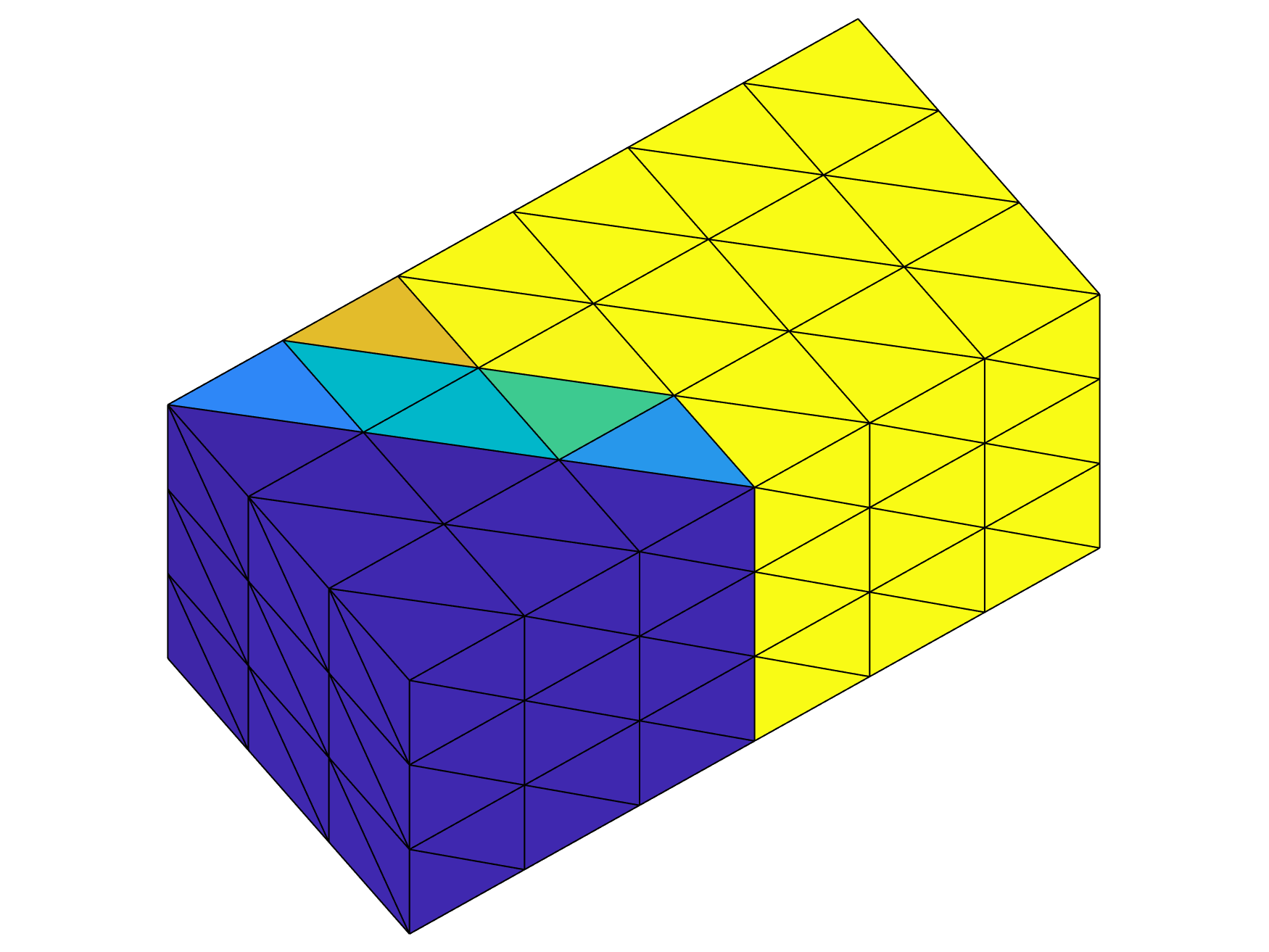}
\includegraphics[width=0.35\textwidth]{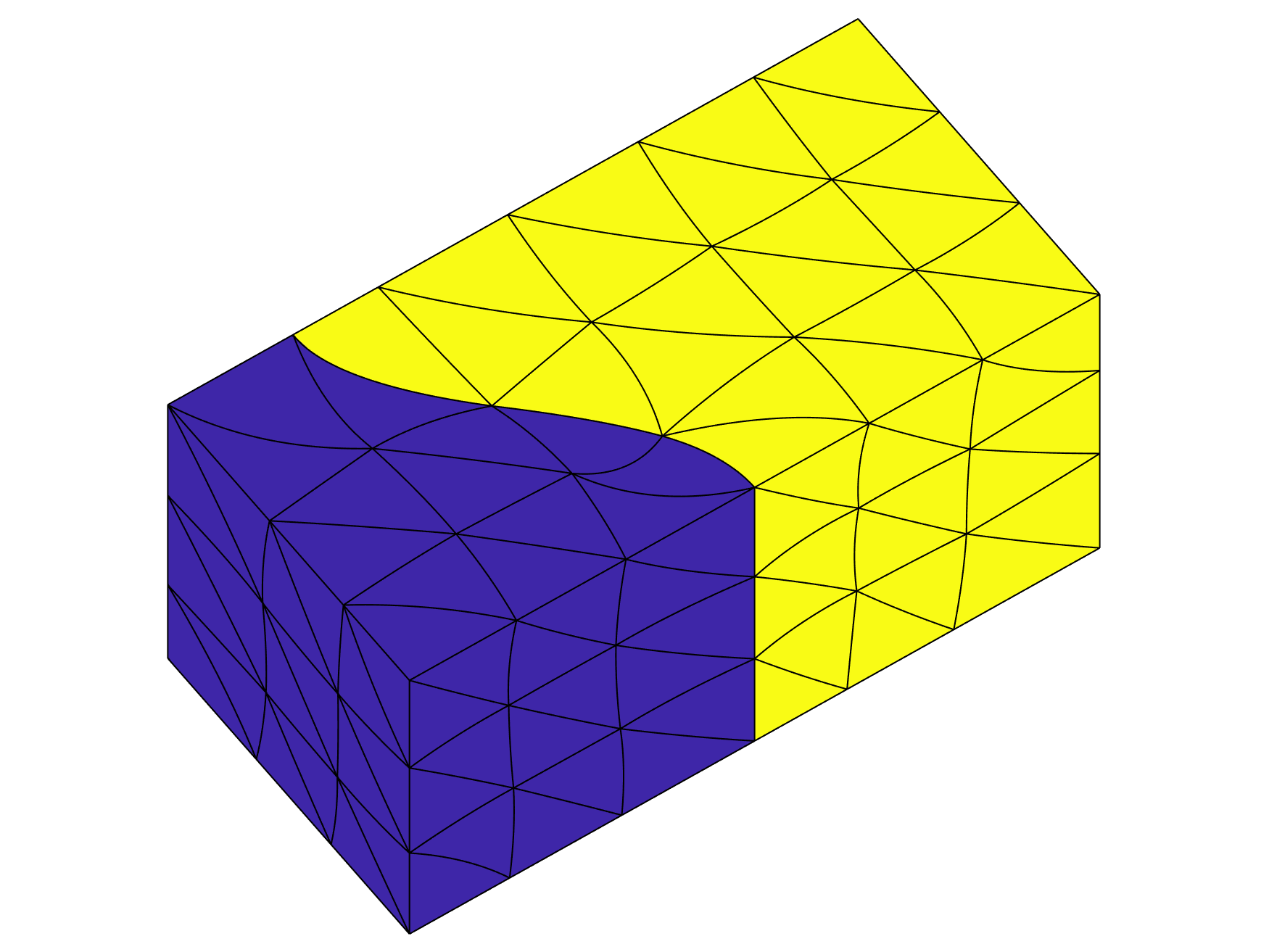}
\includegraphics[width=0.25\textwidth]{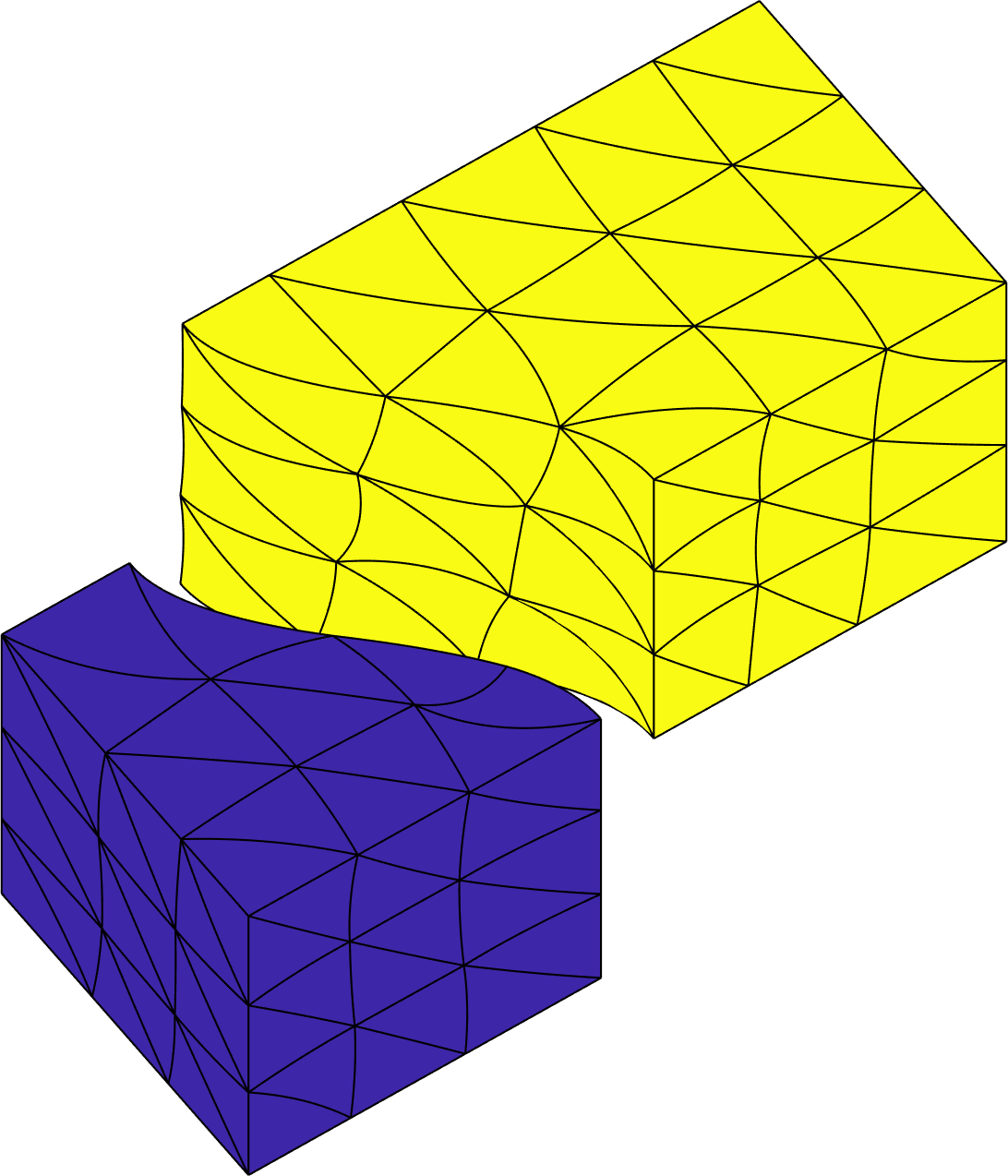}
\caption{The starting point (\textit{left}) and the tracked configuration
 (\textit{middle}) for the \texttt{advec-trig} test case.
 The mesh is separate along the discontinuity surface
 for a clearer view of the tracked mesh (\textit{right}). Colorbar in
 Figure~\ref{fig:advec3d-case0}.}
\label{fig:advec3d-case1}
\end{figure}
}

\begin{figure}[!htbp]
\centering
\ifbool{fastcompile}{}{
\input{_py/advec1v3d_zahr0_case1p0q2_hist.tikz}
}
\caption{SQP convergence history and behavior of adaptive parameters for
 for the \texttt{advec-trig} test case (legend in Table~\ref{tab:legend}).}
\label{fig:advec3d_case1_hist}
\end{figure}

\subsection{Inviscid Burgers' equation}
\label{sec:numexp:iburg}
Next, we consider the time-dependent, inviscid Burgers' equation
that governs nonlinear advection of a scalar quantity
$\func{\phi}{\Omega}{\Rbb}$ through a one-dimensional domain
$\bar\Omega\subset\Rbb$
\begin{equation} \label{eqn:iburg}
  \pder{}{t}\phi(z,t) + \phi(z,t)\pder{}{z}\phi(z,t) = 0, \qquad
  \phi(z,0)=\phi_0(z)
\end{equation}
for all $z\in\bar\Omega$ and $t\in\Tcal$ with the initial condition
$\func{\phi_0}{\Omega}{\Rbb}$. This
leads to a conservation law of the form (\ref{eqn:claw-phys}) over
the space-time domain $\Omega\coloneqq\bar\Omega\times\Tcal$ with
conservative variable $\func{U}{\Omega}{\Rbb}$ and flux function
$\func{F}{\Rbb}{\Rbb^{1\times 2}}$
\begin{equation}
 U : x \mapsto \phi(x_1, x_2), \qquad
 F : W \mapsto \begin{bmatrix} W^2/2 & W \end{bmatrix}.
\end{equation}
that is identical to the advection equation in (\ref{eqn:advec}) with
local space-time flow direction $\beta : x \mapsto (U(x)/2,1)$. For the
DG discretization, we use the smoothed upwind flux described in
\cite{zahr_implicit_2020} as the inviscid numerical flux function.
\subsubsection{Shock acceleration}
\label{sec:numexp:iburg:acc}
First, we consider an accelerating shock
(test case: \texttt{iburg-acc}) in the one-dimensional spatial
domain $\bar\Omega\coloneqq(-0.2,1)$ and temporal domain $\Tcal=(0,1.2)$
produced by the initial condition
\begin{equation}
 \phi_0: z \mapsto \mu_1 H(-z) - \mu_2 H(z)
\end{equation}
and boundary condition $\phi(-0.2,t)=\mu_1$ for $t\in\Tcal$; in this work,
we take $\mu_1=4$ and $\mu_2=3$. The analytical solution is
\begin{equation} \label{eqn:iburg:nrl0}
 \phi : (z,t) \mapsto \mu_1 H(z_\mathrm{s}(t)-z) + \frac{\mu_2(z-1)}{1+\mu_2 t} H(z-z_\mathrm{s}(t)),
\end{equation}
where $\func{z_\mathrm{s}}{\Tcal}{\Rbb}$ is the shock speed
\begin{equation}
 z_\mathrm{s} : t \mapsto (\mu_1/\mu_2+1)(1-\sqrt{1+\mu_2 t}) + \mu_1 t;
\end{equation}
both are shown in Figure~\ref{fig:iburg:nrl0:exact}.
The purpose of this numerical experiment is to provide a detailed
investigation into the proposed HOIST solver for polynomial degrees
$p=1,2,3$ and verify optimal convergence rate of the PDE solution
($\phi$) and shock position ($z_\mathrm{s}$) of the method for a
two-dimensional, nonlinear problem.
\ifbool{fastcompile}{}{
\begin{figure}[!htbp]
\centering
\input{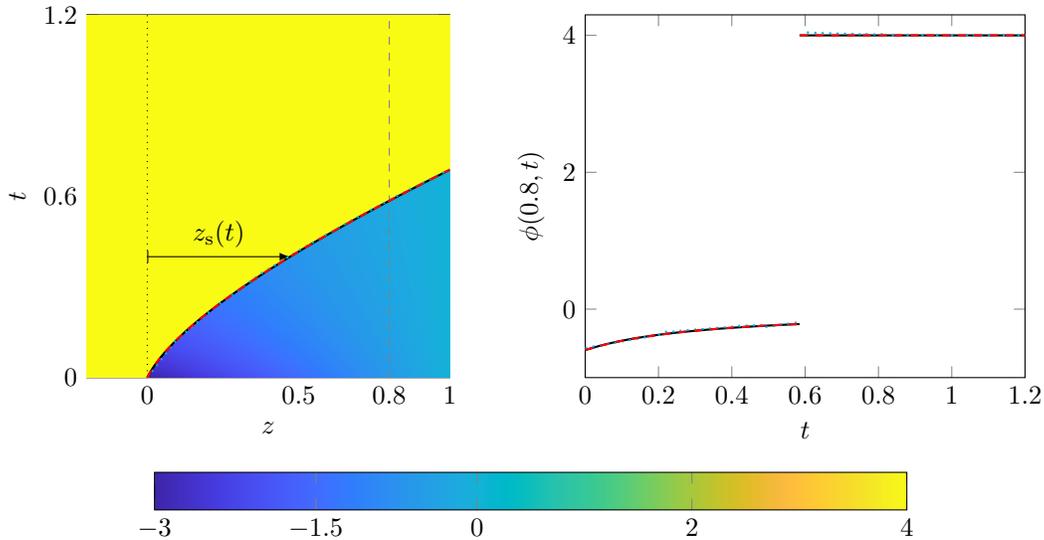}

\colorbarMatlabParula{-3}{-1.5}{0.0}{2}{4}
\caption{Analytical space-time solution to the \texttt{iburg-acc} test case
  with shock surface (\ref{line:iburg-acc:exact})
  (\textit{left}) and the corresponding solution along the line
  $\{(0.8,t)\mid t\in[0,1]\}$ (\ref{line:iburg-acc:slice})
  (\textit{right}). The corresponding quantities for
  the $p=1$ (\ref{line:iburg-acc:p1}) and $p=2$ (\ref{line:iburg-acc:p2})
  HOIST simulations are also shown for comparison. The $p>2$ HOIST solutions
  are not shown for clarity because they are visually indistinguishable
  from the $p=2$ solution, which in turn overlaps the analytical quantity.}
\label{fig:iburg:nrl0:exact}
\end{figure}
}

Due to the curvature of the shock in space-time, we choose to
approximate the mesh with the same polynomial degree as the solution
field, i.e., $q = p$. We discretize the space-time domain with a
structured mesh of 72 right triangles and use this mesh to initialize
all HOIST simulations. Furthermore, the mesh was constructed such that
the hypotenuse of each triangle traverses the discontinuity surface
making the initial mesh far from alignment with the shock
(Figure~\ref{fig:iburg:nrl0:p1to3}). The same HOIST parameters are
used for all polynomial degrees considered (Table~\ref{tab:params}),
except we take $\gamma_0=\gamma_{\min}=1$ in the $p=1$
case because the solution is under-resolved and the more
conservative choice is required.

For all three cases, the proposed solver simultaneously updates the
mesh and solution until the discontinuity is tracked, the fully
discrete PDE is satisfied on the discontinuity-aligned mesh, and
elements away from the discontinuity approach equilateral triangles
(owing to the choice of $K_{\star,e}=K_\star$)
(Figure~\ref{fig:iburg:nrl0:p1to3}). The number of elements removed
in each case is: 1 ($p=1$), 1 ($p=2$), and 7 ($p=3$).
In the $p=1$ case, the solution is somewhat underresolved as evident
by the faceted approximation to the discontinuity, although this is
expected given the coarseness of the mesh. For all other cases, the shock
is tracked using curved, high-order elements and the solution is well-resolved.
Both the space-time discontinuity surface and slice of the conserved
variable from the HOIST solution are visually indistinguishable from
the corresponding analytical quantity for $p\geq 2$
(Figure~\ref{fig:iburg:nrl0:exact}).

The SQP solver converges to the solution the optimization problem
(\ref{eqn:pde-opt}) rapidly considering the starting point is far
from the HOIST solution. The discontinuity is mostly tracked after
only 20 iterations and tracked to the extent possible given the resolution
limits of the discretization by iteration 30
(Figure~\ref{fig:iburg:nrl0:p1to3}).
In all cases, the first-order optimality criteria are driven to
relatively tight tolerances within 100 iterations; after the
element collapses cease and the topology of the shock has
been discovered (around iteration 20), the DG residuals and
optimality condition decrease sharply (Figure~\ref{fig:iburg:nrl0:conv-hist}).
The $p=2$ simulation is an outlier because it is effectively converged by
iteration 10, which is significantly faster than the other polynomial degrees
considered. We attribute this to the $p=1$ simulation being under-resolved
(requires more conservative parameters) and the $p>2$ elements requiring more
frequent re-initializations because they are more susceptible to oscillations.
In our experience, $p=2$ tends to be the ``sweet'' spot for the HOIST
method and is used extensively in the remaining numerical experiments.

\ifbool{fastcompile}{}{
\begin{figure}[!htbp]
\centering
\begin{tikzpicture}
\begin{groupplot}[
  group style={
      group size=4 by 3,
      horizontal sep=0.3cm,
      vertical sep=0.3cm
  },
  width=0.37\textwidth,
  axis equal image,
  xticklabels={,,}, yticklabels={,,},
  xmin=-0.2, xmax=1,
  ymin=0, ymax=1.2
]

\nextgroupplot[ylabel={$p=1$}, title={Iteration 0}]
\addplot graphics [xmin=-0.2, xmax=1, ymin=0, ymax=1.2] {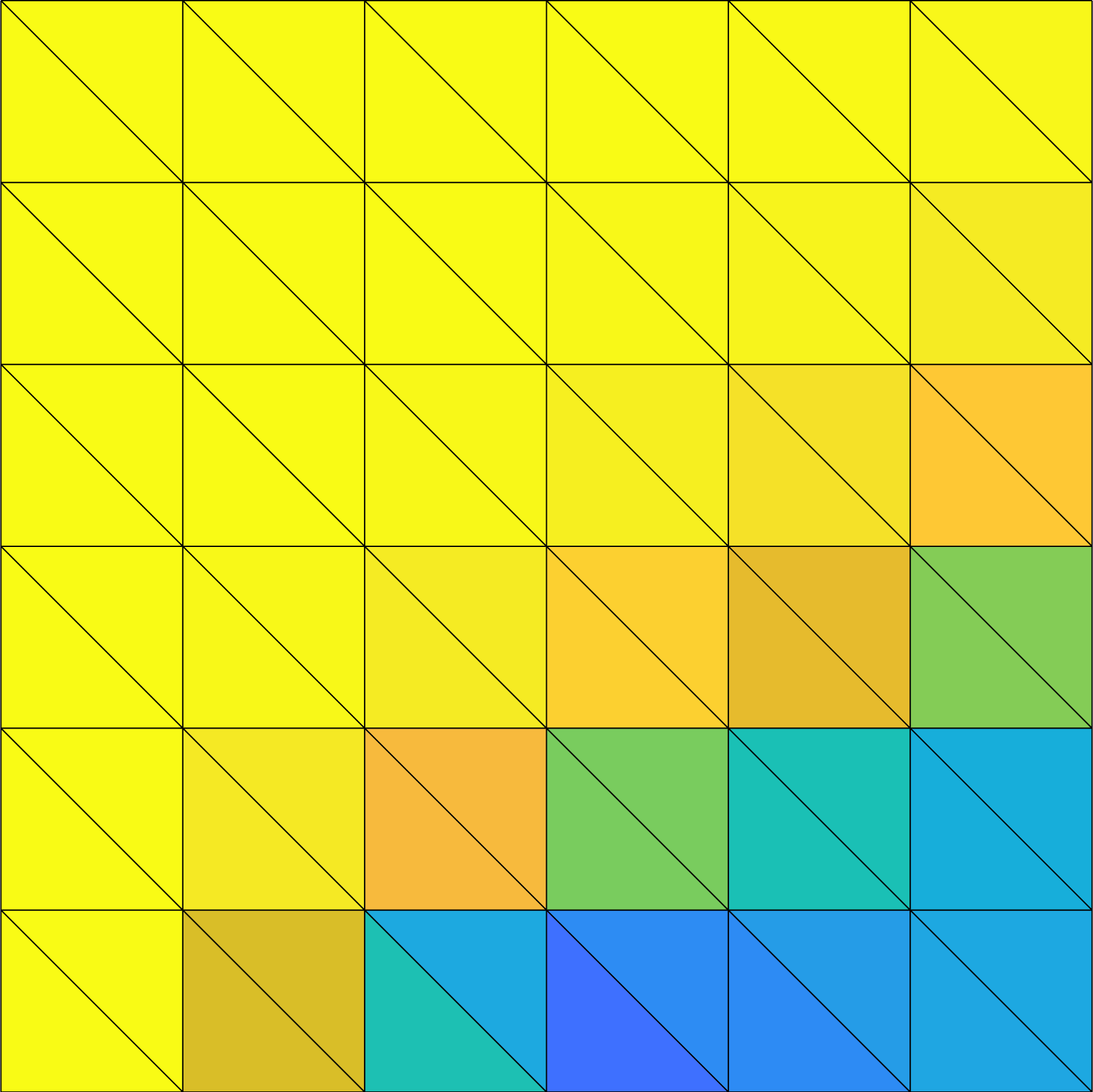};
\nextgroupplot[title={Iteration 10}]
\addplot graphics [xmin=-0.2, xmax=1, ymin=0, ymax=1.2] {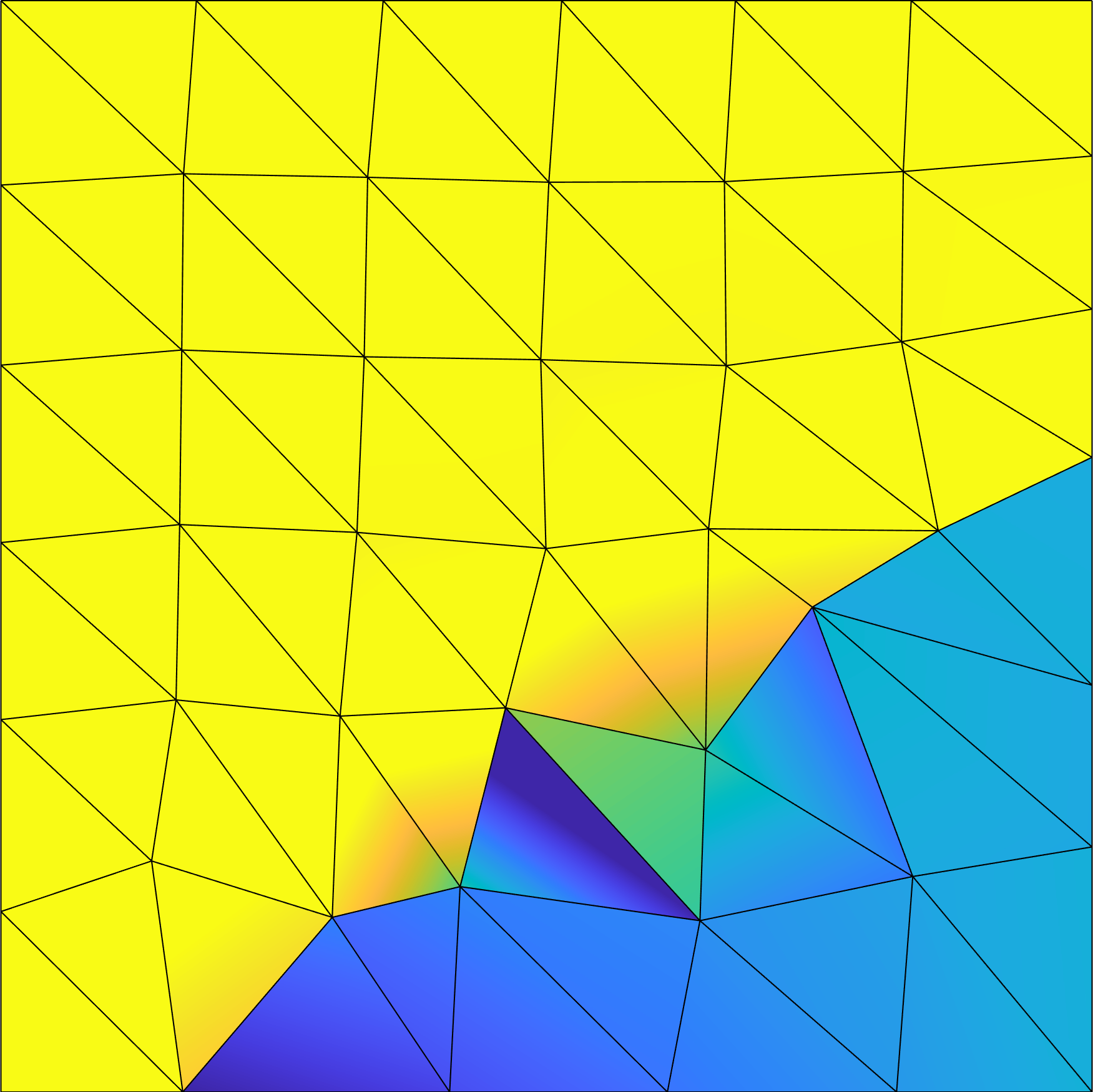};
\nextgroupplot[title={Iteration 20}]
\addplot graphics [xmin=-0.2, xmax=1, ymin=0, ymax=1.2] {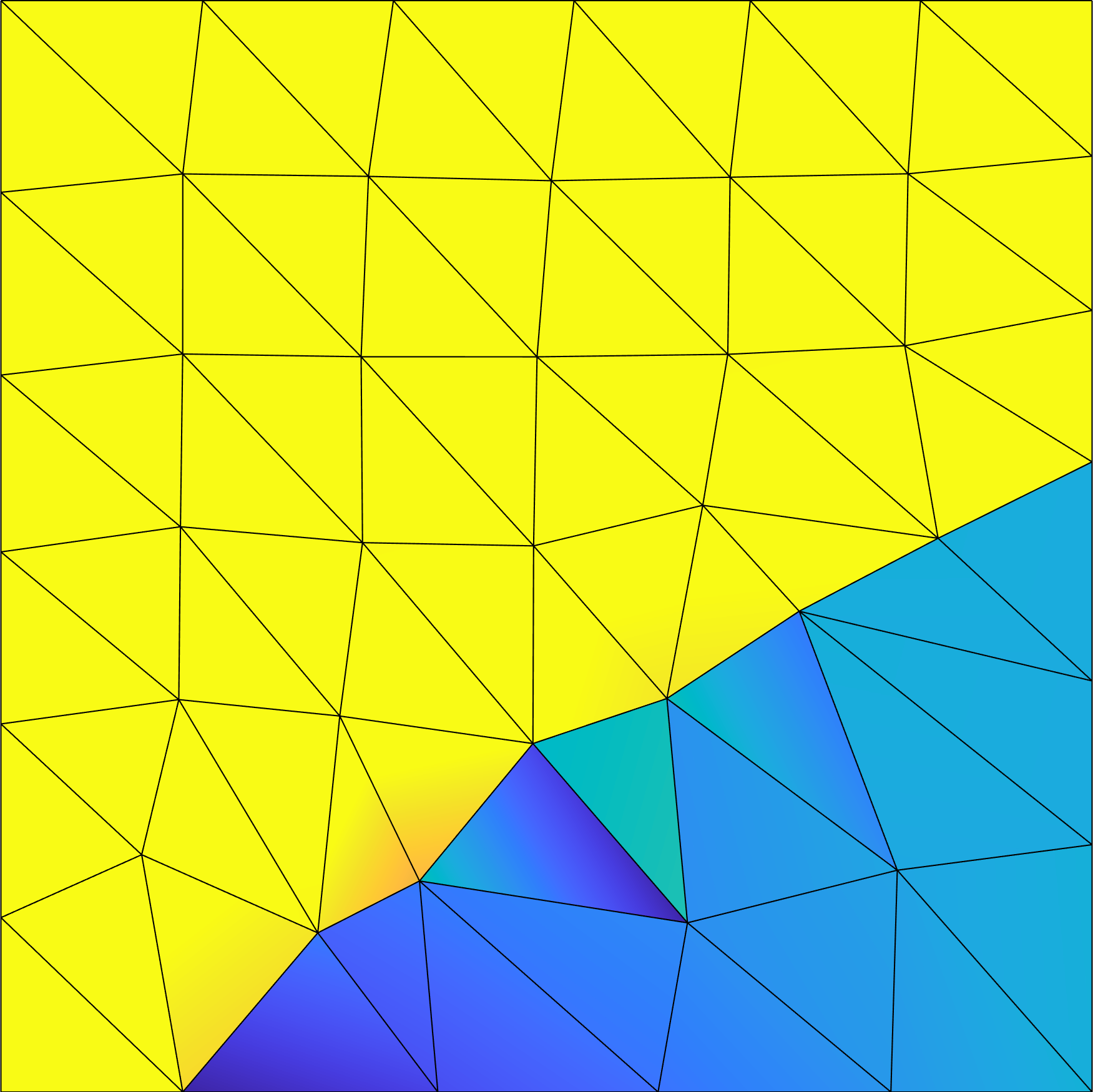};
\nextgroupplot[title={Iteration 30}]
\addplot graphics [xmin=-0.2, xmax=1, ymin=0, ymax=1.2] {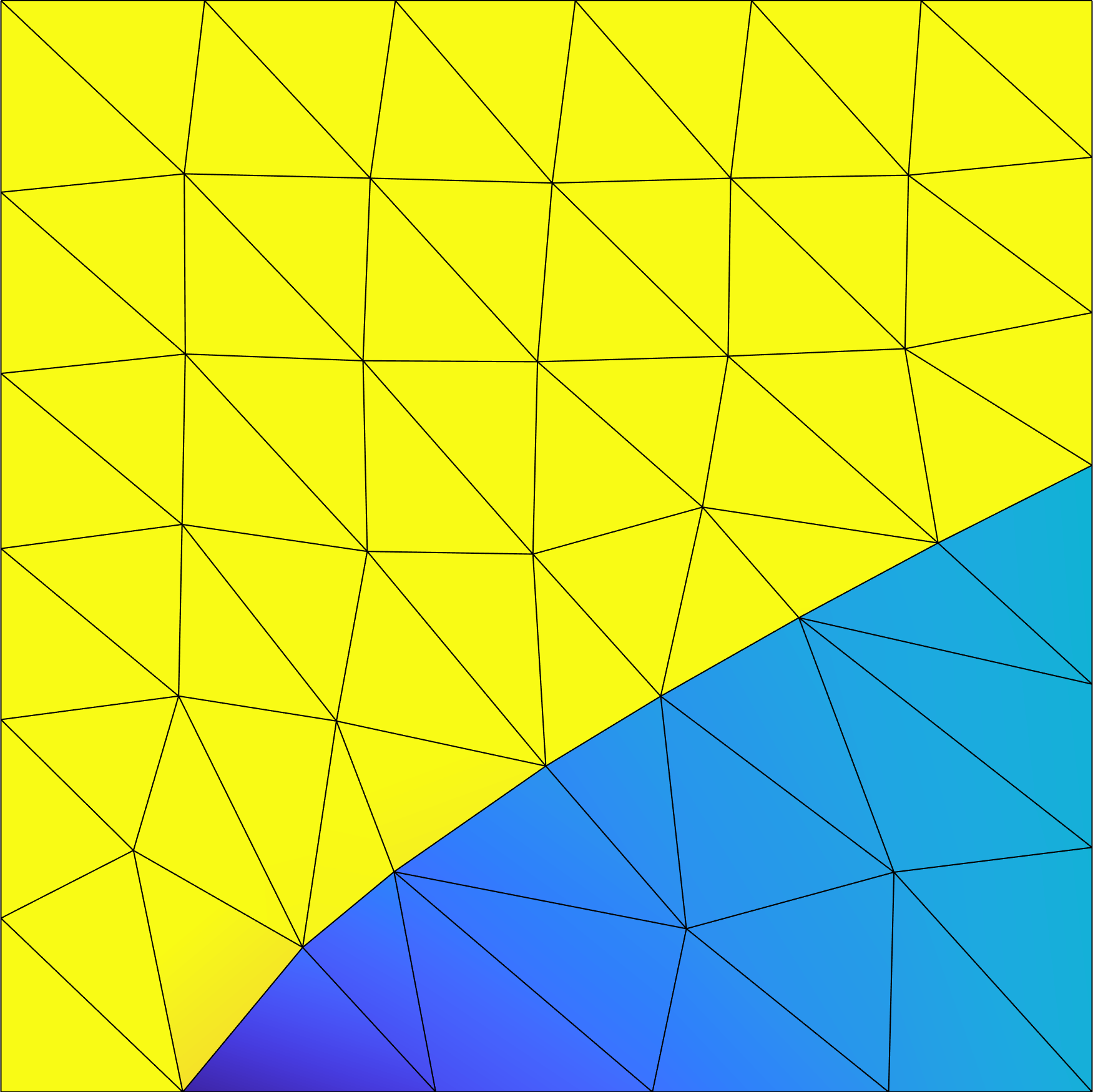};

\nextgroupplot[ylabel={$p=2$}]
\addplot graphics [xmin=-0.2, xmax=1, ymin=0, ymax=1.2] {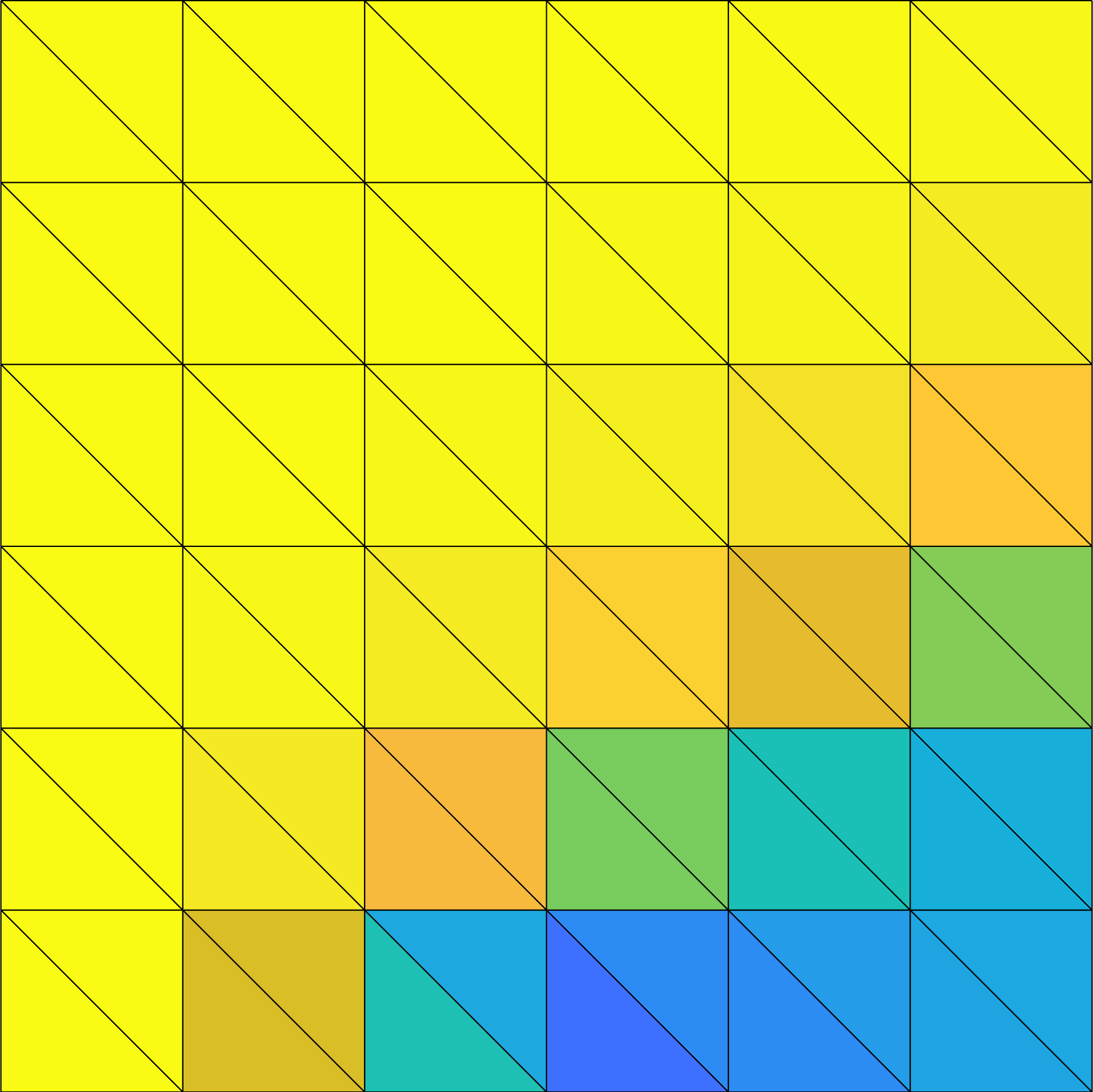};
\nextgroupplot
\addplot graphics [xmin=-0.2, xmax=1, ymin=0, ymax=1.2] {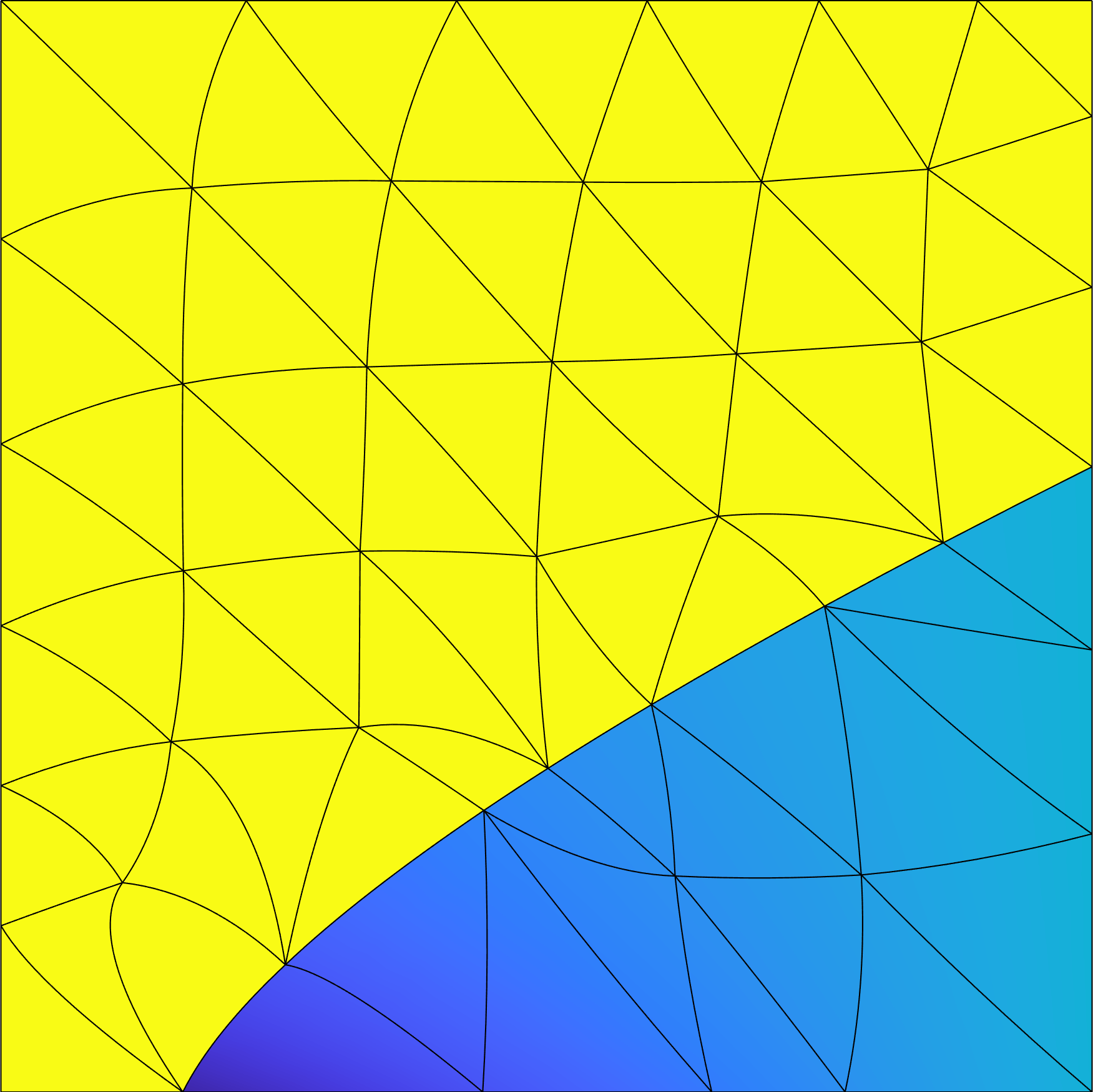};
\nextgroupplot
\addplot graphics [xmin=-0.2, xmax=1, ymin=0, ymax=1.2] {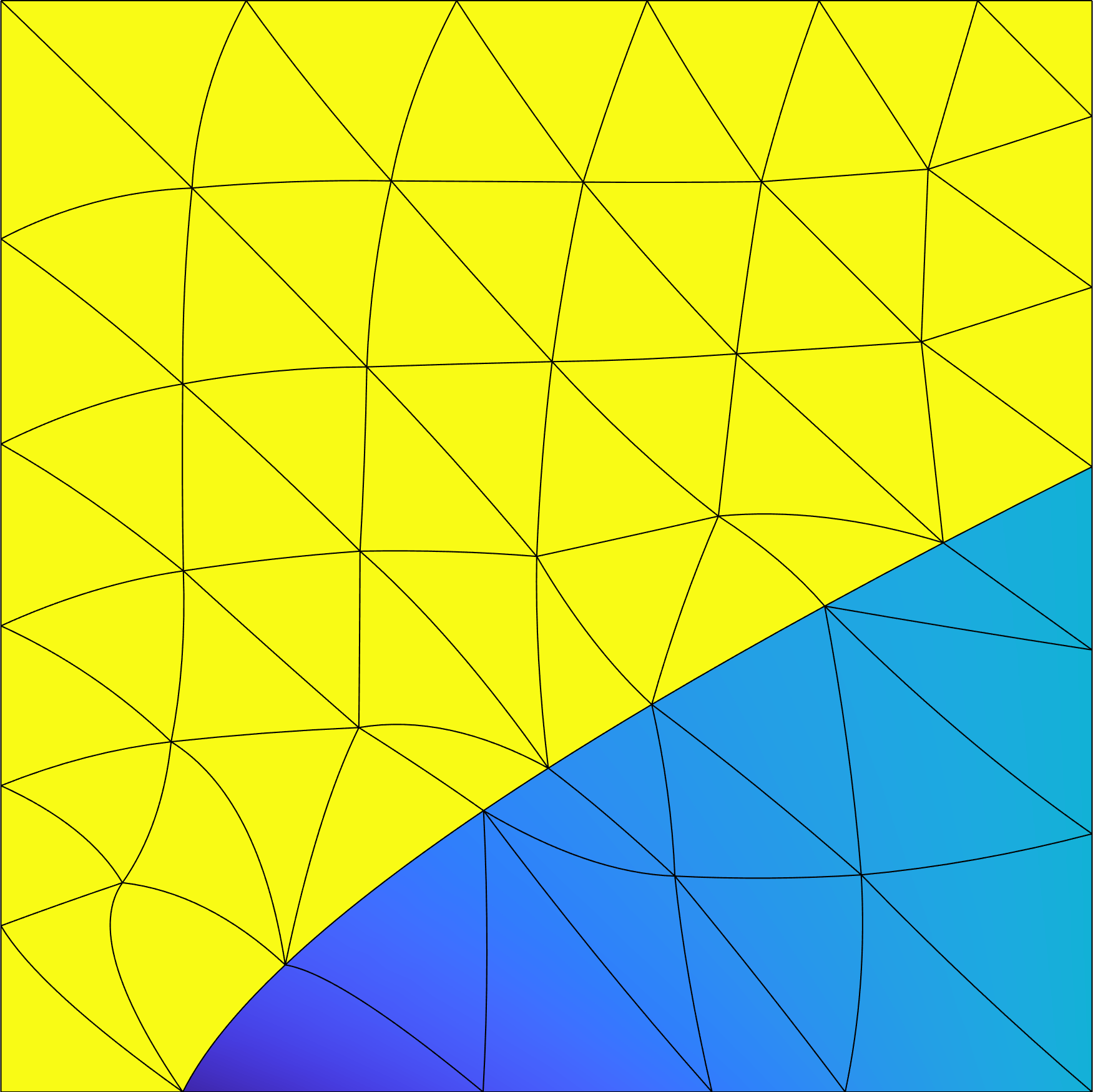};
\nextgroupplot
\addplot graphics [xmin=-0.2, xmax=1, ymin=0, ymax=1.2] {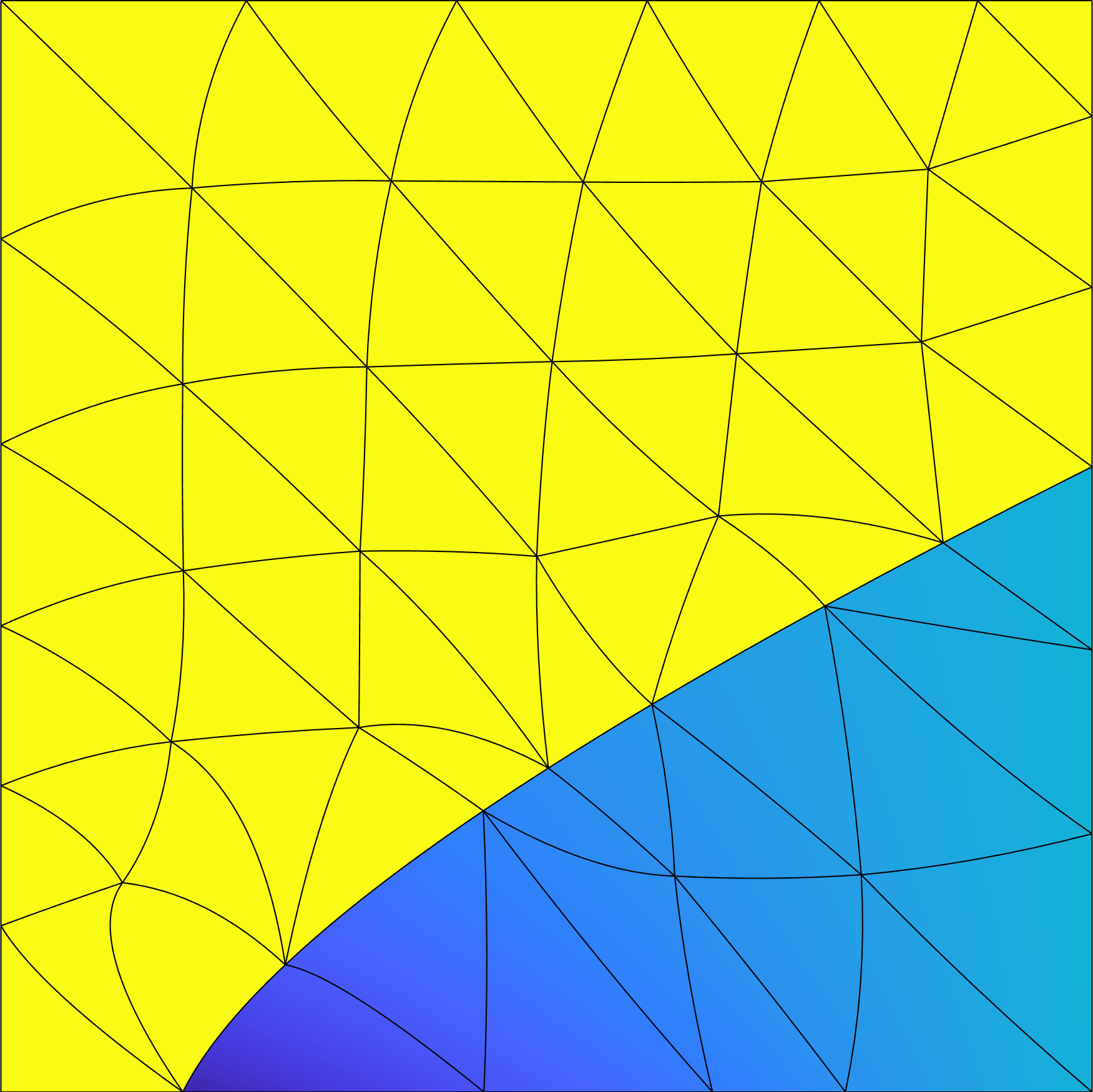};

\nextgroupplot[ylabel={$p=3$}]
\addplot graphics [xmin=-0.2, xmax=1, ymin=0, ymax=1.2] {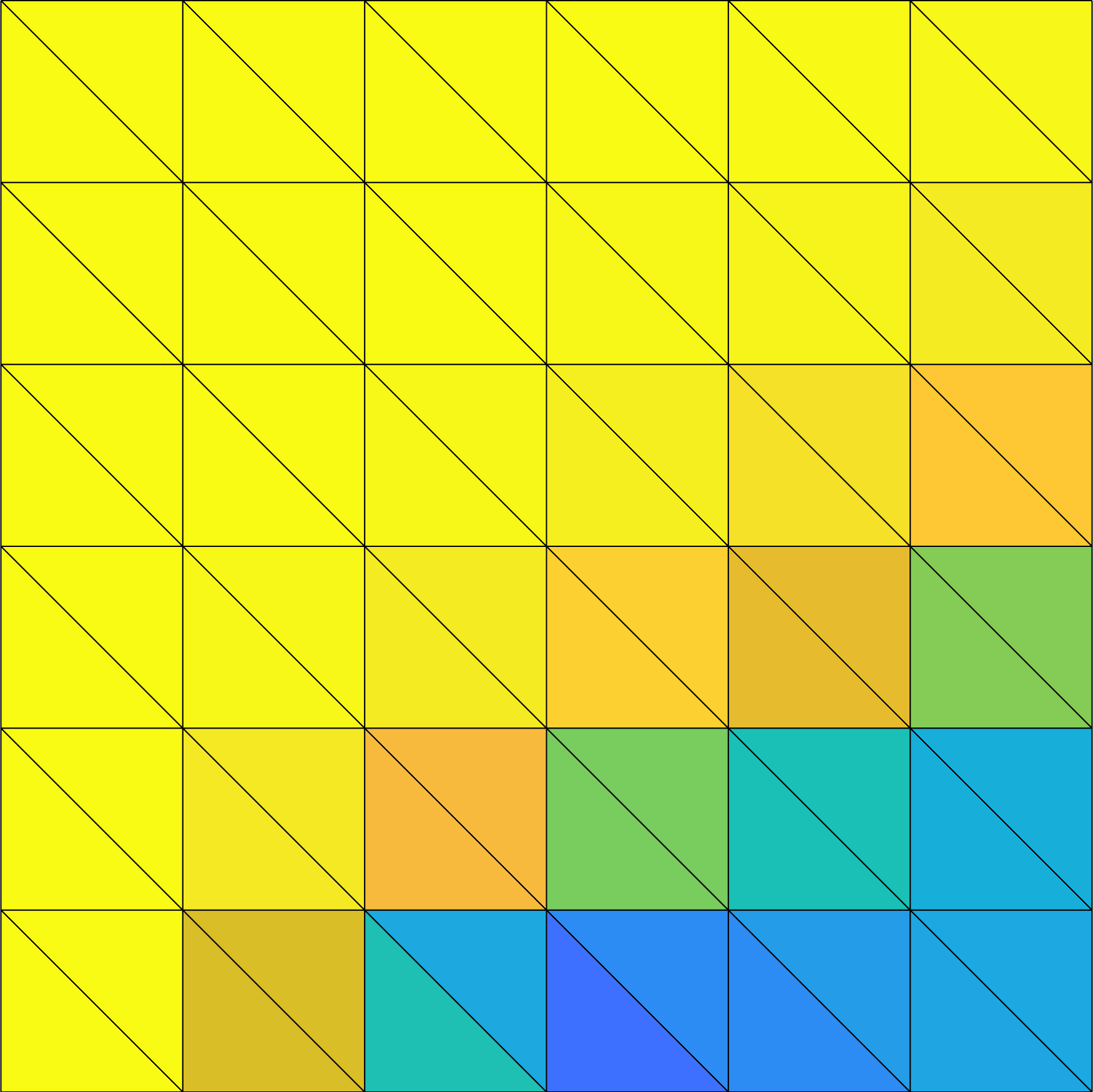};
\nextgroupplot
\addplot graphics [xmin=-0.2, xmax=1, ymin=0, ymax=1.2] {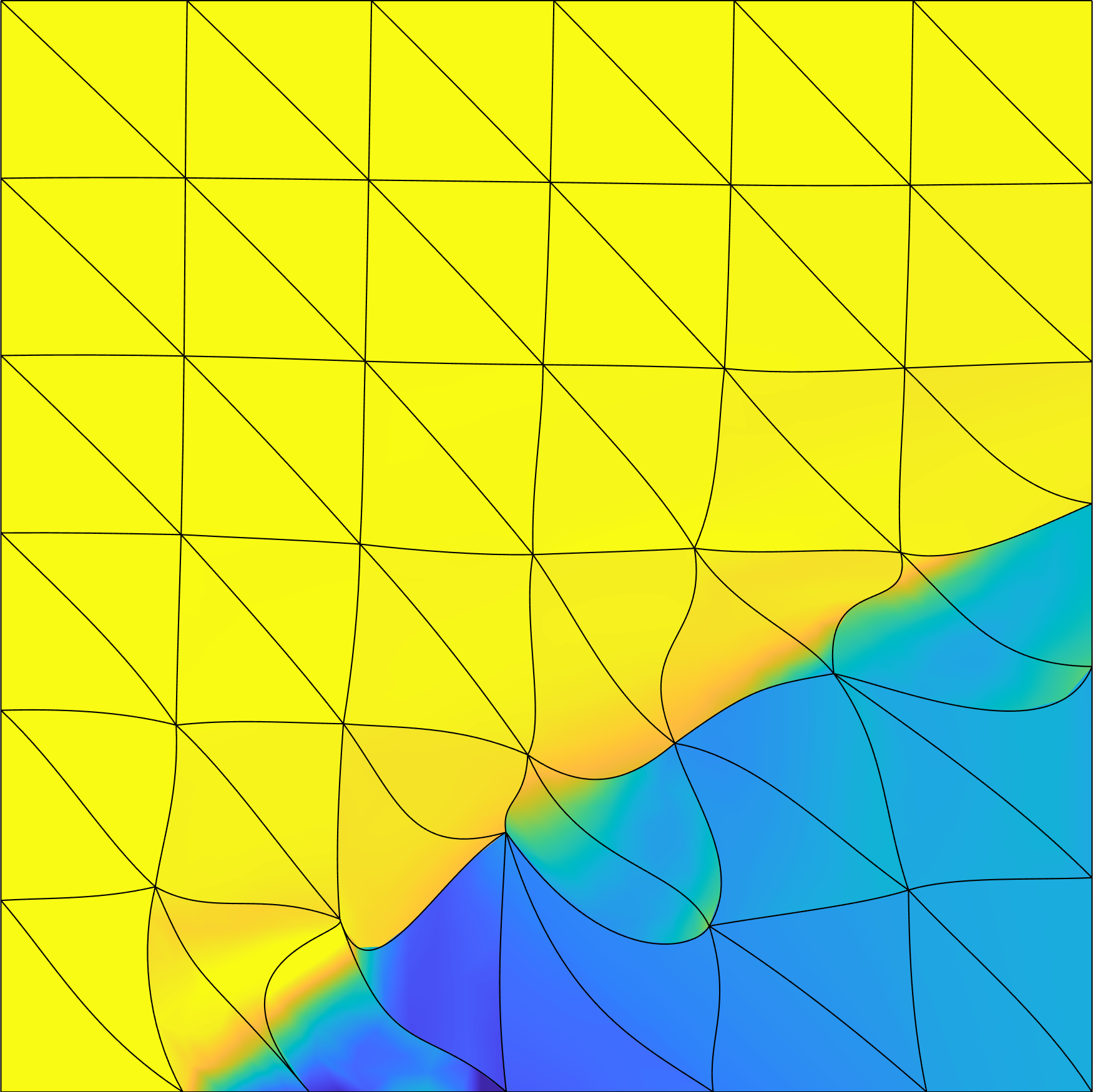};
\nextgroupplot
\addplot graphics [xmin=-0.2, xmax=1, ymin=0, ymax=1.2] {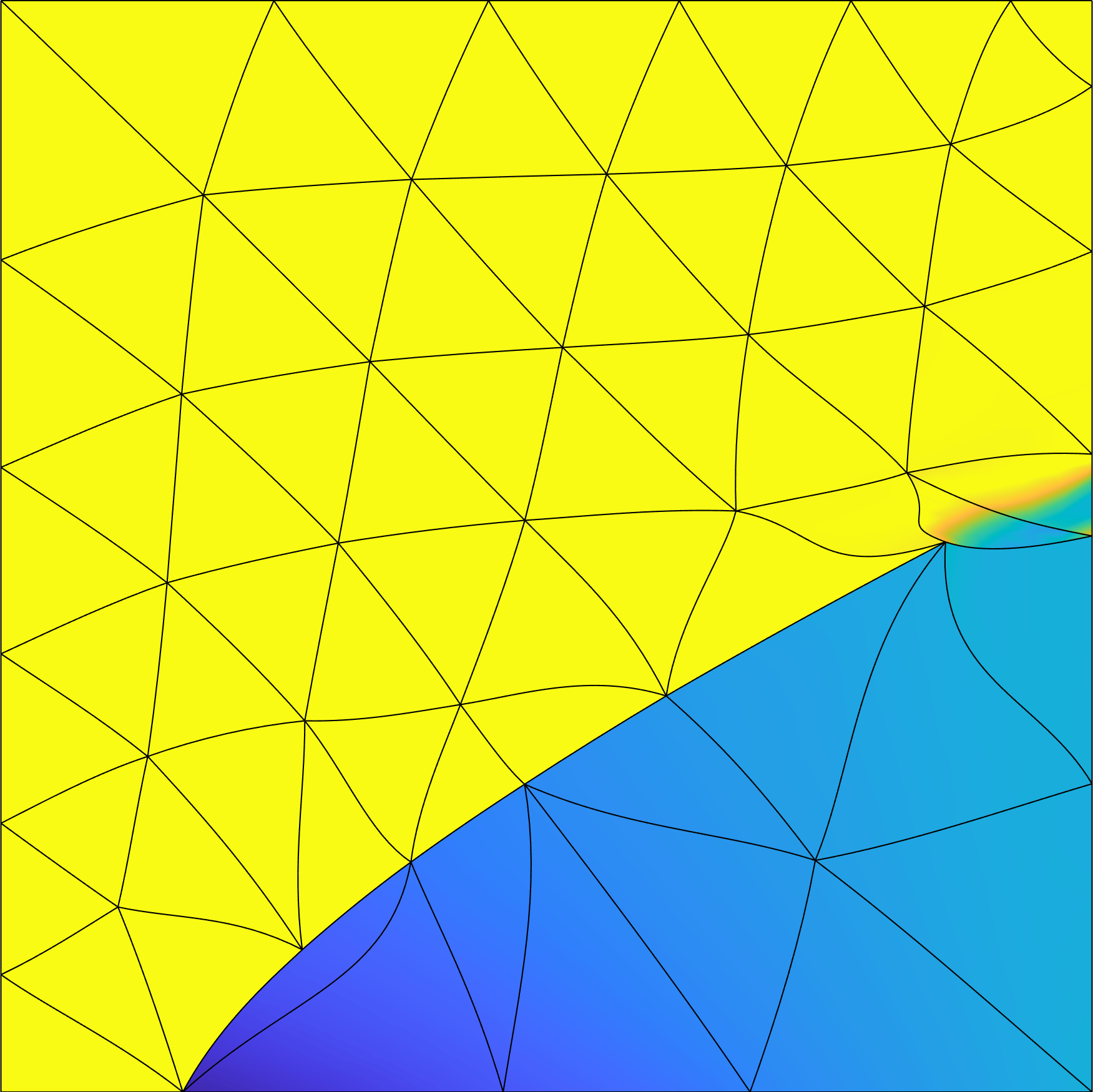};
\nextgroupplot
\addplot graphics [xmin=-0.2, xmax=1, ymin=0, ymax=1.2] {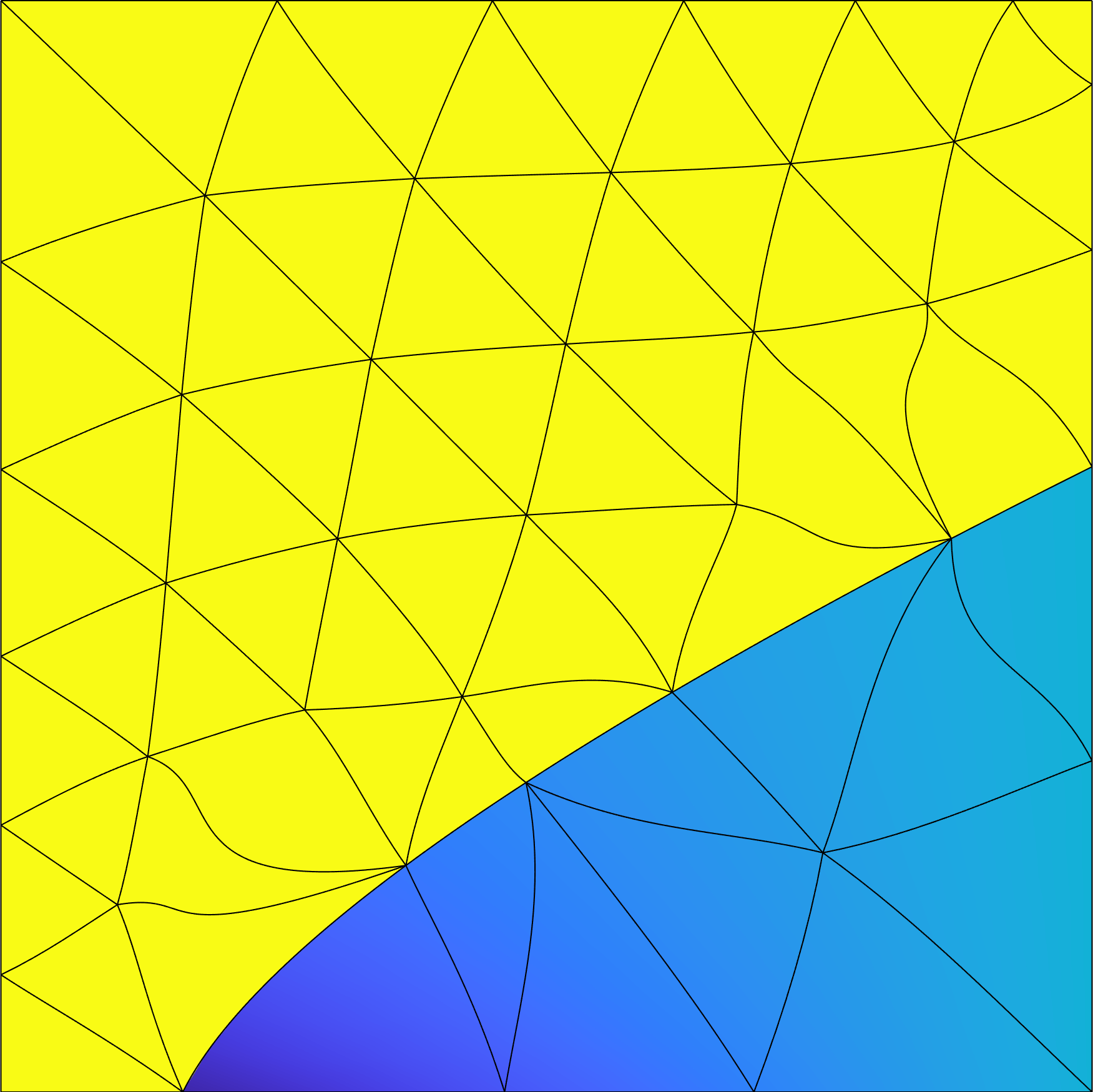};
\end{groupplot}
\end{tikzpicture}

\caption{Selected HOIST iterations for the \texttt{iburg-acc} test case
  for $p=1,2,3$ (\textit{top-to-bottom})
  at iterations $k=0, 10, 20, 30$ (\textit{left-to-right}). Colorbar
  in Figure~\ref{fig:iburg:nrl0:exact}.}
\label{fig:iburg:nrl0:p1to3}
\end{figure}
}

%
%
%

\ifbool{fastcompile}{}{
\begin{figure}[!htbp]
\centering
\input{_py/iburg1v1d_sptm_nrl0_hist.tikz}
\caption{SQP convergence history and behavior of adaptive parameters for
 the \texttt{iburg-acc} test case (legend in Table~\ref{tab:legend})
 for $p=1,2,3$ (\textit{top-to-bottom}).}
\label{fig:iburg:nrl0:conv-hist}
\end{figure}
}

Oscillatory solutions arise in some elements during intermediate iterations,
particularly for the higher polynomial degrees
(Figure~\ref{fig:iburg:nrl0:p1to3}); however, the element-wise
re-initialization (Section~\ref{sec:ist_solve:mod:reinit}) removes them,
which ensures the solver can continue without solution and mesh degradation.
Without re-initializations, the SQP solver can effectively stall because
the oscillatory solution leads to poor search directions
(Section~\ref{sec:ist_solve:mod:reinit}) that eventually cause the
step length ($\alpha_k$) to become very small to satisfy the line
search criteria in (\ref{eqn:steep_descent}). To demonstrate this,
we consider the $p=3$ HOIST simulation without re-initialization
(all other parameters unchanged) and see that while the shock
is mostly tracked, a few oscillatory elements have stalled progress
to the solution (Figure~\ref{fig:iburg:nrl0:reinit-demo}). By contrasting
to the HOIST solution (Figure~\ref{fig:iburg:nrl0:p1to3}) and SQP solver
performance (Figure~\ref{fig:iburg:nrl0:conv-hist}) with re-initialization,
it is clear that re-initialization is important for robust convergence.
\ifbool{fastcompile}{}{
\begin{figure}[!htbp]
\centering
\input{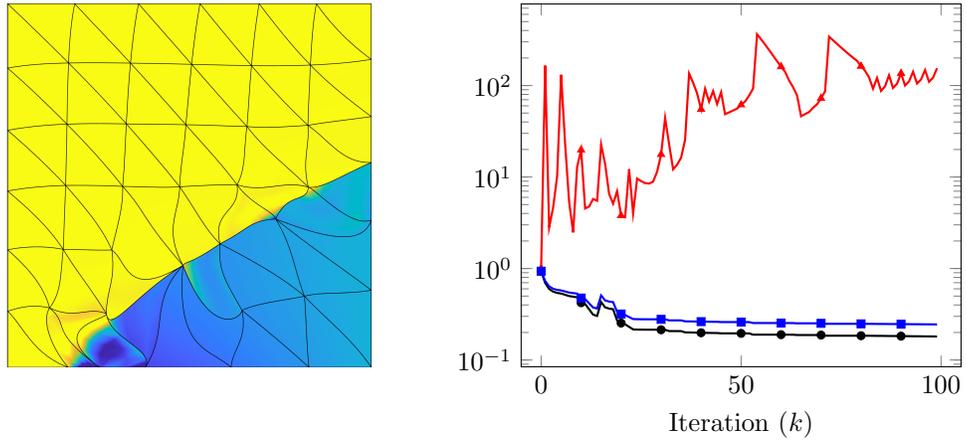}
\caption{HOIST solution ($p=3$) to the \texttt{iburg-acc} test case
  without re-initialization (\textit{left}) (colorbar in
  Figure~\ref{fig:iburg:nrl0:exact}) and the corresponding SQP convergence
  history (\textit{right}) (legend in Table~\ref{tab:legend}).}
\label{fig:iburg:nrl0:reinit-demo}
\end{figure}
}

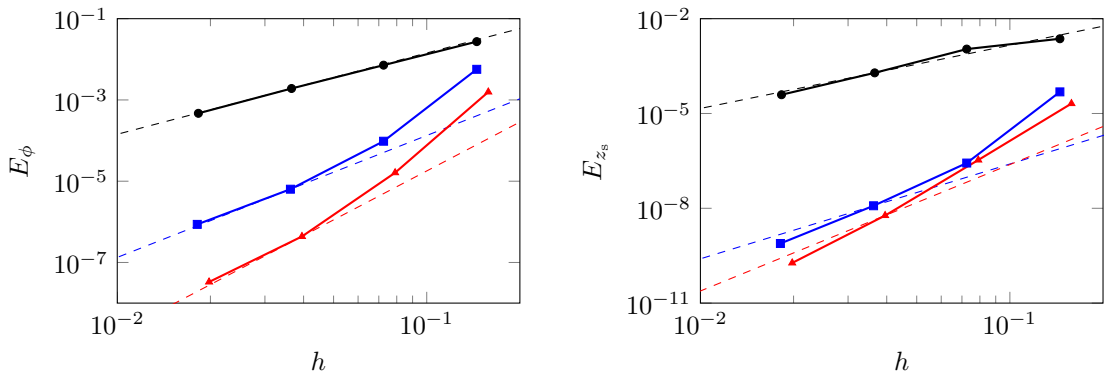
\begin{figure}[!htbp]
\centering
\ifbool{fastcompile}{}{
\begin{tikzpicture}
\begin{groupplot} [
group style={group size = 2 by 1, horizontal sep = 2.4cm}]
\nextgroupplot[width=0.42\textwidth, xlabel={$h$}, ymax=0.1, xmax=0.2, ylabel={$E_\phi$}, xmin=0.01, ymode=log, ymin=1e-08, xmode=log, height=0.325\textwidth]
\addplot [solid, thick, black, mark options={solid, thin}, mark=*, mark size=1.5]
coordinates {
( 1.45095250e-01,  2.72446105e-02)
( 7.25476250e-02,  7.17629824e-03)
( 3.65758478e-02,  1.91283079e-03)
( 1.82879239e-02,  4.68617897e-04)};\label{line:iburg:hconv:p1}

\addplot [dashed, black, forget plot]
coordinates {
( 9.14396195e-03,  1.19551924e-04)
( 2.90190500e-01,  1.20407664e-01)};

\addplot [solid, thick, blue, mark options={solid, thin}, mark=square*, mark size=1.5]
coordinates {
( 1.45095250e-01,  5.67592859e-03)
( 7.25476250e-02,  9.64029259e-05)
( 3.62738125e-02,  6.35782927e-06)
( 1.81369063e-02,  8.65797162e-07)};\label{line:iburg:hconv:p2}

\addplot [dashed, blue, forget plot]
coordinates {
( 9.06845313e-03,  9.93410823e-08)
( 2.90190500e-01,  3.25520859e-03)};

\addplot [solid, thick, red, mark options={solid, thin}, mark=triangle*, mark size=1.5]
coordinates {
( 1.58113883e-01,  1.57083535e-03)
( 7.90569415e-02,  1.61941923e-05)
( 3.95284708e-02,  4.36520693e-07)
( 1.98029509e-02,  3.31050054e-08)};\label{line:iburg:hconv:p3}

\addplot [dashed, red, forget plot]
coordinates {
( 9.90147543e-03,  1.71855897e-09)
( 3.16227766e-01,  1.78798876e-03)};

\nextgroupplot[width=0.42\textwidth, xlabel={$h$}, ymax=0.01, xmax=0.2, ylabel={$E_{z_\mathrm{s}}$}, xmin=0.01, ymode=log, ymin=1e-11, xmode=log, height=0.325\textwidth]
\addplot [solid, thick, black, mark options={solid, thin}, mark=*, mark size=1.5, forget plot]
coordinates {
( 1.45095250e-01,  2.31572916e-03)
( 7.25476250e-02,  1.09128505e-03)
( 3.65758478e-02,  1.93180154e-04)
( 1.82879239e-02,  3.92145310e-05)};

\addplot [dashed, black, forget plot]
coordinates {
( 9.14396195e-03,  1.20737596e-05)
( 2.90190500e-01,  1.21601823e-02)};

\addplot [solid, thick, blue, mark options={solid, thin}, mark=square*, mark size=1.5, forget plot]
coordinates {
( 1.45095250e-01,  4.82980513e-05)
( 7.25476250e-02,  2.70214717e-07)
( 3.62738125e-02,  1.20223406e-08)
( 1.81369063e-02,  7.70118884e-10)};

\addplot [dashed, blue, forget plot]
coordinates {
( 9.06845313e-03,  1.87849072e-10)
( 2.90190500e-01,  6.15543841e-06)};

\addplot [solid, thick, red, mark options={solid, thin}, mark=triangle*, mark size=1.5, forget plot]
coordinates {
( 1.58113883e-01,  2.05959397e-05)
( 7.90569415e-02,  3.37190240e-07)
( 3.95284708e-02,  5.90308890e-09)
( 1.98029509e-02,  1.86635120e-10)};

\addplot [dashed, red, forget plot]
coordinates {
( 9.90147543e-03,  2.32401501e-11)
( 3.16227766e-01,  2.41790521e-05)};

\end{groupplot}\end{tikzpicture}
}
\caption{$h$-convergence of the HOIST method for the
 \texttt{iburg-acc} test case for both error metrics,
 $E_\phi$ (\textit{left}) and $E_{z_\mathrm{s}}$ (\textit{right}),
 for polynomial degrees $p=1$ (\ref{line:iburg:hconv:p1}),
 $p=2$ (\ref{line:iburg:hconv:p2}), and
 $p=3$ (\ref{line:iburg:hconv:p3}). The dashed lines indicate
 the optimal convergence rate ($p+1$).}
\label{fig:iburg:nrl0:hconv}
\end{figure}

Finally, we study the $h$-convergence of the HOIST method for this
problem using the analytical solution to confirm it achieves optimal
convergence rates $\Ocal(h^{p+1})$ for the polynomial degrees $p=1,2,3$.
Because HOIST shock will not be exactly aligned with the analytical
shock position, it is difficult to compute standard integral-based
errors over the space-time domain as one would need to intersect
the HOIST mesh with exact shock location to accurately compute
the integrals using numerical quadrature, which is far from
straightforward for curved meshes, e.g., see \cite{hermes_high-order_2018}.
Instead, we consider two one-dimensional error metrics:
$\func{E_\phi}{\Rbb^{N_\ubm}\times\Rbb^{N_\xbm}}{\Rbb}$, the
$L^1$ error in the solution along the line $\{(0.8,t)\mid t\in[0,1]\}$,
and $\func{E_{z_\mathrm{s}}}{\Rbb^{N_\xbm}}{\Rbb}$, the
$L^2$ error in the shock position, i.e.,
\begin{equation}
 E_\phi : (\ubm,\xbm) \mapsto
 \int_0^1 \left|\phi(\bar{z},t)-\hat\phi(\bar{z},t;\ubm,\xbm)\right|\, dt,
 \qquad
 E_{z_\mathrm{s}} : \xbm \mapsto
 \sqrt{\int_0^1 \left|z_\mathrm{s}(t)-\hat{z}_\mathrm{s}(t;\xbm)\right|^2 \, dt},
\end{equation}
where $\func{\hat\phi}{\Tcal\times\bar\Omega\times\Rbb^{N_\ubm}\times\Rbb^{N_\xbm}}{\Rbb}$ is the HOIST approximation to $\phi$ given by
$\hat\phi: (z,t,\ubm,\xbm)\mapsto\Xi(\ubm)\circ\Gcal_h^{-1}((z,t);\xbm)$
and $\func{\hat{z}_\mathrm{s}}{\Tcal\times\Rbb^{N_\xbm}}{\Rbb}$ is the HOIST
approximation to the shock position. The integral in $E_\phi$ is computed
accurately by slicing the physical mesh to create a one-dimensional mesh,
intersecting it with the exact position of the shock, and integrating using
Gaussian quadrature over each element of the one-dimensional
super-mesh. We use the $L^1$ error
for the solution because it contains a discontinuity and we
expect the $L^1$ error to converge at a given rate provided
the smooth solution and position of the discontinuity converge
at that rate \cite{zahr_implicit_2020}. We consider four levels of
refinement generated by uniformly refining a coarse initial
mesh (each triangular element split into four elements);
due to element collapses, the number of elements in
each refinement level does not necessarily contain
four times the elements of the previous level. We
observe optimal convergence rate for $E_\phi$ and
better-than-optimal convergence rates in $E_{z_\mathrm{s}}$
(Figure~\ref{fig:iburg:nrl0:hconv}, Table~\ref{tab:iburg:nrl0:hconv}).

\begin{table}
\centering
\caption{Tabulated convergence results for the \texttt{iburg-acc} test case
 corresponding to the curves in Figure~\ref{fig:iburg:nrl0:hconv}.
 The segment-wise slope of the curves $(m(\cdot))$ tends to at least the
 optimal rate ($p+1$) as the mesh is refined.}
\label{tab:iburg:nrl0:hconv}
\begin{tabular}{ccc|c|cc|cc|cc}
$p$ & $q$ & $|\Ecal_{h}|$ & $h$ & $E_\phi$ & $m(E_\phi)$ & $E_{z_\mathrm{s}}$ & $m(E_{z_\mathrm{s}})$ \\\hline
1 & 1 & 38 & 1.45e-01 & 2.72e-02 & - & 2.32e-03 & - & \\ 
1 & 1 & 152 & 7.25e-02 & 7.18e-03 & 1.92 & 1.09e-03 & 1.09 & \\ 
1 & 1 & 598 & 3.66e-02 & 1.91e-03 & 1.93 & 1.93e-04 & 2.53 & \\ 
1 & 1 & 2392 & 1.83e-02 & 4.69e-04 & 2.03 & 3.92e-05 & 2.30 & \\ 
 \hline
2 & 2 & 38 & 1.45e-01 & 5.68e-03 & - & 4.83e-05 & - & \\ 
2 & 2 & 152 & 7.25e-02 & 9.64e-05 & 5.88 & 2.70e-07 & 7.48 & \\ 
2 & 2 & 608 & 3.63e-02 & 6.36e-06 & 3.92 & 1.20e-08 & 4.49 & \\ 
2 & 2 & 2432 & 1.81e-02 & 8.66e-07 & 2.88 & 7.70e-10 & 3.96 & \\ 
 \hline
3 & 3 & 32 & 1.58e-01 & 1.57e-03 & - & 2.06e-05 & - & \\ 
3 & 3 & 128 & 7.91e-02 & 1.62e-05 & 6.60 & 3.37e-07 & 5.93 & \\ 
3 & 3 & 512 & 3.95e-02 & 4.37e-07 & 5.21 & 5.90e-09 & 5.84 & \\ 
3 & 3 & 2040 & 1.98e-02 & 3.31e-08 & 3.73 & 1.87e-10 & 5.00 & \\ 

\end{tabular}
\end{table}

\subsubsection{Shock formation and merge}
\label{sec:numexp:iburg:form}
Next, we consider a problem in which a smooth initial condition forms
two distinct shock waves that eventually merge (test case: \texttt{iburg-form})
to demonstrate the ability of the HOIST method to handle shock formation and
intersection. Consider Burgers' equation in (\ref{eqn:iburg}) over the one-dimensional
spatial domain $\bar\Omega\coloneqq(-1,1)$ and temporal domain $\Tcal=(0,1)$
with the initial condition
\begin{equation}
 \phi_0: z \mapsto 1.2~\mathrm{exp}\left(-\frac{(z+0.5)^2}{0.025}\right)-\mathrm{exp}\left(-\frac{(z-0.5)^2}{0.025}\right)
\end{equation}
and boundary conditions $\phi(-1,t)=\phi_0(-1)$ and $\phi(1,t)=\phi_0(1)$
for $t\in\Tcal$.

We discretize the space-time domain using a structured mesh of 400
right triangles and, because shock formation requires substantial
resolution to represent the steepening shock, we only consider
polynomial degrees $p=2,3$. Similar to the accelerating shock
case, the shocks are curved in space-time so we choose the mesh
and solution approximation to be the same polynomial degree
($q = p$). The same HOIST parameters are used for both the
$p=2,3$ simulations (Table~\ref{tab:params}).

\begin{figure}[!htbp]
\centering
\ifbool{fastcompile}{}{
\begin{tikzpicture}
\begin{groupplot}[
  group style={
      group size=3 by 2,
      horizontal sep=0.5cm,
      vertical sep=0.5cm
  },
  width=0.36\textwidth,
  axis equal image,
  xlabel={$x_1$},
  ylabel={$x_2$},
  xtick = {-1.0, -0.5, 0, 0.5, 1.0},
  ytick = {0.0, 0.5, 1.0},
  xmin=-1, xmax=1,
  ymin=0, ymax=1
]

\nextgroupplot[xlabel={}, xticklabels={,,}]
\addplot graphics [xmin=-1, xmax=1, ymin=0, ymax=1] {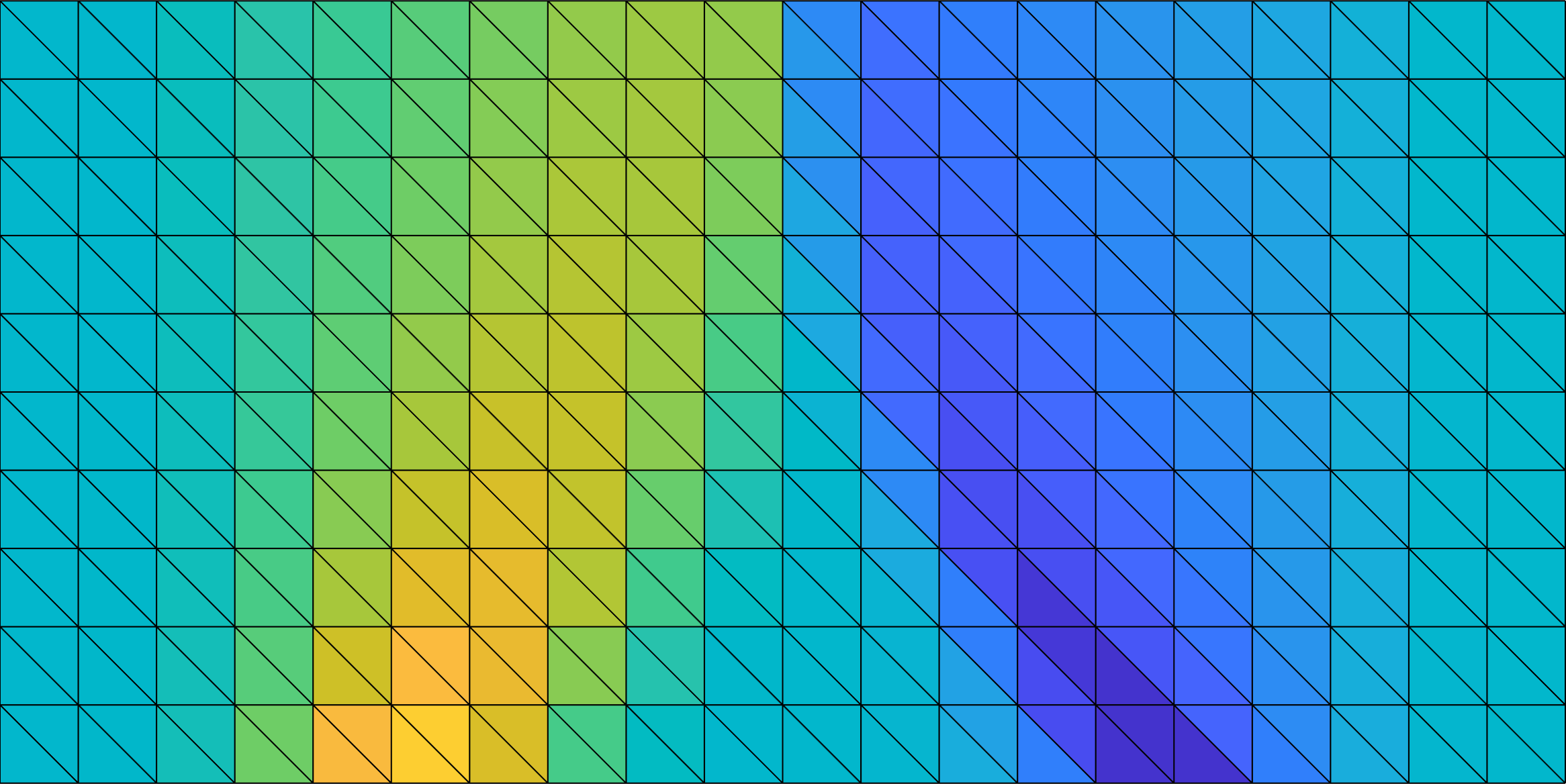};
\nextgroupplot[xlabel={}, ylabel={}, xticklabels={,,}, yticklabels={,,}]
\addplot graphics [xmin=-1, xmax=1, ymin=0, ymax=1] {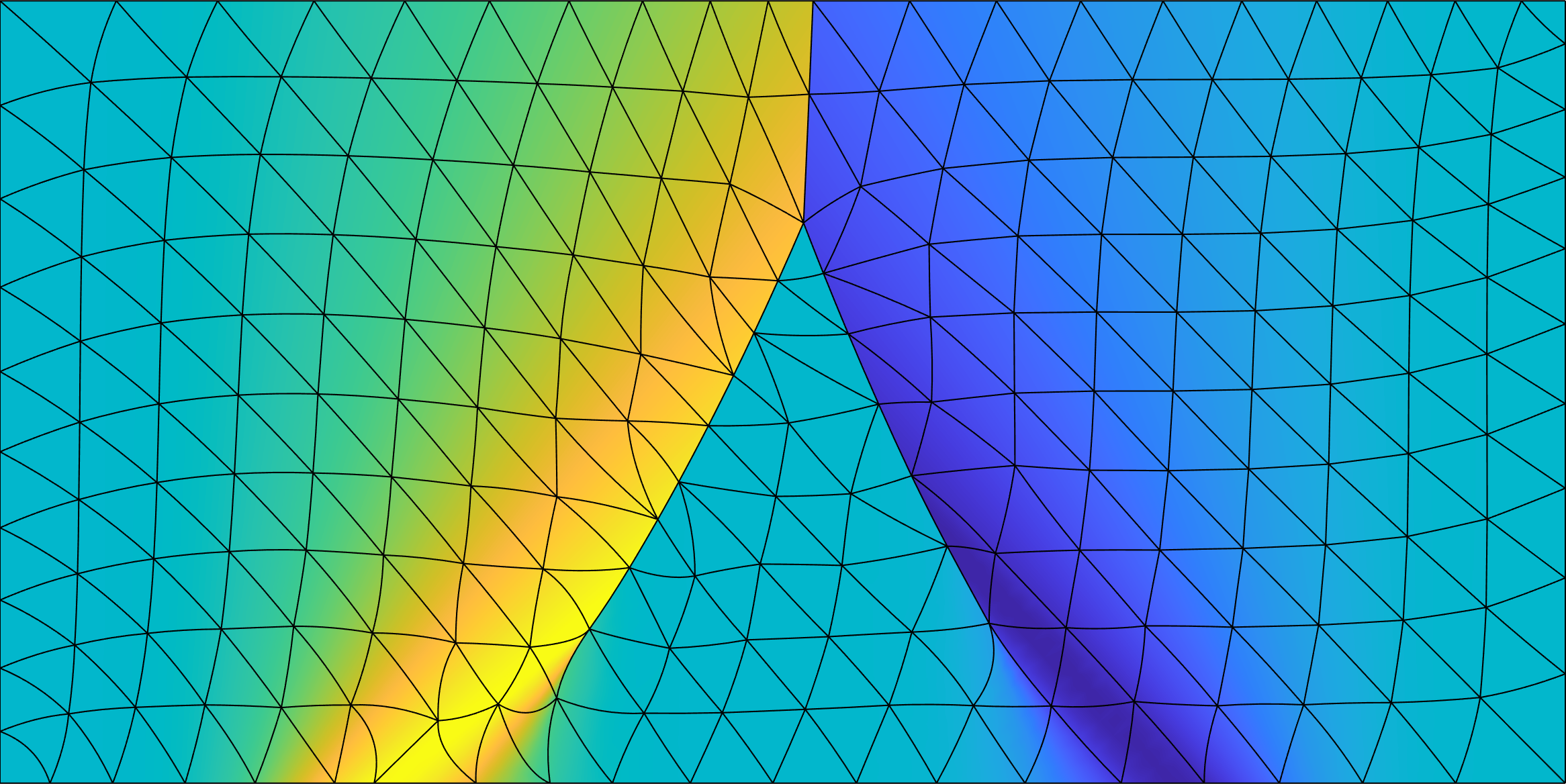};
\nextgroupplot[xlabel={}, ylabel={}, xticklabels={,,}, yticklabels={,,}]
\addplot graphics [xmin=-1, xmax=1, ymin=0, ymax=1] {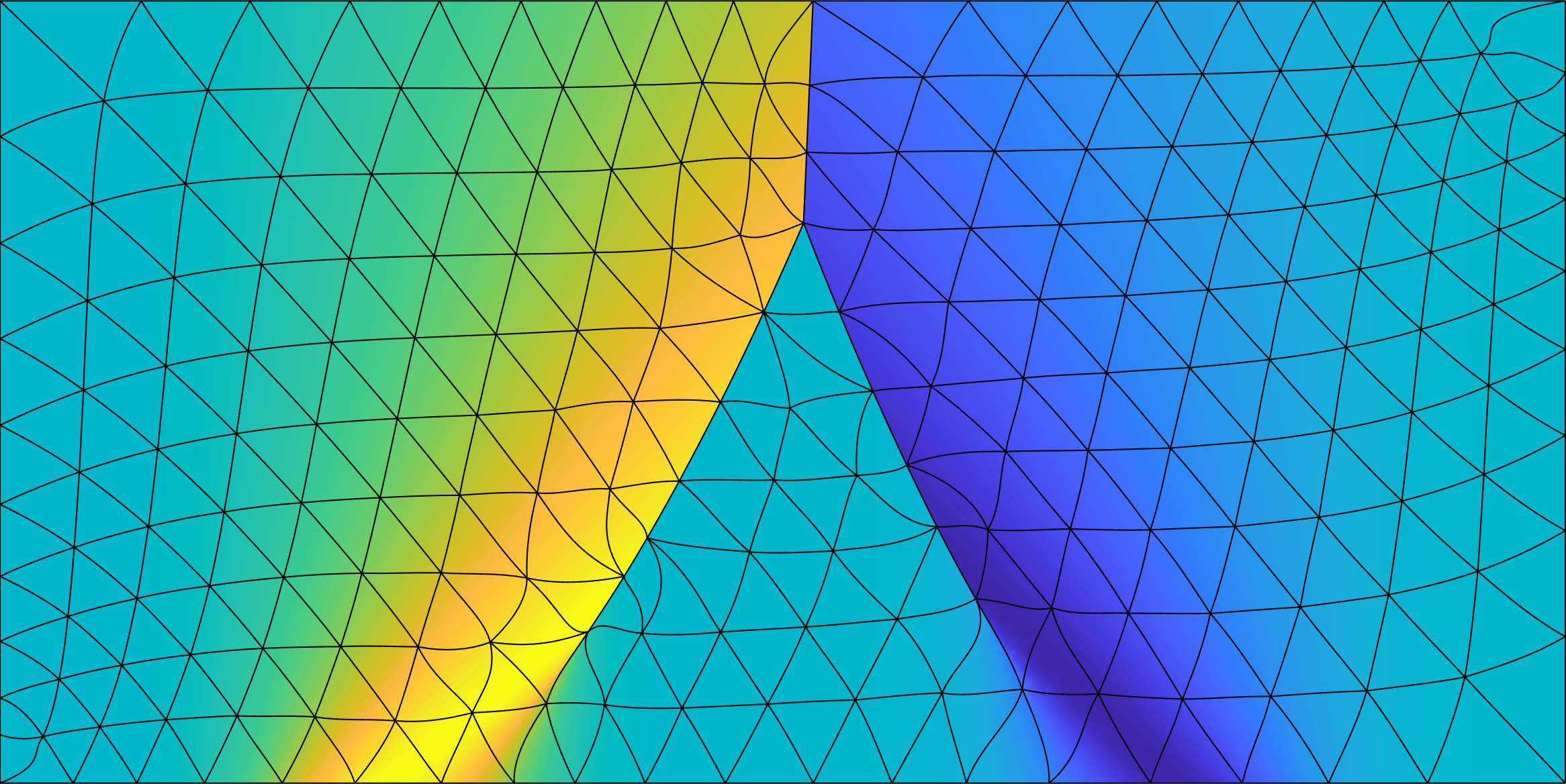};

\nextgroupplot
\addplot graphics [xmin=-1, xmax=1, ymin=0, ymax=1] {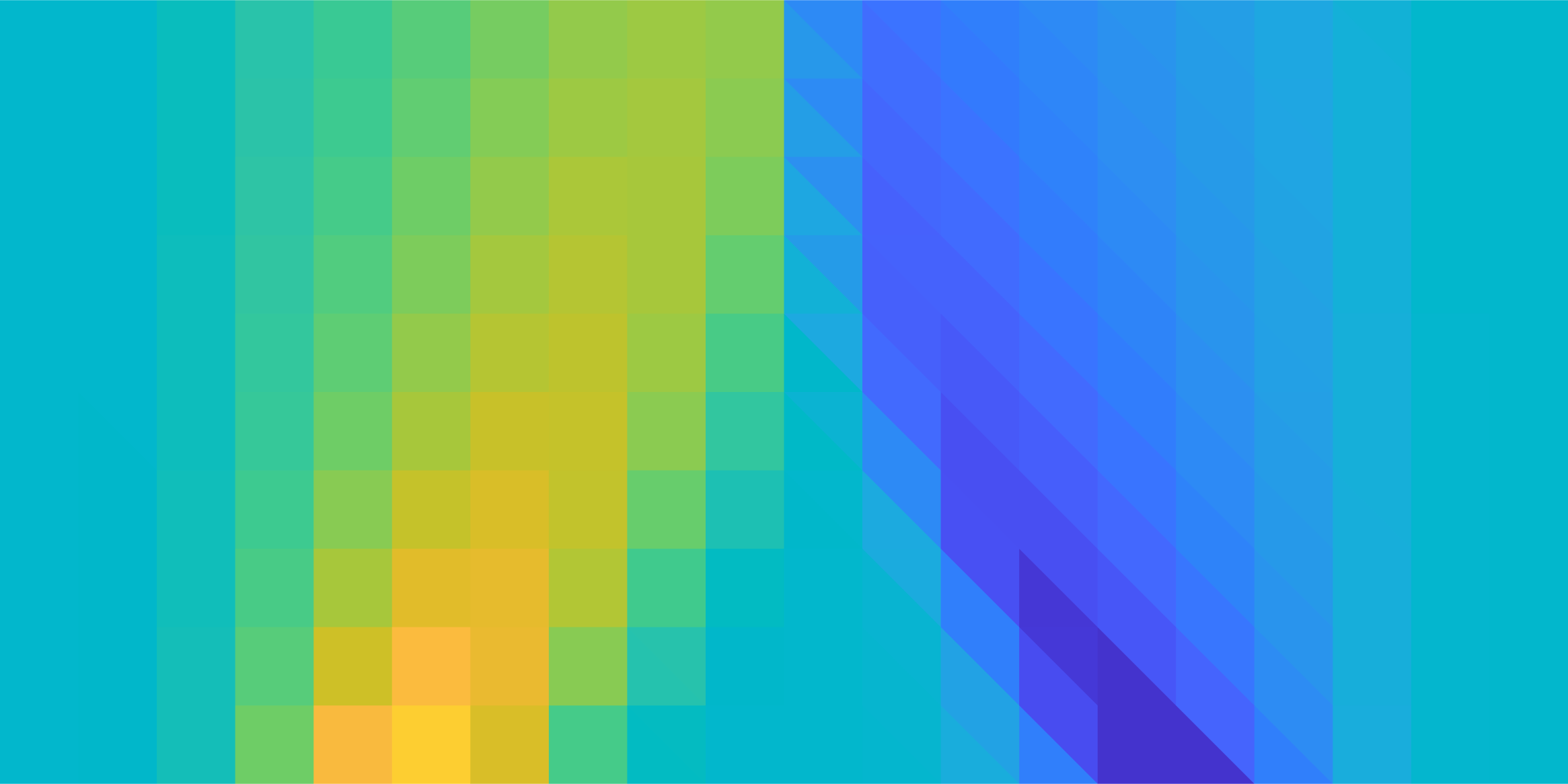};
\nextgroupplot[ylabel={}, yticklabels={,,}]
\addplot graphics [xmin=-1, xmax=1, ymin=0, ymax=1] {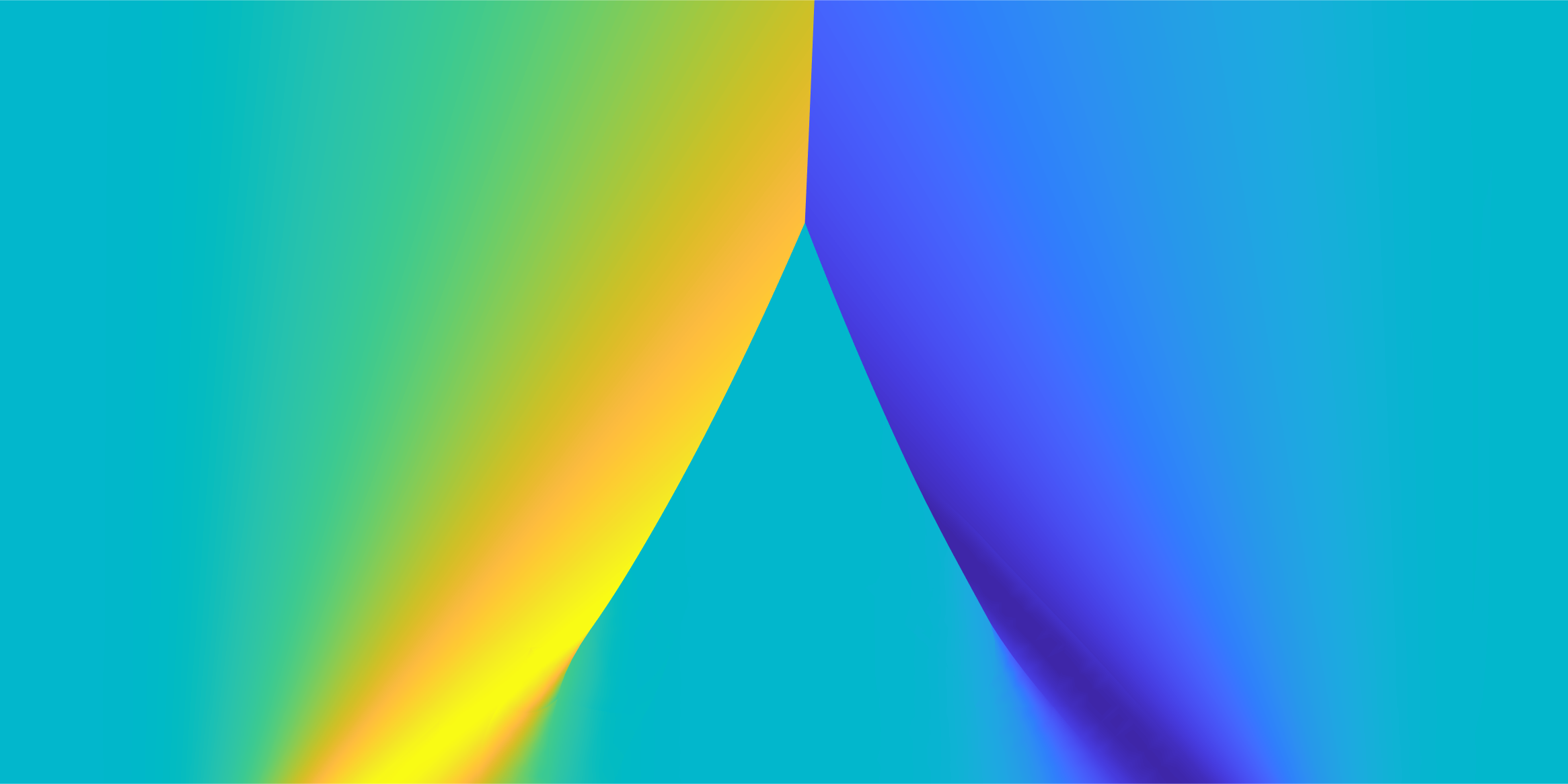};
\nextgroupplot[ylabel={}, yticklabels={,,}]
\addplot graphics [xmin=-1, xmax=1, ymin=0, ymax=1] {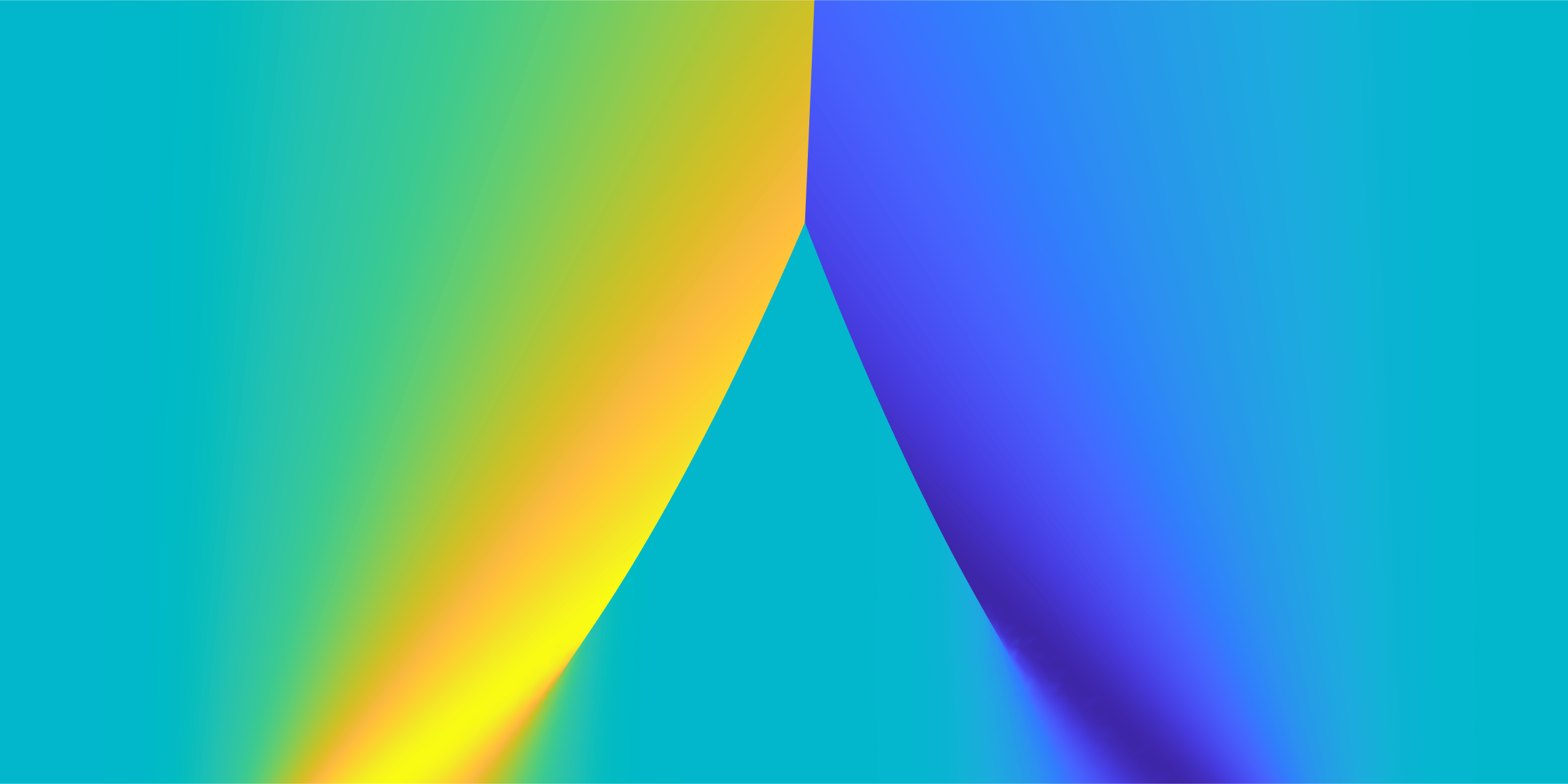};

\end{groupplot}
\end{tikzpicture}

\colorbarMatlabParula{-1}{-0.5}{0.0}{0.6}{1.2}
}
\caption{Starting point (\emph{left}) and HOIST solution for
 polynomial degrees $p=2$ (\textit{middle}) and $p=3$
 (\textit{right}) for the \texttt{iburg-form} test case,
 with (\textit{top row}) and without (\textit{bottom row}) the mesh edges.}
\label{fig:iburg:shkform:sltn}
\end{figure}

\begin{figure}[!htbp]
\centering
\ifbool{fastcompile}{}{
\input{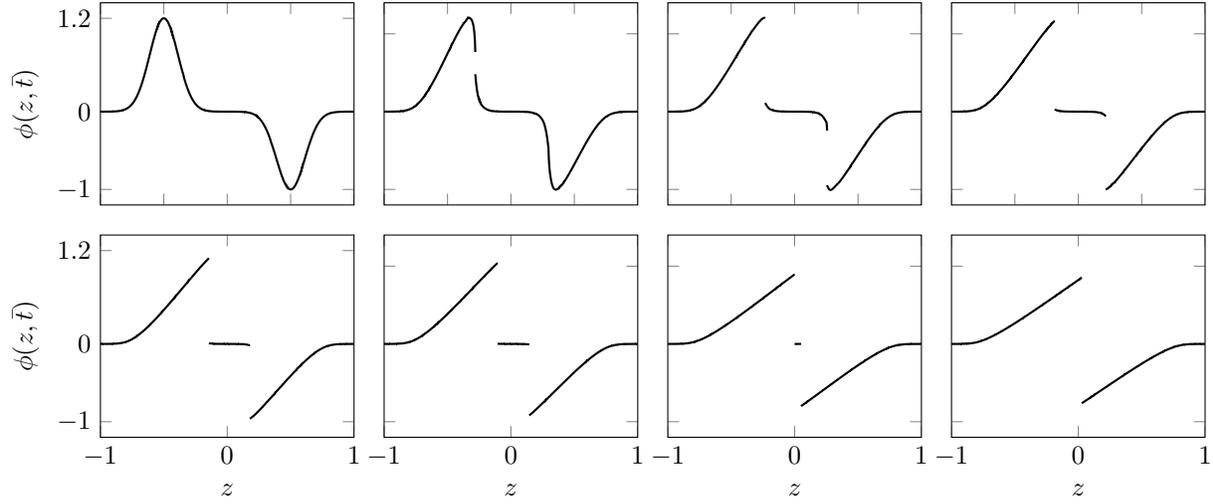}
}
\caption{Temporal slices of the $p=3$ HOIST solution for the
 \texttt{iburg-form} test case at time instances
 $\bar{t}=0, 0.146, 0.218, 0.364, 0.436, 0.655, 0.727$
 (\textit{left-to-right}, \textit{top-to-bottom}). The time instances
 are not uniformly spaced to show relevant features of solution.}
\label{fig:iburg:shkform:shk}
\end{figure}

For both polynomial degrees, the HOIST method is able to track
the shocks during and after formation, including the point at
which they merge (a triple point in space-time) despite
starting from a mesh that incorporates no knowledge of the discontinuity
surface (Figure~\ref{fig:iburg:shkform:sltn}). In space-time, shock
formation essentially requires infinite resolution because the smooth
features continually steepen until the instant at which they become
discontinuities. For this reason, any space-time discretization will
be under-resolved in the vicinity of shock formation. In the region
where the steepness of the forming shock has exceeded the resolution
in the discretization, the HOIST method approximates the steep,
continuous feature with a discontinuity that grows until the
shock is fully formed (Figure~\ref{fig:iburg:shkform:shk}).
Despite the features present in this problem (shock formation
and a triple point), the proposed solver converges rapidly to
a critical point, especially the $p=2$ simulation
(Figure~\ref{fig:iburg:shkform:hist}). The $p=3$ simulation
is slightly slower because the solution becomes more oscillatory
at intermediate iterations, which requires more re-initializations
(jumps in the convergence plots in Figure~\ref{fig:iburg:shkform:hist}).

\begin{figure}[!htbp]
\centering
\ifbool{fastcompile}{}{
\input{_py/iburg1v1d_sptm_shkform_hist.tikz}
}
\caption{SQP convergence history and behavior of adaptive parameters for
 the \texttt{iburg-form} test case (legend in Table~\ref{tab:legend})
 for $p=2$ (\textit{top}) and $p=3$ (\textit{bottom}).}
\label{fig:iburg:shkform:hist}
\end{figure}

\subsection{Inviscid, compressible flow through area variation}
\label{sec:numexp:eulernoz}
Next, we consider inviscid, compressible flow through an area
variation (test case: \texttt{nozzle}) governed by the quasi-one-dimensional
Euler equations
\begin{equation} \label{eqn:eulernoz}
\begin{split}
\pder{}{x}\left(A(x)\rho(x) v(x)\right) &= 0 \\
\pder{}{x}\left(A(x)\left[\rho(x)v(x)^2+P(x)\right]\right) &= P(x)\pder{A}{x}(x) \\
\pder{}{x}\left(A(x)\left[\rho(x)E(x)+P(x)\right]v(x)\right) &= 0,
\end{split}
\end{equation}
for $x\in\Omega\subset(0,10)$, where $\func{\rho}{\Omega}{\Rbb_{>0}}$
is the density of the fluid, $\func{v}{\Omega}{\Rbb}$ is the velocity
of the fluid, and $\func{E}{\Omega}{\Rbb_{>0}}$ is the total energy of
the fluid, implicitly defined as the solution of (\ref{eqn:eulernoz}).
For a calorically ideal fluid, the pressure of the fluid,
$\func{P}{\Omega}{\Rbb_{>0}}$, is related to the energy via
the ideal gas law
\begin{equation}
 P = (\gamma-1)\left(\rho E-\rho v^2/2\right),
\end{equation}
where $\gamma\in\Rbb_{>0}$ is the ratio of specific heats. The
area beneath the nozzle is given by $\func{A}{\Omega}{\Rbb}$,
which we take to be
\begin{equation}
 A : x \mapsto
 \mu_1+(\mu_2-\mu_1)\frac{(10-x)x}{25},
\end{equation}
which has an exit and throat area of $\mu_1$ and $\mu_2$, respectively,
with the throat located midway through the domain ($x=5$); in this work,
we take $\mu_1=3$, $\mu_2=1$. The boundary conditions are $\rho(0)=1$ 
and $\rho(10) = 0.7$, which will lead to a shock at $x_\mathrm{s} = 7.94$.
By combining the area, density, momentum, and energy into a vector
of conservative variables $\func{U}{\Omega}{\Rbb^3}$, defined as
\begin{equation}
 U : x \mapsto \begin{bmatrix} A(x)\rho(x) \\ A(x)\rho(x)v(x) \\ A(x) \rho(x)E(x) \end{bmatrix},
\end{equation}
the quasi-one-dimensional Euler equations in (\ref{eqn:eulernoz}) are
a conservation law of the form (\ref{eqn:claw-phys}). For the DG
discretization, we use a smoothed variant of Roe's flux described
in \cite{zahr_high-order_2020} with Harten-Hyman entropy fix
\cite{harten_self_1983} as the
inviscid numerical flux function. The purpose of this numerical experiment
is to verify optimal convergence rate of the HOIST method for a compressible
flow and provide a direct comparison to popular shock capturing methods
based on artificial viscosity.

\begin{figure}[!htbp]
\centering
\input{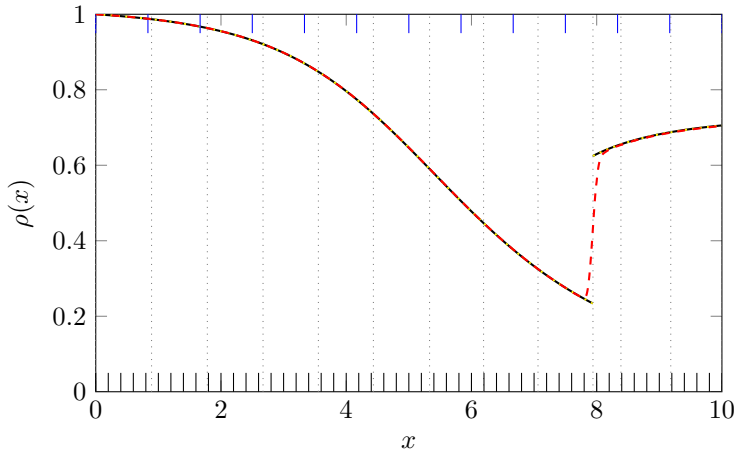}
\caption{The solution (density) of the \texttt{nozzle} test case
 (\ref{line:eulernoz:exact}) and the corresponding HOIST
 (\ref{line:eulernoz:impshktrk-soln}) and shock capturing
 (\ref{line:eulernoz:shkcapt-soln}) approximations. The HOIST simulation
 uses 12 quadratic elements; the element boundaries of the initial
 (\ref{line:eulernoz:shktrk-msh0}) and final (\ref{line:eulernoz:shktrk-msh1})
 mesh are included to show the initial mesh is far from alignment with the
 discontinuity. The shock capturing simulation uses 50 quartic elements
 with the element boundaries indicated with (\ref{line:eulernoz:shkcapt-msh}).}
\label{fig:eulernozii1:sltn}
\end{figure}

\begin{figure}[!htbp]
\centering
\ifbool{fastcompile}{}{
\input{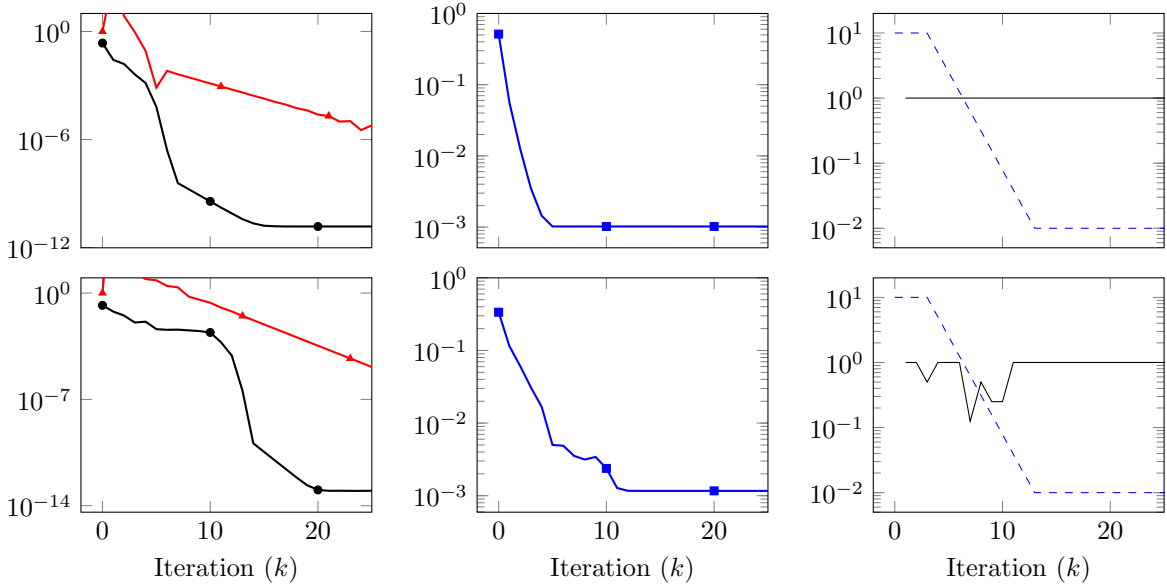}
}
\caption{SQP convergence history and behavior of adaptive parameters for
 \texttt{nozzle} test case (legend in Table~\ref{tab:legend}).
 Only the $p=3$ (\textit{top}) and $p=5$ (\textit{bottom}) simulations are
 included for brevity. The mesh quality term is not used for one-dimensional
 problems and not included in the convergence plots.}
\label{fig:eulernozii1:hist}
\end{figure}

We discretize the one-dimensional domain with $12$ elements of
degree $q=1$ (no benefit to using high-order geometry approximation in
one dimension) and use the HOIST method with polynomial
degrees $p=1,\dots, 5$ (HOIST parameters in Table~\ref{tab:params},
except $c_8=10^{-1}$ in the case $p=1$). 
The $p=2$ HOIST solution on this coarse
discretization (12 elements) provides an excellent approximation
that is indistinguishable from the exact solution to
(\ref{eqn:eulernoz}). The shock capturing method based on
artificial viscosity in \cite{persson_sub-cell_2006} requires
50 elements with a $p=4$ approximation to achieve a similar solution,
although the discontinuity is approximated by a steep gradient, which
is a noticeable deviation from the exact solution
(Figure~\ref{fig:eulernozii1:sltn}). Due to the geometric
simplicity in the one-dimensional setting, the HOIST solvers
converge rapidly, achieving tight tolerances on the optimality
system in only 20 iterations (Figure~\ref{fig:eulernozii1:hist}).

\begin{figure}[!htbp]
\centering
\ifbool{fastcompile}{}{
\begin{tikzpicture}
\begin{groupplot} [
group style={group size = 2 by 1, horizontal sep = 2cm}]
\nextgroupplot[width=0.42\textwidth, ylabel={$E_\rho$}, xlabel={$h$}, xmode=log, ymode=log]
\addplot [solid, thick, black, mark options={solid, thin}, mark=*, mark size=1.5]
coordinates {
( 7.14285714e-01,  2.84445924e-02)
( 3.57142857e-01,  1.05595046e-02)
( 1.78571429e-01,  2.68847047e-03)
( 8.92857143e-02,  6.68178339e-04)
( 4.46428571e-02,  1.67056291e-04)
( 2.23214286e-02,  4.17796878e-05)
( 1.11607143e-02,  1.04369598e-05)};\label{line:eulernoz:hconv:p1}

\addplot [dashed, black, forget plot]
coordinates {
( 5.58035714e-03,  2.61123048e-06)
( 1.42857143e+00,  1.71129601e-01)};

\addplot [solid, thick, blue, mark options={solid, thin}, mark=square*, mark size=1.5]
coordinates {
( 0.00000000e+00,  0.00000000e+00)
( 8.33333333e-01,  1.38840006e-03)
( 4.16666667e-01,  1.45087381e-04)
( 2.08333333e-01,  1.70198464e-05)
( 1.04166667e-01,  2.12526949e-06)
( 5.20833333e-02,  2.97230175e-07)
( 2.60416667e-02,  6.18765737e-08)};\label{line:eulernoz:hconv:p2}

\addplot [dashed, blue, forget plot]
coordinates {
( 1.30208333e-02,  4.64422149e-09)
( 1.66666667e+00,  9.73963838e-03)};

\addplot [solid, thick, red, mark options={solid, thin}, mark=triangle*, mark size=1.5]
coordinates {
( 0.00000000e+00,  0.00000000e+00)
( 0.00000000e+00,  0.00000000e+00)
( 8.33333333e-01,  1.01922369e-04)
( 4.16666667e-01,  6.08733329e-06)
( 2.08333333e-01,  3.52744406e-07)
( 1.04166667e-01,  2.11117676e-08)
( 5.20833333e-02,  1.28902066e-09)};\label{line:eulernoz:hconv:p3}

\addplot [dashed, red, forget plot]
coordinates {
( 2.60416667e-02,  8.24678420e-11)
( 1.66666667e+00,  1.38358080e-03)};

\addplot [solid, thick, magenta, mark options={solid, thin}, mark=pentagon*, mark size=1.5]
coordinates {
( 0.00000000e+00,  0.00000000e+00)
( 0.00000000e+00,  0.00000000e+00)
( 8.33333333e-01,  8.52605642e-06)
( 4.16666667e-01,  2.83387066e-07)
( 2.08333333e-01,  8.95335436e-09)
( 1.04166667e-01,  2.67280388e-10)
( 5.20833333e-02,  8.27699738e-12)};\label{line:eulernoz:hconv:p4}

\addplot [dashed, magenta, forget plot]
coordinates {
( 2.60416667e-02,  2.61016004e-13)
( 1.66666667e+00,  2.80263800e-04)};

\addplot [solid, thick, cyan, mark options={solid, thin}, mark=diamond*, mark size=1.5]
coordinates {
( 0.00000000e+00,  0.00000000e+00)
( 0.00000000e+00,  0.00000000e+00)
( 0.00000000e+00,  0.00000000e+00)
( 1.00000000e+00,  1.23289534e-06)
( 5.00000000e-01,  2.04281963e-08)
( 2.50000000e-01,  3.28501618e-10)
( 1.25000000e-01,  4.81741162e-12)};\label{line:eulernoz:hconv:p5}

\addplot [dashed, cyan, forget plot]
coordinates {
( 6.25000000e-02,  8.02005903e-14)
( 2.00000000e+00,  8.61147281e-05)};

\addplot [solid, thick, gray]
coordinates {
( 2.00000000e-01,  3.55796456e-02)
( 1.00000000e-01,  2.92728555e-02)
( 5.00000000e-02,  1.89001940e-02)
( 2.50000000e-02,  1.90178172e-02)
( 1.25000000e-02,  1.45970828e-02)};\label{line:eulernoz:hconv:shkcapt}

\nextgroupplot[width=0.42\textwidth, ylabel={$E_{x_\mathrm{s}}$}, xlabel={$h$}, xmode=log, ymode=log]
\addplot [solid, thick, black, mark options={solid, thin}, mark=*, mark size=1.5, forget plot]
coordinates {
( 7.14285714e-01,  2.35999746e-02)
( 3.57142857e-01,  9.15325399e-03)
( 1.78571429e-01,  2.41767072e-03)
( 8.92857143e-02,  6.00163677e-04)
( 4.46428571e-02,  1.49447576e-04)
( 2.23214286e-02,  3.73686825e-05)
( 1.11607143e-02,  9.31146813e-06)};

\addplot [dashed, black, forget plot]
coordinates {
( 5.58035714e-03,  2.33554266e-06)
( 1.42857143e+00,  1.53062124e-01)};

\addplot [solid, thick, blue, mark options={solid, thin}, mark=square*, mark size=1.5, forget plot]
coordinates {
( 0.00000000e+00,  0.00000000e+00)
( 8.33333333e-01,  8.80438915e-04)
( 4.16666667e-01,  8.58904049e-05)
( 2.08333333e-01,  1.05090905e-05)
( 1.04166667e-01,  1.36601505e-06)
( 5.20833333e-02,  1.83633153e-07)
( 2.60416667e-02,  2.42294522e-08)};

\addplot [dashed, blue, forget plot]
coordinates {
( 1.30208333e-02,  2.86926802e-09)
( 1.66666667e+00,  6.01729116e-03)};

\addplot [solid, thick, red, mark options={solid, thin}, mark=triangle*, mark size=1.5, forget plot]
coordinates {
( 0.00000000e+00,  0.00000000e+00)
( 0.00000000e+00,  0.00000000e+00)
( 8.33333333e-01,  5.83401967e-05)
( 4.16666667e-01,  3.25033634e-06)
( 2.08333333e-01,  1.75416994e-07)
( 1.04166667e-01,  1.00399316e-08)
( 5.20833333e-02,  6.72247147e-10)};

\addplot [dashed, red, forget plot]
coordinates {
( 2.60416667e-02,  3.92184826e-11)
( 1.66666667e+00,  6.57976954e-04)};

\addplot [solid, thick, magenta, mark options={solid, thin}, mark=pentagon*, mark size=1.5, forget plot]
coordinates {
( 0.00000000e+00,  0.00000000e+00)
( 0.00000000e+00,  0.00000000e+00)
( 8.33333333e-01,  4.18421365e-06)
( 4.16666667e-01,  7.40548538e-08)
( 2.08333333e-01,  1.87963156e-09)
( 1.04166667e-01,  5.49515988e-11)
( 5.20833333e-02,  2.43538523e-12)};

\addplot [dashed, magenta, forget plot]
coordinates {
( 2.60416667e-02,  5.36636707e-14)
( 1.66666667e+00,  5.76209277e-05)};

\addplot [solid, thick, cyan, mark options={solid, thin}, mark=diamond*, mark size=1.5, forget plot]
coordinates {
( 0.00000000e+00,  0.00000000e+00)
( 0.00000000e+00,  0.00000000e+00)
( 0.00000000e+00,  0.00000000e+00)
( 1.00000000e+00,  6.67589684e-07)
( 5.00000000e-01,  4.30910063e-09)
( 2.50000000e-01,  4.50404158e-11)
( 1.25000000e-01,  6.21724894e-14)};

\addplot [dashed, cyan, forget plot]
coordinates {
( 6.25000000e-02,  1.09961953e-14)
( 2.00000000e+00,  1.18070748e-05)};

\end{groupplot}\end{tikzpicture}
}
\caption{$h$-convergence of the HOIST method for the \texttt{nozzle}
 test case for both error metrics, $E_\rho$ (\textit{left})
 and $E_{x_\mathrm{s}}$ (\textit{right}), for polynomial degrees
 $p=1$ (\ref{line:eulernoz:hconv:p1}),
 $p=2$ (\ref{line:eulernoz:hconv:p2}),
 $p=3$ (\ref{line:eulernoz:hconv:p3}),
 $p=4$ (\ref{line:eulernoz:hconv:p4}), and
 $p=5$ (\ref{line:eulernoz:hconv:p5}), as well as the shock capturing method of
 \cite{persson_sub-cell_2006} for polynomial degree $p=2$
 (\ref{line:eulernoz:hconv:shkcapt}). The dashed lines indicate
 the optimal convergence rate ($p+1$).}
\label{fig:eulernozii1:hconv}
\end{figure}
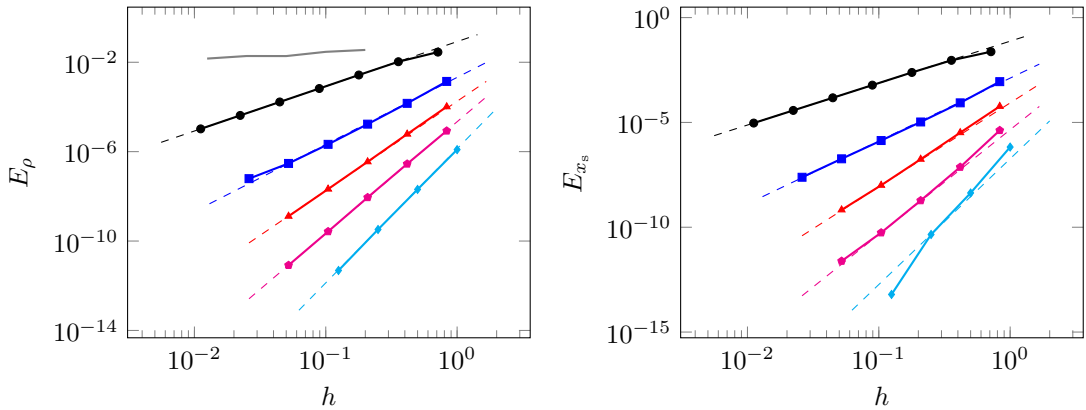

Finally, we study the $h$-convergence of the HOIST method for this
problem and compare to the shock capturing method in
\cite{persson_sub-cell_2006}. We consider two error metrics:
$\func{E_\rho}{\Rbb^{N_\ubm}\times\Rbb^{N_\xbm}}{\Rbb}$, the
$L^1$ error in the density, and $\func{E_{x_\mathrm{s}}}{\Rbb^{N_\xbm}}{\Rbb}$,
the error in the shock position, which are defined as
\begin{equation}
 E_\rho : (\ubm,\xbm) \mapsto
 \int_0^{10} \left|\rho(x)-\hat\rho(x;\ubm,\xbm)\right|\, dx,
 \qquad
 E_{x_\mathrm{s}} : \xbm \mapsto
 \left|x_\mathrm{s}-\hat{x}_\mathrm{s}(\xbm)\right|,
\end{equation}
where $\func{\hat\rho}{\Omega\times\Rbb^{N_\ubm}\times\Rbb^{N_\xbm}}{\Rbb}$
and $\func{\hat{x}_\mathrm{s}}{\Rbb^{N_\xbm}}{\Rbb}$ are the HOIST approximations
to $\rho$ and the shock position, respectively.
Because the shock position is not well-defined for the shock capturing method,
we only consider its $E_\rho$ error. We consider polynomial degrees
$p=1,\dots,5$ with at least four levels of uniform refinement
(up to seven levels for the lower polynomial degrees). The
HOIST method achieves optimal convergence rates,
$\Ocal(h^{p+1})$, for all polynomial degrees considered
for both error metrics, while the shock capturing method
is limited to sub-first-order convergence rate
(Figure~\ref{fig:eulernozii1:hconv}, Table~\ref{tab:eulernozii1:hconv}).
The slow convergence rate of the shock capturing method is
due to the smooth approximation to the discontinuity and
downstream corruption \cite{lee1999osc, bonfiglioli2014conv}. 

\begin{table}
\centering
\caption{Tabulated convergence results for the \texttt{nozzle}
 test case in Figure~\ref{fig:eulernozii1:hconv}. The
 segment-wise slope of the curves $(m(\cdot))$ tends to at
 least the optimal rate ($p+1$) as the mesh is refined.}
\label{tab:eulernozii1:hconv}
\begin{tabular}{ccc|c|cc|cc|cc}
$p$ & $q$ & $|\Ecal_h|$ & $h$ & $E_\rho$ & $m(E_\rho)$ & $E_{x_\mathrm{s}}$ & $m(E_{x_\mathrm{s}})$ \\\hline
1 & 1 & 112 & 8.93e-02 & 6.68e-04 & - & 6.00e-04 & - & \\ 
1 & 1 & 224 & 4.46e-02 & 1.67e-04 & 2.00 & 1.49e-04 & 2.01 & \\ 
1 & 1 & 448 & 2.23e-02 & 4.18e-05 & 2.00 & 3.74e-05 & 2.00 & \\ 
1 & 1 & 896 & 1.12e-02 & 1.04e-05 & 2.00 & 9.31e-06 & 2.00 & \\ 
 \hline
2 & 2 & 48 & 2.08e-01 & 1.70e-05 & - & 1.05e-05 & - & \\ 
2 & 2 & 96 & 1.04e-01 & 2.13e-06 & 3.00 & 1.37e-06 & 2.94 & \\ 
2 & 2 & 192 & 5.21e-02 & 2.97e-07 & 2.84 & 1.84e-07 & 2.90 & \\ 
2 & 2 & 384 & 2.60e-02 & 6.19e-08 & 2.26 & 2.42e-08 & 2.92 & \\ 
 \hline
3 & 3 & 24 & 4.17e-01 & 6.09e-06 & - & 3.25e-06 & - & \\ 
3 & 3 & 48 & 2.08e-01 & 3.53e-07 & 4.11 & 1.75e-07 & 4.21 & \\ 
3 & 3 & 96 & 1.04e-01 & 2.11e-08 & 4.06 & 1.00e-08 & 4.13 & \\ 
3 & 3 & 192 & 5.21e-02 & 1.29e-09 & 4.03 & 6.72e-10 & 3.90 & \\ 
 \hline
4 & 4 & 24 & 4.17e-01 & 2.83e-07 & - & 7.41e-08 & - & \\ 
4 & 4 & 48 & 2.08e-01 & 8.95e-09 & 4.98 & 1.88e-09 & 5.30 & \\ 
4 & 4 & 96 & 1.04e-01 & 2.67e-10 & 5.07 & 5.50e-11 & 5.10 & \\ 
4 & 4 & 192 & 5.21e-02 & 8.28e-12 & 5.01 & 2.44e-12 & 4.50 & \\ 
 \hline
5 & 5 & 10 & 1.00e+00 & 1.23e-06 & - & 6.68e-07 & - & \\ 
5 & 5 & 20 & 5.00e-01 & 2.04e-08 & 5.92 & 4.31e-09 & 7.28 & \\ 
5 & 5 & 40 & 2.50e-01 & 3.29e-10 & 5.96 & 4.50e-11 & 6.58 & \\ 
5 & 5 & 80 & 1.25e-01 & 4.82e-12 & 6.09 & 6.22e-14 & 9.50 & \\ 

\end{tabular}
\end{table}


\subsection{Inviscid, compressible flow}
\label{sec:numexp:euler}
Finally, we consider flow of an inviscid, compressible fluid
through a domain $\bar\Omega\subset\Rbb^d$, which are governed
by the Euler equations
\begin{equation} \label{eqn:euler}
\begin{split}
\pder{}{t}\rho(z,t) + \pder{}{z_j}\left(\rho(z,t) v_j(z,t)\right) &= 0 \\
\pder{}{t}\left(\rho(z,t)v_i(z,t)\right) +
\pder{}{z_j}\left(\rho(z,t) v_i(z,t)v_j(z,t)+P(z,t)\delta_{ij}\right) &= 0 \\
\pder{}{t}\left(\rho(z,t)E(z,t)\right) +
\pder{}{z_j}\left(\left[\rho(z,t)E(z,t)+P(z,t)\right]v_j(z,t)\right) &= 0
\end{split}
\end{equation}
for all $z\in\bar\Omega$ (spatial domain) and $t\in\Tcal\subset\Rbb$
(temporal domain), where $i=1,\dots,d$ and summation is implied over
the repeated index $j=1,\dots,d$. The density of the fluid
$\func{\rho}{\bar\Omega\times\Tcal}{\Rbb_{>0}}$, the velocity of
the fluid $\func{v_i}{\bar\Omega\times\Tcal}{\Rbb}$ in the $z_i$
direction for $i=1,\dots,d$, and the total energy of the fluid
$\func{E}{\bar\Omega\times\Tcal}{\Rbb_{>0}}$ are implicitly defined
as the solution of (\ref{eqn:euler}).
For a calorically ideal fluid, the pressure of the fluid,
$\func{P}{\bar\Omega\times\Tcal}{\Rbb_{>0}}$, is related to the energy
via the ideal gas law
\begin{equation}
 P = (\gamma-1)\left(\rho E - \frac{\rho v_i v_i}{2}\right),
\end{equation}
where $\gamma\in\Rbb_{>0}$ is the ratio of specific heats.
This leads to a conservation law of the form (\ref{eqn:claw-phys})
over the space-time domain
$\Omega\coloneqq\bar\Omega\times\Tcal\subset\Rbb^{d+1}$
by combining the density, momentum, and energy into a vector of
conservative variables $\func{U}{\Omega}{\Rbb^{d+2}}$, defined as
\begin{equation}
 U : x \mapsto
 \begin{bmatrix}
  \rho((x_1,\dots,x_d),x_{d+1}) \\
  \rho((x_1,\dots,x_d),x_{d+1})v((x_1,\dots,x_{d+1}),x_{d+1}) \\
  \rho((x_1,\dots,x_d),x_{d+1})E((x_1,\dots,x_{d+1}),x_{d+1})
 \end{bmatrix}.
\end{equation}
For the DG discretization, we use a smoothed variant of Roe's flux described
in \cite{zahr_high-order_2020} with Harten-Hyman entropy fix
\cite{harten_self_1983} as the inviscid numerical flux function.
For steady problems, the temporal domain is empty ($\Tcal=\emptyset$),
all variables are independent of time, the space-time and spatial domains
are the same $d$-dimensional set ($\Omega=\bar\Omega$), and the vector of
conservative variables is
\begin{equation}
 U : x \mapsto
 \begin{bmatrix}
  \rho(x) \\ \rho(x) v(x) \\ \rho(x) E(x)
 \end{bmatrix}.
\end{equation}
\subsubsection{Sod's shock tube}
\label{sec:numexp:euler:sod}
Sod's shock tube (test case: \texttt{sod}) is a Riemann problem for the
Euler equations that models
an idealized shock tube where the membrane separating a high pressure region
from a low pressure one is instantaneously removed. This is a commonly used
validation problem since it has an analytical solution that features a shock
wave, rarefaction wave, and contact discontinuity. The flow domain is
$\bar\Omega=(0, 1)$, the time domain is $\Tcal=(0,0.2)$, the initial
condition is given in terms of the density, velocity, and pressure as
\begin{equation}
 \rho(z,0) = \begin{cases} 1 & z<0.5 \\ 0.125 & z \geq 0.5 \end{cases}, \quad
 v(z,0) = 0, \quad
 P(z,0) = \begin{cases} 1 & z<0.5 \\ 0.1 & z \geq 0.5, \end{cases}
\end{equation}
and the density, velocity, and pressure are prescribed at $z=0$ and the
velocity is prescribed at $z=1$ (values can be read from the initial condition).
The solution of this problem contains three waves (shock, contact, rarefaction)
that emanate from $z=0.5$ and move at different speeds, which is a generalized
triple point in space-time. The purpose of this numerical experiment
is to demonstrate the HOIST method is able to track various types
of features that are discontinuous (shocks and contacts) or
non-smooth (head and tail of rarefactions) and locate triple
points without seeding with \textit{a priori} information.

\begin{figure}[!htbp]
\centering
\ifbool{fastcompile}{}{
\begin{tikzpicture}
\begin{groupplot}[
  group style={
      group size=2 by 1,
      horizontal sep=0.5cm
  },
  width=0.54\textwidth,
  axis equal image,
  xlabel={$x_1$},
  ylabel={$x_2$},
  xtick = {0, 0.5, 1.0},
  ytick = {0.0, 0.1, 0.2},
  xmin=0, xmax=1,
  ymin=0, ymax=0.2
]
\nextgroupplot
\addplot graphics [xmin=0, xmax=1, ymin=0, ymax=0.2] {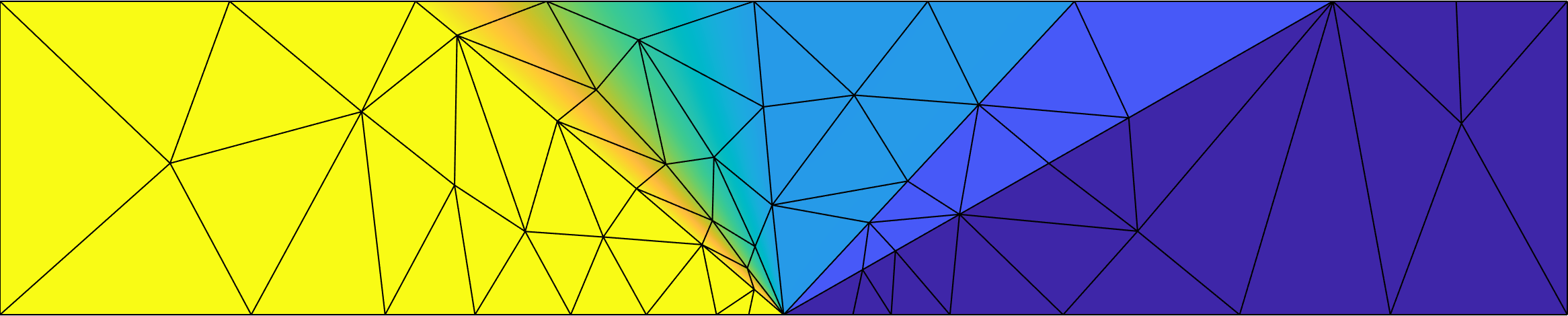};
\nextgroupplot[ylabel={}, yticklabels={,,}]
\addplot graphics [xmin=0, xmax=1, ymin=0, ymax=0.2] {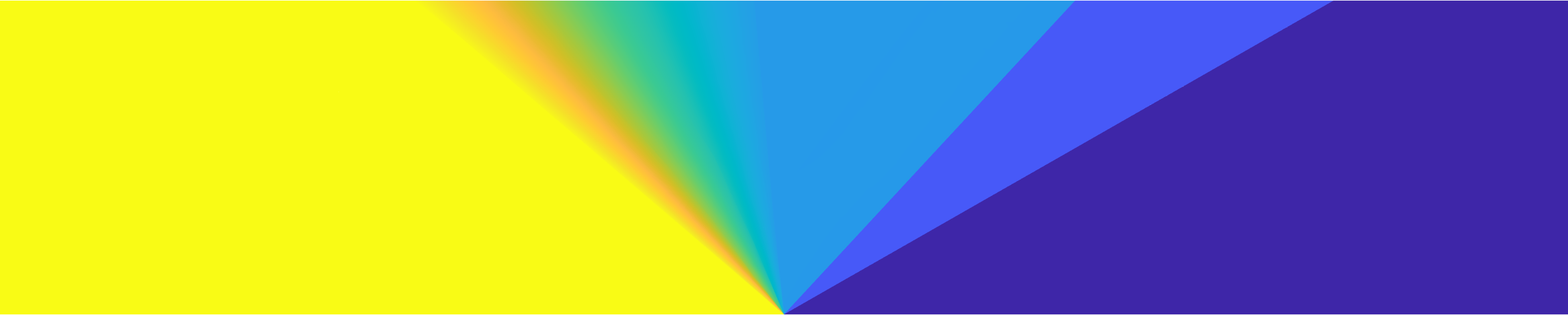};
\end{groupplot}
\end{tikzpicture}
}
\colorbarMatlabParula{0.125}{0.25}{0.5}{0.75}{1}
\caption{HOIST solution to the \texttt{sod} test case with (\textit{left}) and
 without (\textit{right}) mesh. All non-smooth features (shock, contact,
 head and tail of rarefaction) are tracked by the mesh, which leads to
 a well-resolved solution on the coarse mesh.}
\label{fig:euler:sod:final}
\end{figure}

\begin{figure}[!htbp]
\centering
\ifbool{fastcompile}{}{
\input{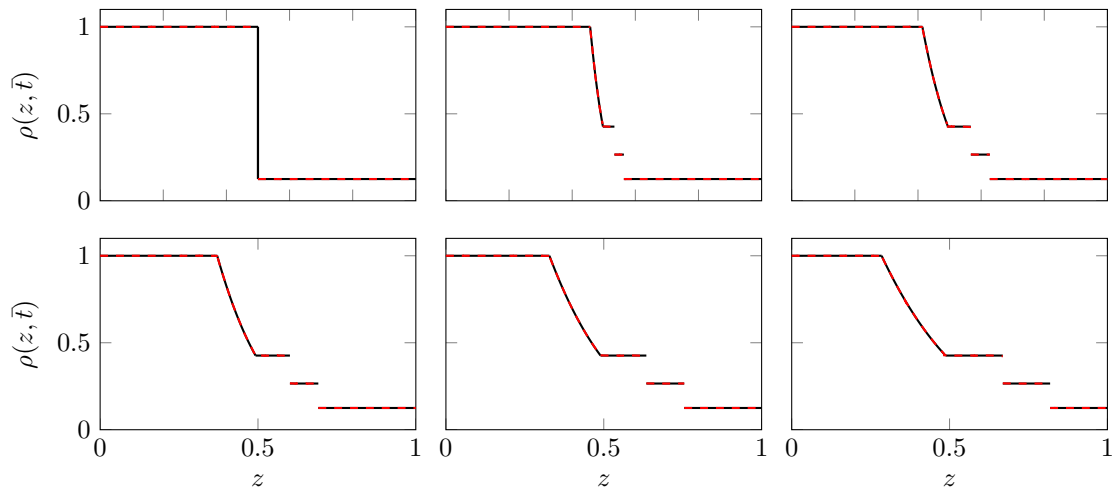}
}
\caption{Temporal slices of the analytical solution (\ref{line:sod:exact})
 and HOIST approximation (\ref{line:sod:hoist}) to the \texttt{sod} test
 case at time instances $\bar{t}=0, 0.0364, 0.0727, 0.109, 0.146, 0.182$
 (\textit{left-to-right}, \textit{top-to-bottom}).}
\label{fig:euler:sod:slice}
\end{figure}

\begin{figure}[!htbp]
\centering
\ifbool{fastcompile}{}{
\begin{tikzpicture}
\begin{groupplot}[
  group style={
      group size=2 by 3,
      horizontal sep=0.4cm,
      vertical sep=0.4cm
  },
  width=0.54\textwidth,
  axis equal image,
  xticklabels={,,},
  yticklabels={,,},
  xmin=0, xmax=1,
  ymin=0, ymax=0.2
]
\nextgroupplot
\addplot graphics [xmin=0, xmax=1, ymin=0, ymax=0.2] {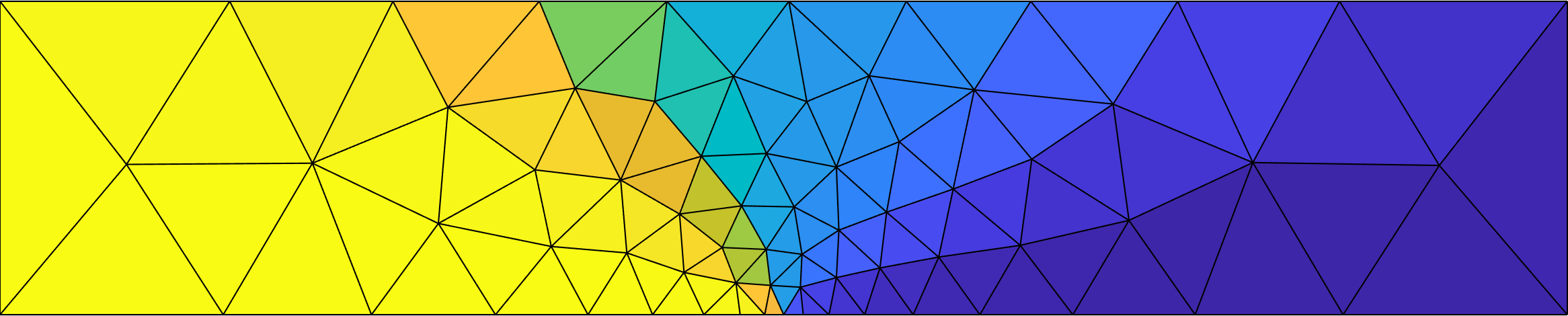};
\nextgroupplot
\addplot graphics [xmin=0, xmax=1, ymin=0, ymax=0.2] {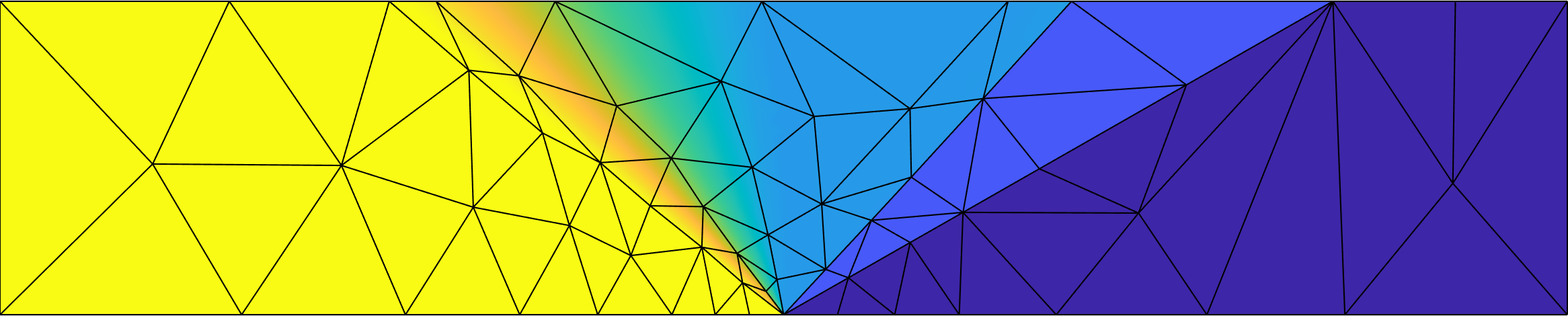};
\nextgroupplot
\addplot graphics [xmin=0, xmax=1, ymin=0, ymax=0.2] {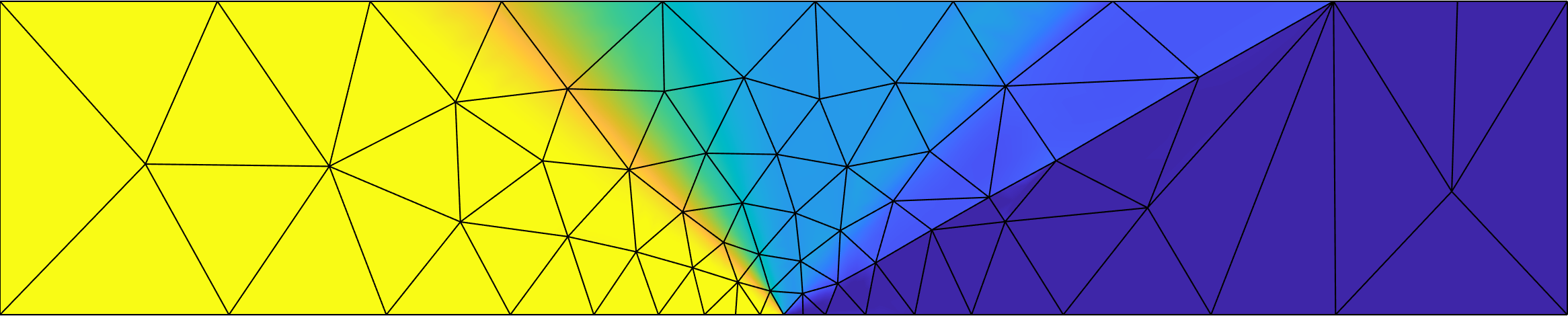};
\nextgroupplot
\addplot graphics [xmin=0, xmax=1, ymin=0, ymax=0.2] {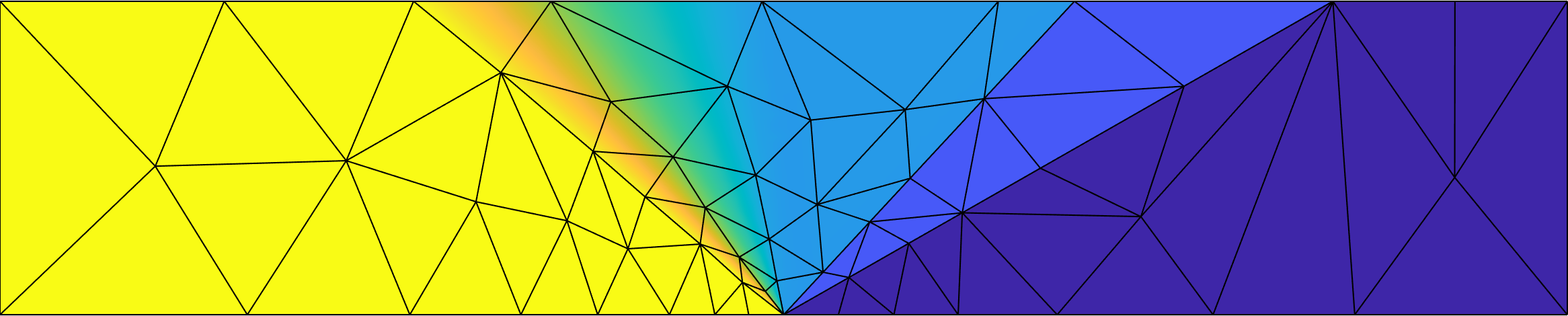};
\nextgroupplot
\addplot graphics [xmin=0, xmax=1, ymin=0, ymax=0.2] {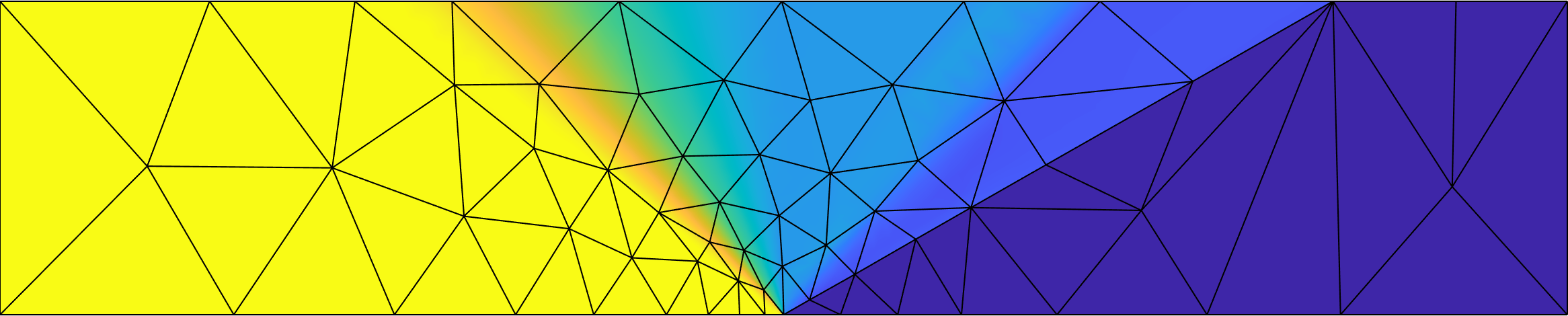};
\nextgroupplot
\addplot graphics [xmin=0, xmax=1, ymin=0, ymax=0.2] {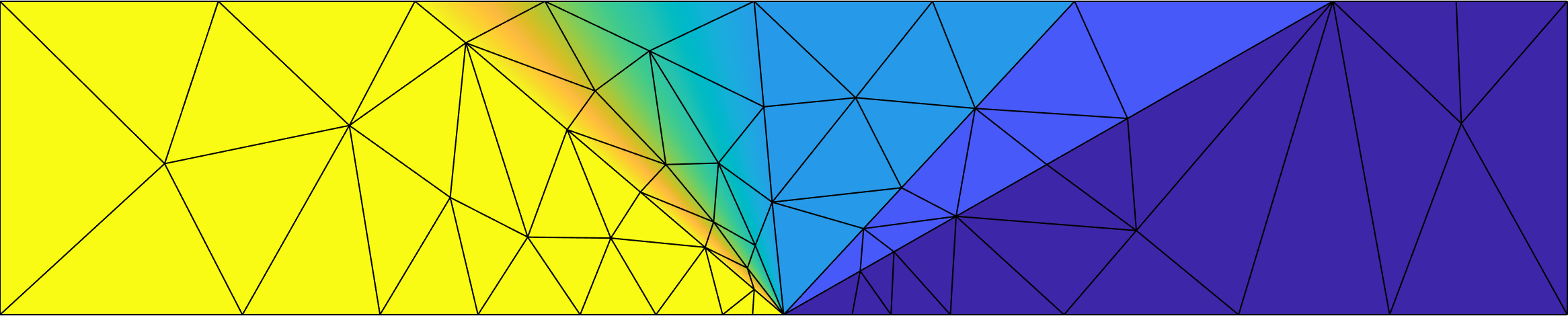};
\end{groupplot}
\end{tikzpicture}
}
\caption{Selected HOIST iterations for \texttt{sod} test case at iterations
  $k=0,10,27,52,66,106$ (\textit{top-to-bottom},
  \textit{left-to-right}). The SQP solver tracks the non-smooth features
  in order based on the strength of the feature (strongest to weakest):
  shock, contact discontinuity, head of rarefaction, and tail of rarefaction.
  Colorbar in Figure~\ref{fig:euler:sod:final}.}
\label{fig:euler:sod:sltn}
\end{figure}

All of the discontinuity surfaces in this problem are straight-sided
so we choose to approximate the geometry using straight-sided elements
($q=1$) and approximate the conservative variables using $p=2$ elements.
We discretize the domain using an unstructured mesh of 108 triangular
elements, generated using DistMesh \cite{persson_simple_2004}, with a
refinement region at the point $x=(0.5,0)$ from which the waves originate
because it separates the flow into multiple regions.
Furthermore, we place a node at $x=(0.5,0)$ and fix
it throughout the HOIST iterations to ensure the discontinuity in the
boundary condition is accurately integrated in the weak formulation of
the conservation law. Using this unstructured mesh and the corresponding
first-order finite volume solution to initialize the HOIST method
(parameters in Table~\ref{tab:params}), it converges
to a mesh that tracks the shock, contact discontinuity, and the head
and tail of the rarefaction (Figure~\ref{fig:euler:sod:final}), which
leads to an accurate approximation to the exact solution
(Figure~\ref{fig:euler:sod:slice}) on a space-time mesh with
only 84 elements (24 elements removed during the solution procedure).
This shows the error-based
HOIST indicator successfully tracks both discontinuous features and
continuous features with discontinuous derivatives. The SQP solver
tracks the features individually based on their relative strength,
i.e., discontinuities are tracked first based on the magnitude of
their jump and then non-smooth features are tracked based on the
magnitude of the derivative jump. The strongest feature, the shock,
is tracked by the mesh first (around iteration 27), next the contact
discontinuity is tracked (around iteration 52), then the head of the
rarefaction is tracked (around iteration 66), and the tail of the
rarefaction, the weakest feature, is tracked (around iteration 106)
(Figure~\ref{fig:euler:sod:sltn}).

\subsubsection{Supersonic flow over two-dimensional diamond in tunnel}
\label{sec:numexp:euler:diamond}
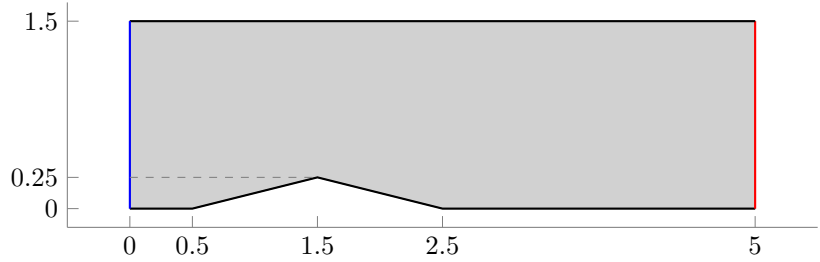
\begin{figure}[!htbp]
\centering
\ifbool{fastcompile}{}{
\begin{tikzpicture}
\begin{axis}[
axis equal image,
axis line style={gray},
axis x line*=bottom,
axis y line*=left,
width=0.7\textwidth,
xtick={0.0, 0.5, 1.5, 2.5, 5.0},
ytick={0.0, 0.25, 1.5},
ymax=1.65,
xmax=5.5,
xmin=-0.5,
ymin=-0.15]
\addplot [opacity=0.6, fill=black!30!white, opacity=0.6, forget plot]
coordinates {
( 0.00000000e+00,  0.00000000e+00)
( 5.00000000e-01,  0.00000000e+00)
( 1.50000000e+00,  2.50000000e-01)
( 2.50000000e+00,  0.00000000e+00)
( 5.00000000e+00,  0.00000000e+00)
( 5.00000000e+00,  1.50000000e+00)
( 0.00000000e+00,  1.50000000e+00)
( 0.00000000e+00,  0.00000000e+00)};

\addplot [thick, color=black]
coordinates {
( 0.00000000e+00,  0.00000000e+00)
( 5.00000000e-01,  0.00000000e+00)
( 1.50000000e+00,  2.50000000e-01)
( 2.50000000e+00,  0.00000000e+00)
( 5.00000000e+00,  0.00000000e+00)};\label{line:diamond:wall}

\addplot [thick, color=red]
coordinates {
( 5.00000000e+00,  0.00000000e+00)
( 5.00000000e+00,  1.50000000e+00)};\label{line:diamond:outlet}

\addplot [thick, color=black]
coordinates {
( 5.00000000e+00,  1.50000000e+00)
( 0.00000000e+00,  1.50000000e+00)};\label{line:diamond:wall}

\addplot [thick, color=blue]
coordinates {
( 0.00000000e+00,  1.50000000e+00)
( 0.00000000e+00,  0.00000000e+00)};\label{line:diamond:inlet}

\addplot [dashed, color=gray, forget plot]
coordinates {
( 0.00000000e+00,  2.50000000e-01)
( 1.50000000e+00,  2.50000000e-01)};

\end{axis}
\end{tikzpicture}
}
\caption{Geometry and boundary conditions for the \texttt{diamond} test case.
 Boundary conditions: slip walls (\ref{line:diamond:wall}), Mach 2 supersonic inflow
(\ref{line:diamond:inlet}), and supersonic outflow (\ref{line:diamond:outlet}).}
\label{fig:euler:diamond:geom}
\end{figure}

Next, we consider steady, supersonic (Mach 2) flow through a two-dimensional
diamond in a tunnel (Figure~\ref{fig:euler:diamond:geom};
test case: \texttt{diamond}) that features
reflecting, intersecting, and curved shocks 
to demonstrate the HOIST method can accurately track these features
with a high-order mesh and resolve the corresponding flow features.

\begin{figure}[!htbp]
\centering
\ifbool{fastcompile}{}{
\begin{tikzpicture}
\begin{groupplot}[
  group style={
      group size=2 by 2,
      horizontal sep=0.4cm,
      vertical sep=0.4cm
  },
  width=0.57\textwidth,
  axis equal image,
  xticklabels={,,},
  yticklabels={,,},
  xmin=0, xmax=5,
  ymin=0, ymax=1.5
]
\nextgroupplot
\addplot graphics [xmin=0, xmax=5, ymin=0, ymax=1.5] {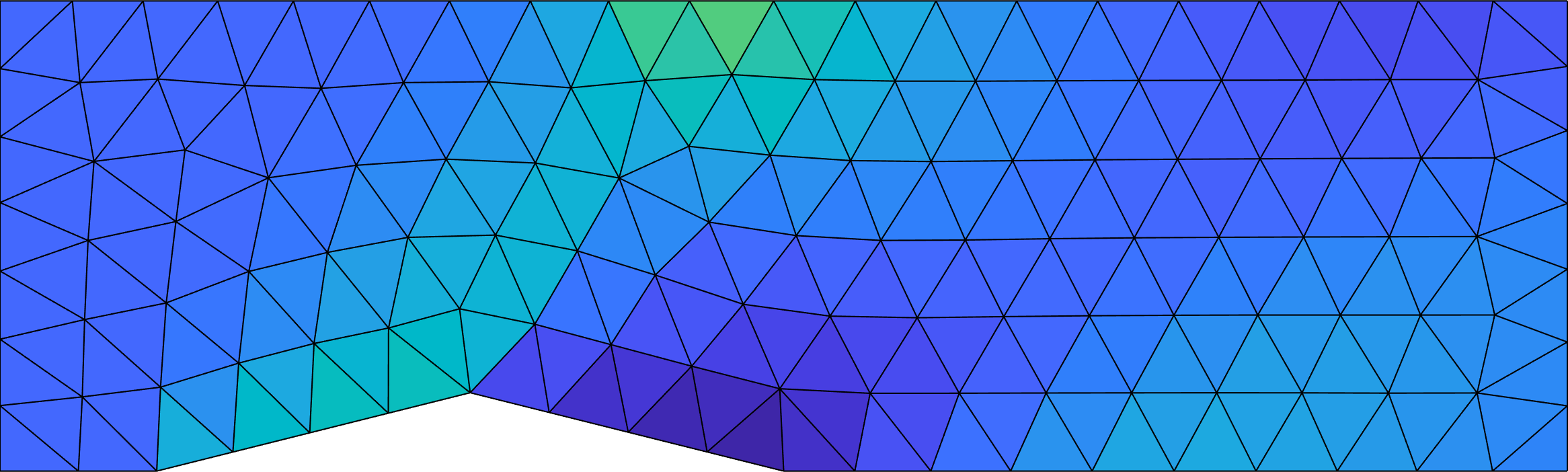};
\nextgroupplot
\addplot graphics [xmin=0, xmax=5, ymin=0, ymax=1.5] {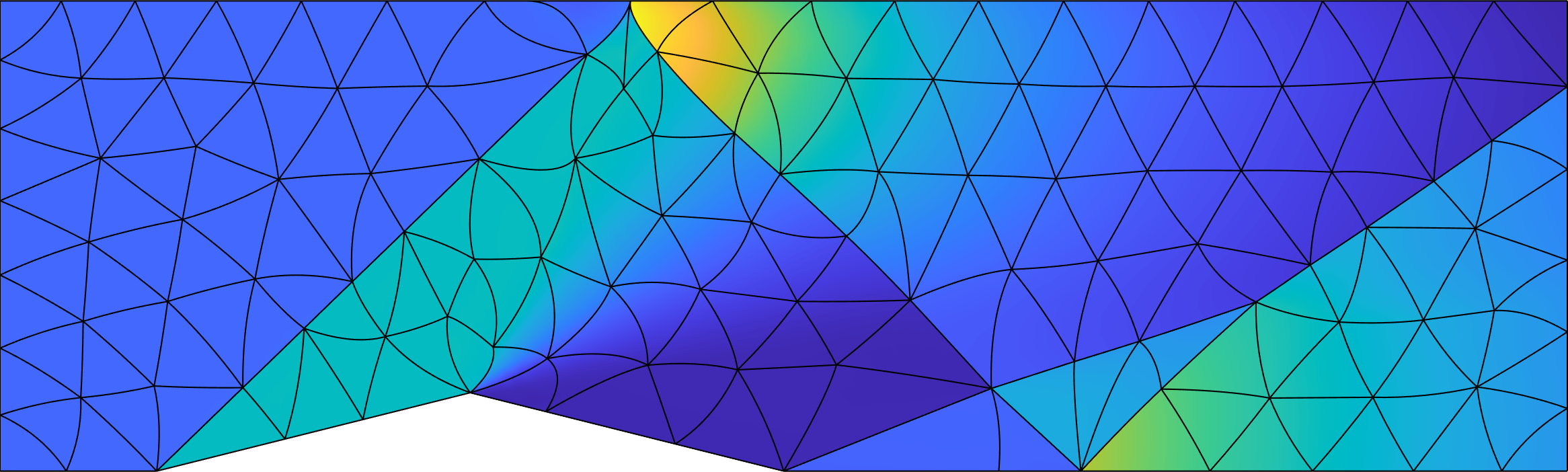};
\nextgroupplot
\addplot graphics [xmin=0, xmax=5, ymin=0, ymax=1.5] {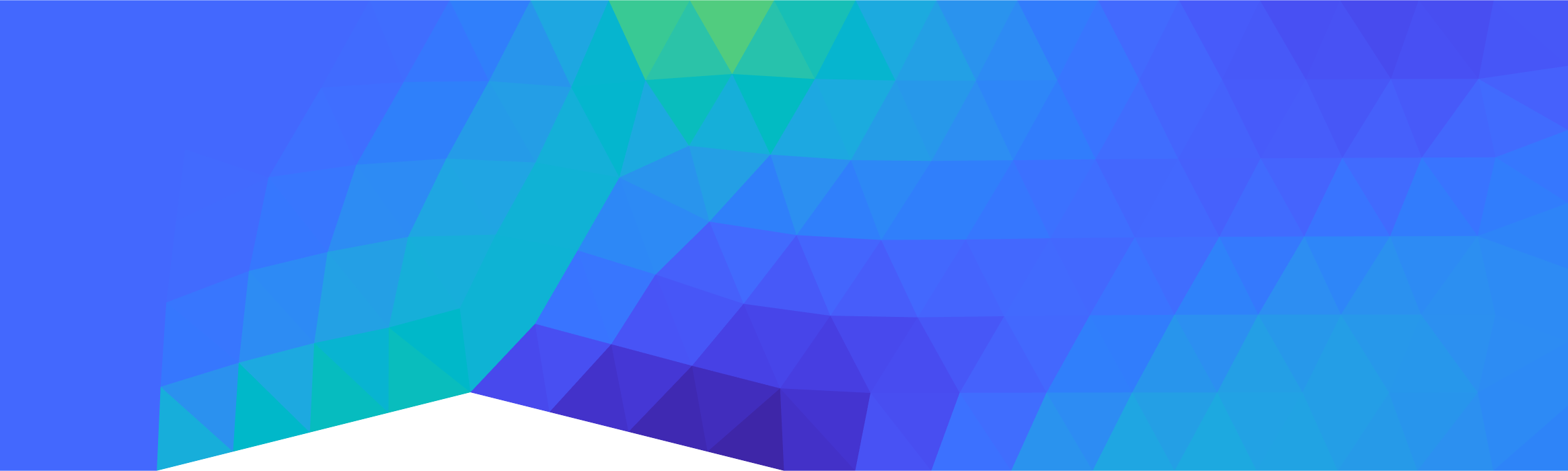};
\nextgroupplot
\addplot graphics [xmin=0, xmax=5, ymin=0, ymax=1.5] {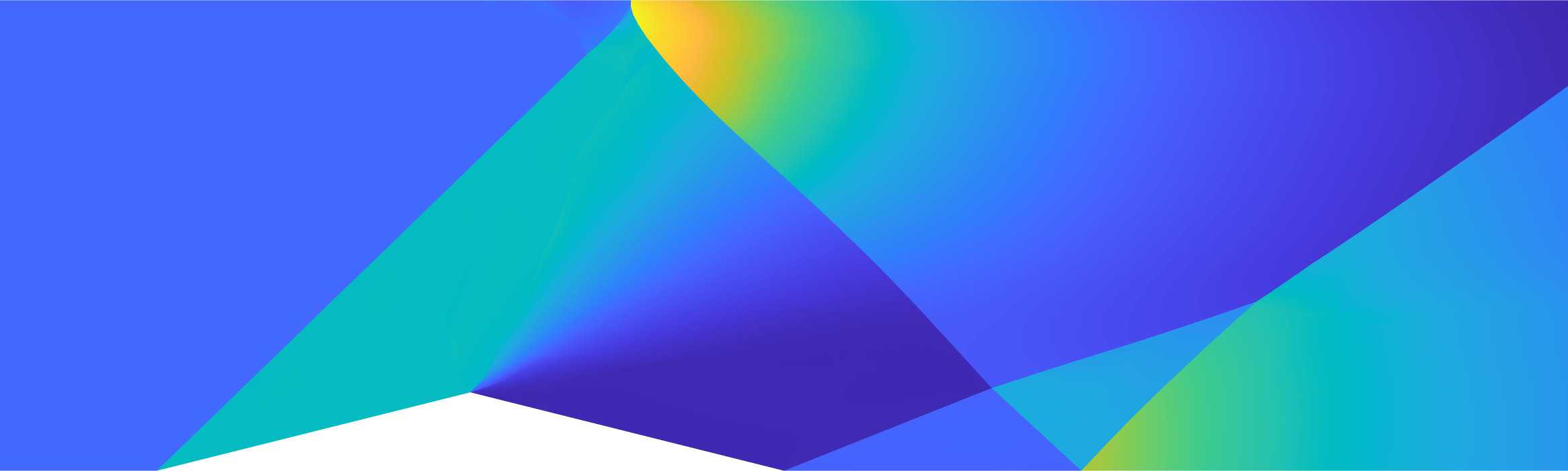};
\end{groupplot}
\end{tikzpicture}
}
\colorbarMatlabParula{0.7}{1}{2}{3}{4.1}
\caption{Starting point (\textit{left}) and HOIST solution (\textit{right})
 for the \texttt{diamond} test case with (\textit{top}) and without (\textit{bottom})
 mesh edges. The high-order mesh tracks the curved shocks, shock-shock interactions,
 and shock reflections and the solution is well-resolved throughout the domain.}
\label{fig:euler:diamond:sltn}
\end{figure}

We discretize the domain using an unstructured triangular mesh, generated
using DistMesh \cite{persson_simple_2004}, consisting of 220 elements and
use a third-order approximation for the geometry and flow variables ($p=q=2$).
Using this unstructured mesh and the corresponding first-order finite volume
solution to initialize the HOIST method (parameters in Table~\ref{tab:params}),
it converges to a mesh that tracks all shocks and their
intersections with the final mesh containing 201 elements (19 elements
removed during the solution procedure). 
The coarse, high-order (curved) elements
conform to the curvature of the discontinuity surfaces and the
four generalized triple points are tracked
(Figure~\ref{fig:euler:diamond:sltn}). A shock
capturing method would require refined elements in the vicinity
of all discontinuity surfaces, particularly the triple points,
which would lead to a mesh with many more elements.

\begin{remark}
There should be a slip line (a line along which there is a weak discontinuity
in the density, but the pressure and velocity direction are continuous)
emerging from the triple point near the trailing edge of the diamond
(at $x = (3.16, 0.26)$) that is
not resolved or tracked because the mesh in that region is too coarse.
As the mesh is refined, the HOIST method will track all slip lines;
we defer a refinement study to the scramjet example in
Section~\ref{sec:numexp:euler:scramjet}.
\end{remark}

\subsubsection{Hypersonic flow through two-dimensional scramjet}
\label{sec:numexp:euler:scramjet}
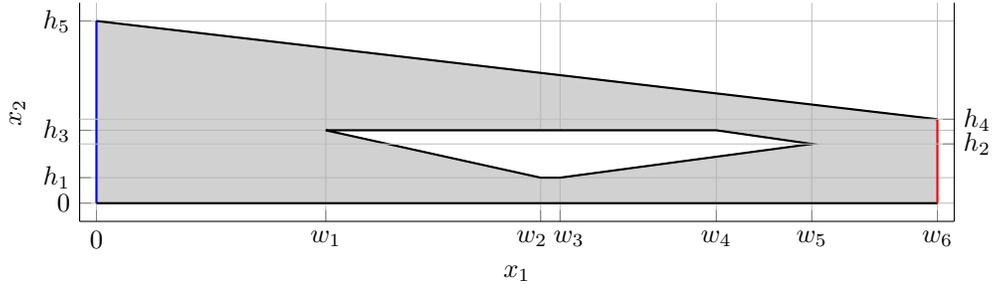
\begin{figure}[!htbp]
\centering
\ifbool{fastcompile}{}{
\begin{tikzpicture}
\begin{axis}[
axis equal image,
axis line style={black},
axis x line*=bottom,
axis y line*=left,
yticklabels={0, $h_1$, $h_3$, $h_5$},
ymax=0.9515,
width=0.8\textwidth,
xtick={0.0, 1.088, 2.1071, 2.1998, 2.941, 3.394, 3.9899},
ytick={0.0, 0.1211, 0.3460, 0.8650},
xlabel={$x_1$},
xticklabels={0, $w_1$, $w_2~~$, $~~w_3$, $w_4$, $w_5$, $w_6$},
xmax=4.069698,
ylabel={$x_2$},
xmin=-0.07979800000000001,
grid=both,
ymin=-0.08650000000000001]
\addplot [opacity=0.6, fill=black!30!white, opacity=0.6, forget plot]
coordinates {
( 0.00000000e+00,  0.00000000e+00)
( 3.98990000e+00,  0.00000000e+00)
( 3.98990000e+00,  3.98600000e-01)
( 0.00000000e+00,  8.65000000e-01)
( 0.00000000e+00,  0.00000000e+00)};

\addplot [opacity=1.0, fill=white, opacity=1.0, forget plot]
coordinates {
( 1.08800000e+00,  3.46000000e-01)
( 2.94100000e+00,  3.46000000e-01)
( 3.39400000e+00,  2.81100000e-01)
( 2.19980000e+00,  1.21100000e-01)
( 2.10710000e+00,  1.21100000e-01)
( 1.08800000e+00,  3.46000000e-01)};

\addplot [thick, color=black]
coordinates {
( 0.00000000e+00,  0.00000000e+00)
( 3.98990000e+00,  0.00000000e+00)};\label{line:scramjet1:wall}

\addplot [thick, color=black, forget plot]
coordinates {
( 3.98990000e+00,  3.98600000e-01)
( 0.00000000e+00,  8.65000000e-01)};

\addplot [thick, color=black, forget plot]
coordinates {
( 1.08800000e+00,  3.46000000e-01)
( 2.94100000e+00,  3.46000000e-01)
( 3.39400000e+00,  2.81100000e-01)
( 2.19980000e+00,  1.21100000e-01)
( 2.10710000e+00,  1.21100000e-01)
( 1.08800000e+00,  3.46000000e-01)};

\addplot [thick, color=blue]
coordinates {
( 0.00000000e+00,  8.65000000e-01)
( 0.00000000e+00,  0.00000000e+00)};\label{line:scramjet1:inflow}

\addplot [thick, color=red]
coordinates {
( 3.98990000e+00,  0.00000000e+00)
( 3.98990000e+00,  3.98600000e-01)};\label{line:scramjet1:outflow}

\end{axis}
\begin{axis}[
axis equal image,
axis line style={black},
axis x line=none,
axis y line*=right,
yticklabels={$h_2$, $h_4$},
width=0.8\textwidth,
ytick={0.2811, 0.3986},
grid=both,
ymax=0.9515,
xmax=4.069698,
xmin=-0.07979800000000001,
ymin=-0.08650000000000001]
\end{axis}
\end{tikzpicture}
}
\caption{Geometry and boundary conditions for the \texttt{scramjet} test case.
Geometry: $w_1=1.088$, $w_2=2.1071$, $w_3=2.1998$, $w_4=2.9410$, $w_5=3.3940$, 
$w_6=3.9899$, $h_1=0.1211$, $h_2=0.2811$, $h_3=0.3460$, $h_4=0.3986$, $h_5=0.8650$.
Boundary conditions:
slip walls (\ref{line:scramjet1:wall}),
Mach 5 supersonic inflow (\ref{line:scramjet1:inflow}), and
supersonic outflow (\ref{line:scramjet1:outflow}).}
\label{fig:euler:scramjet:geom}
\end{figure}

\begin{figure}[!htbp]
\centering
\ifbool{fastcompile}{}{
\begin{tikzpicture}
\begin{groupplot}[
  group style={
      group size=2 by 2,
      horizontal sep=0.4cm,
      vertical sep=0.4cm
  },
  width=0.57\textwidth,
  axis equal image,
  xticklabels={,,},
  yticklabels={,,},
  xmin=0, xmax=4,
  ymin=0, ymax=0.87
]
\nextgroupplot
\addplot graphics [xmin=0, xmax=4, ymin=0, ymax=0.87] {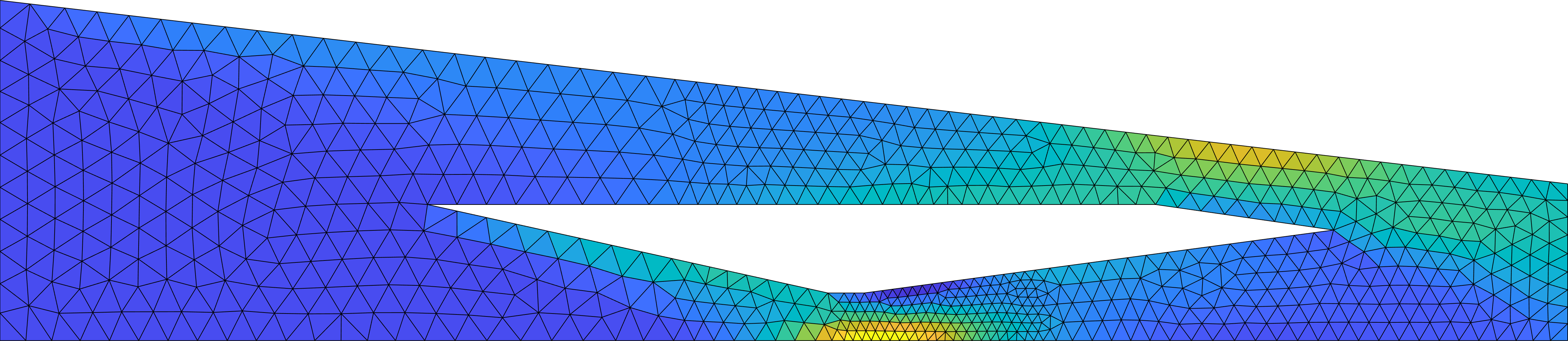};
\nextgroupplot
\addplot graphics [xmin=0, xmax=4, ymin=0, ymax=0.87] {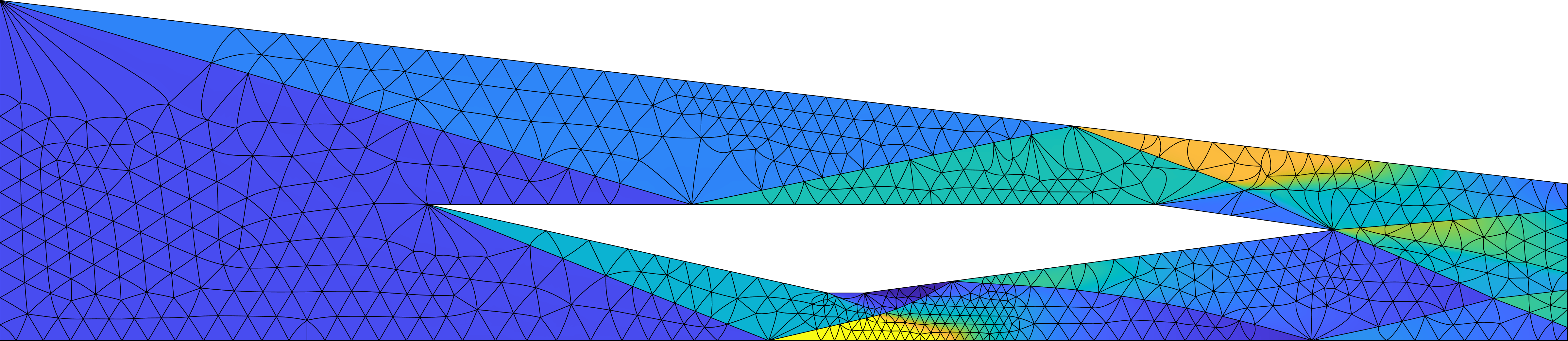};
\nextgroupplot
\addplot graphics [xmin=0, xmax=4, ymin=0, ymax=0.87] {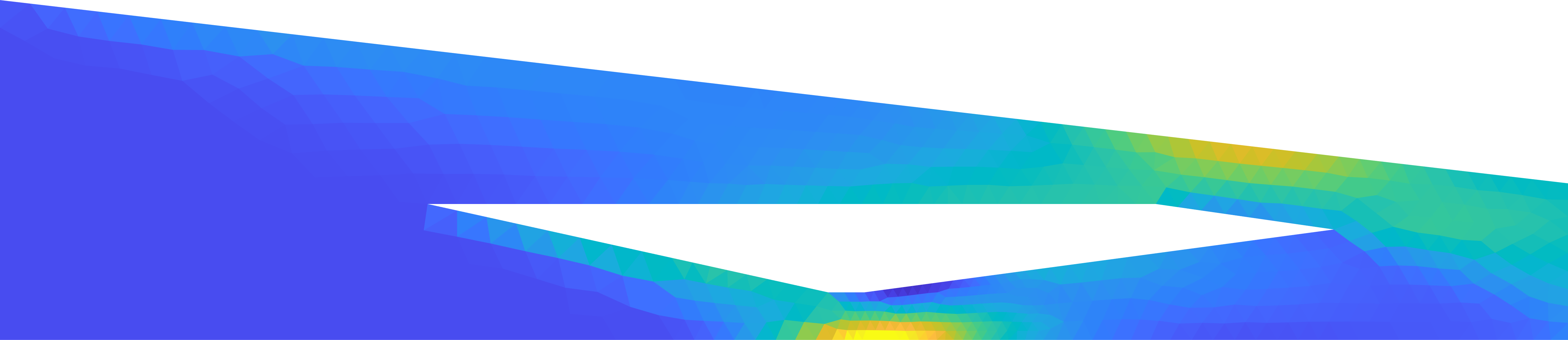};
\nextgroupplot
\addplot graphics [xmin=0, xmax=4, ymin=0, ymax=0.87] {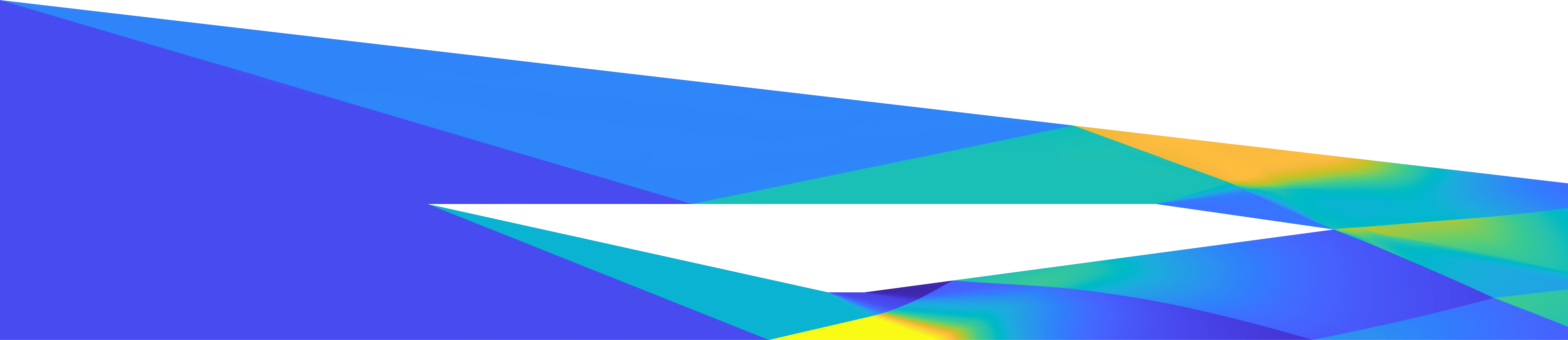};
\end{groupplot}
\end{tikzpicture}
}
\colorbarMatlabParula{0.6}{2}{4}{6}{7}
\caption{Starting point (\textit{left}) and HOIST solution (\textit{right})
 for the \texttt{scramjet} test case on the coarse mesh with (\textit{top}) and
 without (\textit{bottom}) mesh edges. The high-order mesh tracks the curved
 shocks, shock-shock interactions, and shock reflections (except the discontinuity
 emanating from $x=(w_5,h_2)$) and the solution is well-resolved throughout
 the domain.}
\label{fig:euler:scramjet:coarse}
\end{figure}

\begin{figure}[!htbp]
\centering
\ifbool{fastcompile}{}{
\begin{tikzpicture}
\begin{groupplot}[
  group style={
      group size=2 by 2,
      horizontal sep=0.4cm,
      vertical sep=0.4cm
  },
  width=0.57\textwidth,
  axis equal image,
  xticklabels={,,},
  yticklabels={,,},
  xmin=0, xmax=4,
  ymin=0, ymax=0.87
]
\nextgroupplot
\addplot graphics [xmin=0, xmax=4, ymin=0, ymax=0.87] {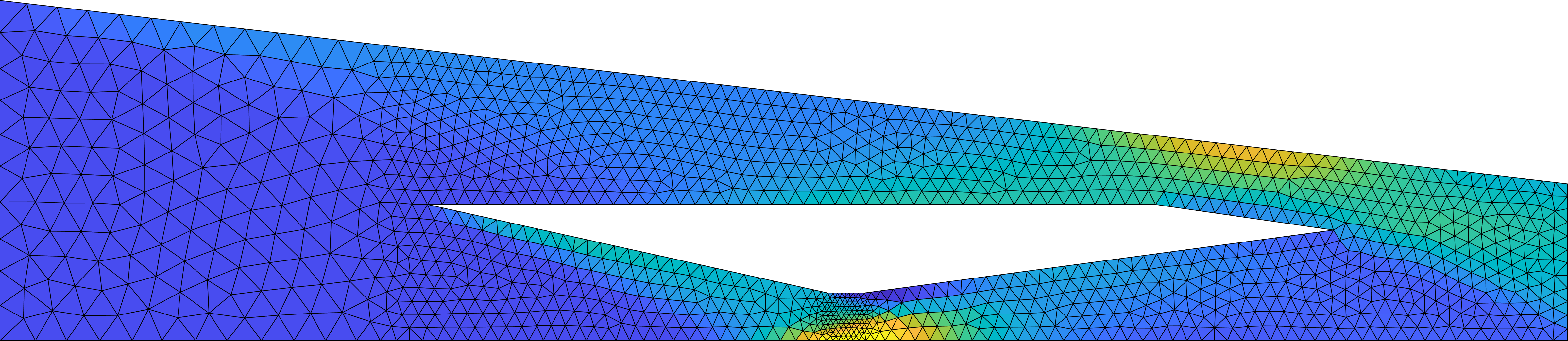};
\nextgroupplot
\addplot graphics [xmin=0, xmax=4, ymin=0, ymax=0.87] {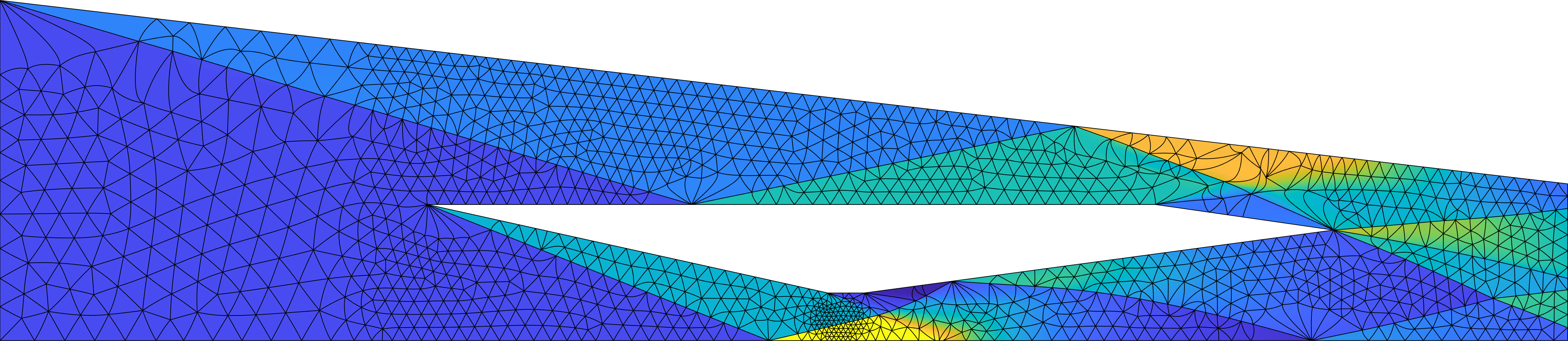};
\nextgroupplot
\addplot graphics [xmin=0, xmax=4, ymin=0, ymax=0.87] {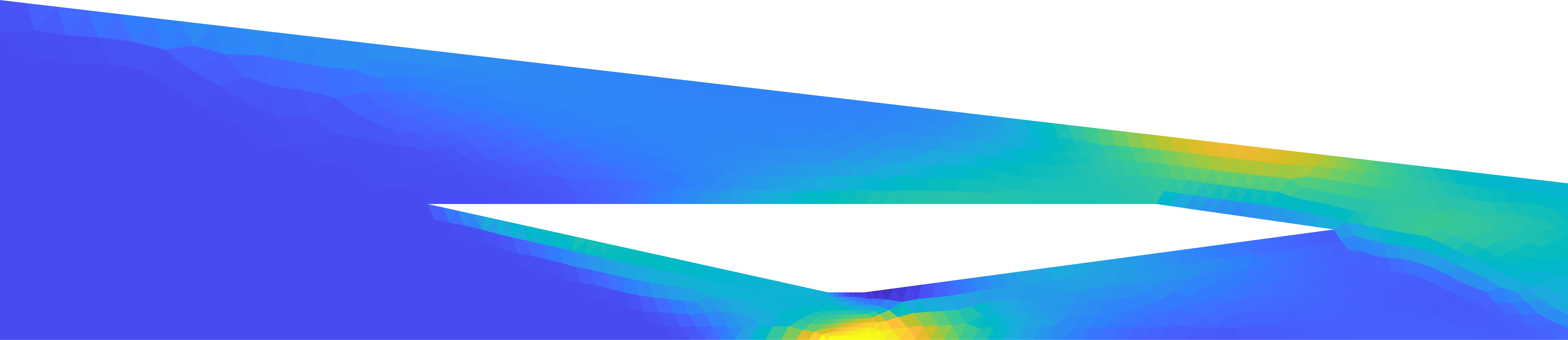};
\nextgroupplot
\addplot graphics [xmin=0, xmax=4, ymin=0, ymax=0.87] {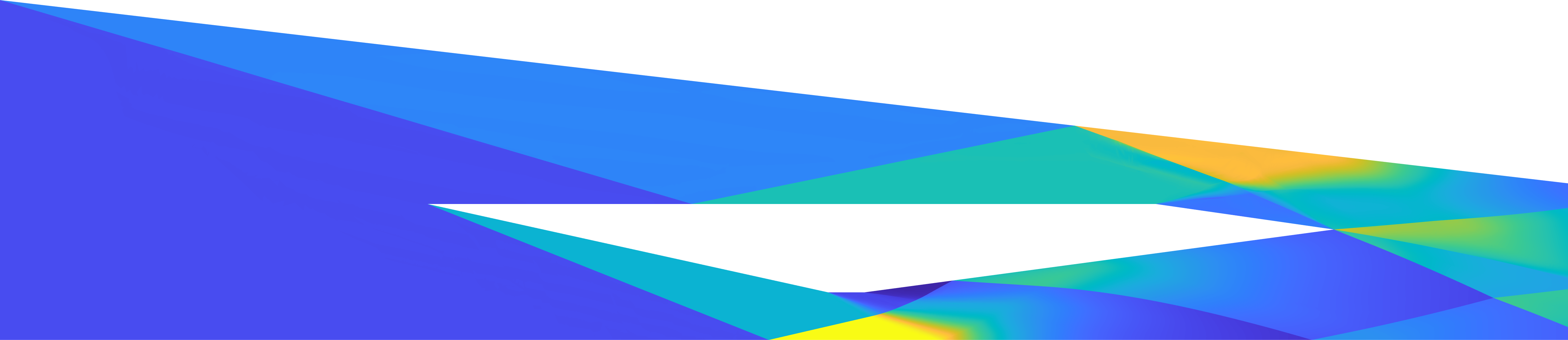};
\end{groupplot}
\end{tikzpicture}
}
\caption{Starting point (\textit{left}) and HOIST solution (\textit{right})
 for the \texttt{scramjet} test case on the fine mesh with (\textit{top}) and without
 (\textit{bottom}) mesh edges. The high-order mesh tracks the curved shocks,
 shock-shock interactions, and shock reflections (including the discontinuity
 emanating from $x=(w_5,h_2)$) and the solution is well-resolved throughout the
 domain. Colorbar in Figure~\ref{fig:euler:scramjet:coarse}.}
\label{fig:euler:scramjet:fine}
\end{figure}

Next, we consider steady, hypersonic (Mach 5) flow through a two-dimensional
scramjet (Figure~\ref{fig:euler:scramjet:geom}; test case: \texttt{scramjet});
it possesses similar features to the diamond case
(Section~\ref{sec:numexp:euler:scramjet}) with more complex discontinuity surfaces.
We discretize the domain using two unstructured triangular meshes, generated
using DistMesh \cite{persson_simple_2004}, consisting of 1442 and 2679 elements
and use a third-order approximation for the geometry and flow variables ($p=q=2$).
Using these unstructured meshes and the corresponding first-order finite volume
solution to initialize the HOIST method (parameters in Table~\ref{tab:params}),
it converges to a mesh that tracks all shocks and their
intersections with the final meshes containing 1299 and 2523 elements
(143 and 156 elements removed during the solution procedure), respectively.
The high-order (curved) elements conform to the curvature of the
discontinuity surfaces and the various triple points and reflection
points are tracked
(Figure~\ref{fig:euler:scramjet:coarse}-\ref{fig:euler:scramjet:fine}).
Both the coarse and fine mesh lead to highly accurate solution as seen
from the flow field without the mesh edges; the main difference between
the two simulations is that the fine mesh is able to completely track all
discontinuities emanating from the corner at $x=(w_5,h_2)$ 
(Figure~\ref{fig:euler:scramjet:fine}), including the slip line
 that is only partially tracked on the coarse mesh (Figure~\ref{fig:euler:scramjet:coarse}).


\subsubsection{Supersonic flow over sphere}
\label{sec:numexp:euler:sph}
\begin{figure}[!htbp]
\centering
\ifbool{fastcompile}{}{
\raisebox{-0.5\height}{\begin{tikzpicture}
\begin{axis}[
axis equal image,
axis line style={gray},
axis x line*=bottom,
axis y line*=right,
width=0.42\textwidth,
xtick={-4, -1, 0},
ytick={-5, -1, 0},
xlabel={$x_1$},
ymax=0.375,
xmax=0.3,
ylabel={$x_3$},
xmin=-4.3,
ymin=-5.375]
\addplot [opacity=0.6, fill=black!30!white, opacity=0.6, forget plot]
coordinates {
(-4.00000000e+00, -5.00000000e+00)
( 0.00000000e+00, -5.00000000e+00)
( 0.00000000e+00, -1.00000000e+00)
(-3.20515776e-02, -9.99486216e-01)
(-6.40702200e-02, -9.97945393e-01)
(-9.60230259e-02, -9.95379113e-01)
(-1.27877162e-01, -9.91790014e-01)
(-1.59599895e-01, -9.87181783e-01)
(-1.91158629e-01, -9.81559157e-01)
(-2.22520934e-01, -9.74927912e-01)
(-2.53654584e-01, -9.67294863e-01)
(-2.84527587e-01, -9.58667853e-01)
(-3.15108218e-01, -9.49055747e-01)
(-3.45365054e-01, -9.38468422e-01)
(-3.75267005e-01, -9.26916757e-01)
(-4.04783343e-01, -9.14412623e-01)
(-4.33883739e-01, -9.00968868e-01)
(-4.62538290e-01, -8.86599306e-01)
(-4.90717552e-01, -8.71318704e-01)
(-5.18392568e-01, -8.55142763e-01)
(-5.45534901e-01, -8.38088105e-01)
(-5.72116660e-01, -8.20172255e-01)
(-5.98110530e-01, -8.01413622e-01)
(-6.23489802e-01, -7.81831482e-01)
(-6.48228395e-01, -7.61445958e-01)
(-6.72300890e-01, -7.40277997e-01)
(-6.95682551e-01, -7.18349350e-01)
(-7.18349350e-01, -6.95682551e-01)
(-7.40277997e-01, -6.72300890e-01)
(-7.61445958e-01, -6.48228395e-01)
(-7.81831482e-01, -6.23489802e-01)
(-8.01413622e-01, -5.98110530e-01)
(-8.20172255e-01, -5.72116660e-01)
(-8.38088105e-01, -5.45534901e-01)
(-8.55142763e-01, -5.18392568e-01)
(-8.71318704e-01, -4.90717552e-01)
(-8.86599306e-01, -4.62538290e-01)
(-9.00968868e-01, -4.33883739e-01)
(-9.14412623e-01, -4.04783343e-01)
(-9.26916757e-01, -3.75267005e-01)
(-9.38468422e-01, -3.45365054e-01)
(-9.49055747e-01, -3.15108218e-01)
(-9.58667853e-01, -2.84527587e-01)
(-9.67294863e-01, -2.53654584e-01)
(-9.74927912e-01, -2.22520934e-01)
(-9.81559157e-01, -1.91158629e-01)
(-9.87181783e-01, -1.59599895e-01)
(-9.91790014e-01, -1.27877162e-01)
(-9.95379113e-01, -9.60230259e-02)
(-9.97945393e-01, -6.40702200e-02)
(-9.99486216e-01, -3.20515776e-02)
(-1.00000000e+00, -6.12323400e-17)
(-4.00000000e+00,  0.00000000e+00)
(-4.00000000e+00, -5.00000000e+00)};

\addplot [black, thick]
coordinates {
( 0.00000000e+00, -1.00000000e+00)
(-3.20515776e-02, -9.99486216e-01)
(-6.40702200e-02, -9.97945393e-01)
(-9.60230259e-02, -9.95379113e-01)
(-1.27877162e-01, -9.91790014e-01)
(-1.59599895e-01, -9.87181783e-01)
(-1.91158629e-01, -9.81559157e-01)
(-2.22520934e-01, -9.74927912e-01)
(-2.53654584e-01, -9.67294863e-01)
(-2.84527587e-01, -9.58667853e-01)
(-3.15108218e-01, -9.49055747e-01)
(-3.45365054e-01, -9.38468422e-01)
(-3.75267005e-01, -9.26916757e-01)
(-4.04783343e-01, -9.14412623e-01)
(-4.33883739e-01, -9.00968868e-01)
(-4.62538290e-01, -8.86599306e-01)
(-4.90717552e-01, -8.71318704e-01)
(-5.18392568e-01, -8.55142763e-01)
(-5.45534901e-01, -8.38088105e-01)
(-5.72116660e-01, -8.20172255e-01)
(-5.98110530e-01, -8.01413622e-01)
(-6.23489802e-01, -7.81831482e-01)
(-6.48228395e-01, -7.61445958e-01)
(-6.72300890e-01, -7.40277997e-01)
(-6.95682551e-01, -7.18349350e-01)
(-7.18349350e-01, -6.95682551e-01)
(-7.40277997e-01, -6.72300890e-01)
(-7.61445958e-01, -6.48228395e-01)
(-7.81831482e-01, -6.23489802e-01)
(-8.01413622e-01, -5.98110530e-01)
(-8.20172255e-01, -5.72116660e-01)
(-8.38088105e-01, -5.45534901e-01)
(-8.55142763e-01, -5.18392568e-01)
(-8.71318704e-01, -4.90717552e-01)
(-8.86599306e-01, -4.62538290e-01)
(-9.00968868e-01, -4.33883739e-01)
(-9.14412623e-01, -4.04783343e-01)
(-9.26916757e-01, -3.75267005e-01)
(-9.38468422e-01, -3.45365054e-01)
(-9.49055747e-01, -3.15108218e-01)
(-9.58667853e-01, -2.84527587e-01)
(-9.67294863e-01, -2.53654584e-01)
(-9.74927912e-01, -2.22520934e-01)
(-9.81559157e-01, -1.91158629e-01)
(-9.87181783e-01, -1.59599895e-01)
(-9.91790014e-01, -1.27877162e-01)
(-9.95379113e-01, -9.60230259e-02)
(-9.97945393e-01, -6.40702200e-02)
(-9.99486216e-01, -3.20515776e-02)
(-1.00000000e+00, -6.12323400e-17)
(-4.00000000e+00,  0.00000000e+00)};\label{line:sph0:wall}

\addplot [blue, thick]
coordinates {
(-4.00000000e+00, -5.00000000e+00)
( 0.00000000e+00, -5.00000000e+00)};\label{line:sph0:inlet}

\addplot [blue, thick, forget plot]
coordinates {
(-4.00000000e+00,  0.00000000e+00)
(-4.00000000e+00, -5.00000000e+00)};

\addplot [red, thick]
coordinates {
( 0.00000000e+00, -5.00000000e+00)
( 0.00000000e+00, -1.00000000e+00)};\label{line:sph0:outlet}

\end{axis}
\end{tikzpicture}} \quad
\raisebox{-0.5\height}{\input{_py/sph1_geom.tikz}} \qquad\quad
\raisebox{-0.425\height}{\includegraphics[width=0.19\textwidth]{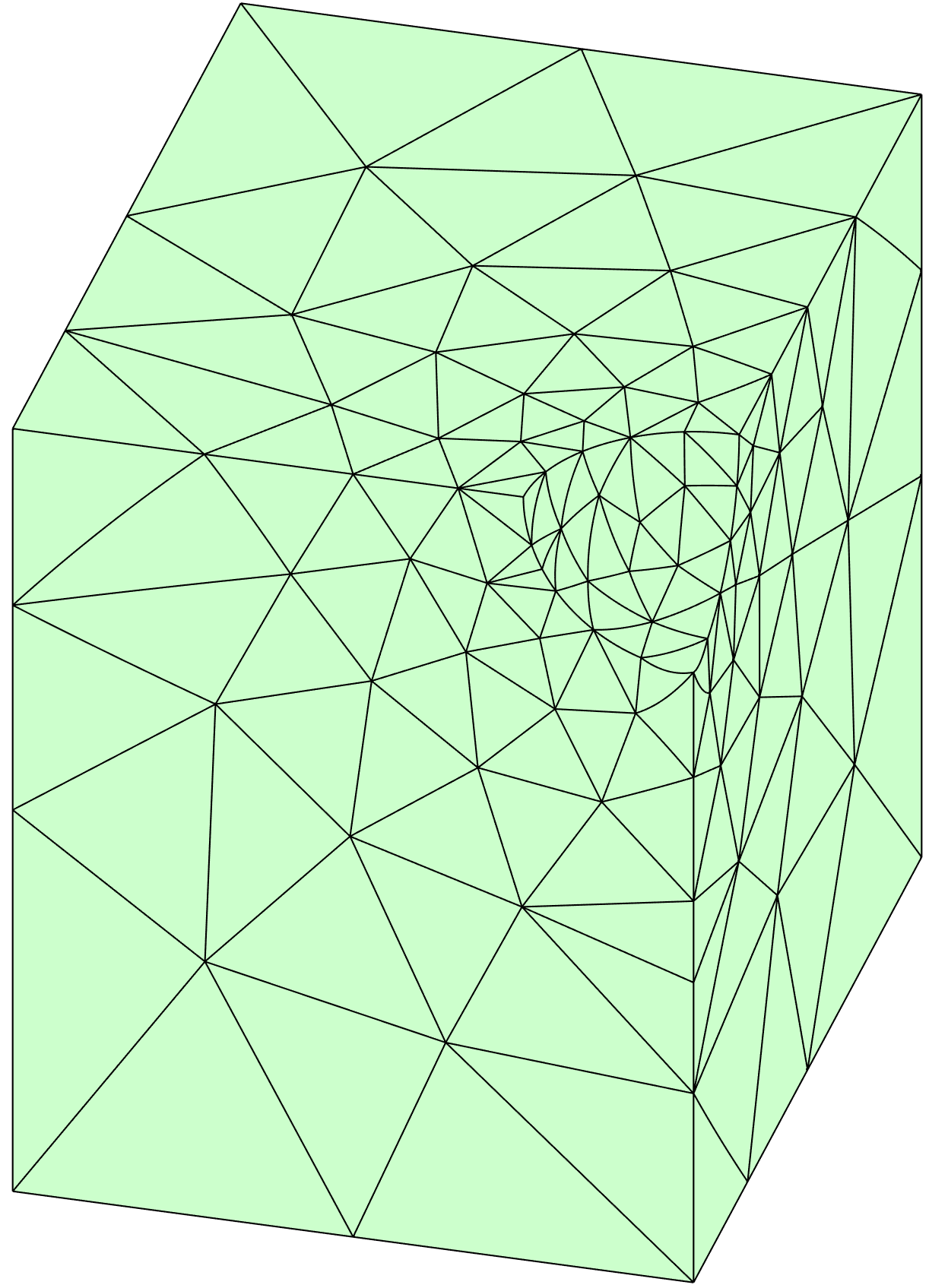}}
}
\caption{Geometry of the \texttt{sphere} test case including a slice along the
 $x_2=0$ plane (\textit{left}) and a three-dimensional view (\textit{middle})
 and the unstructured, high-order mesh with 511 quadratic tetrahedral elements
 (\textit{right}) used to initialize the HOIST simulation.
 Boundary conditions are indicated in the
 left figure: slip walls (\ref{line:sph0:wall}), Mach 2 supersonic
 inflow (\ref{line:sph0:inlet}), and supersonic outflow
 (\ref{line:sph0:outlet}); the $x_2=0$ plane uses slip wall, and the 
 $x_2=5$ plane uses supersonic inflow condition.}
\label{fig:euler:sphere0:geom}
\end{figure}

\begin{figure}[!htbp]
\centering
\ifbool{fastcompile}{}{
\includegraphics[width=0.34\textwidth]{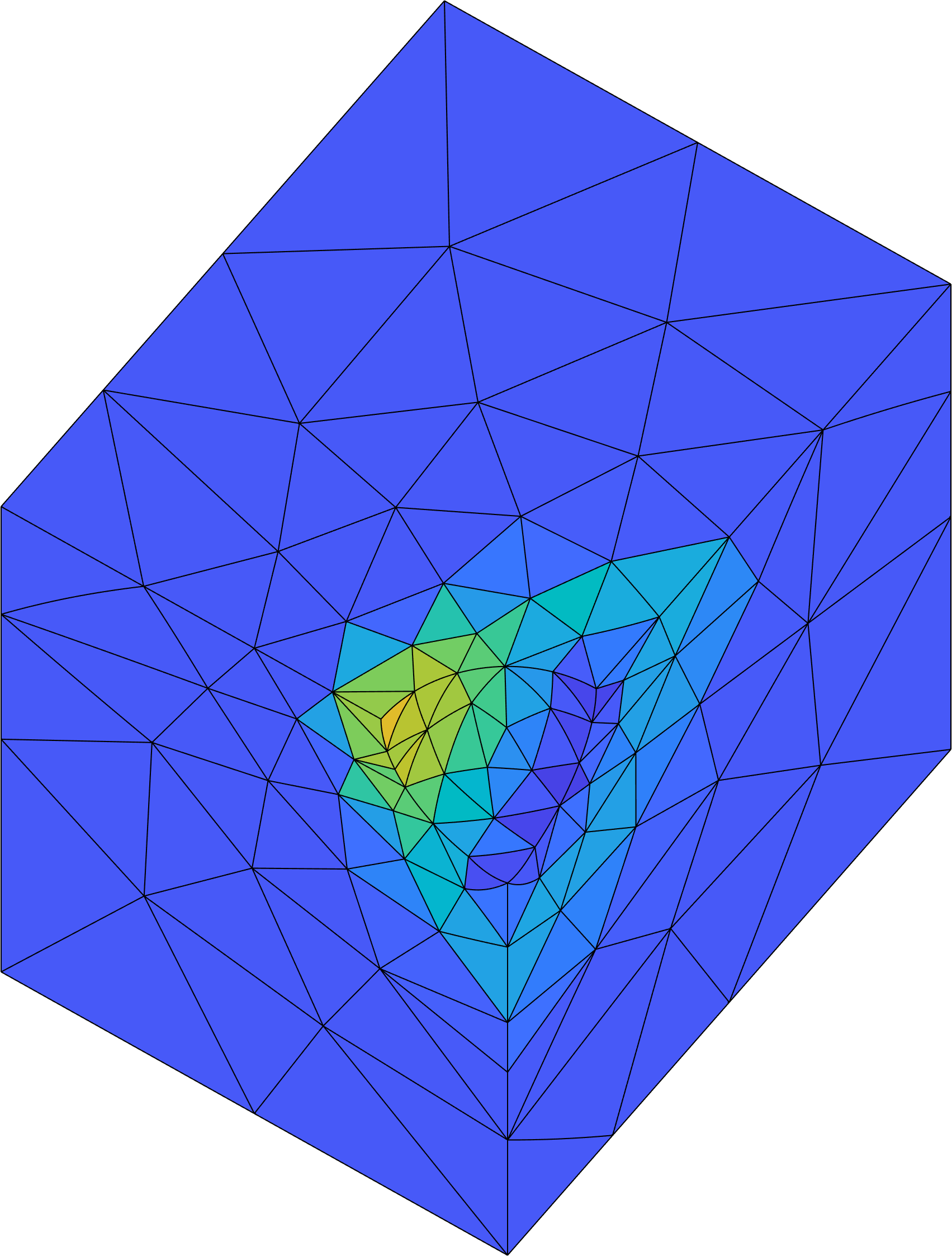} \qquad\qquad
\includegraphics[width=0.34\textwidth]{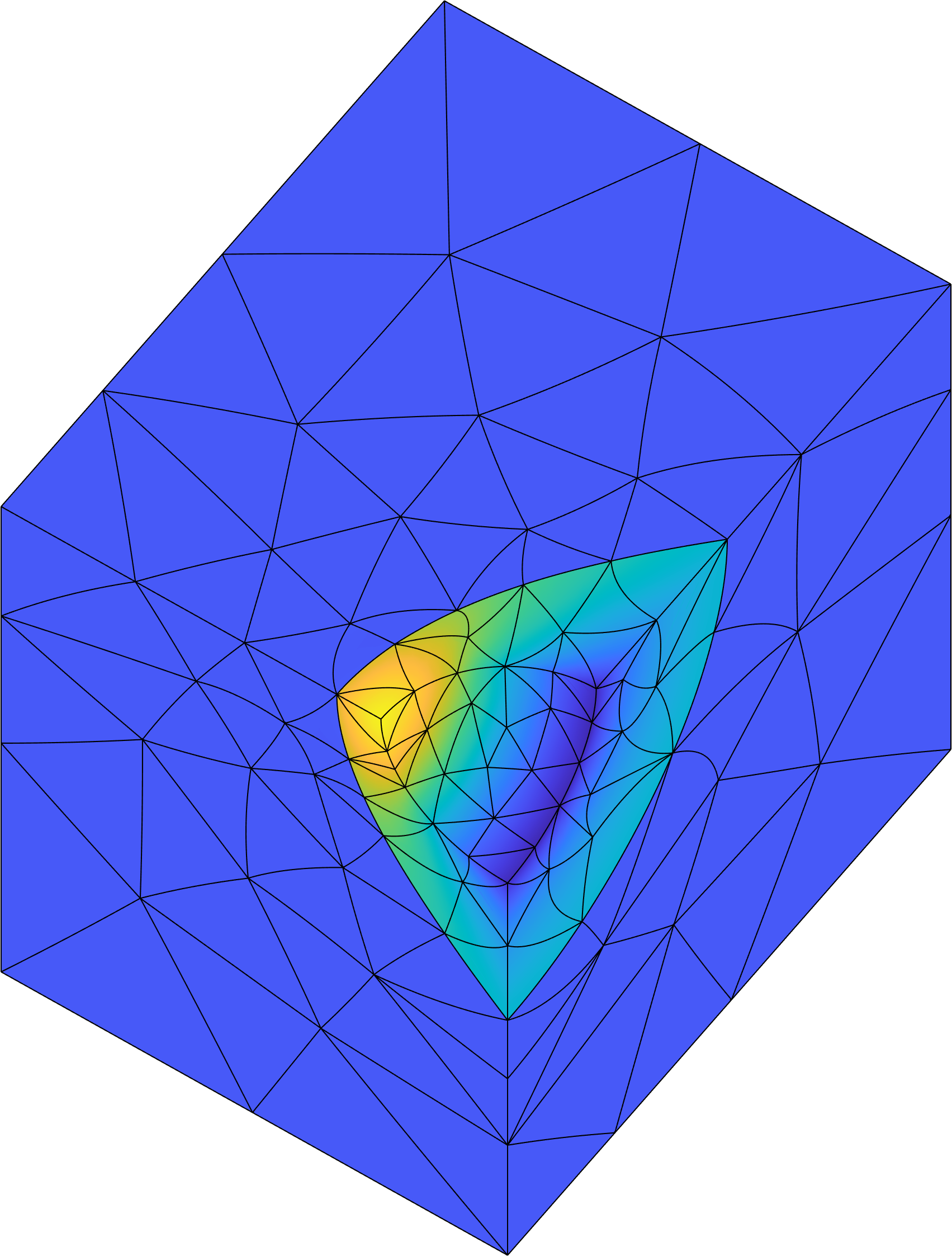}
}
\colorbarMatlabParula{0.8}{2}{3}{4}{4.5}
 \caption{Starting point (\textit{left}) and HOIST solution (\textit{right})
 for the \texttt{sphere} test case (colored by density). The high-order,
 high-quality tetrahedral elements track the curved bow shock and the solution
 is well-resolved throughout the domain.}
\label{fig:euler:sphere0:sltn}
\end{figure}

\begin{figure}[!htbp]
\centering
\ifbool{fastcompile}{}{
\includegraphics[width=0.34\textwidth]{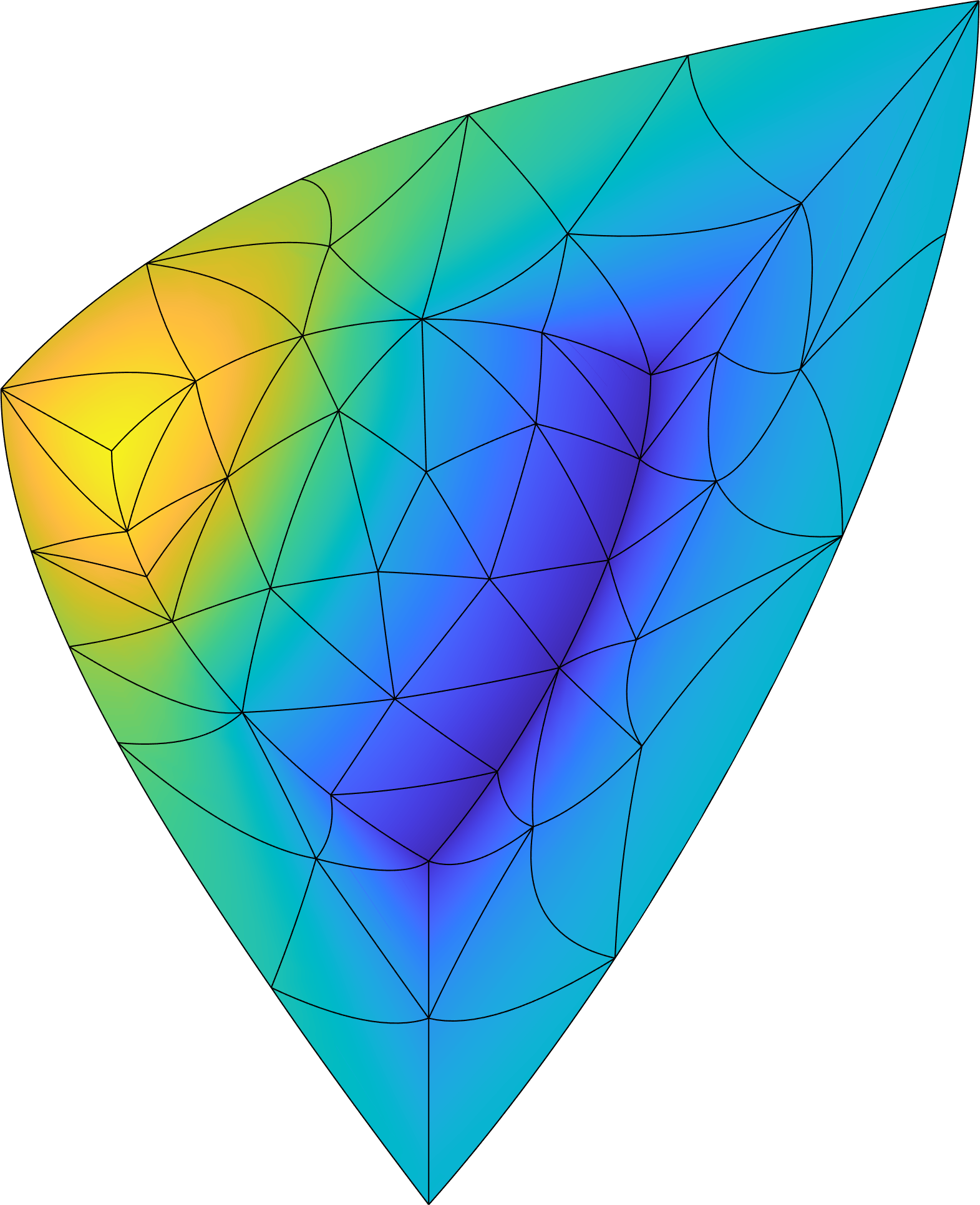} \qquad\qquad
\includegraphics[width=0.34\textwidth]{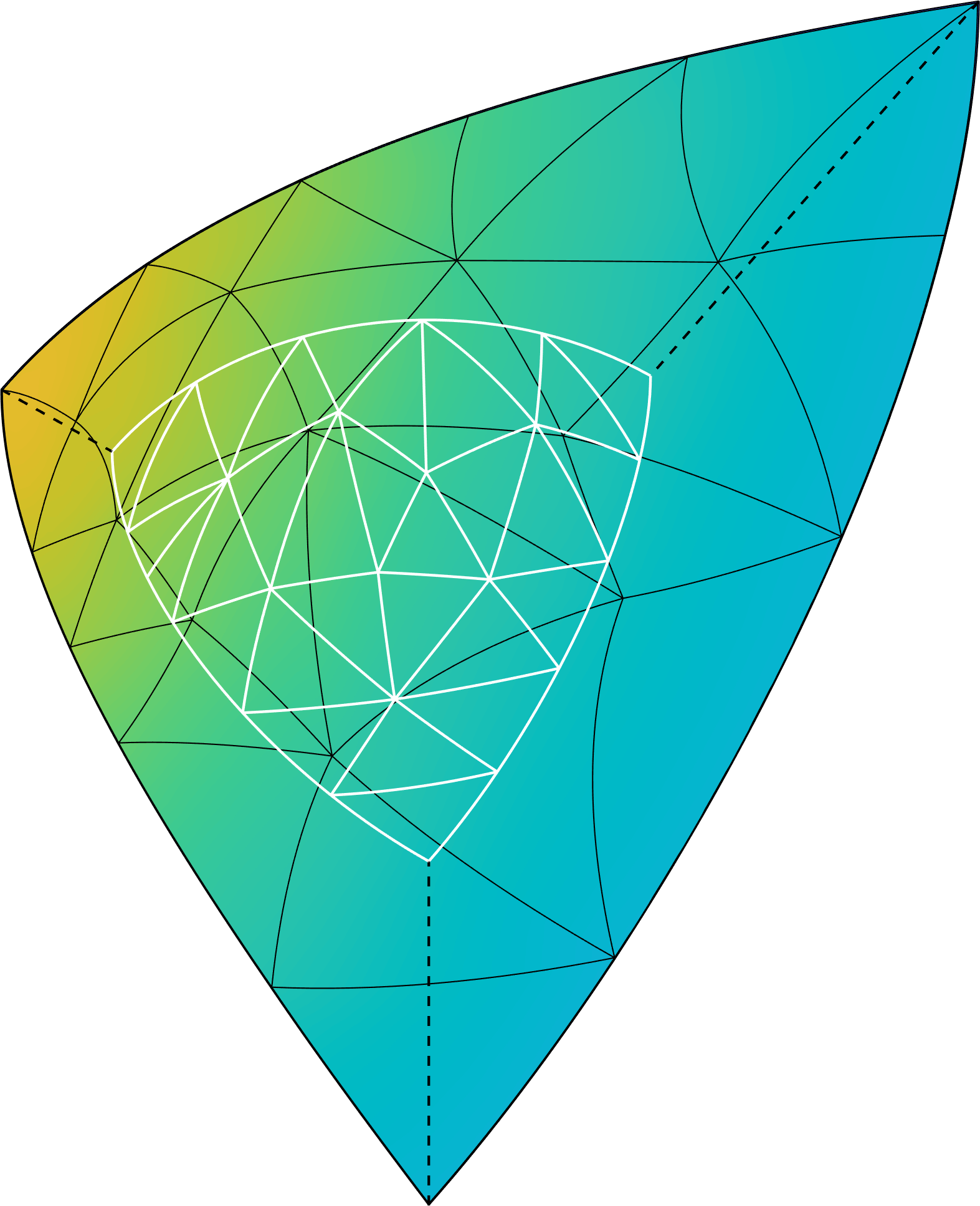}
}
\caption{The HOIST approximation to the shock layer (\textit{left}) and
 the surface of the bow shock (\textit{right}; surface mesh of the sphere
 included in white for reference) for the \texttt{sphere} test case (colored
 by density); the latter is produced
 by extracting the high-order surface triangles of the shock-aligned mesh.
 Both the shock layer and shock surface are represented to high accuracy
 on the extremely coarse mesh due to optimal, automated alignment of the
 high-order mesh with the bow shock. Colorbar in
 Figure~\ref{fig:euler:sphere0:sltn}.}
\label{fig:euler:sphere0:shk}
\end{figure}

Finally, we consider steady, supersonic (Mach 2) flow over a sphere 
(test case: \texttt{sphere}) that features a bow shock to demonstrate
the HOIST method reliably tracks curved shocks in three-dimensional
compressible flows with a coarse, high-order mesh.
To reduce the computational cost, we model 
only a portion of the geometry and use symmetry boundary conditions
(Figure~\ref{fig:euler:sphere0:geom}). We discretize the domain using
a high-order unstructured, tetrahedra mesh consisting of 511 elements
and use a third-order approximation for the geometry and flow variables
($p=q=2$). A straight-sided mesh is generated using DistMesh
\cite{persson_simple_2004}
then high-order nodes of element faces that lie on the sphere are
projected onto the curved boundary. Since this problem contains
a bow shock that does not intersect the sphere, we choose to fix
the nodes on the sphere surface (rather than allow them to slide
along the sphere). Using this unstructured mesh and the corresponding
first-order finite volume solution to initialize the HOIST method
(parameters in Table~\ref{tab:params}), it converges 
to a mesh that tracks the bow shock and a high-quality flow solution
(Figure~\ref{fig:euler:sphere0:sltn}); the final mesh contains
491 elements (20 elements removed during the solution procedure).
The high-order elements conform to the curvature of the bow shock
and the high-order approximation of the flow field leads to
a highly accurate approximation on an extremely coarse mesh.
The analytical stagnation pressure (at $x=(-1,0,0)$)
for this flow configuration is 5.6404 (non-dimensional) and the
corresponding quantity computed from the HOIST simulation is 5.6465,
a relative error of $0.11\%$.
Because the mesh conforms to the bow shock, an accurate
representation of the shock layer (region between bow shock
and sphere surface) and shock surface are directly available
(Figure~\ref{fig:euler:sphere0:shk}).

\begin{remark}
The HOIST parameters used for this simulation are mostly consistent
with those used in previous sections. The main difference is the
reference elements are used as the ideal elements
($K_{\star,e}=\Omega_{0,e}$), which is related to our choice
to fix nodes along the sphere. Tests with our standard approach
that takes $K_{\star,e}=K_\star$ (equilateral tetrahedra) indicated
numerous element collapses occur near the surface of the sphere as
elements in that region are driven toward equilateral tetrahedra,
which can cause the HOIST simulation to fail. 
\end{remark}

%
%

\section{Conclusion}
\label{sec:conclude}
The HOIST method is a high-order numerical method that approximates
solutions of conservation laws using a high-order DG discretization
on a mesh whose elements are aligned with non-smooth features in the
flow. This allows the non-smooth features to be represented
by the inter-element jumps in the DG solution and high-order basis
functions approximate smooth regions of the flow, which eliminates
the need for nonlinear stabilization. As a result, the method
produces highly accurate solutions on coarse meshes and recovers optimal
convergence rates of the DG discretization even for flows with non-smooth
features. While these advantages are shared by explicit shock tracking
approaches \cite{moretti2002thirty, salas2009shock},
HOIST \cite{zahr_optimization-based_2018,zahr_implicit_2020, zahr_high-order_2020}
and other implicit shock tracking methods (e.g., MDG-ICE
\cite{corrigan_moving_2019,kercher_least-squares_2020,kercher_moving_2021})
provide several additional benefits: formulated for
a general, nonlinear conservation law (not tailored to a specific set of
equations) and the formulation is well-suited for problems with intricate
discontinuity surfaces (curved and reflecting shocks, shock formation,
shock-shock interaction). The latter is due to the geometrically complex
problem of generating a feature-aligned mesh being re-cast to solving an
optimization problem over the discrete DG solution and nodal coordinates
of the mesh.

In this work, we introduced an improved solver for the implicit
shock tracking optimization problem based on the SQP solver
originally developed in \cite{zahr_implicit_2020}.
The new solver features a new merit function penalty parameter
(Section~\ref{sec:ist_solve:sqp:lsrch}), an adaptive mesh penalty parameter
in the objective function (Section~\ref{sec:ist_solve:sqp:mshqual}),
and, most importantly, a number of practical robustness measures
(Section~\ref{sec:ist_solve:mod}). The most critical
robustness measures are 1) dimension- and order-independent simplex
element collapses that are boundary-preserving and shock-aware
and 2) element-wise solution re-initialization
that resets the flow solution in oscillatory elements to a constant value.
In addition to the solver developments, we also introduced a general,
automated procedure to parametrize the mesh motion (i.e., nodal coordinates)
that guarantees all planar boundaries will be preserved throughout the implicit
tracking iterations. These advances to the HOIST method were critical to
enable robust convergence for two- and three-dimensional flows with
complex shock surfaces and they eliminated the need for continuation in
the polynomial degree. A series of nine numerical experiments demonstrated
the robustness of the solver even for two- and three-dimensional flows
with complex shock surfaces, the meshes produced are high-quality and track
non-smooth features in the flow (shocks, contacts, rarefactions), and the method
achieves the optimal convergence rate of the DG method even for flows with
non-smooth features. In particular, the diamond and scramjet problems
demonstrated the method can track curved, reflecting, and interacting
shocks, which leads to a highly accurate solution on coarse
meshes, even in the hypersonic flow regime. The sphere problem demonstrated
the robustness and generality of the method for three-dimensional,
compressible flows.

Additional research is required to develop iterative solvers and
preconditioners for the SQP linear system in (\ref{eqn:sqp_sys0}) to
make the approach practical for large-scale problems. Furthermore,
development of a general, automated approach to construct boundary-preserving
parametrizations for curved boundaries is needed for the method
to be useful for complex, three-dimensional geometries.
Other interesting avenues of future research include extension
of the method to viscous conservation laws, investigation into slab-based
space-time approaches to handle complex unsteady flows, and using the
method to study relevant hypersonic flows.

\appendix
\section{Parameters used for numerical experiments}
The HOIST parameters used for the numerical experiments are included in
Table~\ref{tab:params}. The only problems we choose to straighten
re-initialized elements (Remark~\ref{rem:reinit_smth}) are Burgers'
equation with shock formation (\texttt{iburg-form}), the diamond
(\texttt{diamond}), the scramjet (\texttt{scramjet}), and the sphere
(\texttt{sphere}).
\begin{sidewaystable}
\centering
\caption{HOIST parameters used for numerical experiments}
\label{tab:params}
\begin{tabular}{r|c|c|c|c|c}
& $K_{\star,e}$ & $(\gamma_0,\gamma_\mathrm{min},\tau,\sigma_1,\sigma_2)$ & $(\kappa_0,\kappa_\mathrm{min},\upsilon,\xi)$ & $(c_1,c_2,c_3,c_4,c_4')$ & $(c_5,c_6,c_7,c_8)$ \\\hline
\texttt{advec-planar} & $K_\star$ & $(10^{-2},10^{-6},2,10^{-2},10^{-1})$ & $(10^{-10},10^{-10},\text{-},\text{-})$ & $(0.2,10^{-10},0.2,0,\text{-})$ & - \\
\texttt{advec-trig} & $K_\star$ & $(10^{-2},10^{-6},2,10^{-2},10^{-1})$ & $(10^{-10},10^{-10},\text{-},\text{-})$ & $(0.2,10^{-10},0.2,0,0.05)$ & - \\
\texttt{iburg-acc} & $K_\star$ & $(10^{-2},10^{-2},2,10^{-2},10^{-1})$ & $(10^{-2},10^{-10},0.75,1)$ & $(0.2,10^{-10},0.2,0,0.05)$ & $(10^{-1},10^{-2},1,10^{-2})$ \\
\texttt{iburg-form} & $K_\star$ & $(10^{-4},10^{-4},2,10^{-2},10^{-1})$ & $(10^{-2},10^{-10},0.75,1)$ & $(0.2,10^{-10},0.33,0,0.01)$ & $(10^{-1},10^{-2},0.2,10^{-6})$ \\
\texttt{nozzle} & - & $(10,10^{-2},2,10^{-2},10^{-1})$ & - & $(0.2,10^{-10},0.2,0,0.05)$ & $(10^{-2},10^{-2},0.25,10^{-6})$ \\
\texttt{sod} & $K_\star$ & $(10^{-5},10^{-8},1.2,10^{-2},10^{-1})$ & $(10^{-6},10^{-10},0.75,2)$ & $(0.15,10^{-10},0.2,0.05,0.05)$ & $(10^{-2},10^{-2},0,10^{-6})$ \\
\texttt{diamond} & $K_\star$ & $(1,10^{-2},1.2,10^{-2},10^{-1})$ & $(10^{-3},10^{-8},0.5,0.5)$ & $(0.25,10^{-10},0.25,10^{-3},0)$ & $(10^{-2},10^{-2},0.5,10^{-4})$ \\
\texttt{scramjet} & $K_\star$ & $(1,10^{-2},1.1,10^{-2},10^{-1})$ & $(1,10^{-2},0.8,1)$ & $(0.2,10^{-10},0.05,10^{-3},10^{-4})$ & $(10^{-3},10^{-1},0.2,10^{-2})$ \\
\texttt{sphere} & $\Omega_{0,e}$ & $(10^{-3},10^{-3},1.5,10^{-2},10^{-1})$ & $(10^{-1},10^{-4},0.5,1)$ & $(0.2,10^{-10},0.2,10^{-3},10^{-4})$ & $(10^{-2},10^{-1},1,10^{-2})$
\end{tabular}
\end{sidewaystable}

\section*{Acknowledgments}
This material is based upon work supported by the Air Force Office of
Scientific Research (AFOSR) under award numbers FA9550-20-1-0236, FA9550-22-1-0002,
FA9550-22-1-0004. The content of this publication does not necessarily reflect the position
or policy of any of these supporters, and no official endorsement
should be inferred.

\bibliographystyle{plain}
\bibliography{biblio,biblio_intro}

\end{document}